ACADEMY OF SCIENCES OF MOLDOVA
INSTITUTE OF MATHEMATICS AND COMPUTER SCIENCE

Manuscript

**POPA ALEXANDRU**

# ANALYTIC GEOMETRY OF HOMOGENEOUS SPACES

## 111.04 GEOMETRY AND TOPOLOGY

**Doctor Thesis in Mathematics**

| | | |
|---|---|---|
| Scientific Adviser: | ________ | Damian Florin, doctor, associate professor |
| Author: | ________ | Popa Alexandru |

**CHIȘINĂU, 2017**

**ACADEMIA DE ȘTIINȚE A REPUBLICII MOLDOVA**
**INSTITUTUL DE MATEMATICĂ ȘI INFORMATICĂ**

Cu titlu de manuscris

**POPA ALEXANDRU**

# GEOMETRIA ANALITICĂ A SPAȚIILOR OMOGENE

## 111.04 GEOMETRIE ȘI TOPOLOGIE

**Teză de doctor în Matematică**

| | | |
|---|---|---|
| Conducător științific: | ________ | Damian Florin, doctor în matematica, conferențiar universitar |
| Autorul: | ________ | Popa Alexandru |

**CHIȘINĂU, 2017**

## Contents

UDC: 514.742.2:514.120/514.140+514.764.21

## ANNOTATION

**Popa Alexandru, "Analytic geometry of homogeneous spaces", PhD thesis, Chisinau, 2017.**

The thesis is written in English and consists of: introduction, three chapters, general conclusions and recommandations, appendix, 210 bibliography titles, 140 pages of main text, 27 figures, 9 algorithms, 5 tables. The obtained results were published in 9 scientific papers.

**Domain of research:** Geometry of homogeneous spaces.

**Goals and objectives:** The goal of the research is to provide a toolchain that can be used to study of homogeneous spaces by means of linear algebra. The objectives of the research are: introduction of the new concept of the space signature, construction of homogeneous space based on signature concept, construction of the model of homogeneous space with given signature, expression of the measurement of different geometric quantities via signature, different applications of the analytic geometry of homogeneous spaces.

**Scientific innovation of obtained results:**

- Analytic geometry is developed in linear algebra language, even for non–linear spaces.
- One universal theory is developed that uses the elements of space signature as parameters.

**Important scientific problem solved:** The investigation of the homogeneous spaces with linear methods via concept of the signature.

**Theoretical and practical value of the work:** Rezultatele prezentate în teză sunt noi, au un caracter teoretic și cu ajutorul conceptului de signatură prezintă o teorie generală a spațiilor omogene.

**Implementation of scientific rezults:**

- New results can be used in investigation of the problems of differential geometry, in theoretic physics and in other domains where notion of the signature can be applied in the given sense.
- The thesis can be used as the didactic support for optional courses in the university and doctoral studies.

## ADNOTARE

**Popa Alexandru, „Geometria analitică a spațiilor omogene", teză de doctor în științe matematice, Chișinău, 2017.**

Teza este scrisă în engleză și constă din: introducere, trei capitole, concluzii generale și recomandări, apendice, bibliografie din 210 de titluri, 140 de pagini de text de bază, 27 de figuri, 9 algoritmi, 5 tabele. Rezultatele obținute sunt publicate în 9 lucrări științifice.

**Cuvinte-cheie:** Spațiu omogen, spațiu Riemanian, geometrie Klein, metrică proiectivă, geometrie analitică.

**Domeniul de studii:** Geometria spațiilor omogene.

**Scopul și obiectivele tezei:** Scopul cercetării este să se ofere un instrument, care poate fi folosit pentru a studia spații omogene în limbajul algebrei lineare. Obiectivele sunt: argumentarea conceptului de signatură, pe baza lui, construirea spațiilor omogene, construirea modelului spațiului omogen cu signatura dată, expresia măsurării cantităților geometrice via signatura, aplicațiile geometriei analitice a spațiilor omogene.

**Noutatea științifică a rezultatelor obținute:**

- Geometria analitică este dezvoltată folosind limbajul algebrei lineere, chiar și pentru spații nelineere.
- Este dezvoltată o teorie universală, în care elementele signaturii spațiului sunt parametri.

**Problema științifică importantă soluționată:** Cercetarea spațiilor omogene prin metode lineare, aplicând conceptul de signatură.

**Semnificația teoretică și valoarea aplicativă a tezei:** Rezultatele prezentate în teză sunt noi, au un caracter teoretic și cu ajutorul conceptului de signatură prezintă o teorie generală a spațiilor omogene.

**Implementarea rezultatelor:**

- Rezultate noi pot fi folosite în investigarea problemelor în geometria diferențială, în fizica teoretică și în alte domenii unde poate fi aplicat conceptul de signatură în sensul dat.
- Teza poate fi folosită în calitate de suport pentru cursurile opționale universitate și postuniversitare.

УДК: 514.742.2:514.120/514.140+514.764.21

## АННОТАЦИЯ

**Попа Александру, «Аналитическая геометрия однородных пространств», докторская диссертация, Кишинёв, 2017.**

Диссертация написана на английском языке и состоит из: введения, трёх глав, общих выводов и рекомендаций, приложения, библиографии из 210 наименований, 140 страниц основного текста, 27 иллюстраций, 9 алгоритмов, 5 таблиц. Результаты исследований опубликованы в 9 научных работах.

**Ключевые слова:** Однородное пространство, Риманово пространство, геометрия Клейна, проективная метрика, аналитическая геометрия.

**Область исследования:** Геометрия однородных пространств.

**Цели и задачи исследования:** Цель исследования — предоставить инструмент для исследования однородных пространств на языке линейной алгебры. Задачи исследования: введение сигнатуры, построение однородного пространства с его помощью, построение модели однородного пространства по сигнатуре, выражение геометрических мер через сигнатуру, приложения аналитической геометрии однородных пространств.

**Научная новизна и оригинальность:**

- Аналитическая геометрия построена на языке линейной алгебры, даже для нелинейных пространств.
- Разработана одна универсальная теория, в которой элементы сигнатуры пространства являются параметрами.

**Решенная научная проблема:** Исследование однородных пространств с помощью линейных методов, посредством концепцию сигнатуры.

**Теоретическая и прикладная значимость:** результаты данной работы новы, имеют теоретический характер, и с помощью концептии сигнатуры предоставляют общую теорию однородных пространств.

**Внедрение научных результатов:**

- Новые результаты можно использовать в рассмотрении задач дифференциальной геометрии, в физике и других областей, в которых имеет смысл понятие сигнатуры в указанном смысле.
- Диссертацию можно использовать в качестве учебного руководства в факультативных курсах в университете и аспирантуре.

## INTRODUCTION

This thesis is not about analytic geometry in classical sense of this notion. These two topics, analytic geometry and homogeneous spaces, seem incompatible, because homogenous spaces are studied by means of differential geometry. This study isn't about "analytic geometry and homogeneous spaces", but about "analytic geometry of homogeneous spaces". Unlikely somebody will object, that analytic geometry approach, if it is posible, is more preferred than differential geometry approach. The goal of this work is to show this approach isn't only possible, but also remarkably easy.

**Actuality of theme.** The study of homogeneous space gained a great importance in recent years. Because homogeneous spaces arise in different domains of mathematics and physics, there were elaborated different theories that, using different terms and approaches, study essentially the same entity. Depending on domains and methods used, the homogeneous spaces are essentially the same object of research as pseudo–Riemannian, semi–Riemannian spaces, spaces with projective metrics and their geometry is also called "Klein geometry" or "geometry in large". The standard methods used for study of homogeneous spaces is algebraic geometry or differential geometry. The analytic geometry of homogeneous spaces may greatly simplify the non–trivial studies in several domains. Here are some examples:

**Domain 1: Non–linear volume theory in geometry.** It is well known that the volume theory for non–linear spaces is not trivial. Namely, still it is not known symbolic form of volume equation in elementary functions even for the most basic space forms e.g. for orthoscheme. The described theory gives some insights in this domain. The equation is expected have the form:

$$Tr(v) = Alg(Tr(a), Tr(b), Tr(c)),$$

where $a, b, c$ are any three orthoscheme parameters, which are not related among them, $Tr(x)$ is some trigonometric function having the same type as its argument and $Alg(x, y, z)$ is some algebraic function. Moreover, while the type of right–handed equation parameters (and respective trigonometric functions) are known a priori, the type of volume (and respective trigonometric function) is expected to be equal to $K_1 K_2 K_3$. It means that:

$$v = Tr^{-1} Alg(Tr(a), Tr(b), Tr(c)).$$

That is why this equation is so difficult to find based on volume integral. It is much easier to find the algebric function $Alg(x, y, z)$ having the enumerable set of them, and to find the coefficients by the method of unknown coefficients.

For example, this reasoning makes it possible to test to some extent the Rational Volume conjecture: Let $P$ be an elliptic polyhedron whose dihedral angles are in $\pi \cdot \mathbb{Q}$. Then $Vol(P) \in \pi^2 \cdot \mathbb{Q}$. It may be true for angles $\frac{\pi}{2}, \frac{\pi}{3}, \frac{\pi}{4}, \frac{\pi}{5}, \frac{\pi}{6}$ whose elliptic trigonometric functions are simple algebraic expressions of integers, but for e.g. the angle $\frac{\pi}{7}$ there is no such expression so the theory predicts for the polyhedron with such an angle to have not algebraic value of trigonometric function of volume and, the most probably, the volume does not belong to $\pi^2 \cdot \mathbb{Q}$. Unfortunately it is much more difficult to prove some number is not rational than it is so.

**Domain 2: String theory in physics.** As known, the string theory is a promising candidate to Theory of Evething in physics. However, any type of consistent string theory describes a space–time with much more dimensions than physically observable. There are two possible explanations of physical unobservability of higher dimensions: *localization* and *compactification*. For the second one the following example is given: a tube is a two–dimensional space, however when its diameter is rather small, or it is observed from far away, it is percieved as one–dimensional space. This explanation draws the space–time that is homogeneous only at local scale. This seems to not correspond to the most fundamental physical lows that assume globally homogeneous and isotropic physical space–time. This thesis theory describes globally homogeneous two–dimensional space on tube, that enables its physical properties to be fully and globally homogeneous and isotropic. The theory gives also examples of globally homogeneous spaces of 10 and 26 dimensions which reflect observable physical space–time.

At the same time, described theory gives one more possible explanation of physical unabservability of higher dimensions: the 4-dimensional space–time we know is a *limit lineal* in higher dimensional space (the dimension of the limit lineal is always less than the space dimension) and since limit lineals properties may greatly differ from space properties, its inhabitants may be unable to perceive other dimensions.

**Goals and objectives.** The goal of the research is to provide a toolchain that can be used to study of any homogeneous space by means of analytic geometry and linear algebra. The objectives of the research are:

- Introduct of the new concept of the space signature;
- Construction of homogeneous space based on concept of signature;
- Construction of a model of homogeneous space for each given signature;
- Expression of the measurement of different geometric quantities via signature, which reflects their role in analytic geometry of homogeneous spaces;
- Finding of different applications of the analytic geometry of homogeneous spaces.

**Research methodology.** As subject of study, the *model* of homogeneous spaces is constructed. The methods of *linear algebra* is used. In order to see that the model is adequately constructed, its relevant to research properties are verified.

Thesis investigates the notion of homogeneous space and establishes what kind of properties are to be expected. This is achieved by *generalization* of axioms of known geometries. Then, using *duality* method new properties are formulated, dual to known ones.

An important method used in research is *computer modelling*, realized as software project GeomSpace. This project is based on theory results and also is used to test conjuctures, to provide counterexamples and to obtain new results.

**Scientific innovation.** The theory uses methods and language of linear algebra to study non–linear spaces. These technics can be used particularly to describe analytic geometry of non–linear elliptic, hyperbolic, De Sitter and Anti de Sitter spaces.

The main innovation of elaborated theory is space *parameterization* by introduction of space signature. This parameterization allows studying of different homogeneous spaces in one global framework. When the parameters are used as variables in definitions, axioms, equations, theorems, proofs, all these have exactly the same form that describes the reality of all homogeneous spaces simultaneously. When it is necessary to describe some space particularities or to see the difference between two concrete spaces, the concrete values can be put in parameters of each definition, axiom, equation, theorem and proof.

The parameterized approach has many advantages. Firstly, because of the large number of homogeneous spaces, it is next to impossible to describe each space with its geometry one by one. The uniform approach allows to describe any concrete space with all its particularities or to compare two spaces. Secondly, some results are more easy to obtain for non–linear spaces, and then they can be generalized to linear ones. Other results are more easy obtained for linear spaces, and then they can be extended to non–linear ones. The uniform approach drastically simplifies such extensions and generalizations.

**Important scientific problem solved** is the investigation of the homogeneous spaces via linear methods applying the concept of signature.

**Theoretical and practical value of the work.** This work provides one possible universal terminology across different spaces, that facilitates comparison of different spaces properties.

Although the thesis is focused on analytic geometry is has value for differential geometry. The modern differential geometry methods are constructed based on analytic geometry of Euclidean space. This fact has two consequences. First, differential geometry is good to describe the difference between properties of Euclidean space and some non–Euclidean space or to show what Euclidean properties are missing from some non–Euclidean space. However, this approach is useless to describe some non–Euclidean space property that is missing from Euclidean space or to compare such properties of two non–Euclidean spaces. Second, differential geometry is useful for study of locally Euclidean metric spaces. However it has very limited applicability for homogeneous spaces that can't be approximated by Euclidean space in any point. The described theory may be used as toolchain for more universal differential geometry.

The described methods have immediate application also in domains of partial differential equations, general relativity theory, quantum field theory, string theory and M-theory, AdS/CFT correspondence, cosmology and others.

**Obtained scientific results to be defended:**

- Classification of homogeneous spaces based on new form of signature, introduced in this work;
- Elaboration of the universal form of trigonometric equations, common for all homogeneous spaces;
- Introduction of the group of generalized orthogonal matrices and its study in connection to the motion group of homogeneous space;
- Adaptation of the routines and algorithms of linear algebra, that operate on vectors or on vector families, to homogeneous spaces;
- Introduction of decomposition vectors for limit vectors (isotropic vectors) and study of limit vectors with aim of their decomposition vectors;
- Definition and study of volumes in a homogeneous space by corresponding volumes in its metaspase;
- Theorem on isomorphism of the crystallographic groups of dual homogeneous spaces.

**Implementation of scientific results.** Because of theory simplicity, it can be used as facultative course in lyceum or university in order to develop the geometric intuition.

Another implementation of theory is software project Geomepace (http://sourceforge.net/projects/geomspace/). The goal of GeomSpace is to study homogeneous spaces and to make it easy to work with them (visualization, motion, constructions, calculus).

**Approval of work.** The different aspects of the theory were presented at the following scientific events:

**2009** Alba Iulia, România — International Conference on Theory and Applications in Mathematics and Informatics,

**2010** Moscow, Russia — International conference "Metric geometry of surfaces and polyhedra", dedicated to 100th anniversary of N. V. Efimov,

**2010** Moscow, Russia —The International Conference "Geometry, Topology, Algebra and Number Theory, Applications" dedicated to the 120th anniversary of B. N. Delone,

**2014** Chișinău, Moldova — The Third Conference of Mathematical Society of Moldova, IMCS-50,

**2015** Tula, Russia — International conference, dedicated to 85th anniversary of professor S. S. Ryshkov. Algebra, Number Theory and Discrete Geometry: Modern Problems and Applications,

**2015** Iași, România — The 8th Congress of Romanian Mathematicians,

**2016** Chişinău, Moldova — International Conference Mathematics & Information Technologies: Research and Education, MITRE — 2016.

**Publications.** The results of the theory were published in 9 works: [56, 57, 58, 59, 60, 61, 62, 63, 154].

**Summary of thesis contents.** The structure of the thesis is as follows. The first chapter gives an overview of research in domain of homogeneous spaces. Also it defines and classifies the homogeneous spaces. The second chapter constructs the analytic geometry of homogeneous spaces. Reasonong is based on linear algebra apparatus. To make ideas more obvious, there are given also parametrization of some Euclid's and Lobachevsky's axioms, even if description isn't based on them. For similar purpose, the duality principle, described in the first chapter, is often used.

The third chapter gives some ideas regarding the application of the theory in mathematics, that demonstrates its value. The appendix presents more possible applicatons of the theory in physics and presents GeomSpace project. This software project coexist in synergy with the theory. On the one hand, theoretic results are used in GeomSpace. On the other hand, the software is literally the testing ground for theory correctness. If some algorithm doesn't work, it is reviwed.

# 1. ANALYSIS OF SITUATION IN DOMAIN OF HOMOGENEOUS SPACES

The introductory chapter defines basic notions and methods and presents an overview of recent research in homogeneous spaces domain.

## 1.1. The main definitions and notions.

Recall basic algebraic structures.

**Group.** A *group* $G$ is a set equipped with operation

$$a + b = c \in G, \ \forall a, b \in G,$$

that satisfies the following axioms:

1. $a + (b + c) = (a + b) + c$, $\forall a, b, c \in G$,
2. $\exists 0 \in G$, so that $0 + a = a + 0 = a$, $\forall a \in G$,
3. $\forall a \in G \quad \exists(-a) \in G$, so that $a + (-a) = (-a) + a = 0$.

The 0 element is called *neutral* element of the group and $(-a)$ element is called *opposed* or *inverse* element of $a$. The group $G$ is called *commutative* if it additionally satisfies the following property:

4. $a + b = b + a$, $\forall a, b \in G$.

**Field and skew field.** A *field* (*skew field*) $F$ is a set equipped with two operations:

- $u + v = w \in F$, $\forall u, v \in F$
- $u \cdot v = w \in F$, $\forall u, v \in F$,

so that:

1. $F$ forms a commutative group with respect to addition '+' operation,
2. $F \setminus 0$ forms a group with respect to multiplication '·' opetation, where 0 is neutral element of addition '+' operation,

and satisfies the following axioms:

1. $(u + v) \cdot w = (u \cdot w) + (v \cdot w)$, $\forall u, v, w \in F$,
2. $u \cdot (v + w) = (u \cdot v) + (u \cdot w)$, $\forall u, v, w \in F$.

If $F$ is commutative with respect to multiplication '·' operation, it is named field, otherwise it is sometimes named skew field.

**Linear vector space.** A *vector space* over a field $(V, F)$ is a pair composed of a commutative group $V$ and a field (or a skew field) $F$ together with operation:

$$\alpha \cdot u = v \in V,\ \forall u \in V, \forall \alpha \in F,$$

that satisfies the following axioms:

1. $\alpha \cdot (\beta \cdot v) = (\alpha \cdot \beta) \cdot v,\ \forall v \in V, \alpha, \beta \in F,$
2. $1 \cdot v = v, \forall v \in V$, where 1 is neutral element of $F$,
3. $\alpha \cdot (u + v) = (\alpha \cdot u) + (\alpha \cdot v),\ \forall u, v \in V, \alpha \in F,$
4. $(\alpha + \beta) \cdot v = (\alpha \cdot v) + (\beta \cdot v),\ \forall v \in V, \alpha, \beta \in F.$

The elements of the field are called *scalars* and the elements of the group *vectors*.

Important algebric structures that are used in this work are the field of real numbers $\mathbb{R}$, the skew field of square invertible matrices of size $n$, the real $n$-dimensional vector space $\mathbb{R}^n$.

Recall basic geometric properties and operations with linear vector spaces.

**Linear span of vectors.** A *linear span* of a vector set $U = \{u_\omega\}_{\omega \in \Omega}$, where each $u_\omega \in V$ and $\Omega$ is some domain of parameters, is defined as:

$$Span(U) = \sum_{i=1}^{n} \alpha_i u_{\omega_i},\ n \in \mathbb{N}, u_{\omega_i} \in V, \alpha_i \in F.$$

The linear span of some vector family forms a linear vector space.

Vectors $\{v_i\}_{i=\overline{1,n}} \in V$ are called *linearly independent* if no one is present in linear span of all others. Conversely, these vectors are called *linearly dependent* if some vector is present in linear span of others:

$$v_i = \sum_{j \neq i = 1}^{n} \alpha_j v_j.$$

The representation of $v_j$ from above is named *linear combination* of vectors $\{v_j\}_{j \neq i = \overline{1,n}}$. If a vector can be presented as a linear combination of some linearly independent set of vector, this representation is unique.

A minimal subset of vectors $B = \{b_i\}_{i=\overline{1,n}}$ that are linearly independent and any vector $v \in V$ can be presented as linear combination of this vector set, is named *basis* of vector space:

$$V = Span(B).$$

The number $n$ of vectors in space basis doesn't depend on vectors choice and is called space *dimension* and is denoted as $\dim V$. Any space vector can be uniquely represented as linear combination of space basis.

**Scalar product of vectors.** The *scalar product* or *dot product* or *inner product* of vectors from some linear vector space $(V, F)$ is the operation:

$$u \cdot v = \alpha \in F,\ \forall u, v \in V,$$

that satisfies the following properties:

1. $u \cdot v = v \cdot u,\ \forall u, v \in V,$
2. $u \cdot (v + w) = (u \cdot v) + (u \cdot w),\ \forall u, v, w \in V,$
3. $u \cdot (\alpha \cdot v) = (\alpha \cdot u) \cdot v = \alpha \cdot (u \cdot v),\ \forall u, v \in V, \forall \alpha \in F.$

Two vectors are called *orthogonal*, if their scalar product equals to zero:

$$u \perp v \iff u \cdot v = 0_F.$$

Nonzero orthogonal vectors are linearly independent. The mutually orthogonal vectors are usually chosen as the space basis.

The scalar product of vectors is also called the *metric* of the space in sense that the metric in linear vector space is intruduced by means of scalar product. The *norm* of vector is defined as:

$$|v|^2 = v \cdot v,\ \forall v \in V.$$

The scalar product is:

- *Positive definite*, if $v \cdot v > 0_F,\ \forall v \neq 0_V \in V,$
- *Degenerate*, if $v \cdot v \geq 0_F,\ \forall v \in V$ and $\exists v \neq 0_V \in V$, so that $v \cdot v = 0_F$,
- *Indefinite*, if $\exists v \in V$, so that $v \cdot v < 0_F$.

Sylvester's law of inertia states that the number of basis vector with real, zero and imaginary norm doesn't depend on basis choice, so it is property of space metric.

The vector is said to be *normalized* if its norm equals to 1:

$$|v| = 1.$$

**Space linear transformations.** Any linear transformation $\mathfrak{M}$ of linear space can be represented as square matrix $M$ of size equal to space dimension. In this case vector $v = (v_1, ..., v_n)$ is transformed to $v' = \mathfrak{M}v$ as:

$$v'_i = \sum_{j=1}^{n} m_{ij} v_j,\ i = \overline{1, n}.$$

The composition of transformations is represented as product of corresponding transformation matrices:

$$\mathfrak{Q}(\mathfrak{P}v) = (\mathfrak{Q}\mathfrak{P})v$$

where matrix product $M = QP$ is defined as follows:

$$m_{ij} = \sum_{k=1}^{n} q_{ik} p_{kj}.$$

**Projective space.** The $n$-dimensional *projective space* over some field $F$ is denoted as $\mathbb{FP}^n$ and is the set of vectors $\mathbb{F}^{n+1} \setminus O$, where $O = (0, ..., 0)$ is origin, with equivalence defined as:

$$v' \equiv v \Longleftrightarrow \exists \lambda \neq 0, v' = \lambda v.$$

Usually vectors in projective space are written in *homogeneous coordinates* which are chosen from vector coordinates in some unique way, for example normalized against some norm:

$$v = (v_1 : ... : v_{n+1}) \in \mathbb{FP}^n.$$

Note, that a vector in $n$-dimensional projective space has $n + 1$ coordonates.

### 1.2. Definition and classification of the homogeneous spaces.

Informally, the space is homogeneous when it looks equally at each its point.

**Definition 1.2.1** (Homogeneous space). Let $X$ be a non–empty set and $G$ a group that acts on $X$. The structure of $X$ is $\tau : G \to Aut(X)$. The pair $(X, \tau)$ is called *homogeneous space* [67, 84], if:

- $\tau$ is a homomorphism, that is for each $g \in G$ the maping $\tau(g)$ is structure preserving;
- $\tau(G)$ acts transitively on $X$.

The group $G$ is called *simmetry group* or *group of motions* of $X$.

The homogeneity can be viewed in warious ways, resulting in different structures on space: as isometry (giving rigid geometry), as diffeomorphism (giving differential geometry) or as homeomorphism (giving topology). In this thesis the rigid geometry of homogeneous spaces is studied.

**Metric.** Some set $X$ is said to be equipped with *metric* if there is the function $d : X \times X \to \mathbb{R}$, named the *distance*, with the following properties:

1. Non-negativity: $d(X, Y) \geq 0$,
2. Identity of indiscernibles: $d(X, Y) = 0 \Leftrightarrow X = Y$,
3. Symmetry: $d(X, Y) = d(Y, X)$.
4. Triangle inequlity: $d(X, Z) \leq d(X, Y) + d(Y, Z)$.

The homogeneous spaces described in this thesis are complete continuos metric spaces. What does it mean? Firstly, homogeneous space has a rich simmetry group. All objects of the same dimension (points, lines, planes, subspaces) are congruent, that is, map one to another with space motions. The number of degrees of freedom depends only on dimension $n$ and equals to $\frac{n(n-1)}{2}$.

Secondly, homogeneous space is *metric* space. Notably, the metric of homogeneous spaces described in this thesis *does not* always satisfy the triangle inequality, but satisfies its generalized form. The homogeneous space metric means more. Between any two objects (points, lines, planes, subspaces) there exists a measure (distance, angle or inclination) even if these objects are of different dimension. That is, it is possible to find the distance between a point and a plane, or the angle between a line and a subspace. We will introduce space metric by means of motions.

Thirdly, homogeneous space is *continuos* space. In other words, if two objects $X, Y$ are congruent, that is, if exists proper motion $\mathfrak{M}$ that maps $X$ to $Y$ ($Y = \mathfrak{M}X$), then exists also continuos parameterization of motion $\mathfrak{M}(\varphi)$, $\varphi \in [0, 1]$ so that $\mathfrak{M}(0) = \mathfrak{I}$ is identity motion, $\mathfrak{M}(1) = \mathfrak{M}$ is original motion and all $\mathfrak{M}(\varphi)$ for all $\varphi \in [0, 1]$ are space motions, such that imagies $\mathfrak{M}(\varphi)X$ of all objects are space objects.

Fourthly, homogeneous space is *complete* space. That is, all its properties, including described above, are global. The completeness makes the difference betweeen a space and a manifold, which has the same properties at local scale.

*Remark.* Evidently, Euclidean, elliptic and hyperbolic spaces are homogeneous.

*Remark.* Homogeneous spaces arise in different domains of mathematics and physics. Depending on the domain where these spaces arise and on studied properties, the spaces or their geometry are also called: spaces of *projective metric* [104, 128, 193], *pseudo–Euclidean* [111, p. 72], [116, p. 512], [170, p. 737], *semi–Euclidean* [192, p. 1327], *pseudo–Riemannian* or *semi–Riemannian* spaces [14, p.417], [172, p. 193], *Klein geometry* [70, p. 168], [160, p. 103], *geometry in large* [93, 115, 136]. It is important to understand, that all these geometric objects represent the same entity, even if they are obtained in different ways.

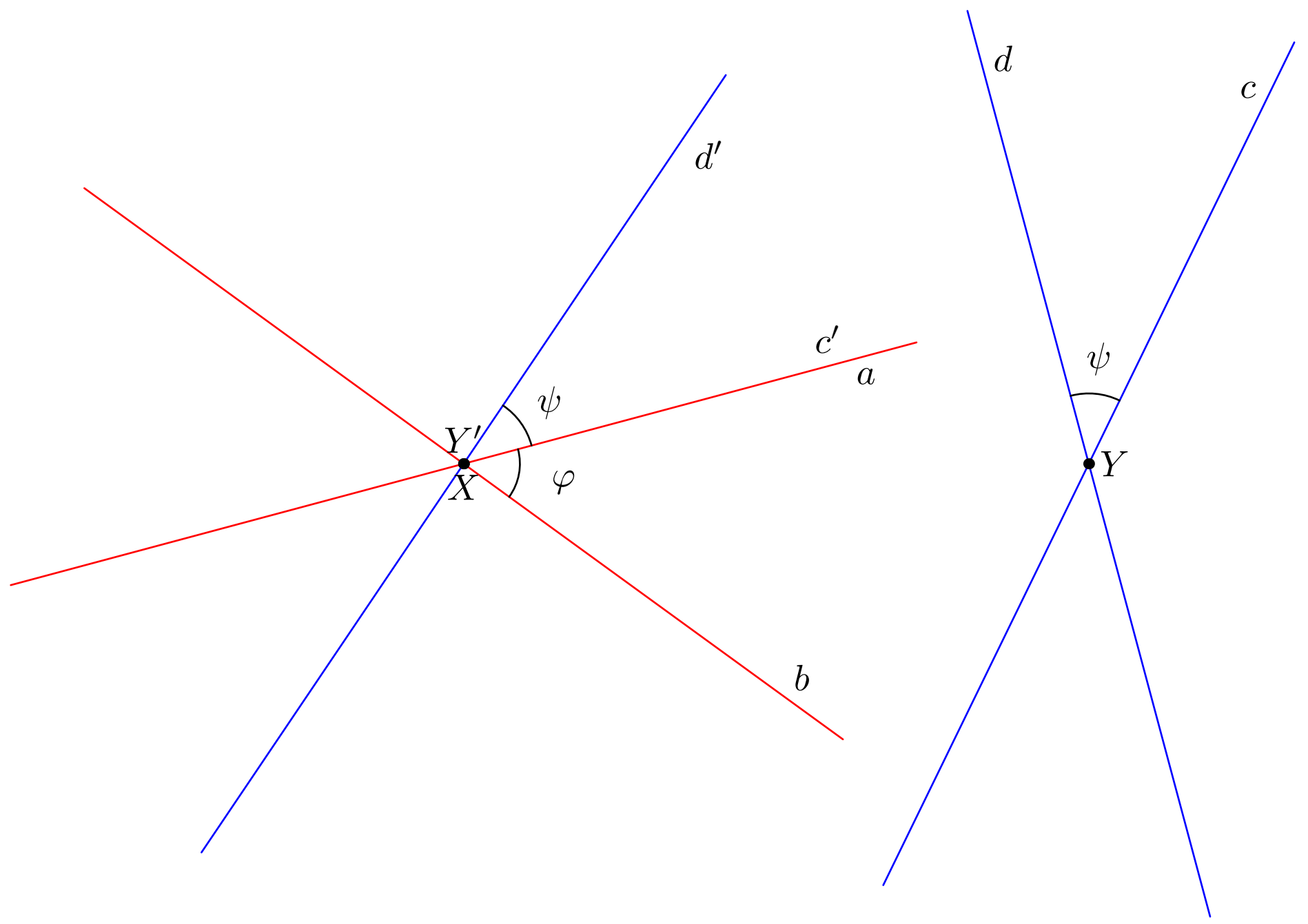

Figure 1.1: Angles $\varphi$ and $\psi$ have the same type because of congruency of points $X, Y$ and lines $a, c$.

In order to classify the homogeneous spaces, consider the following scheme:

- In 0-dimensional space exists one single object, a point. So, there is exactly one zero–dimensional space.

- In 1-dimensional space there are points and single line. The distance between points can be one of these three kinds: *elliptic*, *parabolic*, *hyperbolic*. It means, that motions that map one point to another have elliptic, parabolic or hyperbolic character. All points are congruent. Hence, all distances are of the same kind from these three. It gives us 3 one–dimensional spaces.

- In 2-dimensional space there are points, lines and single plane. All lines can be of some kind from three enumerated above. As any two lines are congruent, all them have the same kind. The angle between two lines can also be of one kind of three: elliptic (this one is present in Euclidean, elliptic and hyperbolic spaces), parabolic and hyperbolic. Consider two figures, each containing two intersecting lines (Figure 1.1, figures $Xab$ and $Ycd$). Because points and lines are congruent, there exists a motion that maps point $Y$ to $Y' = X$ and line $c$ to $c' = a$. Then the angle $\psi$ between lines $c, d$ equals to the angle between $c', d'$, which continues the angle $\varphi$ between lines $a, b$. Therefore, these two angles, and also all angles, have the same kind. It gives us $3 \times 3 = 9$ two–dimensional spaces.

- In 3-dimensional space beside points and lines, there are 2-dimensional planes (and of cource, single 3-dimensional subspace). Each plane can be of one of above enumerated 9 kinds. Moreover, starting with 3-dimensional space the new measure appears, the dihedral angle. As above, dihedral angle between any two planes can be of a single kind from three. It gives us $9 \times 3 = 27$ three–dimensional spaces.

Continuing the scheme, reach the following conclusion:
In $n$-dimensional homogeneous space exist:

- Objects of $n$ different dimensions (points, lines, planes, subspaces of different dimensions). All objects of the same dimension are congruent.
- Measures of $n$ different dimensions, one for each object dimension. Each measure describes the level of departure or deviation between two objects of corresponding dimension and can be of one of three kinds: elliptic, parabolic or hyperbolic.

So there are $3^n$ $n$-dimensional homogeneous spaces.

### 1.3. Space duality.

At different levels of geometry arise different kinds of duality. For example, duality plays an important role in projective geometry. Also, it is easy to observe duality of regular polyhedra of each dimension. Duality is a powerful tool for construction of new figures.

Duality plays fundamental role also in study of homogeneous spaces [2, p. 11]. In this case with the power of duality one can produce not new figures in a space, but whole new spaces with completely new geometry.

**Example.** Consider Poincaré model of hyperbolic geometry (Figure 1.2). As known, points of hyperbolic plane are Euclidean points inside the disk, which we call horizon. Hyperbolic lines are either disk diameters or circles arcs, that intersect horizon orthogonal to it and lie entirely inside the horizon. The angle between intersecting hyperbolic lines is Euclidean angle between these line arcs (betweeen arc tangents at intersection point). The distance between points $E, F$ of line $MN$ (Figure 1.2, $M$ and $N$ lie at horizon) is the relation:

$$d(E, F) = \ln\left(\frac{|EN|}{|FN|}\frac{|FM|}{|EM|}\right),$$

where $|EN|, |FN|, |EM|, |FM|$ are Euclidean lengths of corresponding segments.

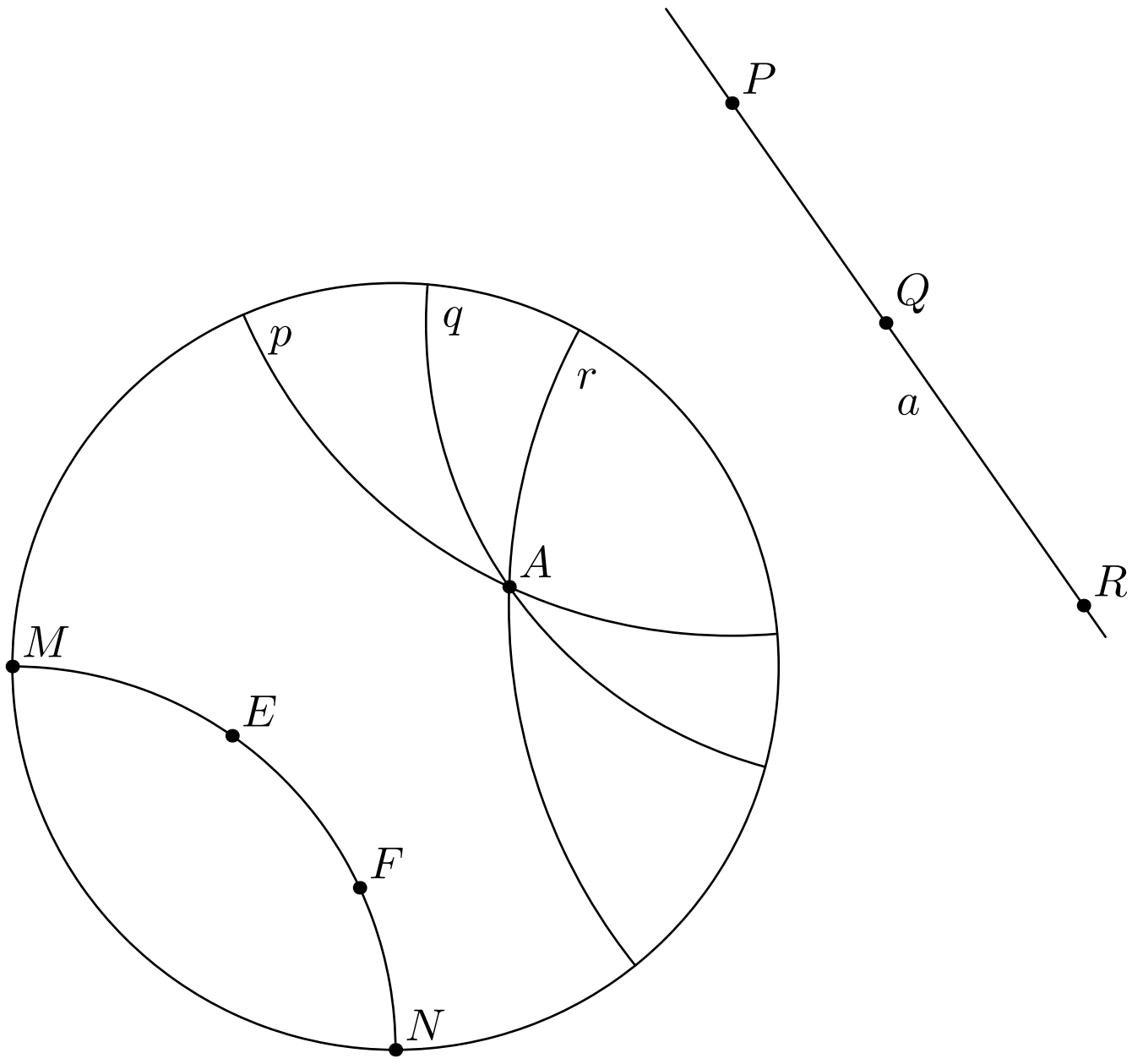

Figure 1.2: Construction of anti–hyperbolic geometry model.

Firstly, observe that in case hyperbolic lines are arcs, their Euclidean centres always lie outside the horizon. Going to the limit, we can consider that horizon diameters are arcs of infinitely large circles and hence their centres are infinitely far (in Euclidean meaning) in direction perpendicular to diameter (opposite directions are the same point at infinity). Thus, each hyperbolic line corresponds to one point from outside horizon and vice versa, each point from outside the horizon (including ones at infinity) corresponds to one hyperbolic line. Let call outside part of horizon *anti–hiperbolic plane.*

Secondly, observe that if it is given a pencil of intersecting lines in hyperbolic plane (in Figure 1.2, lines $p, q, r$ intersect in point $A$), then their Euclidean centres lie on a line (in Euclidean sense) entirely in anti–hyperbolic plane (in Figure 1.2 corresponding points are $\boldsymbol{P}, \boldsymbol{Q}, \boldsymbol{R}$). Let name Euclidean line, that contains all centres of some pencil of intersecting hyperbolic lines, *line* of anti–hyperbolic plane. Since the pencil of intersecting hyperbolic lines is uniquely defined by the intersection point, this hyperbolic point corresponds to the line of anti–hyperbolic plane (in Figure 1.2, point $A$ corresponds to line $a$).

Thirdly, it is easy to define *distances* between anti–hyperbolic points as angles between corresponding hyperbolic lines. Likewise, let define *angle* between anti–hyperbolic lines as distance between centres of two corresponding pencils of intersecting hyperbolic lines.

In this way, we constructed anti–hyperbolic plane and introduced metric in it. It is important, this description is free of contradictions, because each proposition or equation of anti–hyperbolic plane corresponds to true proposition or equation of hyperbolic plane. For example, hyperbolic triangles correspond to anti–hyperbolic triangles (vertices correspond to sides and vice versa).

That is why it isn't necessary to deduce triangle laws in anti–hyperbolic plane. It is sufficient to rewrite corresponding laws of hyperbolic plane. In this case distances become of elliptic kind and angles become of hyperbolic kind.

Describing anti–hyperbolic geometry, we made more then just reformulate known facts from terms of one geometry to another. With duality we can describe absolutely new and unfamiliar notions. As example, one can observe that centres of divergent line arcs of hyperbolic plane also lie in an Euclidean line. However, this line intersects the horizon. So, we can define its part from outside the horizon as new kind of objects in anti–hyperbolic plane, which have no analogy in notions of hyperbolic plane. The distance between divergent hyperbolic planes defines hyperbolic measure of "distance" in new object (the true anti–hyperbolic distance is of elliptic kind).

Moreover, one can observe that centres of parallel line arcs of hyperbolic plane also lie in one Euclidean line. This line is tangent to horizon. It can be considered as one more kind of objects in anti–hyperbolic plane. Parabolic kinded inclination of parallel hyperbolic lines defines parabolic "distance" in this object type. Thus, in anti–hyperbolic plane through each point pass objects of three different kinds: lines (elliptic kind of distance, don't intersect the horizon), hyperbolic kinded "non–lines" (intersect horizon) and only two "limit lines" (parabolic kind of distance, tangent to horizon).

Finally, observe that duality can be applied to any homogeneous plane regardless of the model. This goal is achieved as simple as interchanging the following phrases or notations:

$$\begin{aligned}
\text{Line } l &\longleftrightarrow \text{Point } L,\\
P \in l &\longleftrightarrow p \ni L,\\
P \notin l &\longleftrightarrow p \not\ni L,\\
|AB| = \varphi &\longleftrightarrow \measuredangle ab = \varphi,\\
C = a \cap b &\longleftrightarrow c = AB,\\
a \cap b = \varnothing &\longleftrightarrow \text{Points } A \text{ and } B \text{ don't lie on a line.}
\end{aligned}$$

It affects only obviousness, but not correctness of the description. Of course, the duality can be applied to many–dimensional spaces too. If space dimension is $n$, then $(m-1)$–dimensional planes and measures of initial space correspond to $(n-m)$–dimensional planes and measures of dual space (zero–dimensional planes are points, one–dimensional planes are lines). More precisely, space dual to some $n$-dimensional homogeneous space is also homogeneous $n$-dimensional space. It has the following properties:

- $(m-1)$–dimensional planes of initial space correspond to $(n-m)$–dimensional planes of dual space.

- $(m-1)$–dimensional measures of initial space correspond to $(n-m)$–dimensional measures of dual space and are of the same kind.
- Equal measures of initial space correspond to equal measures of dual space.
- Congruent figures of initial space correspond to congruent figures of dual space.
- Dual space of dual space is initial space.

**1.4. Short history of non–Euclidean geometry.**

Geometry is one of the oldest domains of mathematics. Long time it was considered as natural science. As such, it was considered that there can be the only possible geometry. This geometry was retroactively named Euclidean geometry after Euclid, the author of "Elements", which described the geometry on axiomatic basis in 3rd century BCE.

Curiosly, the Euclid's Elements is the cause of non–Euclidean geometry discovery. Euclid's 4th postulate, that reads "All right–angles are equal to one another", was proven to be theorem shortly after Euclid. It follows from the first common noton, that states "Things equal to the same thing are also equal to one another". The 5th postulate states "If a line segment intersects two straight lines forming two interior angles on the same side that sum to less than two right angles, then the two lines, if extended indefinitely, meet on that side on which the angles sum to less than two right angles". Is has so complex and unintuitive formulation, that generations of geometers tried to prove it as theorem. The attempts to prove the 5th Euclid's postulate by contradiction were in fact theorems of what is now known as hyperbolic geometry.

In 1829, Nikolai Lobachevsky, and independently of him, in 1831, János Bolyai published the axiomatically developed theory where non–Euclidean geometry, presently known as hyperbolic, was considered. Lobachevsky considered the negation of parallel postulate, while Bolyai considered exclusion of it. These ideas were already known to Carl Gauss. The adoption of these ideas has been hampered by common perception of Euclidean geometry as of "natural" or "true" geometry and all other as of "artificial" ones.

Finally, in 1854, Bernhard Riemann demonstrated [19] that Euclidean and hyperbolic geometries are in the same degree real by showing that Euclidean geometry is intrinsic geometry of manifolds with constant zero curvature, while hyperbolic geometry is intrinsic geometry of manifolds with constant negative curvature. The same approach describes also one more geometry, intrinsic for manifolds with constant positive curvature, which is now called elliptic geometry.

Arthur Cayley introduced different ways to define some space metric, Felix Klein in 1871, uses Cayley's metric and constructs models of elliptic, Euclidean, hyperbolic and another 6 plane geometries [197], Table 1.1. Klein proposed "Erlangen Program" [37], a new solution to the problem how to classify and characterize geometries on the basis of projective geometry and group theory. However, 6 of 9 constructed models Klein considered practically not applicable.

Table 1.1: Klein geometries

| Angles \ Lengths | Elliptic | Parabolic | Hyperbolic |
|---|---|---|---|
| Elliptic | Elliptic $\mathbb{S}^2$ | Euclidean $\mathbb{E}^2$ | Lobachevsky $\mathbb{H}^2$ |
| Parabolic | $\mathbb{G}^2_+$ | Galilean $\mathbb{G}^2$ | $\mathbb{G}^2_-$ |
| Hyperbolic | Anti de Sitter $\mathbb{M}^2_+$ | Minkowski $\mathbb{M}^2$ | De Sitter $\mathbb{M}^2_-$ |

In 1892, Hendrik Lorentz presented the transformations of space and time under which the Maxwell's equations are invariant. In 1905 Henri Poincaré observed that the Lorentz transformation is actually a rotation by hyperbolic angle and constructed the Lorentz group of transformations. The idea was used by Albert Einstein for special and then for general theory of relativity, proposed in 1905 [194] and further developed by Hermann Minkowski in 1907. Minkowski geometry found its place in Klein's classification. Also the classical kinematics of Galiley–Newton found its place in Klein's classification, if it is given the geometric interpretation [210]. When Willem de Sitter proposed the cosmologic model of early evolution of universe in 1917, it was observed that also spaces De Sitter and Anti de Sitter have their reserved places in Klein's classification [16]. The geometry proposed by Klein finally became not only practically applicable, but also "natural" in the most strict sense of the word.

Last achievements in theoretical physics (namely, in the String Theory, Green, Schwartz, 1970 — 1980's) [32] generated new homogeneous spaces which leads to necessity of developing new methods of study. The geometrization of physics plays an important role in this direction [87, 208]. The methods of linear algebra provide the universal instruments of research. It follows the study of different mathematical structures by linear methods represents an actual direction of research, which is important for domains of mathematics as well as for its applications.

### 1.5. The present state in the geometry of homogeneous spaces.

Perhaps the success of the differential geometric approach over the axiomatic one in demonstration of non–Euclidean geometry consistence was the key factor why differential geometry was designated as the standard geometric apparatus for study of non–Euclidean geometry. Shortly after Klein constructed models of different homogeneous spaces, the differential geometry of homogeneous spaces was elaborated around 1923 by Élie Cartan [127] and further developed by Charles Ehresmann similarly to what was previously done by Riemann for constructions made by Euclid, Lobachevsky and Bolyai. Riemann introduced the term "metric", Cartan introduced what is now known as "Cartan connection". Still then the differential geometry bacame the main tool to study homogeneous spaces. The differential geometry is beyond the scope of this thesis. The most fruitful domains where the core research is done today.

**Algebraic geometry.** Publications in this domain can be found by the keyword "homogeneous space". This field of research studies homogeneous manifolds rather than spaces, but they are often also named "homogeneous spaces". Due to Lie group — Lie algebra correspondence, it is possible to study geometric objects using algebraic language. The results in this domain answer the question "What does some homogeneous manifold look like?" and mainly enumerate and classify Lie groups with some given property. Because Lie groups can be defined in different ways, this problem is not trivial in most cases. Some recent works in this domain include [5, 6, 7, 8, 15, 20, 23, 26, 64, 69, 78, 80, 81, 82].

**Differential geometry.** Publications in this domain can be found by the keyword "Riemannian space" or "pseudo–Riemannian space". This field of research studies the metric properties of homogeneous spaces in the sense in which classical differential geometry studies manifolds with constant curvature. The results in this domain answer the question "What is the difference between studied space and its linear tangent space?" This domain seems to be more saturated than algebraic geometry. The new approach of the differential geometry related with homogeneous spaces was contribued by Buseman G., Bachmann F., Efimov N. V, Hjelmslev J., Nash L. F., Kallenberg G. W., Borisov Yu. V., Borisenko A. A., Milka A. D., Verner A. L., Schwartz J. T., Naoum A., Roitberg J., Lingenberg R., Karzel H., Struve H., Struve R. etc. Some selected works in this domain are [33, 53, 54, 83, 86, 93, 102, 111, 117, 120, 141, 155, 156]. Notably, there are many results in theory of curves and surfaces in (pseudo–) Riemanninan spaces [1, 12, 22, 50, 55, 65, 74, 89, 90, 101, 112, 113, 116, 119, 130, 135, 139, 142, 177].

The research domain of homogeneous spaces has connections with many compartments of contemporary geometry, from the most recent contributions we mention the monography of authors Bourguignon J. P., Hijazi O., Milhorat J. L., Moroianu A., Moroianu S. [9]. An essential contribution in differential geometry of spaces Finsler, Lagrange, Hamilton and their generalizations, which are closely related with homogeneous geometry, belongs to Miron R. and Anastasiei M. [47], Udriște C. and Balan V. [4].

Important results for some homogeneous spaces of different dimensions were obtained also in the domain of discrete geometry. We mention the works elaborated by Zamorzaev A. [121] on development of the theory of symmetry of homogeneous (Euclidean, Minkowski), results in domain of discrete hyperbolic geometry developed by Macarov V. [43] and in domain of hyperbolic manifolds obtained by Guțul I. [25] and Damian F. [131].

**Metric geometry.** Publications in this domain are difficult to find, because they usually describe geometry of some concrete homogeneous space and the number of spaces is very large.

Sometimes the overview publications can be found by keywords "Klein geometry" and "projective metric spaces". Often the same space is named differently by different authors. The research is done similarly to classical Euclidean geometry applied to geometry of homogeneous spaces. This field of research answers the question "How to compute things in some space?" The active contributrs in this domain in recent years are Hjelmslev J. (mainly in Danish and Dutch sometimes with German summary), Lingenberg R. [66, 40, 198, 199, 200, 201], Bachmann F. [3, 94] and more recently Karzel H. [34], Struve H. and Struve R. [70, 71, 72, 73, 202, 203, 204, 205]. Rosenfeld B., Yaglom I., Yasinskaya E. [165, 166, 167, 190, 192, 193] have contribued to the classification and development of the methods of study of two–dimensional homogeneous spaces. More recently, Pogorelov A. and Poznyak E. studied the rigidity of surfaces in spaces with different metrics [145, 146, 147, 148, 149, 150]. Artykbaev A. studied the metric and differential geometry of Galilean space [88, 89, 90, 91, 92, 93]. Sokolov D. studied the metric properties of manyfolds is pseudo–Euclidean space [172, 173, 174, 175, 176, 177, 178, 179, 180, 181]. Romakina L. described the geometry of several two–dimensional homogeneous spaces, which are seldom used [159, 160, 161, 162, 163, 164]. Several other recent publications are [13, 41].

This thesis research best matches the metric geometry domain. It is about how to compute distances, angles areas and volumes in different homogeneous spaces, what metric properties are dual to other properties, which properties are primary and which are secondary depending on the space. The author's results are [56, 57, 58, 59, 60, 61, 62, 63, 154].

**Applications in physics.** The homogeneous spaces were used by Einstein to elaborate special and then general theory of relativity shortly after their models were constructed by Klein. Different studies with applications in various domains of the physical mathematics were undertook by Hermann Weyl (in 1913), Élie Cartan, George Branceanu, Finsler, Radu Miron and others (1960 — 1970's years). Naturally, there are many publications in domain of homogeneous spaces application to physics, some of them are [10, 17, 27, 29, 42, 45, 47, 52, 68, 74, 85, 100, 114, 122, 137, 185].

### 1.6. The axiomatic method and the modelling methods.

In the thesis two approaches are used: axiomatic and modelling. This section presents these methods and analyzes their advantages and limitations.

**Axiomatic method.** This method operates with *notions* and *propositions.* Each proposition states some properties of the notions or some relations between the notions. New notions are introduced by *definitions* based on existing notions. New propositions are proved based on existing propositions. Such new propositions are named *theorems.* The theory constructed in axiomatic way tries to define as many notions as possible and to prove as many propositions as possible. However each definition and theorem proof is based on already existing notions and proposi-tions. It means that on previous step there are fewer existing notions or propositions. At some stage we reach the limit of provability of propositions and ability to define new notions based on existing. The notions that were not introduced by definitions are named *undefined* notions. The propisitions that are stated without proof are named *axioms*. The undefined notions and axioms are basis of some theory.

The axiomatic system is said to be:

- *Consistent* if it doesn't contain contradictions. That is, it is impossible to prove some statement and its negation based on the axioms.
- *Independent* if no axiom can be proved from other axioms.
- *Complete* if any statement about the undefined notions can be proved or disproved.

The independence of an axiom from some axiomatic system can be reduced to consistency of this system if the negation of the axiom in question is added to it. In 1931, Kurt Gödel demonstrated two incompleteness theorems. The first states that a consistent axiomatic system is never complete. The second, which immediately follows from the first, states that the consistency of the axiomatic system can't be proved within this system. In 1936, Gerhard Gentzen proved that consistency of the axiomatic system can be deduced by other means.

The construction of a "good" axiomatic system, that is consistent and independent, is usually a non-trivial task. The proof of equivalence between two axiomatic systems also pretends much work. All the undefined notions of one system should be either among undefined notions of another system or be definable in it. All axioms of one system should be among either axioms or theorems of another system. And this check should be done in both directions. The extension or generalization of some axiomatic system is also a non-trivial task. Take for example David Hilbert's axiomatic of 3-dimensional Euclidean geometry compared to 2-dimensional case.

**Modelling method.** This method consists of giving the *meaning* (*interpretation*) for notions of some axiomatic system in such a way that all axioms can be checked as true. With a model, the question is not the *provability* of the statements that has to be deduces from axioms, but their *truth* that has to be checked or computed.

In 1936, Alfred Tarski proved the undefinability theorem that states that the truth in a formal (axiomatic) system can't be defined in this system. So the property to be true is not the property of some axiomatic system construction, but the property of some its interpretation (model). This result was also known to Gödel. Tarski also proved that, unlike complete axiomatic systems, there exist complete models in sense there exist models where any proposition can be checked to be true or false. Finally, Tarski proved that the consequence can be proved from its premises if and only if in any model the consequence holds as soon as the premises are met.

The advantages of modelling method is that it is usually much easier to *check* some proposition to be true or false than to *prove* or *disprove* it from axioms. It means that it is usually much easier to check the equivalence of two models than to prove the equivalence of two axiomatic systems and also to extend or generalize some model than to extend or generalize some axiomatic system. Take for example Euclidean 3-dimensional analytic geometry modelled in linear vector space compared to 2-dimensional model.

Today the axiomatic system is considered consistent if *there exists at least one model of this system.* An interesting application of modelling method is proof of axioms independence. If in some model where the axiom in question holds it can be constructed another model where the negation of this axiom holds, leaving all the rest of axioms true, then this axiom is independent from others. Lobachevsky constructed a model of Euclidean space in hyperbolic space and Klein constructed a model of hyperbolic space in Euclidean space. So the parallel axiom was proded to be independent from others.

It should be acknowledged that a model is constructed in some already existing system, which is created, most probably, in axiomatic way. So the construction of a model constructs a *new interpretation* rather than a *new system.* It may seem that the axiomatic approach is more fundamental than the modelling one in that some axiom set defines a model, but not vice versa. This thought is disproved in Löwenheim–Skolem theorem. This theorem states that first–order theories are unable to control the cardinality of their infinite models, and that no first–order theory with an infinite model can have a unique model up to isomorphism. In other words, when someone builds an axiomatic system with a concrete model in mind, always there exist models of this axiomatic system, which are essentially different from the original one. So, the study of an axiomatic system *is distinct subject* from the study of some its model. Still, models are an important subject of research. By convention, the term *geometry* (e.g. Euclidean geometry, hyperbolic geometry) describes geometric research using the axiomatic method, while the term *space* (e.g. Euclidean space, hyperbolic space) describes research using the modelling method.

**Usage of axiomatic and modelling method.** This thesis uses primarly the modelling method of research. However, when the models gives an unusual result, this result is searched also among expected axioms of some homogeneous space. For example *parallel* axiom and its dual *connectability* axiom are explicitely stated as generalization of existing axioms and dual to existing axioms. This is done to ensure that the studied properties are really the properties of geometry of homogeneous spaces, not just of constructed model.

### 1.7. The main results of the thesis.

This thesis will develop an unified theory, different parts contributing to achieve this goal. However, it contains novel results that can be applied to adjacent research domains or present known facts in new light by their exdension or generalization.

The theory will introduce several new concepts, that naturally arise in homogeneous spaces research:

**Type** is a property of vectors, lines and relations between them and describes them: elliptic, parabolic or hyperbolic.

**Measure** is an unification and generalization of the notions of distance, different kinds of angles, areas and volumes equipping them with a type beside the value.

**Signature** is a set of numbers, each being some measure type, that uniquely define the kind of homogeneous spaces and their subspaces and thus can be used for classification of them.

**Unconnectable point pair** refers to points for that there is no straight line that contains them both (and thus violates first axiom of Hilbert axiomatic of Euclidean geometry).

**Unmeasurable angles** refers to the angles whose interior contains points both connectable and unconnectable with the vertex.

**Vector index** is the most natural position of some vector in the space basis that contains this vector.

**Equivalent, interchangeable and non–interchangeable vectors** refers to different "levels of vector congruence".

**Lineal** is linear span of some vector family. Its signature (and thus, also properties) may be very different from subspace properties.

**Limit vectors** is a renaming of *isotropic* vectors. I choose new name because the known properties of isotropic vectors in two–dimensional spaces are very different from their properties in higher dimensional spaces. Also, beside limit vectors there exist also limit lineals.

**Parameterized axioms.** The known Parallel Axiom of Euclidean plane and its variants for hyperbolic and elliptic planes will be generalized. Its dual form is Connectability Axiom. These two axioms have characteristic role in metric properties of homogeneous planes. They also show what properties should we expect from homogeneous space model. This result was published in [57].

**Universal model of homogeneous space.** The universal model of homogeneous spaces will be constructed in linear algebra language. It has two features:

1. The model is linear, even if describes non–linear spaces.

2. The model is parameterized and thus is universal.

The first fact makes it possible to apply linear algebra algorithms to homogeneous spaces with no or little change. The second fact makes it possible to see the whole picture among homogeneous spaces, to compare them and to study their properties in unified manner for all them. These results were published in [56, 59].

**Relations in triangles.** The known facts about relations in triangle will be generalized. The universal forms of triangle equations will be presented. Other triangle metric properties will be described: relation between longest / shortest side and largest / smallest angle as well as inequalities among triangle sides and angles. These relations are important in definition of metric (its triangle inequation axiom). A clear distinction between intrinsic and extrinsic properties of figures will be given, which shows that there are no extrinsic figures properties other than properties of space where the figure resides and which are the same for all figures. This result was published in [60].

**Motions.** Firstly the generalized orthogonal matrices will be introduced by means of orthonormalized vector family constituting the matrix columns. Then full theory of motion group of homogeneous spaces will be elaborated as theory of generalized orthogonal matrix group.

Starting with theorems on multiplication of types and on vector index invariance in space and lemmas on relations among generalized orthogonal matrix rows and on the form of generalized matrix when some element of space signature $k_i = 0$, the structure of generalized orthogonal matrices will be established. These matrices naturally generalize the structure of orthogonal matrices.

Will be proved proposition on main rotation matrices and lemmas on product of generalized orthogonal matrices and on inverse generalized orthogonal matrix as well as theorem stating that generalized orthogonal matrices form the group, that make the link between generalized orthogonal matrix group and homogeneous space motion group. Each homogeneous space motion group or some its subgroup (continuos or discrete) can be described as corresponding generalized orthogonal matrix group / subgroup.

It will be proposed algorithm that describes the decomposition of generalized matrix in product of rotation matrices. It will be demonstrated that given a proper motion $\mathfrak{M}$ it is possible to smoothly parameterize it $\mathfrak{M}(p)$, $p \in [0, 1]$ so that $\mathfrak{M}(0) = \mathfrak{I}$ is trivial motion and $\mathfrak{M}(1) = \mathfrak{M}$ is initial motion. Theorem on coordinate vectors equivalence will establish the relation between type of rotation in plane defined by two coordinate vectors and their equivalence (equivalent, interchangeable or non–interchangeable). Some of these results were published in [63].

**Algorithms.** Thanks to linear nature of elaborated model, it is possible to adapt known algorithms of linear algebra to non–linear homogeneous spaces. These include:

- Vector family orthonormalization,
- Orthonormal vector family completion,
- Canonical form of lineal basis,
- Lineal difference,
- Lineal sum and intersection

toghether with propositions on segment midpoint and on triangle center of gravity and lemma on projection of vector to lineal and to its orthogonal complement. Some of these results were published in [58].

**Limit vectors and lineals.** Theory of limit (isotropic) vectors will be developed. Their study in spaces of higher dimensions is very different from isotropic vectors study in planes in that this thesis does not declare the measure of limit vectors to be zero. Instead their measure will be defined by parameter of motion along them: theorem on type of limit vectors and lemma on measure of limit vectors. It will be proved that this measure is trivial if and only if the motion along the limit vector is trivial. As development of this idea, this thesis will not present limit vectors orthogonal to themselves. Instead new way of more precise orthogonalization is presented. This will be done by introduction of *decomposition* vector pair for each limit vector. It will be stated that indices of decomposition vectors do not depend on space basis choice. Together with lemma that states that if lineal basis contain limit vector, then its decomposition vetors indices are free among indices of this basis, we can declare that limit vectors have double index (indices of decomposition vectors) that is invariant with respect to the space basis choice. This proposition answers the question why isotropic vectors theory on plane is so different from limit vectors theory in space — because their double index takes all available space indices.

The study will continue with lemma which states that limit vectors can only be equivalent or non–interchangeable that leads to algorithm of finding of limit lineal signature. Some results in this area were published in [154].

**Measure between lineals.** Measure between lineals will be defined and computed in the most general case. The measure computing has to provide the information about measure type beside measure value. It will be stated that if the measure value between two lineals equals to zero or they are orthogonal, its type is ambiguous. An algorithm will present measure computation between lineals in the most general case.

**Volume.** The theory of volumes will begin with space volume definition by corresponding volume in metaspace. This definition has the advantage that metaspace itself and the integration bounds are linear except space sphere which is quadric. In these conditions usually it is easier to compose and to compute the integral, at least numerically.

The right triangle area universal equations will be presented. These equations are interesting in that they fulfill statement of proposition on properties of figures and properties of spaces, which leads to the conclusion that the area measure is the same sort of measure as lengths and angles equipped with some type. There will be proved theorem that establishes the area type and the proposition that establishes sufficient and necessary conditions for parabolic volume type and proposed a conjecture about possible type of volume.

**Applications in geometry.** The link between homogeneous space metric and the metric of pseudo–Euclidean, semi–Euclidean, pseudo–Riemannian and semi–Riemannian spaces will be established. It will be shown that the metaspace of some homogeneous space is linear homogeneous space with signature that results from space signature by adding zero element in front, while linear tangent space of some homogeneous space is also homogeneous with the signature where the first element is replaced by zero.

During deduction of examples of crystallographic groups on homogeneous planes, it will be proved lemma on isomorphism between crystallographic groups of dual spaces. This result is to be published.

It will be shown that, generally speaking, the geodesic path in homogeneous space can't be computed based on its extremal length and so there is improper to define geodesic in this way.

**Applications in topology.** The analytic theory of homogeneous spaces leads to some interesting observations in topology. Based exclusively on points connectability notion, which naturally arise here, it is possible to refine the topological notion of *points separability on a line.* As known elliptic line points do not separate the line into two half–lines, while Euclidean and hyperbolic line points do. The notions of *weak separable* and *strong separable* line points will be introduced in this thesis. They make the difference between Euclidean and hyperbolic lines. This result was published in [61].

The *neighborhood* notion will also be generalized with aim of the principle of duality. The generalized neighborhood notion leads to *generalized Hausdorff* spaces notion. While most homogeneous spaces are not classical Housdorff, all them are generalized Hausdorff and this result introduces study of *homogeneous manifolds*, several of them will be presented in this thesis.

**Applications in physics.** Some insights about the structure and the shape of physical space–time and its metric properties will be proposed. These show the theory can be used in physics, including the study of *String Theory* and *M-Theory*.

Hypotetical description of *"non–Euclidean" geometrical optics* will be given. This description is in no way less real then "imaginary geometry" described by Lobachevsky, given corresponding interpretation of described notions.

### 1.8. Domains to which this thesis makes a contribution.

This section describes some domains where this thesis has direct implication. Some of these domains can be improved by this research contribution, others are uncovered by existing today research.

**Elementary geometry.** The Erlangen Program [37] provided the way to classify and fit into hierarchy different geometries by offering the same language of projective geometry and by comparison their symmetry (motion) groups. However it didn't provide "horizontal" classification of geometries situated at the same level in this hierarchy. The results of Isaak Yaglom [190] provide such a classification by means of the generalization of complex numbers, however this formalism limits its applicability to two–dimensional case. This thesis classifies metric geometry of rigid homogeneous spaces of any dimension.

The non-Euclidean geometry is seldom used even when it is the most natural tool to use. For example, long time before elliptic geometry was described, the clear understanding of earth form already existed and the navigation systems of large scale would benefit from elliptic geometry. Still, calculus based on Euclidean geometry were used. The same is the present day reality in matter of space orbit calculus in earth neighborhood. The homogeneous spaces are used about exclusively for needs of Galilean and Einsteinian kinematics. That is why among the spaces that have personal name, more are named after physicists than after mathematicians. Also, many geometric terms have physics origin. The geometries of homogeneous spaces, other than spaces of constant curvature or those used for physics purpose, were never described. This thesis constructs the theory capable of description of each homogeneous space geometry.

**Differential geometry.** The classical differential geometry studied metric properties of non–linear figures in Euclidean space. Because differential geometry is able to make the difference between intrinsic and extrinsic metric properties, when non–Euclidean spaces emerged the existing apparatus of differential geometry was adapted to study of non–Euclidean spaces and their metric properties. It is interesting to recall that Riemann in his research that demonstrated non–Eiclidean geometry consistency used different metrics introduced in Euclidean ball. So the status of Euclidean space changed from single possible one to the main space among all possible spaces.

It should not be forgotten that the tools of differential geometry were developed with Euclidean geometry in mind. Naturally, the basic metric properties are studied with axiomatic or analytic geometry tools. Then the differential geometry is constructed as natural extension of already existing methods. This is the case of differential geometry of Euclidean and non–linear locally Euclidean spaces. However, the existing tools have limited applicability to spaces that can be approximated by Euclidean in no point. The analytic geometry of homogeneous spaces was never developed to become the basis for more universal differential geometry. Instead, considerig existing apparatus of differential geometry as universal tool, geometers try to "adjust" spaces to existing "universal" methods.

As result, some spaces are described in inconsequent and misleading manner, which hides their natural but unexpected properties on the one hand, and tends to suggest to search the properties, these spaces lack off, on the other hand. Not being able to find expected properties, the space is marked as "flawed". Even when some unexpected properties are found, they are often perceived not as space advantages, but as its disadvantages. Even worse is the situation of spaces that doesn't fit in modern differential geometry instrumentarium. Often they are simply declared as having "no right to be called metric space". In this way, instead of studying the unique properties of spaces, often it is limited to study of differences between them and Euclidean space. It is important to not exist such situation, because no science in general and geometry in particular answers the questions nobody asks.

This approach can be observed even in terminology. On elliptic and hyperbolic planes, trianges have "defect", expressing the difference between triangle angles sum and "expected" Euclidean sum of $\pi$, despite the fact that thanks to the "defect", in these spaces edges lengths are fully determined by angles values, unlike Euclidean case. These spaces are "flawed", because they lack length dilation transformation. On the other hand, in semi–Euclidean spaces the number of coordinate vectors, not included in distance bilinear form, is also called "defect", despite the fact that thakns to it exists new transformation of angle dilation. This time "flaw" consists in inability to determine the angles values based on lengths, unlike Euclidean geometry.

This thesis proposes to cover the lack of analytic geometry tools for homogeneous spaces sutable for differential geometry.

### 1.9. Conclusions of chapter 1

It follows from the description above that the theory of homogeneous spaces raised from the solution of the problem of the independency of the Euclid's parallel axiom, followed by the creation of geometries by Klein in XIX century. The Erlangen Program was aimed to elaborate the overview and the classification of the geometric spaces. However this objective was not achieved for rigid geometries of dimension greater than 2. With exception of spaces of constant curvature and spaces needed in physics, no such geometry was even described. That is why the following problem is actual: **to investigate the homogeneous spaces via linear methods applying the concept of signature**. In order to realize this problem, the following objectives are proposed:

1. Introduction of the new concept of the space signature, that depends on space dimension;
2. Construction of homogeneous space based on signature concept;
3. Construction of the model of homogeneous space for each signature;
4. Expression of the measurement of different geometric quantities via signature, which reflects their role in analytic geometry;
5. Different applications of the analytic geometry of homogeneous spaces.

# 2. ANALYTIC GEOMETRY

## 2.1. Types of lines, distances and angles

This chapter analyzes what metric properties are to be expected in homogeneous spaces and constructs the parameterized model of homogeneous spaces. With aim of the model, different elements of analytic geometry of homogeneous spaces are studied.

### 2.1.1. Definition and type of generalized rotations.

Consider real plane $\mathbb{R}^2$ and the transformation in it $\mathfrak{R}$, defined by the matrix:

$$R = \begin{pmatrix} a & -kb \\ b & a \end{pmatrix},$$

$a, b \in \mathbb{R}, k \in \{-1, 0, 1\}$, so that $\det R = a^2 + kb^2 = 1$. Under these conditions, for each fixed value of $k$, the transformation set $\mathfrak{R}$ depends on one parameter, and the parameter $\varphi$ is best chosen in such a way that, for $k = 1$:

$$R_1(\varphi) = \begin{pmatrix} \pm\cos\varphi & -\sin\varphi \\ \sin\varphi & \pm\cos\varphi \end{pmatrix},$$

for $k = 0$:

$$R_0(\varphi) = \begin{pmatrix} \pm 1 & 0 \\ \varphi & \pm 1 \end{pmatrix},$$

and for $k = -1$:

$$R_{-1}(\varphi) = \begin{pmatrix} \pm\cosh\varphi & \sinh\varphi \\ \sinh\varphi & \pm\cosh\varphi \end{pmatrix}.$$

*Remark.* For now we assume that $R(0) = I$ identity matrix, leaving without attention negative sign on main diagonal. We will have a closer look at this case when will analyze the generalized orthogonal matrices.

It is easy to see that $\mathfrak{R}_1$ is the rotation, $\mathfrak{R}_0$ is the Galilean transformation and $\mathfrak{R}_{-1}$ is the Lorentz transformation.

*Properties.* Transformations $\mathfrak{R}_1$, $\mathfrak{R}_0$ and $\mathfrak{R}_{-1}$ have much in common.

- All them have a fixed point, the origin $O = (0, 0)$;

- The rotations mutually commute:

$$\mathfrak{R}(\varphi)\mathfrak{R}(\psi) = \mathfrak{R}(\varphi + \psi) = \mathfrak{R}(\psi)\mathfrak{R}(\varphi); \tag{2.1}$$

- The rotations are invertible:

$$\mathfrak{R}^{-1}(\varphi) = \mathfrak{R}(-\varphi); \tag{2.2}$$

- The orbit of the point $P = (1, 0)$ satisfies the equation (Figure 2.1):

$$x_1^2 + k x_2^2 = 1, \quad k = \{-1, 0, 1\}.$$

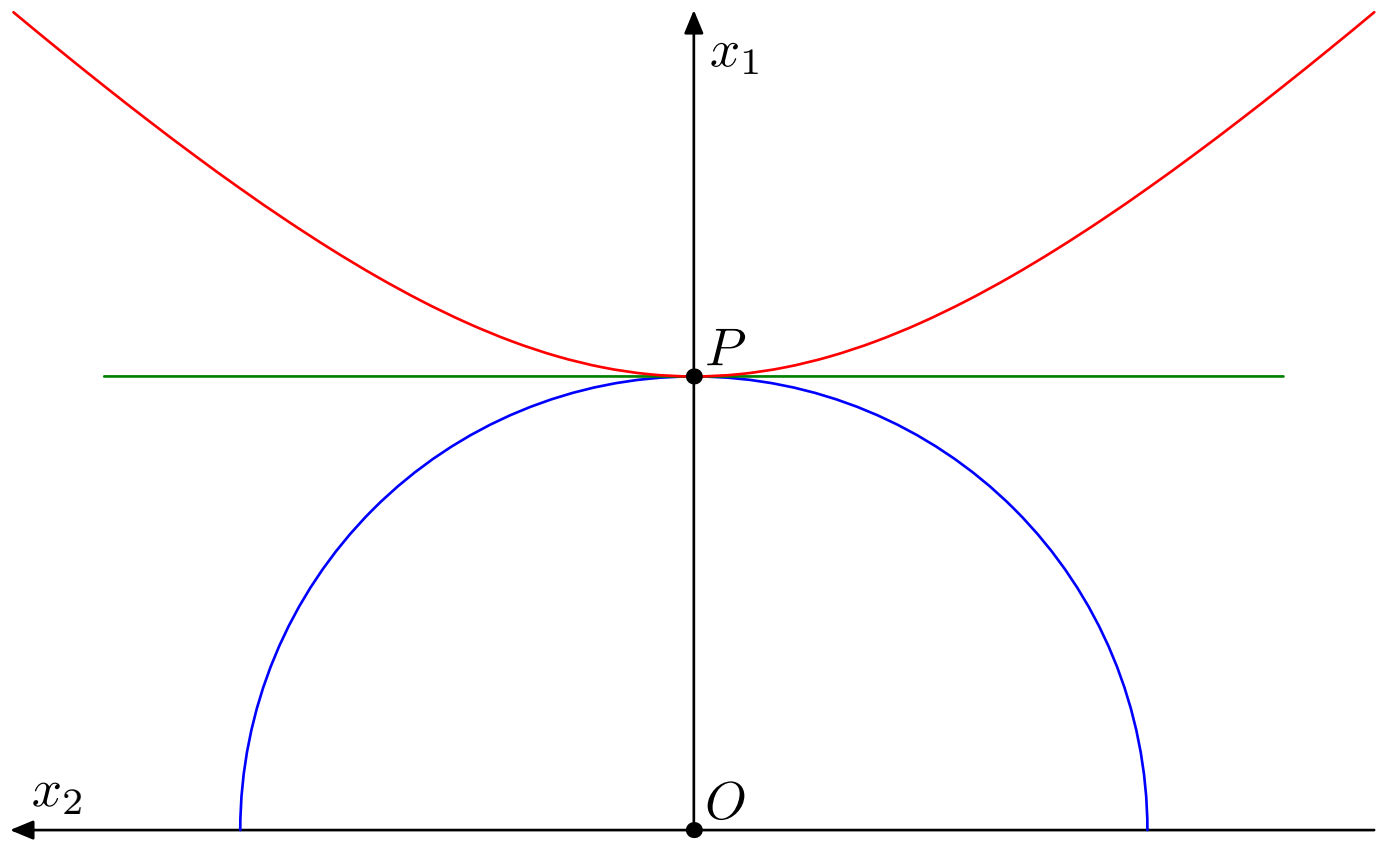

Figure 2.1: The orbit of point $P = (1, 0)$ under the transformations $\mathfrak{R}_1$ (blue), $\mathfrak{R}_0$ (green) and $\mathfrak{R}_{-1}$ (red).

The orbit of point $P$ is circle if $k = 1$ (Figure 2.1, blue line), straight line if $k = 0$ (Figure 2.1, green line) and branch of hyperbola if $k = -1$ (Figure 2.1, red line).

**Definition 2.1.1** (Generalized rotation, the type of a rotation)**.** We call $\mathfrak{R}_1$ *elliptic rotation*, $\mathfrak{R}_0$ *parabolic rotation* and $\mathfrak{R}_{-1}$ *hyperbolic rotation*. Generally, we call transformations $\mathfrak{R}$ *generalized rotations*. We call the parameter $\varphi$ *the angle* of rotation. We call the coefficient $k$ *the type* of rotation, which can be *elliptic* ($k = 1$), *parabolic* ($k = 0$), or *hyperbolic* ($k = -1$).

*Remark.* Generalized rotation, beside angle, has one more parameter, its type.

### 2.1.2. Generalized trigonometric functions.

Let define functions $C(\varphi)$ and $S(\varphi)$ as (for consistency of equation form, assume

$$C(\varphi) = C(\varphi, k) = \sum_{n=0}^{\infty} (-k)^n \frac{\varphi^{2n}}{(2n)!}, \tag{2.3}$$

$$S(\varphi) = S(\varphi, k) = \sum_{n=0}^{\infty} (-k)^n \frac{\varphi^{2n+1}}{(2n+1)!}. \tag{2.4}$$

One can observe that:

$$C(\varphi) = \begin{cases} \cos\varphi, & k = 1; \\ 1, & k = 0; \\ \cosh\varphi, & k = -1. \end{cases}$$

$$S(\varphi) = \begin{cases} \sin\varphi, & k = 1; \\ \varphi, & k = 0; \\ \sinh\varphi, & k = -1. \end{cases}$$

In these notations, the matrix $R(\varphi)$ can be written generally (we chosen positive elements on main diagonal):

$$R(\varphi) = \begin{pmatrix} C(\varphi) & -kS(\varphi) \\ S(\varphi) & C(\varphi) \end{pmatrix} \tag{2.5}$$

Define one more function:

$$T(\varphi) = \frac{S(\varphi)}{C(\varphi)}. \tag{2.6}$$

*Properties.* It is easy to observe that regardless of values $\varphi \in \mathbb{R}, k \in \{-1, 0, 1\}$, the following equalities verify:

$$C^2(\varphi) + kS^2(\varphi) = 1, \tag{2.7}$$

$$C(\varphi \pm \psi) = C(\varphi)C(\psi) \mp kS(\varphi)S(\psi), \tag{2.8}$$

$$S(\varphi \pm \psi) = S(\varphi)C(\psi) \pm C(\varphi)S(\psi), \tag{2.9}$$

**Definition 2.1.2** (Generalized trigonometric functions)**.** We call functions $C(\varphi)$, $S(\varphi)$ and $T(\varphi)$ *generalized cosine*, *sine* and *tangent* respective. All them we call *generalized trigonometric functions.*

### 2.1.3. Translation defined as metarotation. Its type.

Generalized rotation has a fixed point. But translation may have none. We define translation in terms of generalized rotation using additional dimension.

For $n$-dimensional space, consider vector space $\mathbb{R}^{n+1}$. The first coordinate (we will count it 0-th) will be additional. The rotation matrices have two columns different from columns of the identity matrix. When one of them has index 0, we will consider them translations. Otherwise, we will consider them rotations.

Let $n = 1$. We call origin $E$ the vector $e = (1, 0)$. Then, if vector $a = \mathfrak{R}(\varphi)e$, than we call distance $EA$ the angle $\varphi$ between $e$ and $a$.

*Remark.* Different values of $k$ correspond to different kinds of translations. For each translation kind exists corresponding kind of straight line and kind of distance measure: elliptic, parabolic and hyperbolic. So, the type $k$ is a property of:

- translation (and generally, motion),
- distance (and generally, measure),
- straight line (actually, of a vector). For description of higher dimensional objects, more parameters $k$ are used.

The difference between straight lines of different types can be observed from variants of axiom of parallels for elliptic, Euclidean and hyperbolic geometries (Figure 2.2).

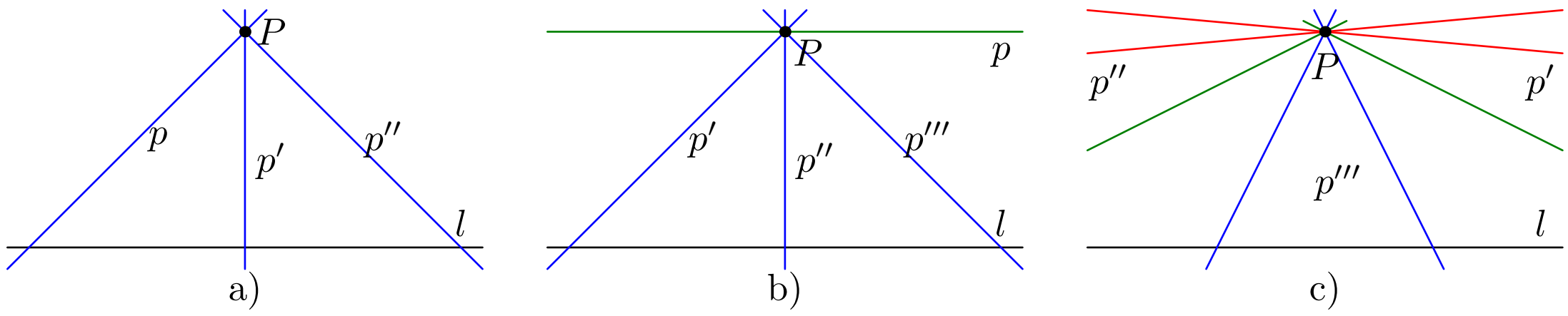

Figure 2.2: Variants of parallel axiom: elliptic a), parabolic b) and hyperbolic c).

**Elliptic axiom:** Figure 2.2 a. For the given line $l$ and the point $P$, not in $l$, there is no line $p$ that pass through $P$ without intersection with $l$. It may be formulated also as: For the given line $l$ and the point $P$, not in $l$, all lines that pass through $P$ intersect $l$.

**Euclidean axiom:** Figure 2.2 b. For the given line $l$ and the point $P$, not in $l$, there is a unique line $p$, that passes through $P$ and doesn't intersect $l$.

**Hyperbolic axiom:** Figure 2.2 c. For the given line $l$ and the point $P$, not in $l$, there are at least two lines $p', p''$ that pass through $P$ and don't intersect $l$.

**Proposition 2.1.1.** *Generally, the axiom of parallels can be formulated as: For the given line $l$ and the point $P$, not in $l$, there are $0^k$ (where $k$ is line type) lines $p$, that pass through $P$ and don't intersect $l$.*

*Remark.* It should be noticed, that $0^k$ is a symbol, not a number used for calculations. It equals to 0 when $k = 1$, to 1 when $k = 0$ and to $\infty$ when $k = -1$.

### 2.1.4. Sequence of unconnectable points.

It is easy to see that rotations in Euclidean, elliptic (Riemann) and hyperbolic (Bolyai–Lobachevsky) geometries, have the type $k = 1$. We can extend the notion of rotation in space to generalized rotation of some type. The best way to explain the difference between different kinds of rotations (as well as between different kinds of angles) is to formulate angular equivalent of parallel axiom — the axiom of points connectability (Figure 2.3). For this, we exchange the following phrases:

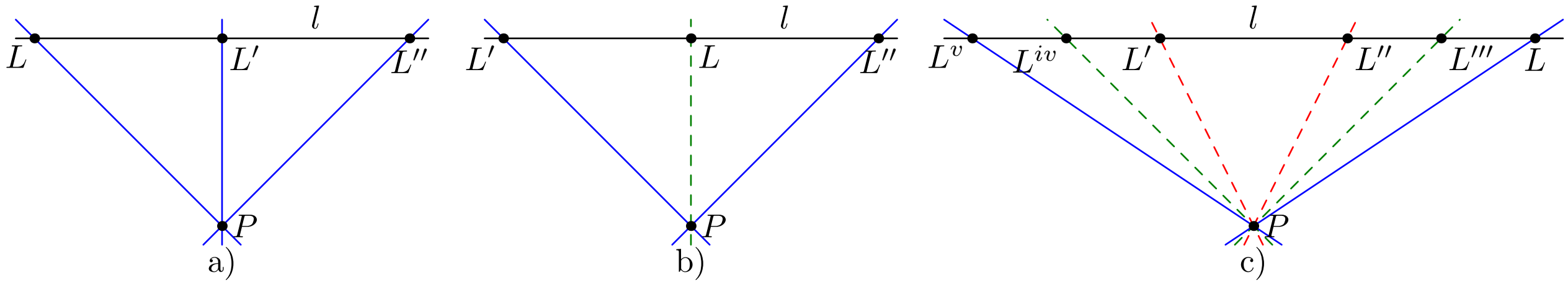

Figure 2.3: Variants of points connectability axiom: elliptic a), parabolic b) and hyperbolic c).

$$\text{Line } l \longleftrightarrow \text{Point } L,$$
$$P \text{ lies in } l \longleftrightarrow p \text{ passes through } L,$$
$$P \text{ doesn't lie in } l \longleftrightarrow p \text{ doesn't pass through } L,$$
$$\text{Distance } AB = \varphi \longleftrightarrow \text{Angle } \measuredangle ab = \varphi,$$
$$C \text{ is intersection point of } a \text{ and } b \longleftrightarrow \text{Line } c \text{ is defined by points } A \text{ and } B,$$
$$\text{Lines } a \text{ and } b \text{ do not intersect } \longleftrightarrow \text{Points } A \text{ and } B \text{ are unconnectable.}$$

*Remark.* The statements above are intentionally written without symbols $P \in l$, $P \notin l$ that make us think of points as of set elements and think of lines as of sets that contain elements. It is done so in order to show complete duality of points and lines as well as of lenghts and angles, opposite to concept that lines consist of points, or that points are elements of lines. For example, two intersected lines define a point (of intersection) as well as two points define a line (that passes through them).

The last phrase, points $A$ and $B$ are unconnectable, is unusual for three aforementioned geometries. It clearly conflicts with Euclidean axiom that through any two points goes a line. This axiom should be changed by one of the following axioms of connectality in order to consider geometries with rotations different from elliptic. The phrase has sense in geometries with angular type 0 or $-1$. The point unconnectability property is dual to the line parallelism property.

**Elliptic type:** Figure 2.3 a. On the line $l$ that doesn't pass through the point $P$, there is no points unconnectable with $P$.

**Parabolic type:** Figure 2.3 b. On the line $l$ that doesn't pass through the point $P$, there is a single point $L$ unconnectable with $P$.

**Hyperbolic type:** Figure 2.3 c. On the line $l$ that doesn't pass through the point $P$, there are at least two points $L'$ and $L''$ unconnectable with $P$.

**Proposition 2.1.2.** *Generally, this axiom can be formulated in this way: On the line $l$ that doesn't pass through the point $P$, there are $0^k$ (where $k$ is angular type) points unconnectable with $P$.*

*Remark.* As in case of parallel axiom, symbol $0^k$ isn't used in calculus.

Similarly to pencils of lines, intersected, parallel or divergent, sequences of points can be defined.

**Definition 2.1.3** (Sequence of points.). Let $X, Y \in \mathbb{R}^{n+1}$. All linear combinations $Z = \alpha X + \beta Y, \alpha, \beta \in \mathbb{R}$ compose a set we call *sequence of points.*

*Remark.* As we will see later, this set has one restriction. So, it has one free parameter.

Just as any two lines define a pencil of lines, any two points ($X$ and $Y$) define a sequence of points. In case $X$ is connectable with $Y$, the sequence is a line (dual to point of intersection in pencil of intersected lines). Lines are shown blue in figures 2.3, 2.4. If $X$ is unconnectable with $Y$, the sequence isn't a line (dual to pencil of parallel or divergent lines). Sequences of unconnectable points are shown green and red in figures 2.3, 2.4.

If we fix a point from a sequence, $X$, and arbitrarily choose another point to define a sequence, $Y$, obtain a pencil of sequences of points (or pencil of lines) defined by point $X$ (Figure 2.4). For any angular type exist unlimited number of sequences of connectable points (Figure 2.4 a, b, c). By connectable axiom (proposition 2.1.2), if angular type is 1, all sequences of points are lines (Figure 2.4 a, blue). If angular type is 0, for any point exists single sequence of unconnectable points (Figure 2.4 b, green). If angular type is $-1$, there are infinitely many sequences of unconnectable points (Figure 2.4 c, red). In this case sequences of points fall in one of two categories, sequences of connectable and unconnectable points. Limit (marginal) sequences (green) of unconnectable points can be viewed as the third category (similarly to difference between parallel and divergent lines). There are two limit sequences of points on plane.

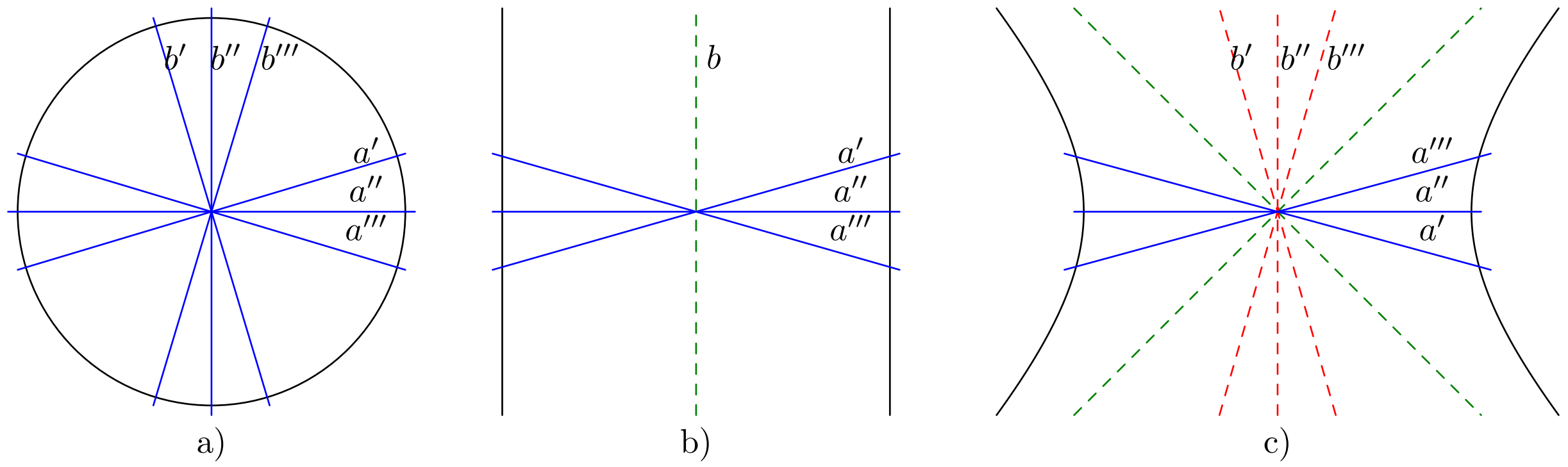

Figure 2.4: Mutual position of different sequences of points and circles: a) elliptic angular type, b) parabolic angular type, c) hyperbolic angular type.

Pencils of connectable points intersect circles with the center in pencil center, no pencil of unconnectable points intersect these circles and limit pencils are asymptotic to circles (Figure 2.4).

If angular type is 1, orthogonal sequence of points defines a line. In this case, when a line rotates counterclockwise, its perpendicular rotates counterclockwise and vice versa (Figure 2.4 a). If angular type is 0, there is only one fixed orthogonal sequence of points (Figure 2.4 b). If angular type is $-1$, orthogonal sequence of points rotates clockwise when line rotates counterclockwise toward the same limit sequence of points and vice versa (Figure 2.4 c).

### 2.2. Homogeneous Space Model

Consider the vector space $\mathbb{R}^{n+1}$ with the collection of numbers $k_1, k_2, ...k_n \in \{-1, 0, 1\}$. Denote $C_i(x) = C(x, k_i)$, $S_i(x) = S(x, k_i)$ and $T_i(x) = \frac{S_i(x)}{C_i(x)}$. Let

$$R_1(\varphi) = \begin{pmatrix} C_1(\varphi) & -k_1 S_1(\varphi) & 0 & \dots & 0 \\ S_1(\varphi) & C_1(\varphi) & 0 & \dots & 0 \\ 0 & 0 & 1 & \dots & 0 \\ \vdots & \vdots & \vdots & \ddots & \vdots \\ 0 & 0 & 0 & \dots & 1 \end{pmatrix},$$

$$R_2(\varphi) = \begin{pmatrix} 1 & 0 & 0 & \dots & 0 \\ 0 & C_2(\varphi) & -k_2 S_2(\varphi) & \dots & 0 \\ 0 & S_2(\varphi) & C_2(\varphi) & \dots & 0 \\ \vdots & \vdots & \vdots & \ddots & \vdots \\ 0 & 0 & 0 & \dots & 1 \end{pmatrix},$$

$\vdots$

$$R_n(\varphi) = \begin{pmatrix} 1 & 0 & \dots & 0 & 0 \\ 0 & 1 & \dots & 0 & 0 \\ \vdots & \vdots & \ddots & \vdots & \vdots \\ 0 & 0 & \dots & C_n(\varphi) & -k_n S_n(\varphi) \\ 0 & 0 & \dots & S_n(\varphi) & C_n(\varphi) \end{pmatrix}.$$

**Definition 2.2.1** (Main space rotations). We call respective motions $\mathfrak{R}_1, ...\mathfrak{R}_n$ *main space rotations.*

### 2.2.1. Meta product of vectors. Invariant biliniar form.

Let:

$$K_0 \equiv 1, \quad K_m = \prod_{i=1}^{m} k_i, \quad m = \overline{1, n} \tag{2.10}$$

*Properties.* It can be observed that $K_m \in \{-1, 0, 1\}$, $\forall m = \overline{0, n}$, like $k_i$, $\forall i = \overline{1, n}$.

**Definition 2.2.2** (Meta product of vectors). We define *meta product of vectors* $\odot$ as

$$x \odot y = \sum_{i=0}^{n} K_i x_i y_i \tag{2.11}$$

**Lemma 2.2.1.** *The meta product of vectors (2.11) is invariant with respect to main space rotations* $\mathfrak{R}_m(\varphi), \forall x, y \in \mathbb{R}^{n+1}, \varphi \in \mathbb{R}$ *and elements* $k_i \in \{-1, 0, 1\}; i, m = \overline{1, n}$*:*

$$x \odot y = \mathfrak{R}_m(\varphi)x \odot \mathfrak{R}_m(\varphi)y. \tag{2.12}$$

*Proof.* Consider vectors $x, y \in \mathbb{R}^{n+1}$:

$$x = (x_0, x_1, ..., x_n),$$
$$y = (y_0, y_1, ..., y_n).$$

Let $x' = \mathfrak{R}_m(\varphi)x$ and $y' = \mathfrak{R}_m(\varphi)y$ $(m = \overline{1,n})$:

$$x' = \{x_0, ..., x_{m-2}, x_{m-1}C_m(\varphi) - k_m x_m S_m(\varphi), x_{m-1}S_m(\varphi) + x_m C_m(\varphi), x_{m+1}, ..., x_n\},$$
$$y' = \{y_0, ..., y_{m-2}, y_{m-1}C_m(\varphi) - k_m y_m S_m(\varphi), y_{m-1}S_m(\varphi) + y_m C_m(\varphi), y_{m+1}, ..., y_n\}$$

Then:

$$\begin{aligned}
x' \odot y' &= \sum_{i=0}^{n} K_i x'_i y'_i = \sum_{i=0}^{m-2} K_i x_i y_i \\
&+ ((x_{m-1}C_m(\varphi) - k_m x_m S_m(\varphi))(y_{m-1}C_m(\varphi) - k_m y_m S_m(\varphi)) \\
&+ k_m(x_{m-1}S_m(\varphi) + x_m C_m(\varphi))(y_{m-1}S_m(\varphi) + y_m C_m(\varphi)))K_{m-1} \\
&+ \sum_{i=m+1}^{n} K_i x_i y_i = \sum_{i=0}^{m-2} K_i x_i y_i \\
&+ (x_{m-1}y_{m-1}C_m^2(\varphi) - k_m(x_{m-1}y_m + x_m y_{m-1})S_m(\varphi)C_m(\varphi) + k_m^2 x_m y_m S_m^2(\varphi) \\
&+ k_m x_{m-1}y_{m-1}S_m^2(\varphi) + k_m(x_{m-1}y_m + x_m y_{m-1})S_m(\varphi)C_m(\varphi) + k_m x_m y_m C_m^2(\varphi))K_{m-1} \\
&+ \sum_{i=m+1}^{n} K_i x_i y_i = \sum_{i=0}^{m-2} K_i x_i y_i \\
&+ (x_{m-1}y_{m-1}(C_m^2(\varphi) + k_m S_m^2(\varphi)) + k_m x_m y_m(C_m^2(\varphi) + k_m S_m^2(\varphi)))K_{m-1} \\
&+ \sum_{i=m+1}^{n} x_i y_i K_i = \sum_{i=0}^{n} K_i x_i y_i = x \odot y
\end{aligned}$$

We used equality (2.7) twice here. It is true for all $m = \overline{1,n}$. Thus, bilinear form $x \odot y$ is invariant with respect to all main rotations $\mathfrak{R}_m$. □

*Properties.* It is easy to see that meta product of vectors (2.11) is symmetric and linear:

$$x \odot y = y \odot x, \tag{2.13}$$
$$(x + y) \odot z = x \odot z + y \odot z, \tag{2.14}$$
$$(\lambda x) \odot y = \lambda(x \odot y) \tag{2.15}$$

*Remark.* Evidently, meta product of vectors (2.11) generalizes dot product of vectors in Euclidean space (for which all $K_i = 1, i = \overline{0,n}$). Yet, product (2.11) will be more generalized further. Such a form is called zero form (or point form, as well).

#### 2.2.2. Space definition by signature.

Consider the vector space $\mathbb{R}^{n+1}$ and let $k_i \in \{-1, 0, 1\}$, $i = \overline{1,n}$. Construct the indicatrix of the $\odot$ product, the $n$-dimensional "unit sphere" $\boldsymbol{B}^n = \{x \in \mathbb{R}^{n+1} \mid x \odot x = 1\}$ (Figure 2.5). Because all main rotations preserve bilinear form defined by the product $\odot$, the sphere is also invariant with respect to these rotations.

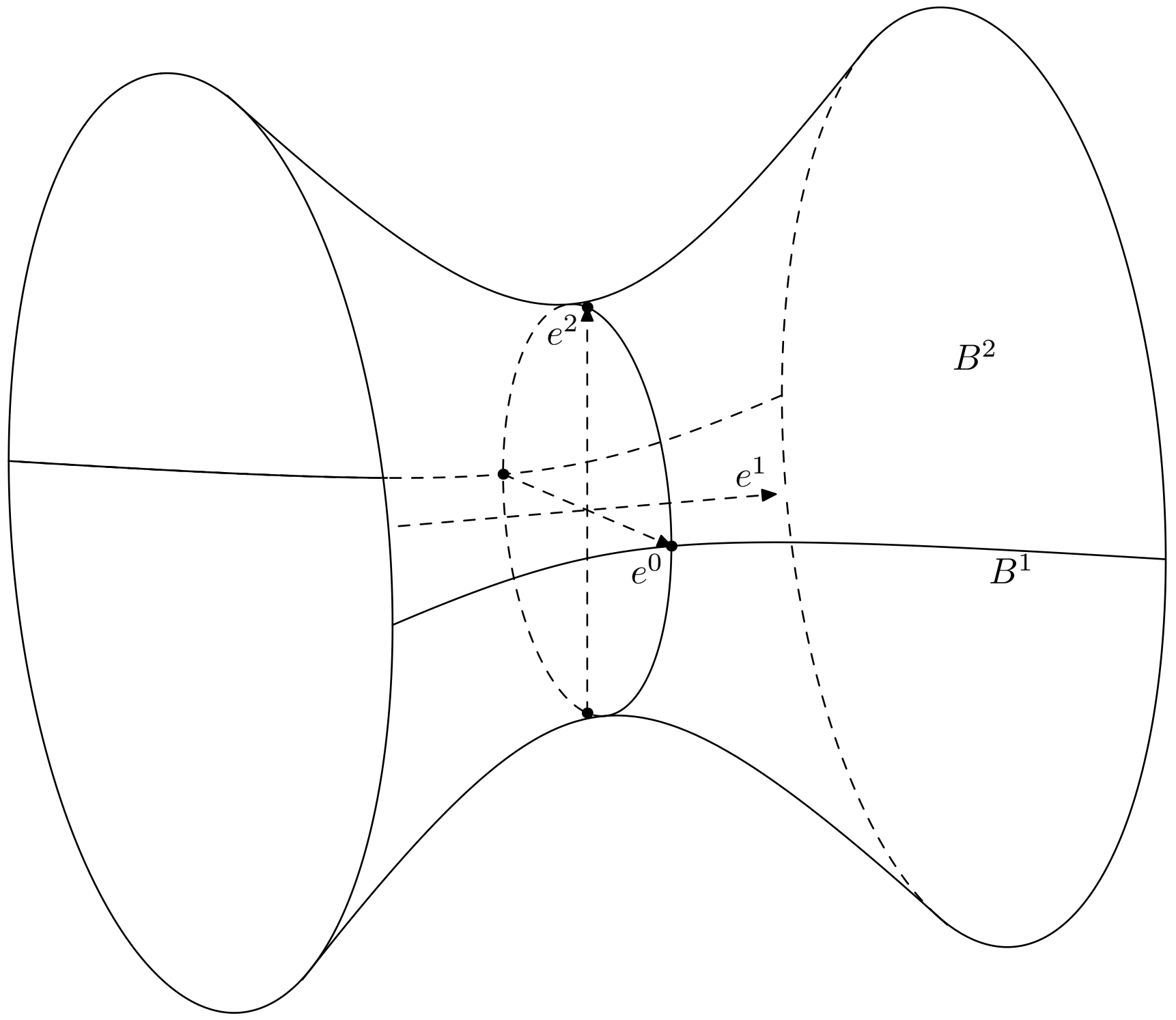

Figure 2.5: The unit sphere of space $\{-1, -1\}$.

*Properties.* Let us examine some properties of the sphere $\boldsymbol{B}^n$:

- In no case origin of $\mathbb{R}^{n+1}$, the vector $o = (0, ..., 0)$ belongs to $\boldsymbol{B}^n$, since $o \odot o = 0 \neq 1$.
- For any $k_i, i = \overline{1,n}$, the vector $e = (1, 0, ..., 0) \in \boldsymbol{B}^n$. Since always $\boldsymbol{K}_0 = 1$, obviously the coordinates of $e$ satisfy the equation of $\boldsymbol{B}^n$.
- If some vector $x \in \boldsymbol{B}^m$, then also $-x \in \boldsymbol{B}^n$.

*Remark.* Generally speaking, bilinear form $\odot$ (2.11) may not be positive definite, therefore it is impossible to define vector $x$ norm as $|x|^2 = x \odot x$. For this purpose we will use $\odot$ generalization.

If for some vector $x \in \mathbb{R}^{n+1}, x \odot x > 0$, then the vector $x' = \frac{x}{\sqrt{x \odot x}} \in \boldsymbol{B}^n$. Because we are interested in direction of vectors rather than in their norm, for sphere vectors it is convenient to use homogeneous coordinates $x = \{x_0 : x_1 : ... : x_n\}$, that have the property $x \equiv \lambda x = \{\lambda x_0 :$

$\lambda x_1 : ... : \lambda x_n\}, \lambda \neq 0$. Because all vectors from the sphere have positive meta product square, we will "norm" them so that $x \odot x = 1$.

**Definition 2.2.3** (Homogeneous space model, signature, metaspace)**.** Define *homogeneous space* $\mathbb{B}^n$ as sphere $\boldsymbol{B}^n \subset \mathbb{R}^{n+1}$. We call the collection of numbers $\{k_1, ..., k_n\}$, each being the type of corresponding main rotation, the *signature* of homogeneous space. We call $\mathbb{R}^{n+1}$ with meta product of vectors $\odot$ (2.11) *metaspace.* Vector $e = (1 : 0 : ... : 0)$ we call *origin* $\boldsymbol{E} = (1 : 0 : ... : 0)$.

*Remark.* Because of homogeneous coordinates, a point in space $\mathbb{B}^n$ is presented by *two* opposite vectors in sphere $\boldsymbol{B}^n$ (from $x \odot x = 1$ follows $(-x) \odot (-x) = 1$). We can define space $\mathbb{B}^n$ by means of real projective space $\mathbb{RP}^n$. It affects the notations, but not the sense of theory. For clarity we will use space sphere here. To be certain, from these two vectors we can choose one, whose first non–zero coordinate is positive (as we already seen, all coordinates of a point in $\mathbb{B}^n$ can't be zero simultaneously).

*Remark.* Defined here signature $\{k_1, ..., k_n\}$ of some homogeneous space $\mathbb{B}^n$ *is not* the same notion as classical signature $(m, n{-}m)$ for pseudo–Euclidean $\mathbb{E}^{m,n-m}$, or classical signature $(m, n{-}m{-}d, d)$ for semi–Euclidean ${}^d\mathbb{E}^m_n$, spaces. Rather, it *plays similar role* in space definition. While classical space signature describes the properties of coordinate axes, defined here signature describes the properties of distances, plane and dihedral angles of different dimensions.

*Remark.* It should be emphasized that a space is defined by its signature, not by the form of meta product of vectors or equation of sphere $\boldsymbol{B}^n$.

**Example.** Homogeneous plane with signature $\{1, -1\}$ has sphere $\boldsymbol{B}^n$ equation:

$$x_0^2 + x_1^2 - x_2^2 = 1,$$

and homogeneous plane with signature $\{-1, -1\}$ has sphere $\boldsymbol{B}^n$ equation:

$$x_0^2 - x_1^2 + x_2^2 = 1.$$

Despite the obvious similarity of these two equations, linear mapping $\mathfrak{F} : \mathbb{R}^3 \to \mathbb{R}^3$ defined as $\mathfrak{F}((x_0 : x_1 : x_2)) = (x_0 : x_2 : x_1)$, that maps the first sphere to the second, is isomorphism of spaces, not automorphism, or even more, a motion. The first plane has line type $k_1 = 1$, the second has $k_1 = -1$. The isomorphism maps lines of the first plane to non–line sequences of points of the second one and vice versa.

**Example.** Homogeneous planes $\{0, 1\}$ (Euclidean), $\{0, 0\}$ (Galilean) and $\{0, -1\}$ (Minkowski) all have space sphere equation:

$$x_0^2 = 1.$$

Obviously, all them are distinct.

*Remark.* One shouldn't forget that meta product $\odot$ is defined in metaspace $\mathbb{R}^{n+1}$, not in space $\mathbb{B}^n$. That is why it plays different role for $\mathbb{B}^n$ than classical dot product of vectors of Euclidean space.

*Properties.* It is easy to see that for any signature $\{k_1, k_2, ...k_n\}$, $E = \mathbb{B}^0 \subset \mathbb{B}^1 \subset ... \subset \mathbb{B}^n$.

**Definition 2.2.4** (Motion, line, plane)**.** Define *motion* of $\mathbb{B}^n$ all possible transformations that are composed of finite product of main rotations. Define *lines* the images $\mathfrak{M}(\mathbb{B}^1)$ under all possible motions $\mathfrak{M} : \mathbb{B}^n \to \mathbb{B}^n$. Similarly, define *m-dimensional planes* as images $\mathfrak{M}(\mathbb{B}^m)$ under all possible motions $\mathfrak{M}$, $m = \overline{0, n-1}$.

For each parameter of type $k_i$ it can be introduced positive scale parameter $r_i \in \mathbb{R}, i = \overline{1, n}$. The value of $\frac{k_1}{r_1^2}$ is gaussian curvature of the space. Other scale parameters have no direct representation, because finite measure of angle doesn't require a scaling. In this case radian measure is natural one. An example of angular scale is degree measure, whose scale is $\frac{180}{\pi}$. However, when angular measure isn't bounded, it makes sense to introduce angular scaling. All scales can be easily embedded in definition of functions $C_i(x)$, $S_i(x)$ and $T_i(x)$ if use instead $C_i\left(\frac{x}{r_i}\right)$, $S_i\left(\frac{x}{r_i}\right)$ and $T_i\left(\frac{x}{r_i}\right)$ respectively, $i = \overline{1, n}$.

**Definition 2.2.5** (Linear space)**.** For theory needs, we call *linear* space, the one which has the first element $k_1 = 0$ in signature.

### 2.2.3. Definition of measure using motions.

Traditionally measure and motion are defined as follows: distance is defined firstly and then motions are defined as all tranformations $\mathfrak{M} : \mathbb{R}^n \to \mathbb{R}^n$ that preserve distances. We use different approach. We define firstly all motions, and then will search for the way to determine measures, so that motions preserve them.

**Definition 2.2.6** (Distance, angle)**.** We say, the *distance* between point $A \in \mathbb{B}^1 \subset \mathbb{B}^n$ and origin $E$ is $\varphi$, if $A = \mathfrak{R}_1(\varphi)E$. We say, the *one–dimensional (plane) angle* between $\mathbb{B}^1$ and some one–dimensional line $\mathbb{B}^0 \subset \mathbb{B}'^1 \subset \mathbb{B}^2$ equals to $\varphi$, if $\mathbb{B}'^1 = \mathfrak{R}_2(\varphi)\mathbb{B}^1$. Similarly, define $(m+1)$*–dimensional angle* $\varphi$ between $\mathbb{B}^m$ and $m$-dimensional plane $\mathbb{B}^{m-1} \subset \mathbb{B}'^m \subset \mathbb{B}^{m+1}$, if $\mathbb{B}'^m = \mathfrak{R}_{m+1}(\varphi)\mathbb{B}^m, \quad \forall m = \overline{0, m-1}$.

*Remark.* It is easy to see that distances, plane and dihedral angles of different dimensions are similarly defined and represent, generally speaking, entities of the same nature.

**Definition 2.2.7** (Measure)**.** We call distances, plane and dihedral angles of different dimensionas generally space *measures.*

**Definition 2.2.8** (Connectable and unconnectable points, measurable measure)**.** Let $X, Y \in \mathbb{B}^n$. If exists a motion that maps $\mathbb{B}^1$ to $XY$, we call points $X$ and $Y$ *connectable* and distance $XY$ *measurable.* If not, we call points $X$ and $Y$ *unconnectable* (dual to parallel lines) and, strictly speaking, the distance $XY$ doesn't exist.

**Lemma 2.2.2.** *Meta product of connectable points of space, $X \odot Y; \quad X, Y \in \mathbb{B}^n$ depends only on distance between them.*

*Proof.* Firstly, consider $\mathbb{B}^1$ case. Let $X = \mathfrak{R}_1(\varphi)E = (C_1(\varphi) : S_1(\varphi))$, $Y = \mathfrak{R}_1(\psi)E = (C_1(\psi) : S_1(\psi))$. Then by (2.8):

$$X \odot Y = \mathfrak{R}_1(\varphi)E \odot \mathfrak{R}_1(\psi)E = C_1(\varphi)C_1(\psi) + k_1 S_1(\varphi)S_1(\psi) = C_1(\varphi - \psi) = C_1(\delta).$$

Let motion $\mathfrak{R}_1(\theta)$ map points $X, Y$ to $X', Y'$. By (2.1):

$$\begin{aligned} X' &= \mathfrak{R}_1(\theta)X = \mathfrak{R}_1(\theta)\mathfrak{R}_1(\varphi)E = \mathfrak{R}_1(\theta + \varphi)E, \\ Y' &= \mathfrak{R}_1(\theta)Y = \mathfrak{R}_1(\theta)\mathfrak{R}_1(\psi)E = \mathfrak{R}_1(\theta + \psi)E; \\ X' \odot Y' &= \mathfrak{R}_1(\theta + \varphi)E \odot \mathfrak{R}_1(\theta + \psi)E = C_1(\theta + \varphi - \theta - \psi) = C_1(\varphi - \psi) = C_1(\delta). \end{aligned}$$

So, in $\mathbb{B}^1$ product $X \odot Y$ depends only on distance $\delta$ between points $X$ and $Y$.

Let now dimension of space $\mathbb{B}^n$ be greater than one. By definition of point connectability, exists a line $l$ that connects $X$ and $Y$. By definition of line, exists a motion $\mathfrak{M}$, so that $\mathfrak{M}\mathbb{B}^1 = l$. Let points $X', Y' \in \mathbb{B}^1$, such that $\mathfrak{M}X' = X$, $\mathfrak{M}Y' = Y$. As it is already shown, product value $X' \odot Y'$ depends only on distance between them. Further, by lemma 2.2.1 motion $\mathfrak{M}$ preserves product of vectors:

$$X \odot Y = \mathfrak{M}X' \odot \mathfrak{M}Y' = X' \odot Y'.$$

So, product value $X \odot Y$ depends only on distance between points $X$ and $Y$. □

**Corollary.** *It is easy to see that motions preserve distances, because by lemma 2.2.1 motions preserve meta product, which by lemma 2.2.2 depends only on distance.*

In case of elliptic, Euclidean and hyperbolic spaces it is sufficient, because all other measures can be calculated based on distances. However, in some spaces exist angular dilation transformations, similar to Euclidean distance dilation transformation. That is why we need a way to find all measures in general case.

### 2.3. Relations in Triangle

Triangle is a plane fugure that is a good illustration of the duality principle. Its three vertices are dual to its three edges, because each edge it defined by two vertices and each vertix is defined by two edges. As we will see, duality in triangle means also metric duality between edges lengths and angles measures.

During deduction of equations and inequations of triangle, consider that all its vertices are mutually connectable. It guarantees that all egdes are lines and then all angles between edges are measurable. Note, that even in this case if $k_2 = 0$ or $k_2 = -1$, one edge contain a point (or points) that is unconnectable with opposed vertex and thus one of internal angles is unmeasurable. Corresponding external angle is always measurable. Triangles with unconnectable vertices are not of much interest. Still they can be considered with care to be taken on edges types and angles measurability (the angle between a line and a non–line sequence of points isn't measurable).

### 2.3.1. Triangle equations.

Consider triangle $\triangle ABC \in \mathbb{B}^2$ with edges $a, b, c$, internal angles $\alpha, \gamma$ and external angle $\beta'$ (Figure 2.6). Let $A = E = (1 : 0 : 0)$ origin, $C = \mathfrak{R}_1(b)E = (C_1(b) : S_1(b) : 0)$ and $B = \mathfrak{R}_2(\alpha)\mathfrak{R}_1(c)E = (C_1(c) : S_1(c)C_2(\alpha) : S_1(c)S_2(\alpha))$.

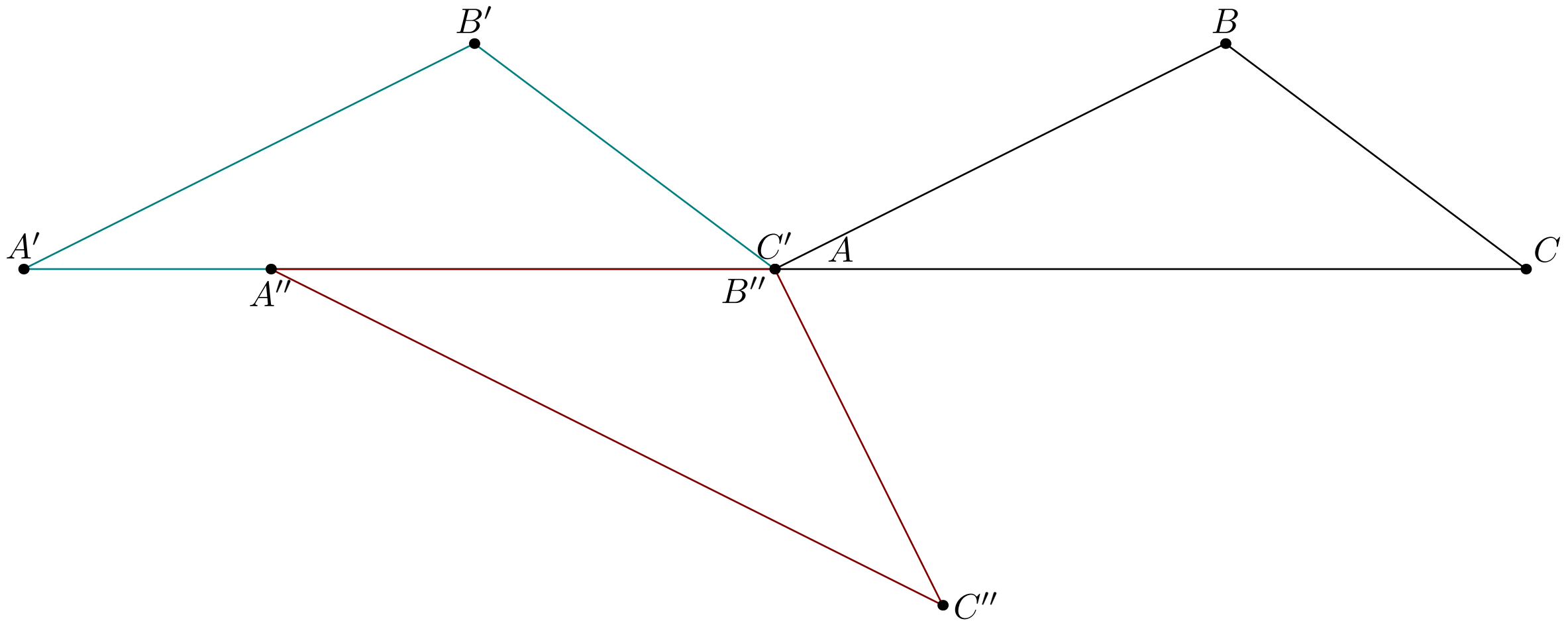

Figure 2.6: Generic triangle equations deduction.

Let now $A'B'C' = \mathfrak{R}_1(-b)(ABC)$ (Figure 2.6, cyan). The point we interested in is $B' = \mathfrak{R}_1(-b)B = (C_1(b)C_1(c) + k_1 S_1(b)S_1(c)C_2(\alpha) : -S_1(b)C_1(c) + C_1(b)S_1(c)C_2(\alpha) : S_1(c)S_2(\alpha))$. On the other hand, $B' = \mathfrak{R}_2(-\gamma)\mathfrak{R}_1(-a)E = (C_1(a) : -S_1(a)C_2(\gamma) : S_1(a)S_2(\gamma))$. This means that:

$$C_1(a) = C_1(b)C_1(c) + k_1 S_1(b)S_1(c)C_2(\alpha), \tag{2.16}$$

$$-S_1(a)C_2(\gamma) = -S_1(b)C_1(c) + C_1(b)S_1(c)C_2(\alpha), \tag{2.17}$$

$$S_1(a)S_2(\gamma) = S_1(c)S_2(\alpha). \tag{2.18}$$

The first equation (2.16) is I cosine low for $\alpha$. Similarly have for $\gamma$:

$$C_1(c) = C_1(a)C_1(b) + k_1 S_1(a)S_1(b)C_2(\gamma). \tag{2.19}$$

The third equation (2.18) is equivalent to:

$$\frac{S_1(a)}{S_2(\alpha)} = \frac{S_1(c)}{S_2(\gamma)},$$

which is a part of sine low.

Now let $A''B''C'' = \mathfrak{R}_1(-c)\mathfrak{R}_2(-\alpha)(ABC)$ (Figure 2.6, brown). Now we interested in $C'' = (C_1(b)C_1(c)+k_1S_1(b)S_1(c)C_2(\alpha) : -C_1(b)S_1(c)+S_1(b)C_1(c)C_2(\alpha) : -S_1(b)S_2(\alpha))$. On the other hand, $C'' = \mathfrak{R}_2(-\beta')\mathfrak{R}_1(a)E = (C_1(a) : S_1(a)C_2(\beta') : -S_1(a)S_2(\beta'))$. It follows:

$$C_1(a) = C_1(b)C_1(c) + k_1 S_1(b)S_1(c)C_2(\alpha), \tag{2.20}$$

$$S_1(a)C_2(\beta') = -C_1(b)S_1(c) + S_1(b)C_1(c)C_2(\alpha), \tag{2.21}$$

$$-S_1(a)S_2(\beta') = -S_1(b)S_2(\alpha). \tag{2.22}$$

The first equation (2.20) is the same I cosine low, the third one (2.22) toghether with (2.18) is equivalent to:

$$\frac{S_1(a)}{S_2(\alpha)} = \frac{S_1(b)}{S_2(\beta')} = \frac{S_1(c)}{S_2(\gamma)}, \tag{2.23}$$

which is sine low. Observe, that in case $k_2 = 1$, we have $\beta = \pi - \beta'$, $\sin\beta = \sin\beta'$. Calculate value of $C_2(\alpha)$ from (2.20) and put it in (2.21):

$$\begin{aligned} S_1(a)C_2(\beta') &= -C_1(b)S_1(c) + S_1(b)C_1(c)\frac{C_1(a) - C_1(b)C_1(c)}{k_1 S_1(b)S_1(c)} \\ &= -C_1(b)S_1(c) + C_1(c)\frac{C_1(a) - C_1(b)C_1(c)}{k_1 S_1(c)}, \\ k_1 S_1(a)S_1(c)C_2(\beta') &= -k_1 S_1(c)^2 C_1(b) + C_1(a)C_1(c) - C_1(b)C_1(c)^2 \\ &= C_1(a)C_1(c) - C_1(b)(C_1(c)^2 + k_1 S_1(c)^2) \\ &= C_1(a)C_1(c) - C_1(b), \\ C_1(b) &= C_1(a)C_1(c) - k_1 S_1(a)S_1(c)C_2(\beta'). \end{aligned} \tag{2.24}$$

The last equation (2.24) is I cosine low for $\beta'$. Remark '$-$' sign in right part of equation. It is so because the angle $\beta'$ is external. For the case $k_2 = 1$ internal angle $\beta = \pi - \beta'$, $\cos\beta = -\cos\beta'$.

Let deduce II cosine low. Consider equations (2.17, 2.21):

$$-S_1(a)C_2(\gamma) = -S_1(b)C_1(c) + C_1(b)S_1(c)C_2(\alpha),$$
$$S_1(a)C_2(\beta') = -C_1(b)S_1(c) + S_1(b)C_1(c)C_2(\alpha)$$

First, using (2.23) change $S_1(b)$ to $\frac{S_1(a)S_2(\beta')}{S_2(\alpha)}$ and $S_1(c)$ to $\frac{S_1(a)S_2(\gamma)}{S_2(\alpha)}$:

$$-S_1(a)C_2(\gamma) = -S_1(a)\frac{S_2(\beta')}{S_2(\alpha)}C_1(c) + C_1(b)S_1(a)\frac{S_2(\gamma)}{S_2(\alpha)}C_2(\alpha),$$
$$-S_2(\alpha)C_2(\gamma) = -C_1(c)S_2(\beta') + C_1(b)S_2(\gamma)C_2(\alpha),$$
$$S_2(\beta')C_1(c) = S_2(\alpha)C_2(\gamma) + C_2(\alpha)S_2(\gamma)C_1(b), \tag{2.25}$$

and

$$S_1(a)C_2(\beta') = -C_1(b)S_1(a)\frac{S_2(\gamma)}{S_2(\alpha)} + S_1(a)\frac{S_2(\beta')}{S_2(\alpha)}C_1(c)C_2(\alpha),$$
$$S_2(\alpha)C_2(\beta') = -C_1(b)S_2(\gamma) + C_1(c)S_2(\beta')C_2(\alpha),$$
$$S_2(\gamma)C_1(b) = -S_2(\alpha)C_2(\beta') + C_2(\alpha)S_2(\beta')C_1(c). \tag{2.26}$$

Now, from the first equation (2.25) calculate $C_1(c)$ and put it in the second one (2.26):

$$S_2(\gamma)C_1(b) = -S_2(\alpha)C_2(\beta') + C_2(\alpha)S_2(\beta')\frac{S_2(\alpha)C_2(\gamma) + C_2(\alpha)S_2(\gamma)C_1(b)}{S_2(\beta')}$$
$$= -S_2(\alpha)C_2(\beta') + C_2(\alpha)S_2(\alpha)C_2(\gamma) + C_2(\alpha)^2S_2(\gamma)C_1(b),$$
$$S_2(\gamma)C_1(b)(1 - C_2(\alpha)^2) = S_2(\alpha)(C_2(\alpha)C_2(\gamma) - C_2(\beta')),$$
$$k_2S_2(\gamma)C_1(b)S_2(\alpha)^2 = S_2(\alpha)(C_2(\alpha)C_2(\gamma) - C_2(\beta')),$$
$$k_2S_2(\alpha)S_2(\gamma)C_1(b) = C_2(\alpha)C_2(\gamma) - C_2(\beta'),$$
$$C_2(\beta') = C_2(\alpha)C_2(\gamma) - k_2S_2(\alpha)S_2(\gamma)C_1(b). \tag{2.27}$$

Equation (2.27) is II cosine low for $b$. When $k_2 = 1$, we have:

$$-\cos\beta = \cos\alpha\cos\gamma - \sin\alpha\sin\gamma C_1(b),$$
$$\cos\beta = -\cos\alpha\cos\gamma + \sin\alpha\sin\gamma C_1(b).$$

Similarly, calculating $C_1(b)$ from (2.26) and putting it in (2.25), obtain:

$$\begin{aligned}
S_2(\beta')C_1(c) &= S_2(\alpha)C_2(\gamma) + C_2(\alpha)S_2(\gamma)\frac{C_2(\alpha)S_2(\beta')C_1(c) - S_2(\alpha)C_2(\beta')}{S_2(\gamma)} \\
&= S_2(\alpha)C_2(\gamma) + C_2(\alpha)^2 S_2(\beta')C_1(c) - C_2(\alpha)S_2(\alpha)C_2(\beta'), \\
S_2(\beta')C_1(c)(1 - C_2(\alpha)^2) &= S_2(\alpha)(C_2(\gamma) - C_2(\alpha)C_2(\beta')), \\
k_2 S_2(\beta')C_1(c)S_2(\alpha)^2 &= S_2(\alpha)(C_2(\gamma) - C_2(\alpha)C_2(\beta')), \\
k_2 S_2(\alpha)S_2(\beta')C_1(c) &= C_2(\gamma) - C_2(\alpha)C_2(\beta'), \\
C_2(\gamma) &= C_2(\alpha)C_2(\beta') + k_2 S_2(\alpha)S_2(\beta')C_1(c). \qquad (2.28)
\end{aligned}$$

Equation (2.28) is II cosine low for $c$. When $k_2 = 1$, obtain as earlier:

$$\cos\gamma = -\cos\alpha\cos\beta + \sin\alpha\sin\beta C_1(c).$$

Similarly to (2.28), have for $a$:

$$C_2(\alpha) = C_2(\beta')C_2(\gamma) + k_2 S_2(\beta')S_2(\gamma)C_1(a). \qquad (2.29)$$

We find a form of I and II cosine lows that don't contain functions $C_1$ or $C_2$ in left part. However, they contain these functions in right part. It makes sense in case $k_1 \neq 0$ (for I cosine low) and $k_2 \neq 0$ (for II cosine low), when it is possible to calculate respective $C^{-1}$ functions. But when $k_1 = 0$, the plane allows distance dilations (II cosine low becomes identity that doesn't contain function $C_1(x)$), while when $k_2 = 0$, the plane allows angles dilations and lengths don't define the angles (I cosine low becomes identity that doesn't contain function $C_2(x)$).

Remark also that we can deduce general form of I and II cosine lows if introduce (may be unmeasurable) angle $\beta$ so that:

$$\begin{aligned}
S_2(\beta) &= S_2(\beta'), \\
C_2(\beta) &= -C_2(\beta'), \\
T_2(\beta) &= -T_2(\beta'),
\end{aligned}$$

then I and II cosine lows take the same form.

Let calculate:

$$\begin{aligned}
k_1 S_1^2(a) &= 1 - C_1^2(a) = (C_1^2(b) + k_1 S_1^2(b))(C_1^2(c) + k_1 S_1^2(c)) \\
&- (C_1(b)C_1(c) + k_1 S_1(b)S_1(c)C_2(\alpha))^2 \\
&= C_1^2(b)C_1^2(c) + k_1 C_1^2(b)S_1^2(c) + k_1 S_1^2(b)C_1^2(c) + k_1^2 S_1^2(b)S_1^2(c) \\
&- C_1^1(b)C_1^2(c) - 2k_1 C_1(b)C_1(c)S_1(b)S_1(c)C_2(\alpha) - k_1^2 S_1^2(b)S_1^2(c)C_2^2(\alpha) \\
&= k_1(C_1^2(b)S_1^2(c) + S_1^2(b)C_1^2(c) - 2C_1(b)C_1(c)S_1(b)S_1(c)C_2(\alpha)) \\
&+ k_1^2 S_1^2(b)S_1^2(c)(1 - C_2^2(\alpha)) \\
&= k_1(C_1^2(b)S_1^2(c) + S_1^2(b)C_1^2(c) - 2C_1(b)C_1(c)S_1(b)S_1(c)C_2(\alpha)) \\
&+ k_1^2 k_2 S_1^2(b)S_1^2(c)S_2^2(\alpha),
\end{aligned}$$

$$\begin{aligned}
S_1^2(a) &= C_1^2(b)S_1^2(c) + S_1^2(b)C_1^2(c) - 2C_1(b)C_1(c)S_1(b)S_1(c)C_2(\alpha) \\
&+ k_1 k_2 S_1^2(b)S_1^2(c)S_2^2(\alpha),
\end{aligned} \tag{2.30}$$

or, dividing (2.30) by (2.16) obtain:

$$T_1^2(a) = \frac{T_1^2(b) + T_1^2(c) - 2T_1(b)T_1(c)C_2(\alpha) + k_1 k_2 T_1^2(b)T_1^2(c)S_2^2(\alpha)}{(1 + k_1 T_1(b)T_1(c)C_2(\alpha))^2} \tag{2.31}$$

Similarly,

$$\begin{aligned}
S_1^2(b) &= C_1^2(a)S_1^2(c) + S_1^2(a)C_1^2(c) + 2C_1(a)C_1(c)S_1(a)S_1(c)C_2(\beta') \\
&+ k_1 k_2 S_1^2(a)S_1^2(c)S_2^2(\beta'),
\end{aligned} \tag{2.32}$$

$$T_1^2(b) = \frac{T_1^2(a) + T_1^2(c) + 2T_1(a)T_1(c)C_2(\beta') + k_1 k_2 T_1^2(a)T_1^2(c)S_2^2(\beta')}{(1 - k_1 T_1(a)T_1(c)C_2(\beta'))^2}, \tag{2.33}$$

and

$$\begin{aligned}
S_1^2(c) &= C_1^2(a)S_1^2(b) + S_1^2(a)C_1^2(b) - 2C_1(a)C_1(b)S_1(a)S_1(b)C_2(\gamma) \\
&+ k_1 k_2 S_1^2(a)S_1^2(b)S_2^2(\gamma),
\end{aligned} \tag{2.34}$$

$$T_1^2(c) = \frac{T_1^2(a) + T_1^2(b) - 2T_1(a)T_1(b)C_2(\gamma) + k_1 k_2 T_1^2(a)T_1^2(b)S_2^2(\gamma)}{(1 + k_1 T_1(a)T_1(b)C_2(\gamma))^2} \tag{2.35}$$

Table 2.1: Triangle equations

$$\frac{S_1(a)}{S_2(\alpha)} = \frac{S_1(b)}{S_2(\beta')} = \frac{S_1(c)}{S_2(\gamma)}$$

$$C_1(a) = C_1(b)C_1(c) + k_1 S_1(b)S_1(c)C_2(\alpha)$$
$$C_1(b) = C_1(a)C_1(c) - k_1 S_1(a)S_1(c)C_2(\beta')$$
$$C_1(c) = C_1(a)C_1(b) + k_1 S_1(a)S_1(b)C_2(\gamma)$$
$$C_2(\alpha) = C_2(\beta')C_2(\gamma) + k_2 S_2(\beta')S_2(\gamma)C_1(a)$$
$$C_2(\beta') = C_2(\alpha)C_2(\gamma) - k_2 S_2(\alpha)S_2(\gamma)C_1(b)$$
$$C_2(\gamma) = C_2(\alpha)C_2(\beta') + k_2 S_2(\alpha)S_2(\beta')C_1(c)$$

$$T_1^2(a) = \frac{T_1^2(b) + T_1^2(c) - 2T_1(b)T_1(c)C_2(\alpha) + k_1k_2T_1^2(b)T_1^2(c)S_2^2(\alpha)}{(1 + k_1T_1(b)T_1(c)C_2(\alpha))^2}$$
$$T_1^2(b) = \frac{T_1^2(a) + T_1^2(c) + 2T_1(a)T_1(c)C_2(\beta') + k_1k_2T_1^2(a)T_1^2(c)S_2^2(\beta')}{(1 - k_1T_1(a)T_1(c)C_2(\beta'))^2}$$
$$T_1^2(c) = \frac{T_1^2(a) + T_1^2(b) - 2T_1(a)T_1(b)C_2(\gamma) + k_1k_2T_1^2(a)T_1^2(b)S_2^2(\gamma)}{(1 + k_1T_1(a)T_1(b)C_2(\gamma))^2}$$
$$T_2^2(\alpha) = \frac{T_2^2(\beta') + T_2^2(\gamma) - 2T_2(\beta')T_2(\gamma)C_1(a) + k_1k_2T_2^2(\beta')T_2^2(\gamma)S_1^2(a)}{(1 + k_2T_2(\beta')T_2(\gamma)C_1(a))^2}$$
$$T_2^2(\beta') = \frac{T_2^2(\alpha) + T_2^2(\gamma) + 2T_2(\alpha)T_2(\gamma)C_1(b) + k_1k_2T_2^2(\alpha)T_2^2(\gamma)S_1^2(b)}{(1 - k_2T_2(\alpha)T_2(\gamma)C_1(b))^2}$$
$$T_2^2(\gamma) = \frac{T_2^2(\alpha) + T_2^2(\beta') - 2T_2(\alpha)T_2(\beta')C_1(c) + k_1k_2T_2^2(\alpha)T_2^2(\beta')S_1^2(c)}{(1 + k_2T_2(\alpha)T_2(\beta')C_1(c))^2}$$

Now calculate:

$$
\begin{aligned}
k_2S_2^2(\alpha) &= 1 - C_2^2(\alpha) = (C_2^2(\beta') + k_2S_2^2(\beta'))(C_2^2(\gamma) + k_2S_2^2(\gamma)) \\
&- (C_2(\beta')C_2(\gamma) + k_2S_2(\beta')S_2(\gamma)C_1(a))^2 \\
&= C_2^2(\beta')C_2^2(\gamma) + k_2C_2^2(\beta')S_2^2(\gamma) + k_2S_2^2(\beta')C_2^2(\gamma) + k_2^2S_2^2(\beta')S_2^2(\gamma) \\
&- C_2^2(\beta')C_2^2(\gamma) - 2k_2C_2(\beta')S_2(\beta')C_2(\gamma)S_2(\gamma)C_1(a) - k_2^2S_2^2(\beta')S_2^2(\gamma)C_1^2(a) \\
&= k_2(C_2^2(\beta')S_2^2(\gamma) + S_2^2(\beta')C_2^2(\gamma) - 2C_2(\beta')S_2(\beta')C_2(\gamma)S_2(\gamma)C_1(a)) \\
&+ k_2^2S_2^2(\beta')S_2^2(\gamma)(1 - C_1^2(a)) \\
&= k_2(C_2^2(\beta')S_2^2(\gamma) + S_2^2(\beta')C_2^2(\gamma) - 2C_2(\beta')S_2(\beta')C_2(\gamma)S_2(\gamma)C_1(a)) \\
&+ k_1k_2^2S_2^2(\beta')S_2^2(\gamma)S_1^2(a),
\end{aligned}
$$

$$
\begin{aligned}
S_2^2(\alpha) &= C_2^2(\beta')S_2^2(\gamma) + S_2^2(\beta')C_2^2(\gamma) - 2C_2(\beta')S_2(\beta')C_2(\gamma)S_2(\gamma)C_1(a) \\
&+ k_1k_2S_2^2(\beta')S_2^2(\gamma)S_1^2(a),
\end{aligned}
\tag{2.36}
$$

or dividing (2.36) by (2.29):

$$
T_2^2(\alpha) = \frac{T_2^2(\beta') + T_2^2(\gamma) - 2T_2(\beta')T_2(\gamma)C_1(a) + k_1k_2T_2^2(\beta')T_2^2(\gamma)S_1^2(a)}{(1 + k_2T_2(\beta')T_2(\gamma)C_1(a))^2}
\tag{2.37}
$$

Similarly,

$$
\begin{aligned}
S_2^2(\beta') &= C_2^2(\alpha)S_2^2(\gamma) + S_2^2(\alpha)C_2^2(\gamma) + 2C_2(\alpha)S_2(\alpha)C_2(\gamma)S_2(\gamma)C_1(b) \\
&+ k_1k_2S_2^2(\alpha)S_2^2(\gamma)S_1^2(b),
\end{aligned}
\tag{2.38}
$$

$$
T_2^2(\beta') = \frac{T_2^2(\alpha) + T_2^2(\gamma) + 2T_2(\alpha)T_2(\gamma)C_1(b) + k_1k_2T_2^2(\alpha)T_2^2(\gamma)S_1^2(b)}{(1 - k_2T_2(\alpha)T_2(\gamma)C_1(b))^2},
\tag{2.39}
$$

and

$$
\begin{aligned}
S_2^2(\gamma) &= C_2^2(\alpha)S_2^2(\beta') + S_2^2(\alpha)C_2^2(\beta') - 2C_2(\alpha)S_2(\alpha)C_2(\beta')S_2(\beta')C_1(c) \\
&+ k_1k_2S_2^2(\alpha)S_2^2(\beta')S_1^2(c),
\end{aligned}
\tag{2.40}
$$

$$
T_2^2(\gamma) = \frac{T_2^2(\alpha) + T_2^2(\beta') - 2T_2(\alpha)T_2(\beta')C_1(c) + k_1k_2T_2^2(\alpha)T_2^2(\beta')S_1^2(c)}{(1 + k_2T_2(\alpha)T_2(\beta')C_1(c))^2}
\tag{2.41}
$$

### 2.3.2. Triangle inequations.

Let deduce triangle inequations. When $k_2 = 1$ and one of angles is obtuse, use corresponding external angle, when $k_1 = 1$, the distance between points $X$ and $Y$ isn't greater then $\pi$, the least distance between $X$ and $Y$ or between $X$ and $-Y$ doesn't exceed $\frac{\pi}{2}$. In these conditions the function $\sin x$ is monotone increasing, even for $k = 1$. From sine low (2.23) follows proposition:

**Proposition 2.3.1.** *The longest edge is opposed to the largest angle and the shortest edge is opposed to the smallest angle.*

Consider separately all 9 homogeneous planes, keeping in mind the properties of $C(x)$ function (for $k = 1, 0 \leq x \leq \frac{\pi}{2}$):

*Properties.*

$$C(x)\begin{cases} \leq 1, \text{ decreasing}, & k = 1; \\ = 1, \text{ constant}, & k = 0; \\ \geq 1, \text{ increasing}, & k = -1. \end{cases}$$

$\mathbb{B}^2 = \{1, 1\}$ — **Riemann plane.** From (2.16) follows:

$$\cos a = \cos b \cos c + \sin b \sin c \cos \alpha < \cos b \cos c + \sin b \sin c = \cos(b - c),$$
$$a > b - c.$$

From (2.24) follows:

$$\cos b = \cos a \cos c - \sin a \sin c \cos \beta' > \cos a \cos c - \sin a \sin c = \cos(a + c),$$
$$b < a + c.$$

From (2.29) follows:

$$\cos \alpha = \cos \beta' \cos \gamma + \sin \beta' \sin \gamma \cos a < \cos \beta' \cos \gamma + \sin \beta' \sin \gamma = \cos(\beta' - \gamma),$$
$$\alpha > \beta' - \gamma.$$

From (2.27) follows:

$$\cos \beta' = \cos \alpha \cos \gamma - \sin \alpha \sin \gamma \cos b > \cos \alpha \cos \gamma - \sin \alpha \sin \gamma = \cos(\alpha + \gamma),$$
$$\beta' < \alpha + \gamma.$$

$\mathbb{B}^2 = \{0, 1\}$ — **Euclidean plane.** From (2.31) follows:

$$a^2 = b^2 + c^2 - 2bc\cos\alpha > b^2 + c^2 - 2bc = (b - c)^2,$$
$$a > b - c.$$

From (2.33) follows:

$$b^2 = a^2 + c^2 + 2ac\cos\beta' < a^2 + c^2 + 2ac = (a + c)^2,$$
$$b < a + c.$$

From (2.29), having $C_1(a) = 1, k_1 = 0, \forall a \in \mathbb{R}$, follows:

$$\cos\alpha = \cos\beta'\cos\gamma + \sin\beta'\sin\gamma \cdot 1 = \cos(\beta' - \gamma),$$
$$\alpha = \beta' - \gamma.$$

Same follows from (2.27):

$$\cos\beta' = \cos\alpha\cos\gamma - \sin\alpha\sin\gamma \cdot 1 = \cos(\alpha + \gamma),$$
$$\beta' = \alpha + \gamma.$$

$\mathbb{B}^2 = \{-1, 1\}$ — **hyperbolic plane.** From (2.16) follows:

$$\cosh a = \cosh b\cosh c - \sinh b\sinh c\cos\alpha > \cosh b\cosh c - \sinh b\sinh c = \cosh(b - c),$$
$$a > b - c.$$

From (2.24) follows:

$$\cosh b = \cosh a\cosh c + \sinh a\sinh c\cos\beta' < \cosh a\cosh c + \sinh a\sinh c = \cosh(a + c),$$
$$b < a + c.$$

From (2.29) follows:

$$\cos\alpha = \cos\beta'\cos\gamma + \sin\beta'\sin\gamma\cosh a > \cos\beta'\cos\gamma + \sin\beta'\sin\gamma = \cos(\beta' - \gamma),$$
$$\alpha < \beta' - \gamma.$$

From (2.27) follows:

$$\cos\beta' = \cos\alpha\cos\gamma - \sin\alpha\sin\gamma\cosh b < \cos\alpha\cos\gamma - \sin\alpha\sin\gamma = \cos(\alpha+\gamma),$$
$$\beta' > \alpha + \gamma.$$

$\mathbb{B}^2 = \{1, 0\}$ — **positive curved Galilean plane.** From (2.16), having $C_2(\alpha) = 1, k_2 = 0, \forall\alpha \in \mathbb{R}$ follows:

$$\cos a = \cos b\cos c + \sin b\sin c \cdot 1 = \cos(b - c),$$
$$a = b - c.$$

Same follows from (2.24):

$$\cos b = \cos a\cos c - \sin a\sin c \cdot 1 = \cos(a + c),$$
$$b = a + c.$$

From (2.37) follows:

$$\alpha^2 = \beta'^2 + \gamma^2 - 2\beta\gamma\cos a > \beta'^2 + \gamma^2 - 2\beta\gamma = (\beta' - \gamma)^2,$$
$$\alpha > \beta' - \gamma.$$

From (2.39) follows:

$$\beta'^2 = \alpha^2 + \gamma^2 + 2\alpha\gamma\cos b < \alpha^2 + \gamma^2 + 2\alpha\gamma = (\alpha + \gamma)^2,$$
$$\beta' < \alpha + \gamma.$$

$\mathbb{B}^2 = \{0, 0\}$ — **Galilean plane.** From (2.31), having $C_2(\alpha) = 0, k_2 = 0, \forall\alpha \in \mathbb{R}$ follows:

$$a^2 = b^2 + c^2 - 2bc \cdot 1 = (b - c)^2,$$
$$a = b - c.$$

Same follows from (2.33):

$$b^2 = a^2 + c^2 + 2ac \cdot 1 = (a + c)^2,$$
$$b = a + c.$$

From (2.37), having $C_1(a) = 1, k_1 = 0, \forall a \in \mathbb{R}$, follows:

$$\alpha^2 = \beta'^2 + \gamma^2 - 2\beta'\gamma \cdot 1 = (\beta' - \gamma)^2,$$
$$\alpha = \beta' - \gamma.$$

Same follows from (2.39):

$$\beta'^2 = \alpha^2 + \gamma^2 + 2\alpha\gamma \cdot 1 = (\alpha + \gamma)^2,$$
$$\beta' = \alpha + \gamma.$$

$\mathbb{B}^2 = \{-1, 0\}$ — **negative curved Galilean plane.** From (2.16), having $C_2(\alpha) = 1, k_2 = 0, \forall \alpha \in \mathbb{R}$, follows:

$$\cosh a = \cosh b \cosh c - \sinh b \sinh c \cdot 1 = \cosh(b - c),$$
$$a = b - c.$$

Same follows from (2.24):

$$\cosh b = \cosh a \cosh c + \sinh a \sinh c \cdot 1 = \cosh(a + c),$$
$$b = a + c.$$

From (2.37) follows:

$$\alpha^2 = \beta'^2 + \gamma^2 - 2\beta'\gamma \cosh a < \beta'^2 + \gamma^2 - 2\beta'\gamma = (\beta' - \gamma)^2,$$
$$\alpha < \beta' - \gamma.$$

From (2.39) follows:

$$\beta'^2 = \alpha^2 + \gamma^2 + 2\alpha\gamma \cosh b > \alpha^2 + \gamma^2 + 2\alpha\gamma = (\alpha + \gamma)^2,$$
$$\beta' > \alpha + \gamma.$$

$\mathbb{B}^2 = \{1, -1\}$ — **Anti de Sitter (positive curved Minkowski) plane.** From (2.16) follows:

$$\cos a = \cos b \cos c + \sin b \sin c \cosh \alpha > \cos b \cos c + \sin b \sin c = \cos(b - c),$$
$$a < b - c.$$

From (2.24) follows:

$$\cos b = \cos a \cos c - \sin a \sin c \cosh \beta' < \cos a \cos c - \sin a \sin c = \cos(a + c),$$
$$b > a + c.$$

From (2.29) follows:

$$\cosh \alpha = \cosh \beta' \cosh \gamma - \sinh \beta' \sinh \gamma \cos b > \cosh \beta' \cosh \gamma - \sinh \beta' \sinh \gamma = \cosh(\beta' - \gamma),$$
$$\alpha > \beta' - \gamma.$$

From (2.27) follows:

$$\cosh \beta' = \cosh \alpha \cosh \gamma + \sinh \alpha \sinh \gamma \cos b < \cosh \alpha \cosh \gamma + \sinh \alpha \sinh \gamma = \cosh(\alpha + \gamma),$$
$$\beta' < \alpha + \gamma.$$

$\mathbb{B}^2 = \{0, -1\}$ — **Minkowski plane.** From (2.31) follows:

$$a^2 = b^2 + c^2 - 2bc \cosh \alpha < b^2 + c^2 - 2bc = (b - c)^2,$$
$$a < b - c.$$

From (2.33) follows:

$$b^2 = a^2 + c^2 + 2ac \cosh \beta' > a^2 + c^2 + 2ac = (a + c)^2,$$
$$b > a + c.$$

From (2.29), having $C_1(a) = 1, k_1 = 0, \forall a \in \mathbb{R}$, follows:

$$\cosh \alpha = \cosh \beta' \cosh \gamma - \sinh \beta' \sinh \gamma \cdot 1 = \cosh(\beta' - \gamma),$$
$$\alpha = \beta' - \gamma.$$

From (2.27) follows:

$$\cosh \beta' = \cosh \alpha \cosh \gamma + \sinh \alpha \sinh \gamma \cdot 1 = \cosh(\alpha + \gamma),$$
$$\beta' = \alpha + \gamma.$$

$\mathbb{B}^2 = \{-1, -1\}$ — **De Sitter (negative curved Minkowski) plane.** From (2.16) follows:

$$\cosh a = \cosh b \cosh c - \sinh b \sinh c \cosh \alpha < \cosh b \cosh c - \sinh b \sinh c = \cosh(b - c),$$
$$a < b - c.$$

From (2.24) follows:

$$\cosh b = \cosh a \cosh c + \sinh a \sinh c \cosh \beta' > \cosh a \cosh c + \sinh a \sinh c = \cosh(a + c),$$
$$b > a + c.$$

From (2.29) follows:

$$\cosh \alpha = \cosh \beta' \cosh \gamma - \sinh \beta' \sinh \gamma \cosh a < \cosh \beta' \cosh \gamma - \sinh \beta' \sinh \gamma = \cosh(\beta' - \gamma),$$
$$\alpha < \beta' - \gamma.$$

From (2.27) follows:

$$\cosh \beta' = \cosh \alpha \cosh \gamma + \sinh \alpha \sinh \gamma \cosh b > \cosh \alpha \cosh \gamma + \sinh \alpha \sinh \gamma = \cosh(\alpha + \gamma),$$
$$\beta' > \alpha + \gamma.$$

Summing all the above up, obtain the following porposition:

**Proposition 2.3.2.** *The shortest triangle edge*

$$a \begin{cases} > b - c, & k_2 = 1; \\ = b - c, & k_2 = 0; \\ < b - c, & k_2 = -1. \end{cases}$$

*The longest triangle edge*

$$b \begin{cases} < a + c, & k_2 = 1; \\ = a + c, & k_2 = 0; \\ > a + c, & k_2 = -1. \end{cases}$$

*Internal angle*

$$\alpha \begin{cases} > \beta' - \gamma, & k_1 = 1; \\ = \beta' - \gamma, & k_1 = 0; \\ < \beta' - \gamma, & k_1 = -1. \end{cases}$$

*External angle*

$$\beta' \begin{cases} < \alpha + \gamma, & k_1 = 1; \\ = \alpha + \gamma, & k_1 = 0; \\ > \alpha + \gamma, & k_1 = -1. \end{cases}$$

**Corollary.** *In cases $k_2 = 0$ and $k_2 = -1$ the triangle has explicit gteatest edge, which is equal to (when $k_2 = 0$) or greater then (when $k_2 = -1$) sum of other two. This edge is opposed to external angle.*

### 2.3.3. Right triangle equations.

Till now the property of vector orthogonality in general case isn't defined (it will be done during study of generalized orthogonal matrix). For now it is sufficient to know that all coordinate vectors are mutually orthogonal.

Strictly speaking, we can't talk about right triangle, as one of its catheti is a line, but another is not (when $k_2 \neq 1$). However, this figure is important for orthogonality analysis. We call it also right triangle, keeping in mind that one cathet may have type $k \neq k_1$.

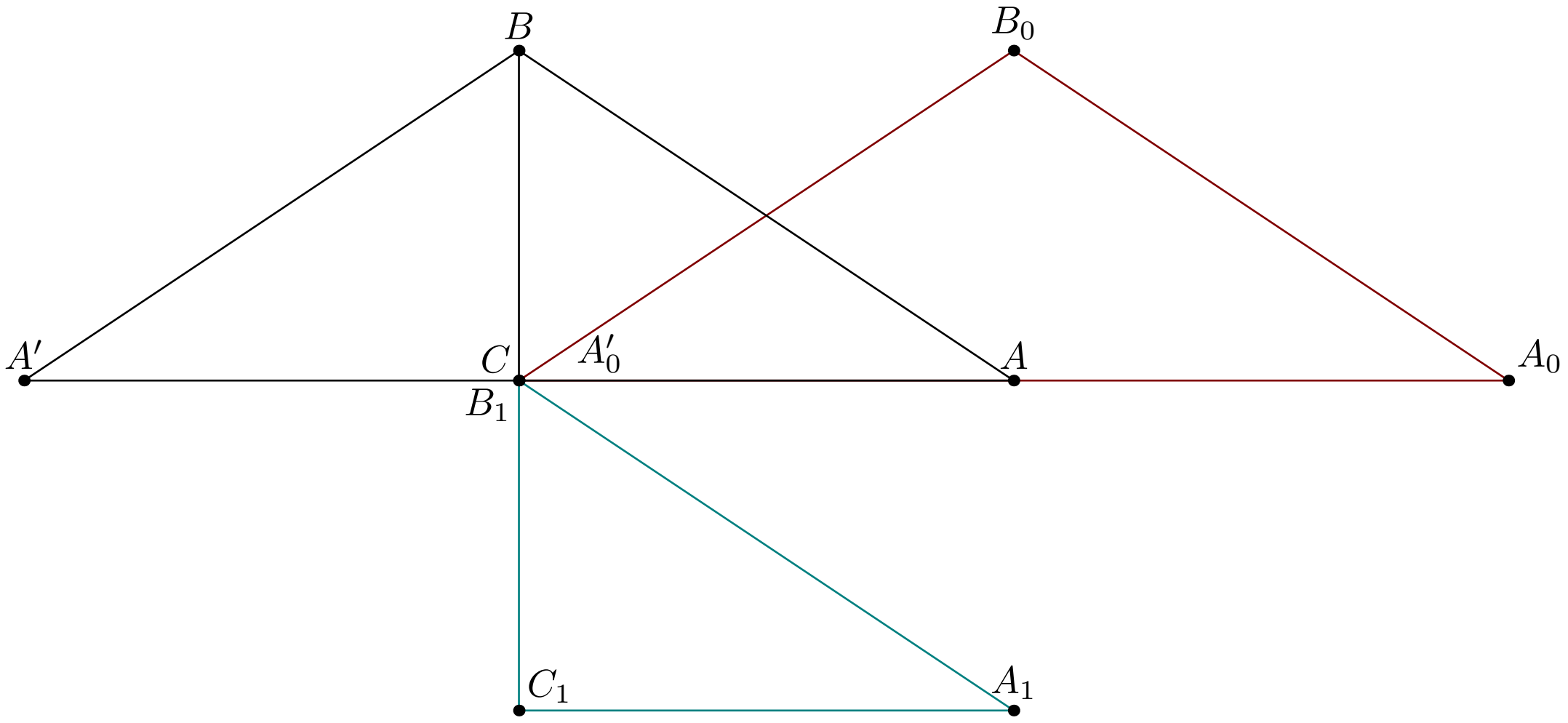

Figure 2.7: Right triangle equations deduction.

Construct right triangle as half of isosceles one (Figure 2.7). Let triangle $\triangle A_0 B_0 A'_0$ have edges $A_0B_0 = A'_0B_0 = c$, $A_0A'_0 = 2b$ and angles $\angle A'_0A_0B_0 = \angle A_0A'_0B_0 = \alpha$ and external angle $\angle A_0B_0A'_0 = 2\beta'$.

Let $A'_0 = E = (1 : 0 : 0)$ be origin, $A_0 = \mathfrak{R}_1(2b)A'_0 = (C_1(2b) : S_1(2b) : 0)$, $B_0 = \mathfrak{R}_2(\alpha)\mathfrak{R}_1(c)A_0 = (C_1(c) : S_1(c)C_2(\alpha) : S_1(c)S_2(\alpha))$.

Let $ABA' = \mathfrak{R}_1(-b)(A_0B_0A'_0)$ (Figure 2.7, black). Now $A' = \mathfrak{R}_1(-b)A'_0 = (C_1(b) : -S_1(b) : 0)$, $A = \mathfrak{R}_1(-b)A_0 = (C_1(b) : S_1(b) : 0)$ and $B = \mathfrak{R}_1(-b)B_0 =$

$$\begin{pmatrix} C_1(b) & k_1S_1(b) & 0 \\ -S_1(b) & C_1(b) & 0 \\ 0 & 0 & 1 \end{pmatrix} \begin{pmatrix} C_1(c) \\ S_1(c)C_2(\alpha) \\ S_1(c)S_2(\alpha) \end{pmatrix}$$

$= (C_1(b)C_1(c) + k_1S_1(b)S_1(c)C_2(\alpha) : -S_1(b)C_1(c) + C_1(b)S_1(c)C_2(\alpha) : S_1(c)S_2(\alpha))$.

Finally, let $C \in AA'$, $AC = A'C = b$. Then $C = (1 : 0 : 0) = E$, the origin. From horizontal symmetry $\triangle A'BA$ and the fact $C$ is origin results that $B$ has the form $B = (x : 0 : y)$. Having in origin ($E = C$) vectors $CA = \{0 : 1 : 0\}$ and $CB = \{0 : 0 : 1\}$, that are both coordinate ones, $CB \perp CA$. So we can consider $\triangle ABC$ as right triangle.

$B = (x : 0 : y)$, that is:

$$C_1(b)C_1(c) + k_1S_1(b)S_1(c)C_2(\alpha) = x, \tag{2.42}$$

$$-S_1(b)C_1(c) + C_1(b)S_1(c)C_2(\alpha) = 0, \tag{2.43}$$

$$S_1(c)S_2(\alpha) = y. \tag{2.44}$$

From the second equality (2.43) follows:

$$T_1(b) = T_1(c)C_2(\alpha). \tag{2.45}$$

Use the value of $C_2(\alpha)$ from this equality and put it to (2.42), obtain:

$$x = C_1(b)C_1(c) + k_1S_1(b)S_1(c)\frac{T_1(b)}{T_1(c)} = \frac{C_1(c)}{C_1(b)}(C_1^2(b) + k_1S_1^2(b)) = \frac{C_1(c)}{C_1(b)}.$$

Now calculate value

$$\begin{aligned}x^2 + k_1k_2y^2 &= \frac{C_1^2(c)}{C_1^2(b)} + k_1k_2S_1^2(c)S_2^2(\alpha) = \frac{C_1^2(c)}{C_1^2(b)} + k_1S_1^2(c)(1 - C_2^2(\alpha))\\ &= \frac{C_1^2(c)}{C_1^2(b)} + k_1S_1^2(c) - k_1S_1^2(c)\frac{T_1^2(b)}{T_1^2(c)} = \frac{C_1^2(c)}{C_1^2(b)} + k_1S_1^2(c) - k_1S_1^2(b)\frac{C_1^2(c)}{C_1^2(b)}\\ &= \frac{C_1^2(c)}{C_1^2(b)}(1 - k_1S_1^2(b)) + k_1S_1^2(c) = \frac{C_1^2(c)}{C_1^2(b)}C_1^2(b) + k_1S_1^2(c)\\ &= C_1^2(c) + k_1S_1^2(c) = 1.\end{aligned}$$

It means that exists $a \in \mathbb{R}$ so that $C_{12}(a) = x = \frac{C_1(c)}{C_1(b)}$, $S_{12}(a) = y = S_1(c)S_2(\alpha)$, that has the type $k = k_1k_2$. The measure $a$ is "distance" parameter of $BC$. Obtain two more equations:

$$C_1(c) = C_{12}(a)C_1(b) \tag{2.46}$$

$$S_{12}(a) = S_1(c)S_2(\alpha) \tag{2.47}$$

Divide (2.47) by (2.45):

$$\frac{S_{12}(a)}{T_1(b)} = \frac{S_1(c)S_2(\alpha)}{T_1(c)C_2(\alpha)} = C_1(c)T_2(\alpha),$$

using value $C_1(c)$ from (2.46), obtain:

$$\frac{S_{12}(a)}{T_1(b)} = C_{12}(a)C_1(b)T_2(\alpha),$$

$$T_{12}(a) = S_1(b)T_2(\alpha). \tag{2.48}$$

The last 6 equations contain the angle between $BC$ and $BA$, which can be unmeasurable. Instead of it consider the angle $\beta'$ between $BA$ and perpendicular to $BC$. In order to deduce next equations, consider transformation $\mathfrak{T}(-a)$ that maps $\mathfrak{T}(-a)B = C$. Having the type $k_1k_2 = K_2$ of distance $a$, obtain:

$$\mathfrak{T}(-a) = \begin{pmatrix} C_{12}(a) & 0 & K_2S_{12}(a)\\ 0 & 1 & 0\\ -S_{12}(a) & 0 & C_{12}(a)\end{pmatrix}.$$

It can be checked this transformation preserves meta product of vectors as well as distances, similar to how it was made in proof of lemma 2.2.1.

Applying $\mathfrak{T}(-a)$, obtain $B_1 = \mathfrak{T}(-a)B = (1 : 0 : 0) = E$, $C_1 = \mathfrak{T}(-a)C = (C_{12}(a) : 0 :$

$-S_{12}(a))$ and $A_1 = \mathfrak{T}(-a)A =$

$$\begin{pmatrix} C_{12}(a) & 0 & K_2 S_{12}(a) \\ 0 & 1 & 0 \\ -S_{12}(a) & 0 & C_{12}(a) \end{pmatrix} \begin{pmatrix} C_1(b) \\ S_1(b) \\ 0 \end{pmatrix}$$

$= (C_{12}(a)C_1(b) : S_1(b) : -S_{12}(a)C_1(b))$. On the other hand, now $A_1 = \mathfrak{R}_2(-\beta')\mathfrak{R}_1(c)E = (C_1(c) : S_1(c)C_2(\beta') : -S_1(c)S_2(\beta'))$ (Figure 2.7, cyan). From here follows:

$$S_1(b) = S_1(c)C_2(\beta'). \tag{2.49}$$

Moreover, having (2.46), obtain:

$$\begin{aligned} S_{12}(a)C_1(b) &= S_1(c)S_2(\beta'), \\ S_{12}(a)\frac{C_1(c)}{C_{12}(a)} &= S_1(c)S_2(\beta'), \\ T_{12}(a) &= T_1(c)S_2(\beta'). \end{aligned} \tag{2.50}$$

Combine last equalities (2.49) and (2.50) with (2.46), obtain:

$$\begin{aligned} &\frac{T_{12}(a)}{S_1(b)} = \frac{T_1(c)S_2(\beta')}{S_1(c)C_2(\beta')} = \frac{T_2(\beta')}{C_1(c)} = \frac{T_2(\beta')}{C_{12}(a)C_1(b)}, \\ &S_{12}(a) = T_1(b)T_2(\beta'). \end{aligned} \tag{2.51}$$

Now, having (2.45) and (2.49):

$$T_1(c)C_2(\alpha) = T_1(b) = \frac{S_1(b)}{C_1(b)} = \frac{S_1(c)C_2(\beta')}{C_1(b)},$$

compute with (2.46):

$$\begin{aligned} C_1(b) &= \frac{S_1(c)C_2(\beta')}{T_1(c)C_2(\alpha)} = C_1(c)\frac{C_2(\beta')}{C_2(\alpha)} \\ &= C_{12}(a)C_1(b)\frac{C_2(\beta')}{C_2(\alpha)}, \\ 1 &= C_{12}(a)\frac{C_2(\beta')}{C_2(\alpha)}, \\ C_2(\alpha) &= C_{12}(a)C_2(\beta'). \end{aligned} \tag{2.52}$$

Now, from (2.48, 2.49, 2.50),

$$
\begin{aligned}
T_{12}(a) = S_1(b)T_2(\alpha) &= T_1(c)S_2(\beta'),\\
S_1(c)C_2(\beta')T_2(\alpha) &= T_1(c)S_2(\beta'),\\
T_2(\beta') &= C_1(c)T_2(\alpha).
\end{aligned} \tag{2.53}
$$

Finally, multiply last equations (2.52) and (2.53), obtain with (2.46):

$$
\begin{aligned}
C_1(c)T_2(\alpha)C_2(\alpha) &= T_2(\beta')C_{12}(a)C_2(\beta'),\\
C_1(c)S_2(\alpha) &= C_{12}(a)C_2(\beta'),\\
C_{12}(a)C_1(b)S_2(\alpha) &= C_{12}(a)S_2(\beta'),\\
S_2(\beta') &= C_1(b)S_2(\alpha).
\end{aligned} \tag{2.54}
$$

It is necessary to change equations (2.46) and (2.52) so that they don't contain $C(x)$ function.

$$
\begin{aligned}
k_1S_1^2(c) = 1 - C_1^2(c) &= (C_{12}^2(a) + k_1k_2S_{12}^2(a))(C_1^2(b) + k_1S_1^2(b)) - C_{12}^2(a)C_1^2(b)\\
&= k_1k_2S_{12}^2(a)C_1^2(b) + k_1C_{12}^2(a)S_1^2(b) + k_1^2k_2S_{12}^2(a)S_1^2(b),\\
S_1^2(c) &= k_2S_{12}^2(a)C_1^2(b) + C_{12}^2(a)S_1^2(b) + k_1k_2S_{12}^2(a)S_1^2(b)
\end{aligned}
$$

Divide last equation by its $C(x)$ form, obtain:

$$
T_1^2(c) = k_2T_{12}^2(a) + T_1^2(b) + k_1k_2T_{12}^2(a)T_1^2(b). \tag{2.55}
$$

Similarly,

$$
\begin{aligned}
k_2S_2^2(\alpha) = 1 - C_2^2(\alpha) &= (C_{12}^2(a) + k_1k_2S_{12}^2(a))(C_2^2(\beta') + k_2S_2^2(\beta')) - C_{12}^2(a)C_2^2(\beta')\\
&= k_2C_{12}^2(a)S_2^2(\beta') + k_1k_2S_{12}^2(a)C_2^2(\beta') + k_1k_2^2S_{12}^2(a)S_2^2(\beta'),\\
S_2^2(\alpha) &= C_{12}^2(a)S_2^2(\beta') + k_1S_{12}^2(a)C_2^2(\beta') + k_1k_2S_{12}^2(a)S_2^2(\beta')
\end{aligned}
$$

Divide last equation by its $C(x)$ form, obtain:

$$
T_2^2(\alpha) = k_1T_{12}^2(a) + T_2^2(\beta') + k_1k_2T_{12}^2(a)T_2^2(\beta') \tag{2.56}
$$

Note, that for $k_2 = 1$ in equations (2.45 — 2.56) we can change external angle $\beta'$ to internal

Table 2.2: Right triangle equations

| | |
|---|---|
| $T_1(b) = T_1(c)C_2(\alpha)$ | $T_{12}(a) = T_1(c)S_2(\beta')$ |
| $S_{12}(a) = S_1(c)S_2(\alpha)$ | $S_1(b) = S_1(c)C_2(\beta')$ |
| $T_{12}(a) = S_1(b)T_2(\alpha)$ | $S_{12}(a) = T_1(b)T_2(\beta')$ |
| $C_2(\alpha) = C_{12}(a)C_2(\beta')$ | $S_2(\beta') = C_1(b)S_2(\alpha)$ |
| $C_1(c) = C_{12}(a)C_1(b)$ | $T_2(\beta') = C_1(c)T_2(\alpha)$ |

| |
|---|
| $T_1^2(c) = k_2T_{12}^2(a) + T_1^2(b) + k_1k_2T_{12}^2(a)T_1^2(b)$ |
| $T_2^2(\alpha) = k_1T_{12}^2(a) + T_2^2(\beta') + k_1k_2T_{12}^2(a)T_2^2(\beta')$ |

$\beta$, which leads to:

$$\beta = \frac{\pi}{2} - \beta',$$
$$\cos\beta = \sin\beta',$$
$$\sin\beta = \cos\beta',$$
$$\tan\beta = \cot\beta',$$
$$\cot\beta = \tan\beta'.$$

### 2.3.4. Properties of figures and properties of spaces.

It is easy to see, that equations of generic and right triangles have the same form for any homogeneous space. It is possible to rigorously describe this fact.

**Proposition 2.3.3.** *For any right–lined figure $\Omega$ in homogeneous space $\mathbb{B}^n$, if there exists equation:*

$$F(p_1, ..., p_n) = 0,$$

*that relates elements $p_1, ..., p_n$ of this figure, then it is possible to find its form:*

$$H(Tr(p_1), ..., Tr(p_n)) = 0,$$

*which is expressed through functions having the properties:*

- *Function $H$ is algebraic and doesn't depend on space $\mathbb{B}^n$;*
- *Function $Tr(p_i)$ is any one from these three: $C_i(p_i)$, $\sqrt{k_i}S_i(p_i)$, $\sqrt{k_i}T_i(p_i)$, which depends only on type $k_i$ of its argument $p_i$.*

*Proof.* Using equations (2.6, 2.7) it is possible to change any possible variant of function $Tr(p_i)$ to any other without violation of equation algebraicity.

Further, it is easy to see, that all equations of generic triangle (2.16 — 2.29) as well as of right triangle (2.45 — 2.56) satisfy conditions of the proposition.

Any equation, that relates elements of any right–lined figure $\Omega \subset \mathbb{B}^n$, can be obtained in the following way:

1. Perform the necessary two–dimensional sections of figure $\Omega$;
2. Devide two–dimensional section into triangles so, that triangles elements coincide either to figure $\Omega$ elements $p_i$ or to their parts;
3. Find the relations among triangles elements (parts of figure $\Omega$ elements) with either generic or right triangle equations;
4. Find equations relating the elements $p_i$ of figure $\Omega$, having equations of these elements parts. Use here equations (2.8, 2.9), which satisfy conditions of the proposition.
5. Bind the equations together, using the method of substitution.

All enumerated above operations take algebraic equations, transform them in algebraic way, and thus, result in algebraic equation. □

*Remark.* The condition, that the figure is right–lined is essential. Non–right–lined figures may not have such equation form.

**Example.** Consider on plane a circle arc with radius $\rho$ having the angle $\phi$. Its length $l$ can be expressed by equation (it is obtained in section 2.8.2):

$$l = S(\rho)\phi.$$

This equation contains non–trigonometric functions of arguments $l, \phi$, unless the case $k_2 = 0$. In this case, $l = S_{12}(l) = T_{12}(l)$ and $\phi = S_2(\phi) = T_2(\phi)$. But exactly in this case ($\mathbb{B}^2 = \{k_1, 0\}$) the circle arc *is right–lined* figure, which coincides with right triangle with catheti $\rho$ and $l$ (see Figure 2.4 b), and the equation satisfies the condition of proposition:

$$\begin{aligned} l &= S_1(\rho)\phi, \\ T_{12}(l) &= S_1(\rho)T_2(\phi), \\ \sqrt{k_1 0} \cdot T_{12}(l) &= \sqrt{k_1}S_1(\rho)\sqrt{0}T_2(\phi). \end{aligned}$$

Proposition 2.3.3 is interesting because allows to delimit universal properties of figures from the properties of space $\mathbb{B}^n$.

**Corollary.** *The properties described by equation $H$ are properties of figure $\Omega$ and don't depend on the space, the properties described by function $Tr$ are properties of space $\mathbb{B}^n$ and don't depend on figure.*

*Proof.* The function $H$ doesn't depend on space. Hence it describes the properties of figure $\Omega$ that don't depend on space. On the other hand, the properties of figure described by function $Tr$ are related to space. However, the function $Tr$ is the same for all figures $\Omega$. It means, that it doesn't depend on the figure and describes the properties of whole space $\mathbb{B}^n$. □

## 2.4. Motion

### 2.4.1. Multiplication of types.

As it was shown during right triangle equations deduction, transformation $\mathfrak{T}(-a)$ preserves product of vectors. Let's show this transformation is a motion.

**Theorem 2.4.1** (Multiplication of characteritsics). *Transformation*

$$\mathfrak{T}(a) = \begin{pmatrix} C_{12}(a) & 0 & -K_2 S_{12}(a) \\ 0 & 1 & 0 \\ S_{12}(a) & 0 & C_{12}(a) \end{pmatrix}$$

*is motion.*

*Proof.* Let $a$, $b$, $c$, $\alpha$ and $\beta'$ be real numbers related by equations (2.45 — 2.54). Let

$$\mathfrak{M} = \mathfrak{R}_2(\beta')\mathfrak{R}_1(c)\mathfrak{R}_2(-\alpha)\mathfrak{R}_1(-b).$$

Because $\mathfrak{M}$ is product of four motions, it is motion by definition. Let calculate elements $m_{ij}, i, j = \overline{0, 2}$ of this motion matrix:

$$\begin{pmatrix} 1 & 0 & 0 \\ 0 & C_2(\beta') & -k_2 S_2(\beta') \\ 0 & S_2(\beta') & C_2(\beta') \end{pmatrix} \begin{pmatrix} C_1(c) & -k_1 S_1(c) & 0 \\ S_1(c) & C_1(c) & 0 \\ 0 & 0 & 1 \end{pmatrix} \begin{pmatrix} 1 & 0 & 0 \\ 0 & C_2(\alpha) & k_2 S_2(\alpha) \\ 0 & -S_2(\alpha) & C_2(\alpha) \end{pmatrix} \begin{pmatrix} C_1(b) & k_1 S_1(b) & 0 \\ -S_1(b) & C_1(b) & 0 \\ 0 & 0 & 1 \end{pmatrix}$$

$$= \begin{pmatrix} C_1(c) & -k_1 S_1(c) & 0 \\ S_1(c)C_2(\beta') & C_1(c)C_2(\beta') & -k_2 S_2(\beta') \\ S_1(c)S_2(\beta') & C_1(c)S_2(\beta') & C_2(\beta') \end{pmatrix} \begin{pmatrix} C_1(b) & k_1 S_1(b) & 0 \\ -S_1(b)C_2(\alpha) & C_1(b)C_2(\alpha) & k_2 S_2(\alpha) \\ S_1(b)S_2(\alpha) & -C_1(b)S_2(\alpha) & C_2(\alpha) \end{pmatrix}$$

$$
\begin{aligned}
m_{00} &= C_1(b)C_1(c) + k_1S_1(b)S_1(c)C_2(\alpha) = C_1(b)C_1(c) + k_1S_1(b)S_1(c)\frac{T_1(b)}{T_1(c)} \\
&= C_1(b)C_1(c) + k_1S_1(b)T_1(b)C_1(c) = C_{12}(a)C_1^2(b) + k_1C_{12}(a)S_1^2(b) \\
&= C_{12}(a)(C_1^2(b) + k_1S_1^2(b)) = C_{12}(a),
\end{aligned}
$$

$$
\begin{aligned}
m_{10} &= C_1(b)S_1(c)C_2(\beta') - S_1(b)C_1(c)C_2(\alpha)C_2(\beta') - k_2S_1(b)S_2(\alpha)S_2(\beta') \\
&= C_2(\beta')(C_1(b)S_1(c) - S_1(b)C_1(c)C_2(\alpha)) - k_2S_1(b)S_2(\alpha)S_2(\beta') \\
&= C_2(\beta')(C_1(b)S_1(c) - C_{12}(a)S_1(b)C_1(b)C_2(\alpha)) - k_2S_1(b)S_2(\alpha)S_2(\beta') \\
&= C_1(b)(S_1(b) - C_{12}(a)S_1(b)C_2(\alpha)C_2(\beta')) - k_2S_1(b)S_2(\alpha)S_2(\beta') \\
&= C_1(b)S_1(b)(1 - C_{12}(a)C_2(\alpha)C_2(\beta')) - k_2S_1(b)S_2(\alpha)S_2(\beta') \\
&= C_1(b)S_1(b)(1 - C_2^2(\alpha)) - k_2S_1(b)S_2(\alpha)S_2(\beta') \\
&= k_2C_1(b)S_1(b)S_2^2(\alpha) - k_2S_1(b)S_2(\alpha)S_2(\beta') \\
&= k_2S_1(b)S_2(\alpha)(C_1(b)S_2(\alpha) - S_2(\beta')) = k_2S_1(b)S_2(\alpha)\cdot 0 = 0,
\end{aligned}
$$

$$
\begin{aligned}
m_{20} &= C_1(b)S_1(c)S_2(\beta') - S_1(b)C_1(c)C_2(\alpha)S_2(\beta') + S_1(b)S_2(\alpha)C_2(\beta') \\
&= S_2(\beta')(C_1(b)S_1(c) - S_1(b)C_1(c)C_2(\alpha)) + S_1(b)S_2(\alpha)C_2(\beta') \\
&= S_2(\beta')\left(C_1(b)S_1(c) - S_1(b)C_1(c)\frac{T_1(b)}{T_1(c)}\right) + S_1(b)S_2(\alpha)C_2(\beta') \\
&= \frac{S_2(\beta')}{C_1(b)S_1(c)}(C_1^2(b)S_1^2(c) - S_1^2(b)C_1^2(c)) + S_1(b)S_2(\alpha)C_2(\beta') \\
&= \frac{S_2(\beta')}{C_1(b)S_1(c)}(S_1^2(c) + k_1S_1^2(b)S_1^2(c) - S_1^2(b) - k_1S_1^2(b)S_1^2(c)) + S_1(b)S_2(\alpha)C_2(\beta') \\
&= \frac{S_1(c)S_2(\beta')}{C_1(b)} - \frac{T_1(b)S_1(b)S_2(\beta')}{S_1(c)} + S_1(b)S_2(\alpha)C_2(\beta') \\
&= \frac{S_1(b)S_2(\beta')}{C_1(b)C_2(\beta')} - \frac{T_1(b)S_1(b)S_2(\beta')C_2(\beta')}{S_1(b)} + S_1(b)C_2(\beta')\frac{S_2(\beta')}{C_1(b)} \\
&= T_1(b)T_2(\beta') - T_1(b)S_2(\beta')C_2(\beta') + T_1(b)S_2(\beta')C_2(\beta') = T_1(b)T_2(\beta') = S_{12}(a),
\end{aligned}
$$

$$\begin{aligned}
m_{01} &= k_1 S_1(b)C_1(c) - k_1 C_1(b)S_1(c)C_2(\alpha) = k_1 S_1(b)C_1(c) - k_1 C_1(b)S_1(c)\frac{T_1(b)}{T_1(c)} \\
&= k_1 S_1(b)C_1(c) - k_1 S_1(b)C_1(c) = 0,
\end{aligned}$$

$$\begin{aligned}
m_{11} &= k_1 S_1(b)S_1(c)C_2(\beta') + C_1(b)C_1(c)C_2(\alpha)C_2(\beta') + k_2 C_1(b)S_2(\alpha)S_2(\beta') \\
&= C_2(\beta')(k_1 S_1(b)S_1(c) + C_1(b)C_1(c)C_2(\alpha)) + k_2 C_1(b)S_2(\alpha)S_2(\beta') \\
&= C_2(\beta')\left(k_1 S_1(b)S_1(c) + C_1(b)C_1(c)\frac{T_1(b)}{T_1(c)}\right) + k_2 C_1(b)S_2(\alpha)S_2(\beta') \\
&= C_2(\beta')\left(k_1 S_1(b)S_1(c) + \frac{S_1(b)C_1^2(c)}{S_1(c)}\right) + k_2 C_1(b)S_2(\alpha)S_2(\beta') \\
&= \frac{S_1(b)C_2(\beta')}{S_1(c)}(k_1 S_1^2(c) + C_1^2(c)) + k_2 C_1(b)S_2(\alpha)S_2(\beta') \\
&= C_2(\beta')C_2(\beta') + k_2 C_1(b)S_2(\beta')\frac{S_2(\beta')}{C_1(b)} = C_2^2(\beta') + k_2 S_2^2(\beta') = 1,
\end{aligned}$$

$$\begin{aligned}
m_{21} &= k_1 S_1(b)S_1(c)S_2(\beta') + C_1(b)C_1(c)C_2(\alpha)S_2(\beta') - C_1(b)S_2(\alpha)C_2(\beta') \\
&= S_2(\beta')(k_1 S_1(b)S_1(c) + C_1(b)C_1(c)C_2(\alpha)) - C_1(b)S_2(\alpha)C_2(\beta') \\
&= S_2(\beta')\left(k_1 S_1(b)S_1(c) + C_1(b)C_1(c)\frac{T_1(b)}{T_1(c)}\right) - C_1(b)C_2(\beta')\frac{S_2(\beta')}{C_1(b)} \\
&= \frac{S_1(b)S_2(\beta')}{S_1(c)}(k_1 S_1^2(c) + C_1^2(c)) - C_2(\beta')S_2(\beta') = C_2(\beta')S_2(\beta') - C_2(\beta')S_2(\beta') = 0,
\end{aligned}$$

$$m_{02} = -k_1 k_2 S_1(c)S_2(\alpha) = -K_2 S_{12}(a),$$

$$\begin{aligned}
m_{12} &= k_2 C_1(c)S_2(\alpha)C_2(\beta') - k_2 C_2(\alpha)S_2(\beta') = k_2\frac{T_2(\beta')}{T_2(\alpha)}S_2(\alpha)C_2(\beta') - k_2 C_2(\alpha)S_2(\beta') \\
&= k_2 C_2(\alpha)S_2(\beta') - k_2 C_2(\alpha)S_2(\beta') = 0,
\end{aligned}$$

$$m_{22} = k_2 C_1(c) S_2(\alpha) S_2(\beta') + C_2(\alpha) C_2(\beta') = k_2 \frac{T_2(\beta')}{T_2(\alpha)} S_2(\alpha) S_2(\beta') + C_2(\alpha) C_2(\beta')$$
$$= k_2 \frac{C_2(\alpha) S_2^2(\beta')}{C_2(\beta')} + C_2(\alpha) C_2(\beta') = \frac{C_2(\alpha)}{C_2(\beta')} (k_2 S_2^2(\beta') + C_2^2(\beta') = \frac{C_2(\alpha)}{C_2(\beta')} = C_{12}(a).$$

The motion:

$$\mathfrak{M} = \begin{pmatrix} m_{00} & m_{01} & m_{02} \\ m_{10} & m_{11} & m_{12} \\ m_{20} & m_{21} & m_{22} \end{pmatrix} = \begin{pmatrix} C_{12}(a) & 0 & -K_2 S_{12}(a) \\ 0 & 1 & 0 \\ S_{12}(a) & 0 & C_{12}(a) \end{pmatrix} = \mathfrak{T}(a).$$

So, transformation $\mathfrak{T}(a)$, and its inverse, $\mathfrak{T}(-a)$, are motions. □

**Corollary.** *Transformations:*

$$\mathfrak{R}_{ij}(\varphi) = \begin{pmatrix} 1 & \dots & 0 & \dots & 0 & \dots & 0 \\ \vdots & \ddots & \vdots & \ddots & \vdots & \ddots & \vdots \\ 0 & \dots & C_{i+1\dots j}(\varphi) & \dots & -K_{ij} S_{i+1\dots j}(\varphi) & \dots & 0 \\ \vdots & \ddots & \vdots & \ddots & \vdots & \ddots & \vdots \\ 0 & \dots & S_{i+1\dots j}(\varphi) & \dots & C_{i+1\dots j}(\varphi) & \dots & 0 \\ \vdots & \ddots & \vdots & \ddots & \vdots & \ddots & \vdots \\ 0 & \dots & 0 & \dots & 0 & \dots & 1 \end{pmatrix} \tag{2.57}$$

*where $K_{ij} = k_{i+1} \cdot \dots \cdot k_j$ are motions.*

**Definition 2.4.1** (Rotation, translation)**.** We call motions $\mathfrak{R}_{ij}(\varphi)$ *rotations* in space $\mathbb{B}^n$. Particularly, $\mathfrak{R}_{0i}(\varphi)$ we call *translations* $\mathfrak{T}_i(\varphi)$ in space.

**Corollary.** *The order of coordinate vectors in space $\mathbb{B}^n$ is, generally speaking, important.*

*Proof.* As it was shown, translation $\mathfrak{T}_i(\varphi), i = \overline{0,n}$ along coordinate vector $x^i$ has type $K_i = k_1 \cdot \dots \cdot k_i$. If some type $k_m = 0 \quad (1 \le m \le i)$, then the types of all coordinate vectors, $K_i = 0, \quad i = \overline{m,n}$. Otherwise, if $K_{m-1} \neq 0$, then all types $K_i \neq 0, \quad i = \overline{0,m-1}$. It means that two groups of coordinate vectors $\{x^i\}_{i=\overline{0,m-1}}, \{x^i\}_{i=\overline{m,n}}$ are strictly ordered. □

**Corollary.** *In ordered coordinate vector family $\{x^i\}_{i=\overline{0,n}}$, if $0 \le i < j \le n$ then type $K_i$ of vector $x^i$, type $K_j$ of vector $x^j$ and type $K_{ij}$ of rotation $\mathfrak{R}_{ij}(\varphi)$ in plane $x^i x^j$ are related by equation:*

$$K_j = K_i \cdot K_{ij}. \tag{2.58}$$

### 2.4.2. Vector index. Natural product of vectors.

We defined $K_{ij}$ for the case $i < j$. Let generalize this definition. Define:

$$K_{ij} = \begin{cases} 1, & i = j; \\ \prod_{p=i+1}^{j} k_p, & i < j; \\ \frac{1}{K_{ji}}, & i > j. \end{cases}$$

**Definition 2.4.2** ($i$-th product)**.** Define *$i$-th product*, $i = \overline{0, n}$ of vectors $x$ and $y$ as follows:

$$x \odot_i y = \sum_{j=0}^{n} K_{ij} x_j y_j. \tag{2.59}$$

*Remark.* $K_{ij} = \infty$, when $i > j$ and some type $k_p = 0, \quad j < p \leq i$. In this case in order to have finite value, we have to require that respective coordinate $x_j = 0$ or $y_j = 0$. Then the summand $K_{ij} x_j y_j = 0$. For formal correctness we can define $K_{ij}$ for case $i > j$ as $K_{ij} = K_{ji}$. Really, if $K_{ji} = \pm 1$, then $K_{ij} = \frac{1}{K_{ji}} = K_{ji} = \pm 1$. If $K_{ji} = 0$, we can consider that also $K_{ij} = 0$, keeping in mind that respective coordinate $x_j$ or $y_j$ should be zero (it isn't necessary to require that $x_j = 0$ or $y_j = 0$ if $K_{ij} = 0, \quad i < j$). It is more convenient to use the definition for theory and this note for practical calculus.

*Properties.*

- It is easy to see that meta product $\odot$ is zero product $\odot_0$.
- It is easy to check that $\odot_i$ product is invariant with respect to the motions:

  $$x \odot_i y = \mathfrak{M} x \odot_i \mathfrak{M} y$$

  It is sufficient to check (absolutely the same way it was done for $\odot$) invariancy with respect to the main rotations $\mathfrak{R}_p(\varphi)$.
- If $e^i, \quad i = \overline{0, n}$ is some coordinate vector ($e^i_i = 1$, all other coordinates are zero), than for any signature:

  $$e^i \odot_i e^i = 1,$$

  because $K_{ii} = 1$.

**Definition 2.4.3** (Vector index)**.** We say, the number $i$, $\quad 0 \le i \le n$ is *index of vector* $x \in \mathbb{R}^{n+1}$, if:

$$x \odot_i x > 0$$

In this case the notation is $x^i$.

*Remark.* Tere is no guaranty that each vector $x \in \mathbb{R}^{n+1}$ has an index. On the other hand, some vectors may have several indices.

**Lemma 2.4.2.** *If $K_{ij} = 1$ and some vector $x$ has index $i$, then this vector has also index $j$. The converse is also true. If $i \neq j$ are indices of some vector $x$, then $K_{ij} = 1$.*

*Proof.* Having the remark after $\odot_i$ definition, we can write:

$$K_{ij} = \frac{K_j}{K_i}$$

Then:

$$x \odot_i y = \sum_{p=0}^{n} K_{ip} x_p y_p = \sum_{p=0}^{n} \frac{K_p}{K_i} x_p y_p = \frac{K_j}{K_i} \sum_{p=0}^{n} \frac{K_p}{K_j} x_p y_p = K_{ij} x \odot_j y$$

That is, $x \odot_i y = x \odot_j y$ if and only if $K_{ij} = 1$. It is true also for $x = y$. It means that for any two indices $i \neq j$, $x \odot_i x = x \odot_j x > 0$ if and only if $K_{ij} = 1$. □

At the same time we proved the following

**Corollary.** *If some vector $x$ has more than one index, for example $i \neq j$, then its "scalar square" $x \odot_i x = x \odot_j x$ has the same value regardless of product index choice from vector indices.*

**Theorem 2.4.3.** *Vector index doesn't depend on space basis choice.*

*Proof.* Let coordinate vector family $\{f^i\}_{i=\overline{0,n}}$ result from $\{e^i\}_{i=\overline{0,n}}$ on motion $\mathfrak{M}$: $f^i = \mathfrak{M} e^i, i = \overline{0,n}$. Let vector $x$ have coordinates $x = (a_0 : \dots : a_n)$ in basis $\{e^i\}$ and coordinates $x = (b_0 : \dots : b_n)$ in basis $\{f^i\}$. Then vectors $a$ and $b$ of coordinates of $x$ are related by equality: $b = Ma$, where $M$ is the matrix of motion $\mathfrak{M}$. In basis $\{e^i\}$ scalar square $x \odot_p x = a \odot_p a$. In basis $\{f^i\}$ scalar square $x \odot_p x = b \odot_p b = Ma \odot_p Ma$. Since $\mathfrak{M}$ is the motion, it preserves vector product, that is $b \odot_p b = Ma \odot_p Ma = a \odot_p a$. In other words, in both bases, $\{e^i\}$ and $\{f^i\}$, the value $x \odot_p x$ is the same. □

**Proposition 2.4.4.** *Unless otherwise is specified, consider index of vector $x$ as the least index (if there exist any).*

**Corollary.** *Vector index is the property of vector, not of its representation in the coordinate system.*

The product indices of interest are those, which correspond to indices of multiplied vectors themselves.

**Definition 2.4.4** (Natural product)**.** Define *natural* product of vectors $x^i$ and $y^j$, the one that has index $k = \min(i, j)$:

$$x^i \odot y^j = x^i \odot_{\min(i,j)} y^j.$$

*Remark.* For space points, the natural product is zero product.

Further, if no product index is specified, the natural product is referred, not zero product.

**Proposition 2.4.5.** *Scalar square of any vector, calculated by means of natural product, is nonnegative.*

*Proof.* Consider vector $x$. If $x \odot_0 x \geq 0$, it satisfies the proposition. Let $x \odot_0 x = \sum_{i=0}^{n} K_{0i} x_i^2 < 0$. In this case there exists at least one negative summand. Let the first negative summand be $K_{0p} x_p^2 < 0$. $K_{0p} = -1 < 0$, because $x_p^2 \geq 0$. In this case,

$$x \odot_0 x < 0,$$
$$x \odot_p x = \frac{1}{K_{0p}} x \odot_0 x = -x \odot_0 x > 0.$$

It means that index of vector $x$ equals to $p$ and $x^p \odot x^p > 0$. □

### 2.4.3. Generalized orthogonal matrix.

**Definition 2.4.5** (Normalized and orthogonal vectors, generalized orthogonal matrix)**.** We call vector $x^i$ in $\mathbb{R}^{n+1}$ with given signature, *normalized*, if $x^i \odot x^i = 1$. We call vectors $x^i$, $y^j$ *orthogonal*, if $x^i \odot y^j = 0$. Define square matrix $M$ of size $n + 1$ *generalized orthogonal* if all its columns $m^j$ have index $j$, are normalized and any two columns are orthogonal. Further, for the generalized orthogonal matrices acting in metaspace the shorter term is also used: *GM-orthogonal.*

**Proposition 2.4.6.** *All main rotation matrices are GM-orthogonal.*

*Proof.* Consider matrix $M$ of rotation $\mathfrak{R}_j(\varphi)$. It has only two columns $m^{j-1} = (0 : \ldots : 0 : C_j(\varphi) : S_j(\varphi) : 0 : \ldots : 0)$ and $m^j = (0 : \ldots : 0 : -k_j S_j(\varphi) : C_j(\varphi) : 0 : \ldots : 0)$, different from the identity matrix. Check normalization and orthogonality condition for these columns.

Other cases are trivial:

$$\begin{aligned}
m^{j-1} \odot_{j-1} m^{j-1} &= \sum_{i=0}^{n} K_{j-1\,i} m_{i\,j-1}^2 = K_{j-1\,j-1} C_j^2(\varphi) + K_{j-1\,j} S_j^2(\varphi) \\
&= K_{j-1\,j-1}(C_j^2(\varphi) + k_j S_j^2(\varphi)) = 1 \cdot 1 = 1, \\
m^j \odot_j m^j &= \sum_{i=0}^{n} K_{ji} m_{ij}^2 = K_{j\,j-1} k_j^2 S_j^2(\varphi) + K_{jj} C_j^2(\varphi) \\
&= K_{jj}(k_j S_j^2(\varphi) + C_j^2(\varphi)) = 1 \cdot 1 = 1, \\
m^{j-1} \odot_{j-1} m^j &= \sum_{i=0}^{n} K_{j-1\,i} m_{i\,j-1} m_{ij} = K_{j-1\,j-1} C_j(\varphi)(-k_j S_j(\varphi)) + K_{j-1\,j} S_j(\varphi) C_j(\varphi) \\
&= -K_{j-1\,j} C_j(\varphi) S_j(\varphi) + K_{j-1\,j} C_j(\varphi) S_j(\varphi) = 0.
\end{aligned}$$

□

**Lemma 2.4.7.** *Product of two GM-orthogonal matrices is GM-orthogonal matrix.*

*Proof.* Let $X, Y$ be two GM-orthogonal matrices. It means that $X$ is composed of columns $\{x^0, ..., x^n\}$ and $Y$ is composed of columns $\{y^0, ..., y^n\}$, for which $x^i \odot x^j = y^i \odot y^j = \delta_{ij}$ for all $i, j = \overline{0, n}$, where $\delta_{ij} = 1, \quad i = j$ and $\delta_{ij} = 0, \quad i \neq j$. Let $Z = XY$ with elements $z_{ij} = \sum_{p=0}^{n} x_{ip} y_{pj}$. Let $z^i$ and $z^j$ be two columns of $Z$ and $i \leq j$. Calculate:

$$\begin{aligned}
z^i \odot z^j &= \sum_{p=0}^{n} K_{ip} z_{pi} z_{pj} = \sum_{p=0}^{n} K_{ip} \left( \sum_{m_1=0}^{n} x_{pm_1} y_{m_1 i} \right) \left( \sum_{m_2=0}^{n} x_{pm_2} y_{m_2 j} \right) \\
&= \sum_{p=0}^{n} K_{ip} \sum_{m_1=0}^{n} \sum_{m_2=0}^{n} x_{pm_1} x_{pm_2} y_{m_1 i} y_{m_2 j} = \sum_{m_1=0}^{n} \sum_{m_2=0}^{n} y_{m_1 i} y_{m_2 j} \sum_{p=0}^{n} K_{ip} x_{pm_1} x_{pm_2} \\
&= \sum_{m_1=0}^{n} \sum_{m_2=0}^{n} y_{m_1 i} y_{m_2 j} x^{m_1} \odot_i x^{m_2} = \sum_{m_1=0}^{n} \sum_{m_2=0}^{n} y_{m_1 i} y_{m_2 j} \frac{K_{\min(m_1, m_2)}}{K_i} x^{m_1} \odot_{\min(m_1, m_2)} x^{m_2} \\
&= \sum_{m_1=0}^{n} \sum_{m_2=0}^{n} y_{m_1 i} y_{m_2 j} K_{i\,\min(m_1, m_2)} \delta_{m_1 m_2} = \sum_{m=0}^{n} y_{mi} y_{mj} K_{im} = \delta_{ij}.
\end{aligned}$$

□

**Proposition 2.4.8.** *The type $K_{ij}$ divides element $m_{ij}$ of some GM-orthogonal matrix $M$ (from $K_{ij} = 0$ follows that $m_{ij} = 0$).*

*Proof.* The column $m^j$ of GM-orthogonal matrix is normalized. It means that:

$$m^j \odot m^j = \sum_{i=0}^{n} K_{ji} m_{ij}^2 = \sum_{i=0}^{n} \frac{m_{ij}^2}{K_{ij}} = 1.$$

Each summand $\frac{m_{ij}^2}{K_{ij}}$ is finite, that is $K_{ij}$ divides $m_{ij}^2$, and thus also $K_{ij}$ divides $m_{ij}$. □

**Lemma 2.4.9.** *Rows $m^i$ and $m^j$ of a GM-orthogonal matrix $M$ satisfy equation (row normalization when $i = j$ or orthogonality when $i \neq j$):*

$$K_{\max(i,j)} \sum_{p=0}^{n} \frac{m_{ip} m_{jp}}{K_p} = \delta_{ij}$$

*Proof.* Let $M$ be some GM-orthogonal matrix. Construct matrix $W$ with elements $w_{ij} = \frac{m_{ij}}{\sqrt{K_{ij}}}$. Matrix $W$ is orthogonal (it may be complex, in this case $W$ isn't unitary, but exactly orthogonal). Really, product of columns $w^i \cdot w^j = \delta_{ij}$. Let $i < j$:

$$w^i \cdot w^i = \sum_{p=0}^{n} w_{pi}^2 = \sum_{p=0}^{n} \left( \frac{m_{pi}}{\sqrt{K_{pi}}} \right)^2 = \sum_{p=0}^{n} \frac{m_{pi}^2}{K_{pi}} = \frac{1}{K_i} \sum_{p=0}^{n} K_p m_{pi}^2 = m^i \odot m^i = 1,$$

$$w^i \cdot w^j = \sum_{p=0}^{n} w_{pi} w_{pj} = \sum_{p=0}^{n} \frac{m_{pi}}{\sqrt{K_{pi}}} \frac{m_{pj}}{\sqrt{K_{pj}}} = \sum_{p=0}^{n} \frac{m_{pi}}{\sqrt{K_{pi}}} \frac{m_{pj}}{\sqrt{K_{pi} K_{ij}}} = \frac{1}{\sqrt{K_{ij}}} \sum_{p=0}^{n} \frac{m_{pi} m_{pj}}{K_{pi}}$$

$$= \frac{1}{\sqrt{K_{ij}}} \sum_{p=0}^{n} K_{ip} m_{pi} m_{pj} = \frac{1}{\sqrt{K_{ij}}} (m^i \odot m^j) = \frac{1}{\sqrt{K_{ij}}} \cdot 0 = 0,$$

if $K_{ij} \neq 0$. Otherwise, $K_{ij} = 0$ follows $k_q = 0, \quad i < q \leq j$. In this case:

$$\begin{cases} w_{pj} = \frac{m_{pj}}{\sqrt{K_{pj}}} = \frac{0}{\sqrt{K_{pj}}} = 0, & 0 \leq p < q, \\ w_{pi} = \frac{m_{pi}}{\sqrt{K_{pi}}} = \frac{m_{pi}}{\infty} = 0, & q \leq p \leq n. \end{cases}$$

So,

$$w^i \cdot w^j = \sum_{p=0}^{n} w_{pi} w_{pj} = \sum_{p=0}^{q-1} w_{pi} \cdot 0 + \sum_{p=q}^{n} 0 \cdot w_{pj} = 0.$$

In orthogonal matrix $W$ also product of rows $w^i \cdot w^j = \delta_{ij}$. It means $(i < j)$:

$$K_i \sum_{p=0}^{n} \frac{m_{ip}^2}{K_p} = \sum_{p=0}^{n} \frac{m_{ip}^2}{K_{ip}} = \sum_{p=0}^{n} \left( \frac{m_{ip}}{\sqrt{K_{ip}}} \right)^2 = \sum_{p=0}^{n} w_{ip}^2 = w^i \cdot w^i = 1,$$

$$K_j \sum_{p=0}^{n} \frac{m_{ip} m_{jp}}{K_p} = \sum_{p=0}^{n} \frac{m_{ip} m_{jp}}{K_{jp}} = \sum_{p=0}^{n} \frac{m_{ip}}{\sqrt{K_{jp}}} \frac{m_{jp}}{\sqrt{K_{jp}}} = \sum_{p=0}^{n} \frac{m_{ip}}{\sqrt{K_{ji} K_{ip}}} \frac{m_{jp}}{\sqrt{K_{jp}}}$$

$$= \frac{1}{\sqrt{K_{ji}}} \sum_{p=0}^{n} \frac{m_{ip}}{\sqrt{K_{ip}}} \frac{m_{jp}}{\sqrt{K_{jp}}} = \sqrt{K_{ij}} \sum_{p=0}^{n} w^i \cdot w^j = \sqrt{K_{ij}} \cdot 0 = 0.$$

□

**Lemma 2.4.10.** *For GM-orthogonal matrix $M$ with elements $m_{ij}$ inverse matrix $W = M^{-1}$ is GM-orthogonal and has elements $w_{ij} = \frac{m_{ji}}{K_{ji}}$.*

*Proof.* Let calculate elements $z_{ij}$ of matrix $Z = WM$:

$$z_{ii} = \sum_{p=0}^{n} w_{ip} m_{pi} = \sum_{p=0}^{n} \frac{m_{pi} m_{pi}}{K_{pi}} = \frac{1}{K_i} \sum_{p=0}^{n} K_p m_{pi}^2 = m^i \odot m^i = 1,$$

$$z_{ij} = \sum_{p=0}^{n} w_{ip} m_{pj} = \sum_{p=0}^{n} \frac{m_{pi} m_{pj}}{K_{pi}} = \frac{1}{K_i} \sum_{p=0}^{n} K_p m_{pi} m_{pj} = m^i \odot m^j = 0.$$

□

**Theorem 2.4.11.** *GM-orthogonal matrices form a group. The isometry group (group of motions) of space is subgroup of this group.*

*Proof.* Evidently, identity matrix $I$ is GM-orthogonal. By lemma 2.4.7 product of two GM-orthogonal matrices is also GM-orthogonal matrix. By lemma 2.4.10, inverse matrix of some GM-orthogonal matrix is also GM-orthogonal. It follows, that GM-orthogonal matrices form a group.

By proposition 2.4.6 all main rotations has GM-orthogonal matrices. Any motion can be obtained as product of main rotations, which by lemma 2.4.7 has GM-orthogonal matrix. □

*Remark.* GM-orthogonal matrices for metaspace of a homogeneous spaces play the same role as orthogonal matrices for Euclidean space. Moreover, orthogonal matrices are particular case of GM-orthogonal matrices when signature is $\{1, 1, ..., 1\}$.

The way of finding inverse matrix is inacceptable if some element of signature $k_p = 0$. In this case, by proposition 2.4.8 elements $m_{ij} = 0, \quad i < p \leq j$. Then the equality from lemma 2.4.10, although true, is indeterminate $\frac{0}{0}$. We need another way of finding inverse matrix.

**Lemma 2.4.12.** *If some element $k_m = 0$ of signature $\{k_1, ..., k_n\}$, then GM-orthogonal matrix has form:*

$$W = \begin{pmatrix} A & O \\ B & C \end{pmatrix},$$

*where $A$ is GM-orthogonal matrix of size $m$ with signature $\{k_1, ..., k_{m-1}\}$, $C$ is GM-orthogonal matrix of size $n - m + 1$ with signature $\{k_{m+1}, ..., k_n\}$, $O$ is zero matrix of size $(n - m + 1) \times m$ and $B$ is arbitrary rectangular matrix of size $m \times (n - m + 1)$.*

*Proof.* For columns $w^i, w^j; 0 \le i \le j < m$ have:

$$w^i \odot w^j = \frac{1}{K_i} \sum_{p=0}^{n} K_p w_{pi} w_{pj} = \frac{1}{K_i} \sum_{p=0}^{m} K_p w_{pi} w_{pj} + \frac{1}{K_i} \sum_{p=m+1}^{n} K_p w_{pi} w_{pj}$$
$$= \frac{1}{K_i} \sum_{p=0}^{m} K_p a_{pi} a_{pj} + \frac{1}{K_i} \sum_{i=m+1}^{n} 0 \cdot b_{pi} b_{pj} = \delta_{ij}.$$

It means that matrix $A$ is GM-orthogonal with signature $\{k_1, ..., k_{m-1}\}$ and matrix $B$ is arbitrary.

For columns $w^i, w^j; m \le i \le j \le n$ have:

$$w^i \odot w^j = \frac{1}{K_i} \sum_{p=0}^{n} K_p w_{pi} w_{pj} = \frac{1}{0} \sum_{p=0}^{m} K_p 0 \cdot 0 + \frac{1}{K_{mi}} \sum_{p=m+1}^{n} K_{mp} c_{pi} c_{pj} = \delta_{ij}.$$

It means that matrix $O$ is necessarily zero and matrix $C$ is GM-orthogonal with signature $\{k_{m+1}, ..., k_n\}$.

□

**Lemma 2.4.13.** *If some element $k_m = 0$ of signature $\{k_1, ..., k_n\}$, then inverse matrix $W^{-1}$ has the form:*

$$W^{-1} = \begin{pmatrix} A^{-1} & O \\ -C^{-1} B A^{-1} & C^{-1} \end{pmatrix},$$

*where $A^{-1}$ and $C^{-1}$ are calculated by lemma 2.4.10, if their signature don't contain zero elements, or recursively by lemma 2.4.13 otherwise.*

*Proof.*

$$W^{-1}W = \begin{pmatrix} A^{-1} & O \\ -C^{-1}BA^{-1} & C^{-1} \end{pmatrix} \begin{pmatrix} A & O \\ B & C \end{pmatrix}$$
$$= \begin{pmatrix} A^{-1}A + OB & A^{-1}O + OC \\ -C^{-1}BA^{-1}A + C^{-1}B & -C^{-1}BA^{-1}O + C^{-1}C \end{pmatrix}$$
$$= \begin{pmatrix} I_m & O \\ -C^{-1}B + C^{-1}B & I_{n-m+1} \end{pmatrix} = I_{n+1}.$$

□

**2.4.4. GM-orthogonal matrix decomposition in product of rotations.**

Consider an GM-orthogonal matrix $X$. Search for rotation matrices, whose product gives $X$. The matrix $X\mathfrak{R}_{ij}(\varphi)$ has all columns $x^p$ as $X$, except $i$-th and $j$-th. These columns equal $x'^i = x^i C_{ij}(\varphi) + x^j S_{ij}(\varphi)$ and $x'^j = -K_{ij} x^i S_{ij}(\varphi) + x^j C_{ij}$.

**Algorithm 2.1** (GM-orthogonal matrix decomposition in product of rotations).

1. For rows $r$ from $n$ to 1 do:

   (a) Divide elements of the row $x_{ri}$, $i = \overline{0, r}$ in three categories: having the type $K_{ri}$ equal to 1, 0 and $-1$. Note, that in $i$-th row $i$-th element is from category 1, because its type $K_{ii} = 1$. We will right multiply $X$ by $\mathfrak{R}_{ir}(\varphi)$, $i = \overline{0, r}$ so that in $r$-th row remains one element of category 1 and one element of category $-1$, different from 0. All these rotations are elliptic. For elements $x_{ri}$ and $x_{rj}$ with the same type we can use $\cos\varphi = \frac{x_{ri}}{\sqrt{x_{ri}^2 + x_{rj}^2}}$ and $\sin\varphi = \frac{x_{rj}}{\sqrt{x_{ri}^2 + x_{rj}^2}}$. Moreover, always $x_{rr} \neq 0$. In this case:

$$x'_{ri} = x_{ri}\cos\varphi + x_{rj}\sin\varphi = x_{ri}\frac{x_{ri}}{\sqrt{x_{ri}^2 + x_{rj}^2}} + x_{rj}\frac{x_{rj}}{\sqrt{x_{ri}^2 + x_{rj}^2}}$$
$$= \frac{x_{ri}^2 + x_{rj}^2}{\sqrt{x_{ri}^2 + x_{rj}^2}} = \sqrt{x_{ri}^2 + x_{rj}^2},$$
$$x'_{rj} = -x_{ri}\sin\varphi + x_{rj}\cos\varphi = -x_{ri}\frac{x_{rj}}{\sqrt{x_{ri}^2 + x_{rj}^2}} + x_{rj}\frac{x_{ri}}{\sqrt{x_{ri}^2 + x_{rj}^2}}$$
$$= \frac{-x_{ri}x_{rj} + x_{rj}x_{ri}}{\sqrt{x_{ri}^2 + x_{rj}^2}} = 0.$$

(b) Now have one element of types 1 and $-1$ different from zero (let they be $r$-th and $p$-th), and element of category 1 has greater absolute value than element of category $-1$, because this normalized row satisfies equality $x_{rr}^2 - x_{rp}^2 = 1$. It means that exists $\varphi \in \mathbb{R}$ so that $\cosh\varphi = \frac{x_{rr}}{\sqrt{x_{rr}^2 - x_{rp}^2}}$ and $\sinh\varphi = \frac{-x_{rp}}{\sqrt{x_{rr}^2 - x_{rp}^2}}$. In this case:

$$
\begin{aligned}
x'_{rp} &= x_{rp}\cosh\varphi + x_{rr}\sinh\varphi = x_{rp}\frac{x_{rr}}{\sqrt{x_{rr}^2 - x_{rp}^2}} + x_{rr}\frac{-x_{rp}}{\sqrt{x_{rr}^2 - x_{rp}^2}} \\
&= \frac{x_{rp}x_{rr} - x_{rr}x_{rp}}{\sqrt{x_{rr}^2 - x_{rp}^2}} = 0, \\
x'_{rr} &= x_{rp}\sinh\varphi + x_{rr}\cosh\varphi = x_{rp}\frac{-x_{rp}}{\sqrt{x_{rr}^2 - x_{rp}^2}} + x_{rr}\frac{x_{rr}}{\sqrt{x_{rr}^2 - x_{rp}^2}} \\
&= \frac{-x_{rp}^2 + x_{rr}^2}{\sqrt{x_{rr}^2 - x_{rp}^2}} = \sqrt{x_{rr}^2 - x_{rp}^2} = \pm 1.
\end{aligned}
$$

(c) For category 0 exist parabolic rotations that preserve element of category 1 ($x_{rr}$) and elements of category 0 transform to 0. For this, if one such element is in $q$-th columns, $\varphi = -\frac{x_{rq}}{x_{rr}}$:

$$
\begin{aligned}
x'_{rq} &= x_{rq} \cdot 1 + x_{rr}\varphi = x_{rq} + x_{rr} \cdot \left(-\frac{x_{rq}}{x_{rr}}\right) = 0, \\
x'_{rr} &= 0 \cdot x_{rq}\varphi + x_{rr} \cdot 1 = x_{rr} = \pm 1.
\end{aligned}
$$

Remaining nonzero element $x_{rr}$ equals to 1 or $-1$, because the last row is normalized.

(d) We can consider the first $r$ columns as having $r$ elements (the last $n - r + 1$ equal to zero). They form GM-orthogonal matrix of size $r$. The last $(n-r+1)$ columns (without their last elements) are GM-orthogonal to the first $r$ columns. Because these columns have $r$ elements, all columns from $(r + 1)$-th to $(n - r + 1)$-th are necessarily vanish (except their last $n - r + 1$ elements).

2. At this phase we can consider the obtained matrix as having the size $r$ instead of $r + 1$ and repeat the process for it. At the end obtain matrix $E$ whose elements on main diagonal equal to 1 or $-1$ and all other elements equal to 0. It is reflexion matrix on point or line or plane or hyperplane. We get equality: $X \prod_{j=1}^{q} M_j = E$. It's easy to see that $X = E \prod_{j=q}^{1} M_j^{-1}$.

The matrix $E$ can't result from product of rotations (the described algorithm would detect these rotations). In order to consider motions of $\mathbb{B}^n$ equivalent to GM-orthogonal matrices, we need to declare transformations that have matrices $E$ also motions (they preserve product of vectors). These motions, however, are not proper (there is no continous parametrization $M(\varphi)$ on segment $[0, 1]$ so that $M(0) = I$, $M(1) = E$ and all $M(\varphi)$ are GM-orthogonal for any $\varphi \in [0, 1]$).

*Remark.* Having the fact, that for any vector $x \in \mathbb{R}^{n+1}$, vector $-x = x$, the matrix $-I = I$. Really, $-Ix = -x = x$. It means that for any motion matrix $M$, the matrix $-M$ represents the same motion. To be certain, from two matrices $M$ and $-M$ choose one, that in column zero, $m^0$, has positive the first nonzero element. It is always possible, because its normality condition $m^0 \odot m^0 = 1$ means it can't have all zero elements at the same time.

**Proposition 2.4.14.** *Sinse all matrices $M_j$ in decomposition of $X$ have determinant equal to* 1 *and determinant of $E$ equal to $\pm 1$, determinant of any motion matrix $X$ equals to $\pm 1$.*

**Definition 2.4.6** (Imporper motion). Define *motion of genus II* as motions that have matrix $X$ with determinant $-1$. Define *improper motions* as motions that have matrix $X = E \prod_{j=q}^{1} M_j^{-1}$ and matrix $E$ contains at least one negative element.

Evidently, all motions of genus II are improper, however these two notions are, generally, not equivalent. Speaking about homogeneous spaces, the notion of motions of genus II plays little role unlike the improper motions notion.

**Lemma 2.4.15.** *Proper motions $\mathfrak{X}$ have continuos parameterization $\mathfrak{X}(p)$, $p \in [0, 1]$, so that $\mathfrak{X}(0) = \mathfrak{I}$ is identity motion, $\mathfrak{X}(1) = \mathfrak{X}$ is initial motion and all $\mathfrak{X}(p)$ are proper motions.*

*Proof.* Using algorithm 2.1 of GM-orthogonal matrix decomposition in product of rotations, represent motion matrix as $X = E \prod_{j=q}^{1} M_j^{-1} = \prod_{j=q}^{1} M_j^{-1}$, since $E = I$ for proper motions. Each matrix $M_j$ is rotation with some angle $\varphi_j$. These rotations can be parameterized $M_j(p), p \in [0, 1]$ if, instead of angle $\varphi_j$, use angle $p\varphi_j$. The motion matrix parameterization $X(p) = \prod_{j=q}^{1} M_j^{-1}(p), p \in [0, 1]$.

Because $0 \cdot \varphi_j = 0$, all $M_j(0) = I$, and then $X(0) = \prod_{j=q}^{1} M_j^{-1}(0) = \prod_{j=q}^{1} I = I$. When $p = 1$ all rotation angles equal to initial ones, therefore $X(1) = X$. Finally, note that for any value of parameter $p$ the multiplied matrices are rotations, it follows that their product is motion. This motion is proper, because for any value of $p$ the matrix $E = I$. □

**Theorem 2.4.16.** *The number of freedom degrees in homogeneous space $\mathbb{B}^n$ equals to $\frac{n(n+1)}{2}$.*

*Proof.* In decomposition algorithm 2.1 of GM-orthogonal matrix of size $n + 1$ (that is, of motion matrix) in rotations, we constructed as many rotations as many elements has the matrix below the main diagonal (at each step we constructed one rotation for each element from diagonal left side). So, dimension of motion space, that is the number of freedom degrees of $\mathbb{B}^n$, equals to $\frac{n(n+1)}{2}$. □

**Corollary.** *The number of freedom degrees of space* $\mathbb{B}^n$ *depend on its dimension and doesn't depent on its signature.*

### 2.4.5. Equivalence of coordinate axes.

As follows from theorem 2.4.1 on product of types, the order of coordinate axes matters. We can strictly formalize this.

**Definition 2.4.7** (Equivalent, interchangeable and non-interchangeable vectors)**.** We say two vectors with indices $i, j$ are *non-interchangeable*, if there is no homogeneous space isomorphic to given in which corresponding two vectors have indices $j, i$. We say vectors are *interchangeable*, if such isomorphic homogeneous space exists. We say vectors are *equivalent*, if space contains a motion that interchanges them (up to sign).

*Remark.* Space isomorphism in case of interchangeable vectors may be automorphism. It is important that the automorphism is not motion.

**Example.** On Minkowski plane with signature $\{0, -1\}$ vectors $e^1$ and $e^2$ are interchangeable. Corresponding automorphism: $Aut : \mathbb{R}^3 \to \mathbb{R}^3$ is defined as $Aut((x_0 : x_1 : x_2)) = (x_0 : x_2 : x_1)$. This automorphism matrix is

$$\begin{pmatrix} 1 & 0 & 0 \\ 0 & 0 & 1 \\ 0 & 1 & 0 \end{pmatrix}$$

This matrix isn't GM-orthogonal, the automorphism isn't motion, and vectors $e^1$ and $e^2$ are not equivalent.

**Theorem 2.4.17** (coordinate vectors equivalence)**.** *Two vectors* $e^i, e^j, \quad i < j$ *are:*

- *non-interchangeable, if* $K_{ij} = 0$*;*
- *interchangeable, if* $K_{ij} = -1$*;*
- *equivalent, if* $K_{ij} = 1$.

*Proof.* As it was shown, vector $e^i$ type equals to $K_i$, $i = \overline{0, n}$. Each next type $K_i$ results from previous $K_{i-1}$ by multiplication by corresponding type element $k_i$ of space signature: $K_i = K_{i-1}k_i$. In this way, elements of space signature indicate how coordinate vector type changes on switching to the next dimension. Each zero element of signature $k_i = 0$ divides coordinate vectors in non-interchangeable groups. In this case any vector $e^i$ of previous group has lower index than any vector $e^j$ of the next group. This is true regardless of space and follows that $e^i$ and $e^j$ are non-interchangeable. It happens if and only if between groups exists a zero element $k_p = 0, i < p \leq j$, and thus $K_{ij} = 0$.

Inside some group, between any two vectors $e^i, e^j$ the type $K_{ij} = \pm 1$. If $K_{ij} = -1$, then $K_i = -K_j$. In this case exists isomorphic space, where corresponding vectors have indices $j, i$. This space signature differs from initial space signature only in this group, for which elements can be calculated as $k_i = \frac{K_i}{K_{i-1}}$. It is always possible inside a group, because for any $i, j$, $K_{ij} \neq 0$, and thus exists $K_{ji} = \frac{1}{K_{ij}} = K_{ij} = \pm 1$. It means that vectors $e^i$ and $e^j$ are interchangeable. When $K_i = -K_j = 0$, then isomorphism is actually automorphism, when $K_i = -K_j = \pm 1$, then isomorphism is not automorphism.

Finally, if $K_{ij} = 1$, then exists elliptic motion $\mathfrak{R}_{ij}\left(\frac{\pi}{2}\right)$, that transforms $e^i \to e^j$, and $e^j \to -e^i$. In this case vectors $e^i$ and $e^j$ are equivalent. □

*Remark.* When vector $x^i$ is equivalent (interchangeable) with vector $y^j$, then, evidently, also $y^j$ is equivalent (interchangeable) with $x^i$, because from $K_{ij} = 1$ ($K_{ij} = -1$) follows $K_{ji} = \frac{1}{K_{ij}} = 1$ ($K_{ji} = \frac{1}{K_{ij}} = -1$). But when vectors $x^i, y^j$ are non-interchangeable, their relation is "in one sense", because in this case their order matters.

## 2.5. Lineal

### 2.5.1. Planes and lineals. Their signature.

Having $\mathbb{B}^m \subset \mathbb{R}^{m+1}$, all $m$-dimensional planes of $\mathbb{B}^n$ ($m < n$) lie in $(m+1)$–dimensional plane of $\mathbb{R}^{n+1}$ with one global restriction for points $X \in \mathbb{R}^{m+1} \subset \mathbb{R}^{n+1}$: $X \odot X = 1$. Leaving this restriction aside (since it doesn't change under the motions), we can consider $m$-dimensional planes $\mathbb{B}^m$ as $(m+1)$–dimensional planes $\mathbb{R}^{m+1}$.

By definition, $m$-dimensional planes $L^m$ result from subspace $\mathbb{B}^m$ using all possible motions. The basis of $\mathbb{R}^{m+1}$ is composed of the first $m+1$ coordinate vectors $e^i$, each of dimension $n+1$. Multiplying left the basis matrix of $\mathbb{B}^m$ by all possible GM-orthogonal matrices, obtain basis matrix of $L^m$ as the first $m+1$ columns of GM-orthogonal matrix. Being a subspace, the signature of $\mathbb{B}^m$ contains the first $m$ elements of $\mathbb{B}^n$ signature.

What happens if we take arbitrary $m+1$ columns of GM-orthogonal matrix as the basis? Let column indices be $i_0, i_1, ...i_m$ and $i_p = \overline{0, n}$. Evidently, the motions that preserve this figure affect either only these columns (internal motions of figure) or affect none of them (motions of $\mathbb{B}^n$ that preserve all points of figure). In this way, the type of $p$-th basis vector of figure equals to $K'_p = K_{i_p}, p = \overline{0, m}$. Then the elements $k'_p = \frac{K'_p}{K'_{p-1}} = \frac{K_{i_p}}{K_{i_{p-1}}} = K_{i_{p-1} i_p}$, $p = \overline{1, m}$ form figure signature. This signature, generally, is not the subspace signature. These figures are, generally, not planes.

**Definition 2.5.1** (Lineal). We call *lineals* any intersection between sphere $\mathbb{B}^n$ and the linear span of a family of vectors of the metaspace.

*Remark.* All planes are lineals. But not all lineals are planes.

*Remark.* Obviously, congruent lineals have equal signatures. However, not all lineals with equal signatures are congruent.

**Example.** Consider Minkowski plane with signature $\{0, -1\}$. Its first coordinate line signature (the lineal constructed on vectors $\{e^0, e^1\}$) is $\{K_1 = k_1 = 0\}$. The signature of its second coordinate line (constructed on vectors $\{e^0, e^2\}$) also is $\{K_2 = k_1 k_2 = 0 \cdot (-1) = 0\}$. However, there is no motion that transforms the first lineal to the second. Thus, they are not congruent.

It may happen, that some lineal have $K_0' \neq 1$. In this case, in coordinate matrix of its basis there is no 0-th column, which in turn means, that lineal has no space points. This lineal doesn't intersect space sphere and has no image.

**Definition 2.5.2** (Improper lineal)**.** Define a lineals *improper*, if it has no points.

Because improper lineals have no points, they have no information about their position, still they have information about their orientation.

### 2.5.2. Projection of vector on lineal and on its orthogonal complement.

**Definition 2.5.3** (Projection of vector on lineal)**.** Define $v'$ *projection* of vector $v$ on linal $L^m$, if $v' \in L^m, v'' = v - v' \perp L^m$. Define $v''$ *projection* of vector $v$ on orthogonal complement of lineal $L^m$.

**Lemma 2.5.1.** *Projection $v'$ of vector $v$ on lineal $L^m$ can be computed as:*

$$v' = \sum_{i=0}^{m} (v \odot_i l^i) l^i,$$

*where vector family $\{l^i\}_{i=\overline{0,m}}$ is orthonormal basis of $L^m$.*

*Proof.* Obviously, $v' \in L_m$.

Consider:

$$\begin{aligned} v'' \odot_j l^j &= (v - v') \odot l^j = v \odot_j l^j - \sum_{i=0}^{m} (v \odot_i l^i)(l^i \odot_j l^j) \\ &= v \odot_j l^j - \sum_{i=0}^{m} (v \odot_i l^i) \delta_{ij} = v \odot_j l^j - v \odot_j l^j = 0 \end{aligned}$$

for all $j = \overline{0, m}$. □

*Remark.* If index $p$ of vector $v''$ is less than $j$, then:

$$v'' \odot_p l^j = K_{pj} v'' \odot_j l^j = 0,$$

and orthogonality really has place.

*Remark.* When some $K_p = 0$, then for $i < p \leq j$, expression:

$$l^i \odot_j l^j = \frac{l^i \odot_i l^j}{K_{ij}} = \frac{0}{0}$$

has indeterminate value. This happens when one vector is orthogonal to all other vectors. In this case the perpendicular is determined ambiguously. Yet, any value of this expression, for example 0 gives correct result, which corresponds to some perpendicular vector.

### 2.5.3. Orthonormalization of the vector family.

Consider a vector family $\{x^i\}_{i=\overline{0,m}}$. In order to orthonormalize it, consider the following algorithm:

**Algorithm 2.2** (Vector family orthonormalization)**.**

1. Take the vector with lowest index $x$ and normalize it:

   $$y = \frac{1}{\sqrt{x \odot x}} x.$$

   Add $y$ to family of orthonormal vectors $\{y^i\}$ and remove $x$ from initial vector family $\{x^i\}$.

2. While family $\{x^i\}$ contains at least one vector do:

   (a) Choose from family $\{x^i\}$ the vector with lowest index $x$ and find its projection on orthogonal complement of $\{y^i\}$:

   $$x' = x - \sum_i y^i (x \odot y^i).$$

   (b) Remove vector $x$ from family $\{x^i\}$, and if $x'$ is nonzero, normalize it:

   $$y = \frac{1}{\sqrt{x' \odot x'}} x'$$

   and add to family $\{y^i\}$.

### 2.5.4. Completion of the orthonormal vector family.

Consider orthonormal vector family $\{x^i\}_{i=\overline{0,m}} \in \mathbb{B}^n$, $m < n$. The objective of the following algorithm is to complete this family to contain $n+1$ of vectors so that all vectors to be orthonormal.

*Remark.* Orthonormal family of coordinate vectors $\{e^i\}$ is complete, because it contains $n+1$ of vectors.

**Algorithm 2.3** (Orthonormal vector family completion).

1. For each coordinate vector $\{e^j\}_{j=\overline{0,n}}$ find a vector $y$ orthogonal to family $\{x^i\}$. If vector $y$ isn't zero, normalize it and add to family $\{x^i\}$. All coordinate vectors can't be linear independent with family $\{x^i\}$, because they form complete vector family.

### 2.5.5. Basis change of the lineal. Canonical form of the lineal.

Consider some space lineal $L^m \subset \mathbb{B}^n$ defined by the matrix of size $(n+1) \times (m+1)$. This matrix columns $\{l^i\}_{i=\overline{0,m}}$ form vectors of lineal basis.

Consider vector $v = (v_0 : \dots : v_n) \in L^m \subset \mathbb{B}^n$. Let vector $v$ coordinates in basis $\{l^i\}$ be $a = (a_0 : \dots : a_m)$. Then $v = L^m a$. Let $\mathfrak{M}$ be internal motion of lineal $L^m$, defined by matrix $M$ of size $(m+1)$. Let coordinates of $v$ in new basis $L'^m$ be $b = (b_0 : \dots : b_m)$. Then $b = Ma$. At the same time, $v = L'^m b$. Because coordonates of vector $v$ in $\mathbb{B}^n$ don't change, the following matrix equality is true:

$$L^m a = v = L'^m b = L'^m (Ma) = (L'^m M) a.$$

This equality doesn't depend on vector $a$ choice, so

$$L^m = L'^m M \tag{2.60}$$

is equation of lineal basis change.

It is necessary to find the unique form of lineal definition. In order to find such a basis, the following algorithm can be used:

**Algorithm 2.4** (Canonical form of lineal basis).

1. Consider coordinate vector family $\{e^p\}_{p=\overline{0,n}}$ of space $\mathbb{R}^{n+1}$. Start with empty basis of $L$.
2. Until new basis has $m+1$ elements, find projection $e'^p$ of the next $e^p$ on $L^m$.
   (a) If $e'^p$ is nonzero, find projection $e''^p$ of $e'^p$ on orthogonal complement to existing basis $\{l^i\}$ of lineal $L$.
   (b) If $e''^p$ isn't zero, normalize it and add to existing basis $\{l^i\}$ of lineal $L$.

## 2.6. Limit Vectors and Lineals

### 2.6.1. Limit vectors. Their decomposition vectors.

Until now we considered only vectors that have an index. However, there exist vectors that have no index.

**Definition 2.6.1** (Limit vector)**.** Define vector with no index *limit vector.*

*Properties.* Limit vectors have the following properties:

- Limit vector is orthogonal to itself. As it was already shown, natural product of vectors has the property $x \odot x \geq 0, \forall x \in \mathbb{B}^n$. And if for some vector $x$ there exists index $i$ so that $x \odot_i x > 0$, then index $i$ is vector index, and thus the vector is not limit. It follows that for some limit vector $x$, $x \odot_i x = 0, i = \overline{0, n}$, and thus $x$ is orthogonal to itself.

- Motions transform limit vectors to limit vectors. Consider a limit vector $x$ and some motion $\mathfrak{M}$. Let vector $y = \mathfrak{M}x$. Vector $y$ is also limit. Really, if it would have index, then the motion $\mathfrak{M}^{-1}$ should preserve this index and vector $x = \mathfrak{M}^{-1}y$ would have index. However, vector $x$ is limit one. It means that vector $y$ is also limit.

- If two limit vectors $x, y$ are not collinear and not orthogonal, then not collinear vectors $x + y, x - y$ have indices. The fact $x + y, x - y$ are not collinear follows from the fact $x, y$ are not collinear. Calculate:

  $$(x + y) \odot (x + y) = x \odot x + 2x \odot y + y \odot y = 2x \odot y \neq 0,$$

  because vectors $x, y$ are not orthogonal. Similarly,

  $$(x - y) \odot (x - y) = x \odot x - 2x \odot y + y \odot y = -2x \odot y \neq 0.$$

  It means vectors $x + y, x - y$ have indices.

*Remark.* Some authors use term *isotropic* vectors to describe limit vectors. This name suggests the idea of immutability of vectors direction. However, this property of limit vectors on two–dimensional homogeneous planes, in general, is not true for limit vectors of larger dimensions.

Consider a limit vector $x \in \mathbb{B}^n$.Construct vectors $a, b \in \mathbb{B}^n$ so that vector $a$ has all coordinates $a_i$ equal to those of $x$ (coordinates that enter in bilinear form with sign "+" or coefficient 0), the rest equal to 0. Vector $b$ has all coordinates $b_j$ equal to those of $x$ (coordinates that enter in bilinear form with sign "–"), the rest equal to 0.

**Definition 2.6.2** (Decomposition vectors of limit vector)**.** We call vector pair $a, b$ *decomposition vectors* of vector $x$.

*Properties.* Decomposition vectors have the following properties:

- Vector $x = a + b$ (this property justifies the name of decomposition vectors).
- Vectors $a, b$ are indexed. Really, since vector $x$ isn't zero, among $a_i$ and $b_j$ there are nonzero coordinates, it means:

$$a \odot a > 0,$$
$$b \odot b > 0.$$

- Vectors $a, b$ have equal measure value. Really:

$$x \odot x = a \odot a - b \odot b = 0,$$
$$a \odot a = b \odot b.$$

- Decomposition vectors $a, b$ are interchangeable (the rotation type in plane $a, b$ equals to $-1$). Really, if in same bilinear form the vectors enter with opposite signs and the value of the bilinear form equals to 0, then the types of $a$ and $b$ are opposed, and thus the type of motion from $a$ to $b$ equals to $-1$.
- Vectors $a$ and $b$ are orthogonal. Really, because all nonzero coordinates of vector $a$ have different position from nonzero coordinates of vector $b$.
- No decomposition vector is orthogonal to limit vector:

$$a \odot x = a \odot (a + b) = a \odot a + a \odot b = a \odot a \neq 0,$$
$$b \odot x = b \odot (a + b) = b \odot a + b \odot b = b \odot b \neq 0.$$

*Remark.* Decomposition of limit vector depends on space basis choice and thus is not unequally possible.

**Example.** Consider Minkowski plane with signature $\{0, -1\}$ and limit vector $x = \{1, 1\}$. Its decomposition vectors are $a = \{1, 0\}$ and $b = \{0, 1\}$. The basis change

$$\mathfrak{R}_2(\varphi) = \begin{pmatrix} \cosh\varphi & \sinh\varphi \\ \sinh\varphi & \cosh\varphi \end{pmatrix}$$

maps $x \to x' = \{\cosh\varphi + \sinh\varphi, \cosh\varphi + \sinh\varphi\}, a \to a' = \{\cosh\varphi, \sinh\varphi\}, b \to b' = \{\sinh\varphi, \cosh\varphi\}$. In new basis decomposition vectors are also $a'' = \{\cosh\varphi + \sinh\varphi, 0\}, b'' = \{0, \cosh\varphi + \sinh\varphi\}$.

**Lemma 2.6.1.** *Decomposition vector indices do not depend on space basis choice.*

*Proof.* All possible pairs of decomposition vectors result from some given pair by all possible basis changes. And by theorem 2.4.3, basis change doesn't change vector indices. □

### 2.6.2. Type of the limit vector.

Because for limit vector $x \odot x = 0$, it may seem that its measure equals to zero, However, the distance is defined as measure of motion. The motion along limit vector is not trivial, and thus its measure should not be zero.

**Theorem 2.6.2** (Mean characteistic)**.** *Limit vector type always equals to* 0.

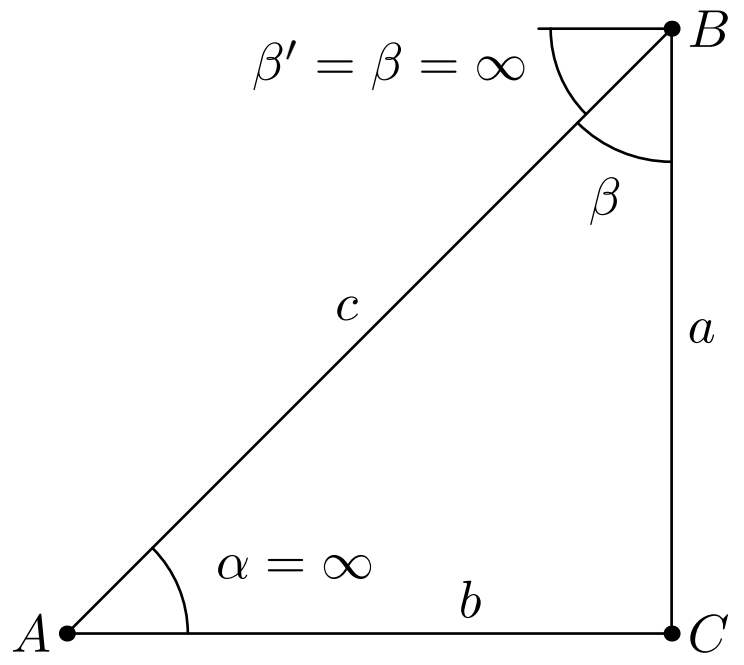

Figure 2.8: Right triangle with hypotenuse lying on limit vector.

*Proof.* Consider a plane $\mathbb{B}^2$ with signature $\{k_1, k_2 = -1\}$ (the requirement $k_2 = -1$ is necessary for existence of limit vectors, see decomposition vectors properties). Construct the right triangle $\triangle ABC$ with catheti $a, b$, hypotenuse $c$, internal angle $\alpha = \infty$ and external angle $\beta'$ (Figure 2.8). By (2.48):

$$\begin{aligned}
T_{12}(a) &= S_1(b) \tanh \alpha = S_1(b) \tanh \infty = S_1(b),\\
1 + k_1 k_2 T_{12}^2(a) &= 1 - k_1 T_{12}^2(a) = 1 - k_1 S_1^2(b),\\
\frac{1}{C_{12}^2(a)} &= C_1^2(b),\\
C_{12}(a) C_1(b) &= C(c) = 1.
\end{aligned}$$

We used (2.46) here. This equality is true regardless of $c$, it follows that its type equals to 0.

Moreover:

$$
\begin{aligned}
C(c) &= C_{12}(a)C_1(b) = 1,\\
C_{12}^2(a) &= \frac{1}{C_1^2(b)},\\
1 - k_1k_2S_{12}^2(a) &= 1 + k_1S_{12}^2(a) = 1 + k_1T_1^2(b),\\
S_{12}(a) &= T_1(b) = T_1(b)\tanh\beta',\\
\tanh\beta' &= 1,\\
\beta' &= \infty.
\end{aligned}
$$

We used (2.51) here. In this case exists also internal angle $\beta = \infty$, and hypotenuse lies on limit vector with type 0. □

*Remark.* While deduction the $T(x)$ form (2.55) of Pythagoras theorem (2.46) the type of $c$ was considered equal to $k_1$. In constructed triangle hypotenuse $c$ has the type 0. So, equation (2.55) is not applicable in this case. However, its modification for this case is true:

$$
\begin{aligned}
0 = 0 \cdot T^2(c) &= k_1(k_2T_{12}^2(a) + T_1^2(b) + k_1k_2T_{12}^2(a)T_1^2(b))\\
&= k_1T_1^2(b) + k_1k_2T_{12}^2(a)(1 + k_1T_1^2(b)) = k_1T_1^2(b) + k_1k_2\frac{T_{12}^2(a)}{C_1^2(b)}\\
&= \frac{k_1(S_1^2(b) + k_2T_{12}^2(a))}{C_1^2(b)} = \frac{k_1(S_1^2(b) - T_{12}^2(a))}{C_1^2(b)} = \frac{k_1 \cdot 0}{C_1^2(b)} = 0.
\end{aligned}
$$

**Corollary.** *Rotation from indexed to limit vector, orthogonal to it, has type equal to* 0.

*Remark.* The type of catheti $a$ and $b$ have opposite values while the type of hypotenuse $c$ equals to 0, which is always the mean value between the types of $a$ and $b$.

### 2.6.3. Measure of the limit vector.

Construct the motion along limit vector $p$ as limit case of motions along indexed vectors $e$ and $h$ (Figure 2.9).

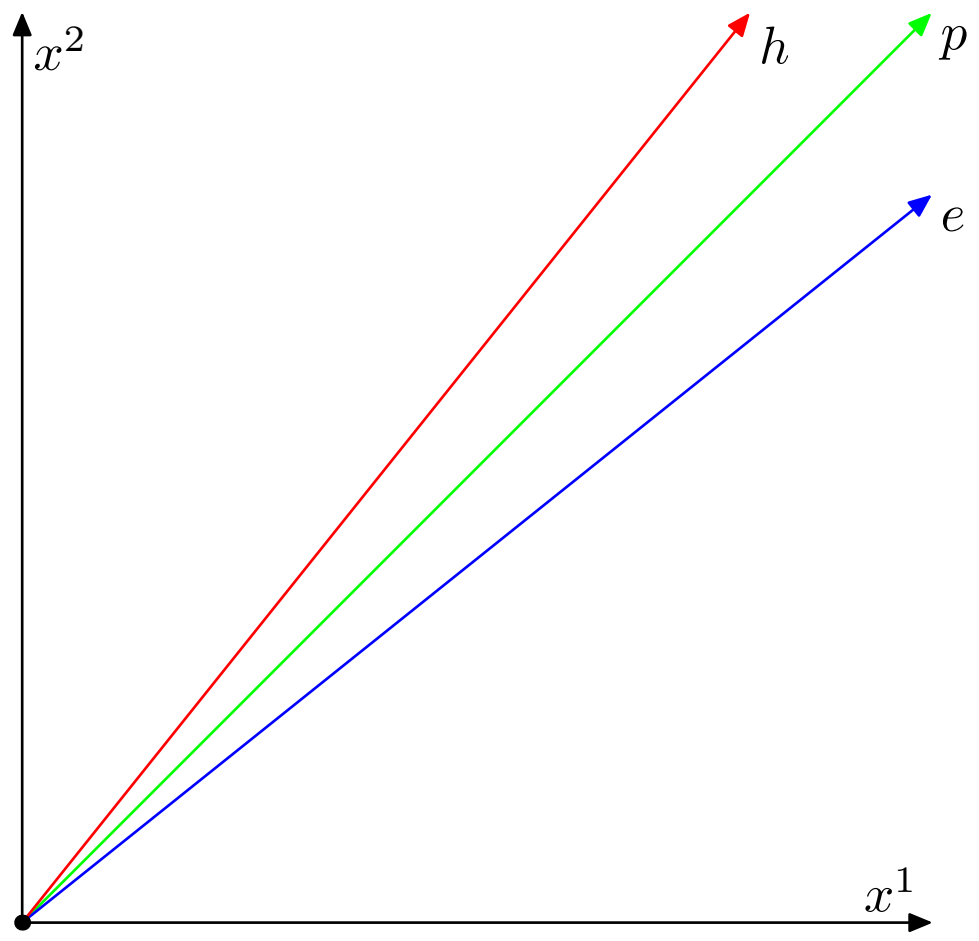

Figure 2.9: Approximation of motion along limit vector by motions along indexed vectors.

The motion along $e$ equals $\mathfrak{E} = \mathfrak{R}_2(\varphi)\mathfrak{T}_1(a)\mathfrak{R}_2(-\varphi)$:

$$\begin{aligned}
\mathfrak{E} &= \begin{pmatrix} 1 & 0 & 0 \\ 0 & \cosh\varphi & \sinh\varphi \\ 0 & \sinh\varphi & \cosh\varphi \end{pmatrix} \begin{pmatrix} C_1(a) & -k_1 S_1(a) & 0 \\ S_1(a) & C_1(a) & 0 \\ 0 & 0 & 1 \end{pmatrix} \begin{pmatrix} 1 & 0 & 0 \\ 0 & \cosh\varphi & -\sinh\varphi \\ 0 & -\sinh\varphi & \cosh\varphi \end{pmatrix} \\
&= \begin{pmatrix} C_1(a) & -k_1 S_1(a) & 0 \\ S_1(a)\cosh\varphi & C_1(a)\cosh\varphi & \sinh\varphi \\ S_1(a)\sinh\varphi & C_1(a)\sinh\varphi & \cosh\varphi \end{pmatrix} \begin{pmatrix} 1 & 0 & 0 \\ 0 & \cosh\varphi & -\sinh\varphi \\ 0 & -\sinh\varphi & \cosh\varphi \end{pmatrix} \\
&= \begin{pmatrix} C_1(a) & -k_1 S_1(a)\cosh\varphi & k_1 S_1(a)\sinh\varphi \\ S_1(a)\cosh\varphi & C_1(a)\cosh^2\varphi - \sinh^2\varphi & -(C_1(a)-1)\cosh\varphi\sinh\varphi \\ S_1(a)\sinh\varphi & (C_1(a)-1)\cosh\varphi\sinh\varphi & -C_1(a)\sinh^2\varphi + \cosh^2\varphi \end{pmatrix}.
\end{aligned}$$

Note, that by (2.7, 2.8):

$$\begin{aligned}
C_1(a) - 1 = C_1^2\left(\frac{a}{2}\right) - k_1 S_1^2\left(\frac{a}{2}\right) - C_1^2\left(\frac{a}{2}\right) - k_1 S_1^2\left(\frac{a}{2}\right) &= -2k_1 S_1^2\left(\frac{a}{2}\right), \\
C_1(a)\cosh^2\varphi - \sinh^2\varphi &= \cosh^2\varphi - \sinh^2\varphi + (C_1(a)-1)\cosh^2\varphi \\
&= 1 - 2k_1 S_1^2\left(\frac{a}{2}\right)\cosh^2\varphi, \\
-C_1(a)\sinh^2\varphi + \cosh^2\varphi &= -\sinh^2\varphi + \cosh^2\varphi - (C_1(a)-1)\cosh^2\varphi \\
&= 1 + 2k_1 S_1^2\left(\frac{a}{2}\right)\cosh^2\varphi.
\end{aligned}$$

Let now $\varphi \to \infty, a \to 0$, so that $S_1\left(\frac{a}{2}\right)\cosh\varphi = \frac{\lambda}{2}$ is constant. Then:

$$\begin{aligned}
1 + 2k_1 S_1^2\left(\frac{a}{2}\right)\cosh^2\varphi &= 1 + k_1\frac{\lambda^2}{2},\\
1 - 2k_1 S_1^2\left(\frac{a}{2}\right)\cosh^2\varphi &= 1 - k_1\frac{\lambda^2}{2},\\
(C_1(a) - 1)\cosh\varphi\sinh\varphi &= -2k_1 S_1^2\left(\frac{a}{2}\right)\cosh^2\varphi\tanh\varphi = -k_1\frac{\lambda^2}{2},\\
S_1(a)\cosh\varphi &= 2S_1\left(\frac{a}{2}\right)C_1\left(\frac{a}{2}\right)\cosh\varphi = \lambda,\\
S_1(a)\sinh\varphi &= S_1(a)\cosh\varphi\tanh\varphi = \lambda,\\
C_1(a) &= 1.
\end{aligned}$$

We used (2.9) here. Thus:

$$\mathfrak{P}(\lambda) = \begin{pmatrix} 1 & -k_1\lambda & k_1\lambda \\ \lambda & 1 - k_1\frac{\lambda^2}{2} & k_1\frac{\lambda^2}{2} \\ \lambda & -k_1\frac{\lambda^2}{2} & 1 + k_1\frac{\lambda^2}{2} \end{pmatrix}$$

Now, construct the motion along $h$, $\mathfrak{H} = \mathfrak{R}_2(\varphi)\mathfrak{T}_2(b)\mathfrak{R}_2(-\varphi)$:

$$\begin{aligned}
\mathfrak{H} &= \begin{pmatrix} 1 & 0 & 0 \\ 0 & \cosh\varphi & \sinh\varphi \\ 0 & \sinh\varphi & \cosh\varphi \end{pmatrix}\begin{pmatrix} C_{12}(b) & 0 & -K_2 S_{12}(b) \\ 0 & 1 & 0 \\ S_{12}(b) & 0 & C_{12}(b) \end{pmatrix}\begin{pmatrix} 1 & 0 & 0 \\ 0 & \cosh\varphi & -\sinh\varphi \\ 0 & -\sinh\varphi & \cosh\varphi \end{pmatrix}\\
&= \begin{pmatrix} C_{12}(b) & 0 & -K_2 S_{12}(b) \\ S_{12}(b)\sinh\varphi & \cosh\varphi & C_{12}(b)\sinh\varphi \\ S_{12}(b)\cosh\varphi & \sinh\varphi & C_{12}(b)\cosh\varphi \end{pmatrix}\begin{pmatrix} 1 & 0 & 0 \\ 0 & \cosh\varphi & -\sinh\varphi \\ 0 & -\sinh\varphi & \cosh\varphi \end{pmatrix}\\
&= \begin{pmatrix} C_{12}(b) & K_2 S_{12}(b)\sinh\varphi & -K_2 S_{12}(b)\cosh\varphi \\ S_{12}(b)\sinh\varphi & \cosh^2\varphi - C_{12}(b)\sinh^2\varphi & (C_{12}(b) - 1)\cosh\varphi\sinh\varphi \\ S_{12}(b)\cosh\varphi & -(C_{12}(b) - 1)\cosh\varphi\sinh\varphi & C_{12}(b)\cosh^2\varphi - \sinh^2\varphi \end{pmatrix}.
\end{aligned}$$

As before, let $\varphi \to \infty, b \to 0$, so that $S_{12}\left(\frac{b}{2}\right)\cosh\varphi = \frac{\mu}{2}$. Then:

$$\mathfrak{P}'(\mu) = \begin{pmatrix} 1 & K_2\mu & -K_2\mu \\ \mu & 1 + K_2\frac{\mu^2}{2} & -K_2\frac{\mu^2}{2} \\ \mu & K_2\frac{\mu^2}{2} & 1 - K_2\frac{\mu^2}{2} \end{pmatrix}$$

What represent the parameters $\lambda$ and $\mu$? Because $a \to 0$:

$$S_1\left(\frac{a}{2}\right)\cosh\varphi = T_1\left(\frac{a}{2}\right)C_1\left(\frac{a}{2}\right)\cosh\varphi \to T_1\left(\frac{a}{2}\right)\cosh\varphi = \frac{\lambda}{2},$$
$$S_{12}\left(\frac{b}{2}\right)\cosh\varphi = T_{12}\left(\frac{b}{2}\right)C_{12}\left(\frac{b}{2}\right)\cosh\varphi \to T_{12}\left(\frac{b}{2}\right)\cosh\varphi = \frac{\mu}{2},$$

where by (2.45) $\frac{\lambda}{2}$ and $\frac{\mu}{2}$ are projections of $\frac{a}{2}$ and $\frac{b}{2}$ on decomposition vectors $x^1$ and $x^2$ respectively (having zero types of $\lambda, \mu$, $T(\lambda) = \lambda, T(\mu) = \mu$).

From $K_2 = k_1 k_2 = -k_1$ follows that $\mathfrak{P}(\lambda) = \mathfrak{P}'(\mu)$. And thus, the parabolic parameter $\lambda$ or $\mu$ represents the length of vector $p$.

**Lemma 2.6.3.** *Value of limit vector measure $x \in \mathbb{B}^n$ equals to values of its decomposition vectors measurs, $a$ and $b$:*

$$|x| = |a| = |b|$$

*Proof.* Consider a limit vector $x \in \mathbb{B}^n$, that has more then two coordinates different from zero. Having orthogonal vectors $a$ and $b$, we can normalize them $a' = \frac{1}{|a|}a, b' = \frac{1}{|b|}b$ and include in orthonormal basis of $\mathbb{B}^n$. In this new basis all coordinates of vector $x = \lambda a' + \lambda b'$ equal to 0, except $p$-th and $q$-th ($p$ and $q$ are indices of vectors $a'$ and $b'$ respectively), which equal to $\lambda = |a| = |b|$. As it was shown, the value of vector measure is

$$|x|^2 = \lambda^2 = |a|^2 = |b|^2.$$

□

So, the type of limit vector always equals to 0, while the types of its decomposition vectors are opposite: either 1 and $-1$, or both 0. The values of limit and its decomposition vectors measures are all equal.

*Remark.* The measure of limit vector depend on space basis choise and *is not* invariant of motions. Thus, it can't be considered a measure in the strict sense.

**Example.** Consider limit vector of unite measure $p = \{0 : 1 : 1\}$. Let

$$q = \mathfrak{R}_2(\varphi)p = \begin{pmatrix} 1 & 0 & 0 \\ 0 & \cosh\varphi & \sinh\varphi \\ 0 & \sinh\varphi & \cosh\varphi \end{pmatrix}\begin{pmatrix} 0 \\ 1 \\ 1 \end{pmatrix} = \begin{pmatrix} 0 \\ \cosh\varphi + \sinh\varphi \\ \sinh\varphi + \cosh\varphi \end{pmatrix}$$

So, measure $q = \mathfrak{R}_2(\varphi)p$ equals to $\cosh\varphi + \sinh\varphi > 1$.

**Proposition 2.6.4.** *The motion $\mathfrak{P}(\lambda) = I$ (is trivial motion) if and only if $\lambda = 0$.*

**Proposition 2.6.5.** *The measure ratio of collinear limit vectors is preserved on motions.*

*Proof.* Consider two collinear limit vectors $x, y \in \mathbb{B}^n$ with measures $\lambda \neq 0, \mu \neq 0$ respectively. It means that $y = \frac{\mu}{\lambda}x$. Consider a motion $\mathfrak{M}$. Then, $y' = \mathfrak{M}y = \mathfrak{M}\left(\frac{\mu}{\lambda}x\right) = \frac{\mu}{\lambda}\mathfrak{M}x = \frac{\mu}{\lambda}x'$. So the measure ratio of vectors $x$ and $y$ equals to measure ratio of vectors $x'$ and $y'$ equals to $\frac{\mu}{\lambda}$ (the proof doesn't rely on limit property of vectors, so the proposition is not specific to limit vectors only). □

After all, the measure of limit vectors is not, strictly speaking, a measure. However it is usefull from two points of view. Firstly, it helps to estimate the measure ratio of collinear vectors instead of simply declare all these measure zero. Secondly, this measure definition corresponds to the principle that only zero vector has zero measure.

### 2.6.4. Orthogonalization of the limit vectors.

Consider two vectors, one limit $x \in \mathbb{B}^n$ and one indexed $y \in \mathbb{B}^n$. Let them be orthogonal $x \odot y = 0$. And let decomposition vectors be $a, b \in \mathbb{B}^n, x = a + b$. Choose space basis in a way that vectors $a$ and $b$ are collinear to some coordinate vectors (and vector $x$ has two coordinates, $p$-th $q$-th, different from zero). If at least one of coordinates $y_p$ or $y_q$ is different from zero, then in this space basis:

$$x \odot x = a_p^2 - b_q^2 = 0,$$
$$x \odot y = a_p y_p - b_q y_q = 0,$$
$$y \odot y \neq 0.$$

From here follows that $a_p = \pm b_q, y_p = \pm y_q$ (either both signs are "+" or both signs are "−").

Formally, vectors $x$ and $y$ are orthogonal, however decomposition vectors $a$ and $b$ are not orthogonal to vector $y$. We can "orthogonalize" the vectors $x$ and $y$ even more. We can construct indexed vector $z = y - \frac{x_p}{a_p}x$:

$$x \odot z = x \odot \left(y - \frac{x_p}{a_p}x\right) = x \odot y - \frac{x_p}{a_p}x \odot x = 0,$$
$$z \odot z = \left(y - \frac{x_p}{a_p}x\right) \odot \left(y - \frac{x_p}{a_p}x\right) = y \odot y - 2\frac{x_p}{a_p}x \odot y + \left(\frac{x_p}{a_p}\right)^2 x \odot x = y \odot y \neq 0.$$

In this case vector $z$ is orthogonal to both decomposition vector $a$ and $b$, because all their not zero coordinates have different positions (for vectors $x$ and $y$ these positions are $p$ and $q$).

*Remark.* If some orthogonal vector family has limit vector $x$ and indexed vector $y$, then replacing $y$ by $z$, the orthogonality isn't affected, because having both vectors $x, y$ orthogonal to all others, vector $z$, as their linear combination, is also orthogonal to all other vectors.

Now, consider two limit vectors $x, y \in \mathbb{B}^n$. And let them be orthogonal $x \odot y = 0$. Let their decomposition vectors be $a, b, c, d \in \mathbb{B}^n, x = a + b, y = c + d$. As earlier, choose space basis in which vectors $a$ and $b$ are collinear to some coordinate vectors (and vector $x$ has two coordinate, $p$-th and $q$-th, different from zero). In this basis:

$$x \odot x = a_p^2 - b_q^2 = 0,$$
$$y \odot y = 0.$$

If index $p$ is among indices $i$, then index $q$ have to be among indices $j$, and:

$$a_p = \pm b_q,$$
$$x \odot y = a_p c_p - b_q d_q = 0,$$
$$c_p = \pm d_q \neq 0.$$

In these equations either both signs are "+", or both signs are "−".

Although in this case vectors $x, y$ are orthogonal, their decomposition vectors are not. Vectors $x, y$ can be more "orthogonalized":

$$z = y - \frac{c_p}{a_p} x.$$

Because $\frac{a_p}{c_p} = \frac{b_q}{d_q}$, vectors $x$ and $z$ are free of common non-zero coordinates (the common positions of non-zero coordinates of vectors $x, y$ are $p, q$). Then:

$$x \odot z = x \odot \left( y - \frac{c_p}{a_p} x \right) = a \odot y - \frac{c_p}{a_p} x \odot x = 0.$$

Vectors $x$ and $z$ are orthogonal with their decomposition vectors ($z = e + f$):

$$a \odot e = 0,$$
$$b \odot f = 0.$$

Obviously, because $p \notin \{j\}$ and $q \notin \{i\}$, also:

$$a \odot f = 0,$$
$$b \odot e = 0.$$

**Proposition 2.6.6.** *We will consider limit vectors orthogonal to indexed vectors or orthogonal among them, if all their decomposition vectors are orthogonal to indexed vectors or are mutually orthogonal.*

*Otherwise, we will "orthogonalize" them using described above rules.*

*Remark.* In case of such orthogonalization, a limit vector is not more orthogonal to itself. By orthogonalizing it with itself, obtain zero vector, which is similar to orthogonalization of indexed vectors.

*Remark.* Because decomposition vectors depend on space basis choice, this orthogonalization is also ambiguous. This is similar to ambiguity of non-interchangeable vectors orthogonalization, which also depends on space basis choice.

### 2.6.5. Limit lineals. Double index of limit vector.

**Definition 2.6.3** (Limit lineal)**.** Define *limit lineal* as lineal with orthonormal basis that contains limit vectors.

*Remark.* Limit vectors may lie also in non-limit lineals. In this case they are not part of orthonormal basis.

**Lemma 2.6.7.** *If orthonormal basis of some lineal contain limit vectors, their decomposition vectors indices are free in lineal.*

*Proof.* Consider by proposition 2.6.6 that all decomposition vectors of all basis limit vectors are mutually orthogonal and are orthogonal to indexed basis vectors. It means that indexed basis vectors don't take indices of decomposition vectors.

Moreover, decomposition vectors do not lie in lineal. Really, if one of decomposition vectors $a$ of limit basis vector $x$ lies in lineal, then another one also lies in lineal, because it is linear combination of lineal vectors $b = x - a$. And if both decomposition vectors lie in lineal, they can be chosen as basis vectors, while the limit vector $x = a + b$ is not more basis one, because it is linear combination of basis vectors. Thus, no one vector $a$ or $b$ lies in lineal if $x$ is basis vector. It follows that indices of vectors $a, b$ are free in lineal. □

**Proposition 2.6.8.** *Consider double index of limit vectors as its decomposition vectors index pair (decomposition vectors are always interchangeable, and thus their indices are different). By lemma 2.6.1 this index doesn't depend on space basis choice and by lemma 2.6.7 it is free in lineal. Denote* $x^{ij} = a^i + b^j$*, where vectors* $a^i, b^j$ *are decomposition vectors of limit vector* $x^{ij}$*.*

**Corollary.** *The limit lineal dimension can't be greater then* $n - m$*, where* $n$ *is space dimension and* $m$ *is number of limit vectors in lineal basis.*

*Proof.* All indices, including subindex of double indices, are unique in a lineal. Their number can't be greater then the maximum number of indices $n + 1$. As each limit basis vector takes double index, the number of basis vectors can't be greater then $n - m + 1$, where $m$ is the number of basis limit vectors. And thus, the limit lineal dimension can't be greater then $n - m$. □

**Corollary.** *The number of vectors in lineal basis can't be greater then $\frac{n+1}{2}$, where n is space dimension.*

*Proof.* Each basis limit vector takes two unique indices. The greatest possible number of indices is $n+1$. It follows, the number of basis limit vectors can't be greater then $\frac{n+1}{2}$. □

### 2.6.6. Signature of the limit lineal.

The limit lineal signature can't be deduced from the space signature and basis vector indices.

**Lemma 2.6.9.** *Orthogonal limit vectors can be either equivalend (the type of motion from one to another equals to* 1*), or non-interchangeable (the type of motion from one to another equals to* 0*).*

*Proof.* Consider two orthogonal limit vectors $x^{ij}, y^{pq} \in \mathbb{B}^n$ with their decomposition vectors $a^i, b^j, c^p, d^q \in \mathbb{B}^n, x^{ij} = a^i + b^j, y^{pq} = c^p + d^q$. Having interchangeability of decomposition vectors ($K_{ij} = K_{pq} = -1$), the following cases are possible:

1. $K_{ip} = K_{jq} = 1$. In this case vector pairs $a^i, c^p$ and $b^j, d^q$ are equivalent and thus, vectors $x^{ij} = a^i + b^j, y^{pq} = c^p + d^q$ are equivalent.

2. $K_{ip} = K_{jq} = -1$. In this case:

$$K_{iq} = K_{ip}K_{pq} = -1 \cdot -1 = 1,$$
$$K_{jp} = \frac{K_{jq}}{K_{pq}} = \frac{-1}{-1} = 1.$$

   That is, vector pairs $a^i, d^q$ and $b^j, c^p$ are equivalent and thus, vectors $x^{ij} = a^i + b^j, y^{pq} = c^p + d^q$ are equivalent.

3. $K_{ip} = K_{jq} = 0$. In this case vector pairs $a^i, c^p$ and $b^j, d^q$ are non-interchangeable and thus, vectors $x^{ij} = a^i + b^j, y^{pq} = c^p + d^q$ are non-interchangeable.

4. $K_{ip} = K_{jq} = \infty$. In this case:

$$K_{pi} = \frac{1}{K_{ip}} = \frac{1}{\infty} = 0,$$
$$K_{qj} = \frac{1}{K_{jq}} = \frac{1}{\infty} = 0$$

   and vector pairs $c^p, a^i$ and $d^q, b^j$ are non-interchangeable and thus, vectors $y^{pq} = c^p + d^q, x^{ij} = a^i + b^j$ are non-interchangeable.

□

**Algorithm 2.5** (Limit lineal signature).

1. Divide all basis vectors of lineal in equivalence groups so that in each group exist only equivalent and / or interchangeable vectors, and each vector of one group is non-interchangeable to all vectors of the next group. Instead of limit vectors use corresponding decomposition vector pairs.
2. In each group, as earlier, find the elements of signature as types of motion from one indexed vector to another.
3. In each group add to signature the number 0 as the type of motion from the last indexed vector to the first limit vector (by theorem 2.6.2 and corollary from it), if such vectors exist.
4. In each group add to signature the number 1 as the type of motion between limit vectors (by lemma 2.6.9), if there exist more than one.
5. As earlier, insert the number 0 between groups of non-interchangeable vectors, from the first to the last.

### 2.7. Constructions and Calculus

Homogeneous space $\mathbb{B}^n$ exactly coincides with metaspace $\mathbb{R}^{n+1}$ up to point normalization. In turn, metaspace is linear vector space. In this way, many constructions, calculus and algorithms for linear vector spaces are applicable to $\mathbb{B}^n$ space. It should not be forgotten, that result have always to be normalized.

**Proposition 2.7.1.** *If the lineal $A^m$ is defined by some points in general positon $\{A_i\}_{i=\overline{0,m}}$, then the linearly independent family of vectors $\{A_i\}$ determines the basis of the lineal $A^m$.*

This basis can be orthonormalized and, if desired, canonized. It follows, that $m$-dimensional lineal is defined by $m + 1$ points in general position.

#### 2.7.1. Midpoint of a segment and center of gravity of a triangle.

**Proposition 2.7.2.** *Up to normalization, the midpoint of segment $AB$ is $C = A + B$.*

*Proof.* Points $A, B$ form basis of lineal $AB$. It means that any point $X$ of lineal $AB$ can be presented as $X = \alpha A + \beta B$. Then $A \odot X = A \odot (\alpha A + \beta B) = \alpha A \odot A + \beta A \odot B = \alpha + \beta A \odot B$; $B \odot X = B \odot (\alpha A + \beta B) = \alpha A \odot B + \beta B \odot B = \alpha A \odot B + \beta$. Having $C(d(X, Y)) = X \odot Y$, where $d(X, Y)$ is distance between $X$ and $Y$:

$$\begin{aligned} d(A, X) &= d(B, X), \\ C(d(A, X)) &= C(d(B, X)), \\ A \odot X &= B \odot X, \\ \alpha + \beta A \odot B &= \alpha A \odot B + \beta, \\ (\beta - \alpha)(A \odot B - 1) &= 0. \end{aligned}$$

Next, assuming that $A \odot B \neq 1$, obtain $\alpha = \beta$. For nonlinear spaces always $A \odot B \neq 1$ when $A \neq B$. For linear spaces the result is also true, in this case $\alpha = \beta = \frac{1}{2}$. □

**Proposition 2.7.3.** *The center of gravity of triangle $\triangle ABC$ (the intersection point of medians) is $M = A + B + C$, up to normalization.*

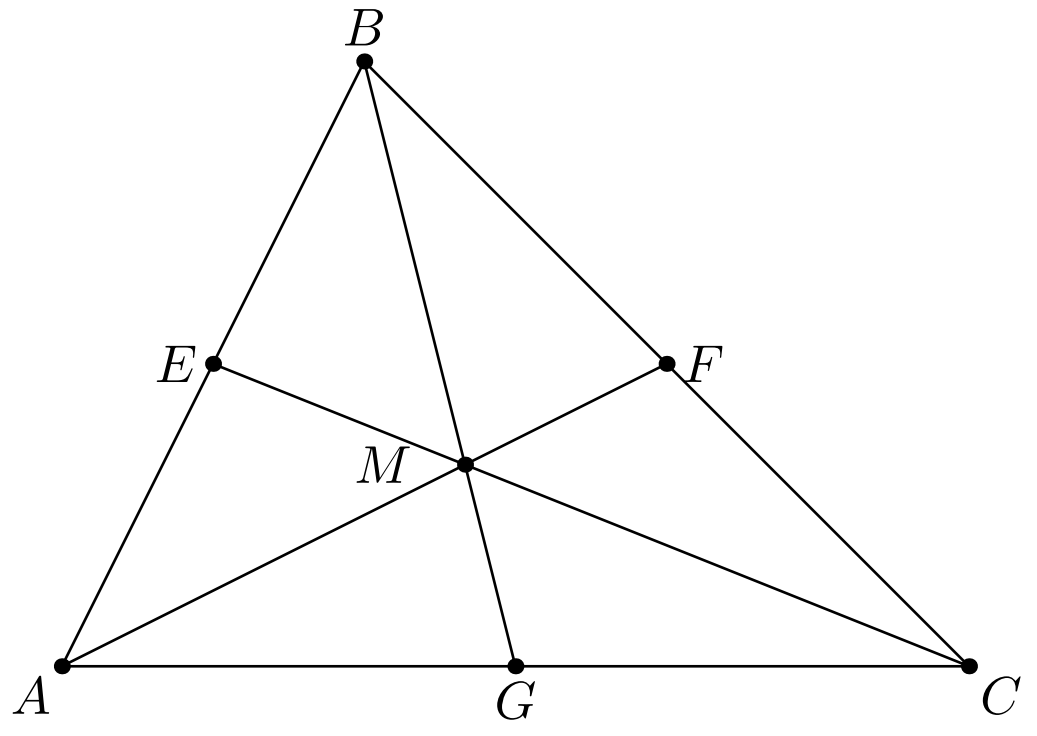

Figure 2.10: Up to normalization, the center of gravity of triangle $\triangle ABC$ is $M = A + B + C$.

*Proof.* All plane points can be presented as $X = \alpha A + \beta B + \gamma C$ in basis $\{A, B, C\}$. Edge $AB$ midpoint is $E = \theta A + \theta B$ (Figure 2.10). Each point of median $CE$ has the form $\varphi E + \psi C = \theta\varphi A + \theta\varphi B + \psi C$, it follows that $\alpha = \beta = \theta\varphi$. Similarly, edge $BC$ midpoint is $F = \xi B + \xi C$. Each point of median $AF$ has the form $\varphi A + \psi F = \varphi A + \psi\xi B + \psi\xi C$, it follows that $\beta = \gamma = \psi\xi$. The point $M$ belongs to both medians, that is $\alpha = \beta = \gamma$ and $M = A + B + C$ up to normalization. □

**Corollary.** *The triangle center of gravity (the median intersection point) exists in all homogeneous spaces.*

### 2.7.2. Sum, intersection and difference of lineals.

**Definition 2.7.1** (Sum, intersection and difference of lineals)**.** Define *sum* of lineals $A^p$, $B^q$ the lineal $C^r$ (denote $C^r = A^p + B^q$), constructed on union of lineal $A^p$ and $B^q$ bases. Define *intersection* of lineals $A^p$, $B^q$ lineal $C^r$ (denote $C^r = A^p \cap B^q$), composed of all vectors lying in $A^p$ and $B^q$ simultaneously. Define *difference* of lineals $A^p$, $B^q$ lineal $C^r$ (denote $C^r = A^p - B^q$) composed of all vectors that lie in $A^p$, and are orthogonal to $B^q$.

To be certain, consider lineals $A^p$, $B^q$ with orthonormal bases $\{a^i\}_{i=\overline{0,p}}$ and $\{b^j\}_{j=\overline{0,q}}$.

**Algorithm 2.6** (Lineal difference)**.**

1. Find projections $a'^i$ of vectors $a^i$ to orthogonal complement of $B^q$:

$$a'^i = a^i - \sum_{j=0}^{q} (a^i \odot b^j) b^j$$

2. Orthonormalize family $\{a'^i\}$.

**Algorithm 2.7** (Lineal sum and intersection)**.**

1. Denote the basis of sum as $\{w^i\}$ and the basis of intersection as $\{v^i\}$. Use one more vector family $\{h^i\}$.
2. Copy all vectors $\{a^i\}$ to $\{w^i\}$ and $\{h^i\}$, fill with zero vectors all free indices of $\{w^i\}$ and $\{h^i\}$, so that each family contains $n + 1$ of vectors.
3. For each vector $b^i$ do:

   (a) Find vector:

   $$b'^i = b^i - \sum_j (b^i \odot w^j) w^j,$$

   orthogonal to family $\{w^i\}$. Using the same coefficients, calculate:

   $$h = -\sum_j (b^i \odot w^j) h^j.$$

   (b) If $b'^i$ isn't zero, it is linear independent with $\{w^i\}$. Normalize it:

   $$w = \frac{1}{\sqrt{b' \odot b'}} \cdot b'$$

   and add to $\{w^i\}$. Using the same factor, add vector:

   $$h' = \frac{1}{\sqrt{b' \odot b'}} \cdot h$$

to $\{h^i\}$.

(c) If $b'^i$ is zero, a equality is found, in which non-trivial linear combination composed of vectors $\{a^i\}$ and $\{b^i\}$, equals to zero vector. That is, non-trivial linear combination composed of vectors $\{a^i\}$ equals to non-trivial linear combination composed of $\{b^i\}$. Evidently, these linear combinations represent a vector from intersection. This linear combination precisely equals to $h$, because all operations performed with $b^i$ were performed also with $h$, however it is composed only of vectors $\{a^i\}$. Vector $h$ isn't zero, because the linear combination isn't trivial. It is also linear independent with $\{v^i\}$, because it includes non-zero component of new vector $b^i$, which is linear independent with $\{v^i\}$. Orthonormalize $h$ and add to $\{v^i\}$.

4. Exclude zero vectors from family $\{w^i\}$. The family $\{v^i\}$ is intersection basis and the family $\{w^i\}$ is sum basis.

**Proposition 2.7.4.** *It is easy to see, that dimension of $A^p + B^q$ is not less then* $\max(p, q)$ *and is not greater then $p + q$, and dimension of $A^p \cap B^q$ is not greater then* $\min(p, q)$.

**Theorem 2.7.5.** *The sum of lineals $A^p$ and $B^q$ dimensions equals to the sum of lineals $A^p + B^q$ and $A^p \cap B^q$ dimensions.*

*Proof.* In described algorithm each vector from $\{a^i\}$ and $\{b^i\}$ was taken once and each time strictly one vector was added to sum $\{w^i\}$ or intersection $\{v^i\}$. It means that the sum $p+q+2$ of number of vectors in $\{a^i\}$ and $\{b^i\}$ equals to the sum of the number of vectors in $\{w^i\}$ and $\{v^i\}$, which in turn is greater then the sum of dimensions $A^p + B^q$ and $A^p \cap B^q$ by 2. □

**Corollary.** *Dimension of lineal $A^p + B^q$ equals to sum of dimensions of lineals $A^p \cap B^q$, $A^p - B^q$ and $B^q - A^p$.*

*Proof.* If choose bases of lineals $A^p$ and $B^q$ so that they contain basis $\{v^i\}$ of lineal $A^p \cap B^q$, then remaining basis vectors belong to lineals $A^p - B^q$ and $B^q - A^p$ bases respectively. In this situation vector number doesn't change and all vectors from all bases are linear independent. So, the assertion is proved. □

### 2.7.3. Coordinate matrix and state matrix.

Consider a space $\mathbb{B}^n$ with $m+1$ vectors $\{v^i\}_{i=\overline{0,m}}$. Let $v^i$ coordinates be $(v_{0i} : \dots : v_{ni}), i = \overline{0, m}$. And let vectors $v^i$ be ordered and linear independent. Compose the matrix $V$ with elements $\{v_{ij}\}, i = \overline{0, n}, j = \overline{0, m}$.

**Definition 2.7.2** (Coordinate matrix, state matrix)**.** Define *coordinate matrix* of vector family $\{v^i\}$ as rectangular matrix composed of this vector family coordinates. Define *state matrix* of vector family as matrix composed of elements $v^i \odot_i v^j$.

The state matrix shows how orthonormal the vector family is. It is as closer to identity matrix as more normalized and orthogonal to each other are the vectors. If space signature has zero entries, some status matrix elements can have any value even in orthonormalized vector family.

**Lemma 2.7.6.** *If the number of vectors in family $\{v^i\}$ equals to $n+1$ and $V$ and $W$ are respectively coordinate and state matrices of this family, then the volume of parallelepiped constructed on vectors $\{v^i\}$ equals to $|\det V|$ (in this case $V$ is square matrix) and*

$$\det W = (\det V)^2 \,. \tag{2.61}$$

*Proof.* First, let $\{v^i\}_{i=\overline{0,n}}$ be orthonormalized. Then the parallelepiped (cube) volume equals to 1, the matrix $V$ is GM-orthogonal and $\det V = \pm 1$. That is, $|\det V| = 1$ equals to parallelepiped volume. All elements on main diagonal of matrix $w_{ii} = v^i \odot_i v^i = 1$, because all vectors $v^i$ are normalized. All elements above main diagonal $w_{ij} = v^i \odot_i v^j = 0$, $i < j$, because all vector pairs $v^i, v^j$ are orthogonal (elements below main diagonal may be different from 0). It means that matrix $W$ is lower triangular with all elements on main diagonal equal to 1 and $\det W = 1 = (\det V)^2$.

Next, use linear property of vector product and note, that matrix determinant equals to zero if it contains proportional columns or rows. If some row or column is multiplyed by $\alpha$, the determinant of resulting matrix grows by $\alpha$ factor. When some row or column is sum of two rows / columns, the matrix determinant equals to sum of determinants of matrices containing first and second row / column.

Second, if instead of some $v^i$ use $v'^i = \alpha v^i$, the parallelepiped volume grows by $\alpha$ times and $|\det V'| = |\alpha||\det V|$. In this case $\det W' = \alpha^2 \det W$ (the factor $\alpha$ results from row $i$ and column $i$). Again, $\det W' = \alpha^2 \det W = \alpha^2 (\det V)^2 = (\alpha \det V)^2 = \left(\det V'\right)^2$.

Third, if instead of $v^i$ use $v'^i = v^i + \alpha v^j$, then the parallelepiped volume doesn't change as well as determinants $V$ and $\det W' = \det W = (\det V)^2 = (\det V')^2$.

Finally note, that all matrices $V$ can be obtained from GM-orthogonal one by operations from the second and the third point. It means that equation (2.61) holds for all matrices and volume of parallelepiped constructed on vectors $v^i$ always equals to $|\det V|$. □

When we need to determine the parallelepiped volume, it is convenient to use state matrix. Its elements don't depend on space basis choise, don't change on motions and it is always square, even when the number of vectors is less then space dimension (the matrix $V$ isn't square in this case).

### 2.7.4. Measure between lineals.

The measure between congruent lineals was defined as parameter of motion that maps the first one to another. This motion sometimes isn't easy to find. Additionally, there exists measure between non-congruent lineals (they may have even different dimensions). So, we need to find the way to compute the measure between lineals in general case.

**Definition 2.7.3** (Measure between lineals)**.** Define *measure between lineals* $A^p$ and $B^q$, $p \leq q$ as measure between $A^p$ and its projection on lineal $B^q$.

*Remark.* This definition generalizes previous one, issued for lineals $\mathbb{B}^m$ and $\mathfrak{R}_m(\varphi)\mathbb{B}^m$.

The cube constructed on basis vectors $\{a^i\}$ of lineal $A^p$ has $(p+1)$–dimensional volume in $\mathbb{R}^{n+1}$ equal to one. Let vectors $\{a'^i\}$ be projections of vectors $\{a^i\}$ on lineal $B^q$. Evidently, projection of $A^p$ on $B^q$ is based on vectors $\{a'^i\}$. Because the metaspace $\mathbb{R}^{n+1}$ is linear, by (2.45) the parallelepiped volume constructed on vectors $\{a'^i\}$, equals to $C(\varphi)$, where $\varphi$ is angle between linear planes in metaspace, which is equal to measure between lineals in $\mathbb{B}^n$ (Figure 2.11).

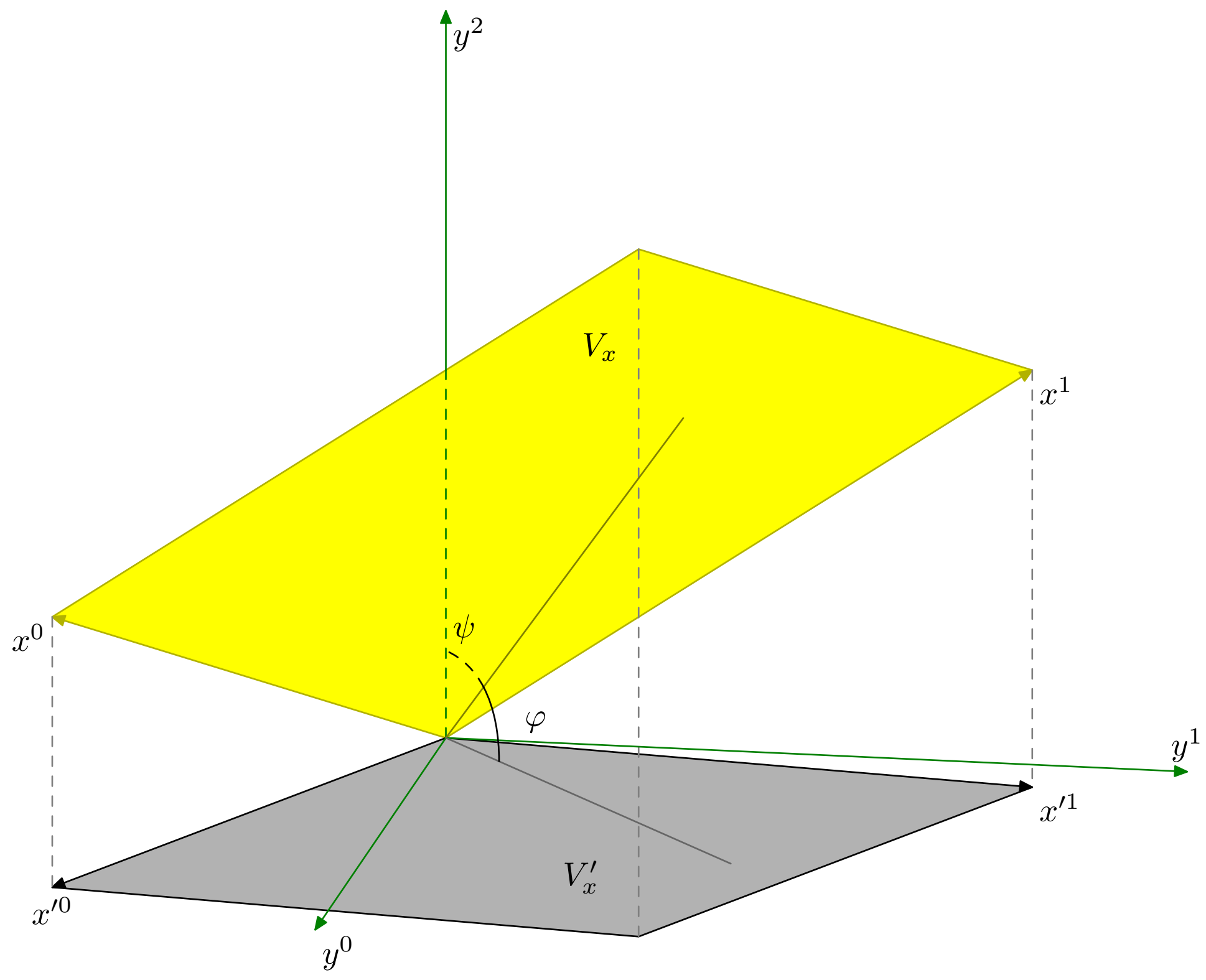

Figure 2.11: Computing the measure between lineals $X^2$ and $Y^2$.

We can use state matrix of vector family $\{a'^i\}$ in order to find the volume of parallelepiped constructed on them. In this way we have to be sure, that dimension of $B^q$ is not less then dimension of $A^p$, otherwise vector family $\{a'^i\}$ contains more vectors then dimension of lineal constructed on them, from which follows that vectors are linear dependent and status matrix determinant always equals to zero. Note that if $p \leq q$ and vectors $\{a'^i\}$ are linear dependent, it results that lineals are orthogonal.

If the type of measure between lineals equals to 0 or $-1$, it may happen that direct measure isn't measurable. In this case helpful information is complementary measure, which is now measurable. For this purpose we need to project each vector from family $\{a^i\}$ to orthogonal complement of $B^q$ resulting the family $\{a''^i\}$, and find the volume of parallelepiped constructed on this family. This volume by (2.50) equals to $S(\varphi) = C(\psi)$, where $\psi$ is complementary angle between planes in metaspace, that equals to complementary measure between lineals. For its correct computing we need to require that lineal $A^p$ dimension is not greater than dimension of orthogonal complement of $B^q$. In order to define the algorithm of measure computing between lineals, use the following lemmas.

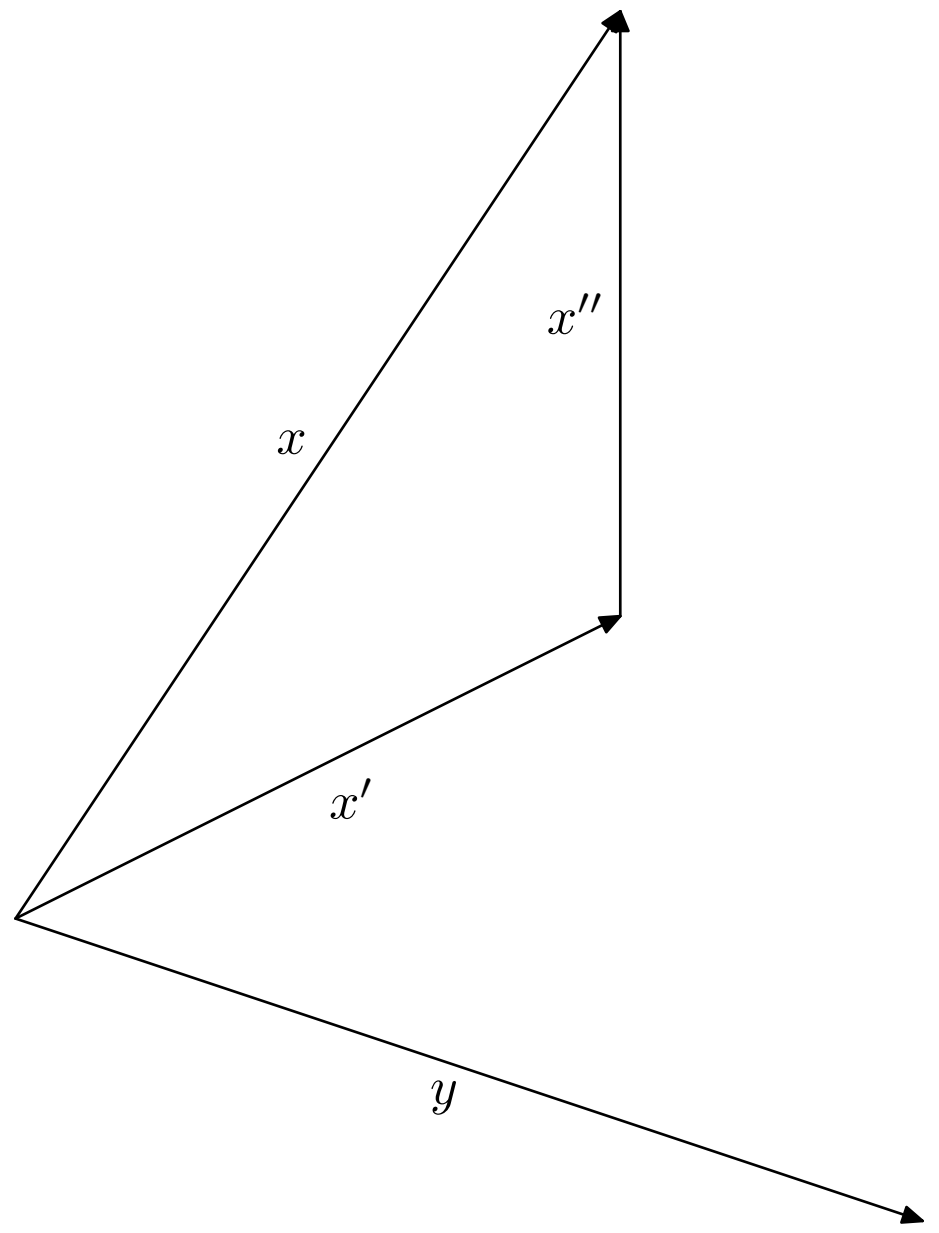

Figure 2.12: From $x \perp y$ follows $x' \perp y$.

**Lemma 2.7.7** (Three orthogonalities rule)**.** *If some vector $x$, that doesn't lie in a lineal is orthogonal to vector $y$, that lies in that lineal ($x \odot y = 0$), then projection $x'$ of vector $x$ on that lineal is also orthogonal to vector $y$ ($x' \odot y = 0$, Figure 2.12).*

*Proof.* By projection definition, $x = x' + x''$, where $x'$ lies in lineal and $x''$ is orthogonal to lineal, and thus, to all vectors from the lineal, it follows that $x'' \odot y = 0$. Having two different vectors, $x$ and $x''$, orthogonal to vector $y$, their linear combination $x' = x - x''$ is orthogonal to $y$:

$$x' \odot y = (x - x'') \odot y = x \odot y - x'' \odot y = 0 - 0 = 0.$$

□

**Lemma 2.7.8.** *In linear space $\mathbb{R}^{n+1}$ the angle $\varphi$ between two lineals $A^p$ and $B^q$ equals to the angle $\psi$ between $A^p - B^q$ and $B^q - A^p$.*

*Proof.* Let $p \leq q$. Choose lineal $A^p$ basis $\{a^i\}$ so that it contains all basis vectors $\{u^i\}$ of intersection $A^p \cap B^q$. Remaining vectors of $A^p$ basis form the basis $\{x^i\}$ of lineal $A^p - B^q$. Projections of vectors $\{u^i\}$ on lineal $B^q$ equal to $\{u^i\}$, because they also lie in $B^q$. Let projections of vectors $\{x^i\}$ on $B^q$ equal to $\{x'^i\}$. Then projection of unite cube constructed on vectors $\{a^i\}$ is parallelepiped constructed on $\{u^i\} \cup \{x'^i\}$. And for each vector triplet $u^i, x^j, x'^j, \quad u^i \perp x'^j$ by the three orthogonalities rule (lemma 2.7.7, $u^i \perp x^j$ and $x'^j$ is projection of $x^j$). That is, the volume of projected parallelepiped equals to product of volumes of cube constructed on vectors $\{u^i\}$, and parallelepiped constructed on vectors $\{x'^i\}$. But volume of cube constructed on $\{u^i\}$ equals to one. On the other hand, because all vectors $\{x'^i\}$ are orthogonal to $\{u^i\}$, they lie in $B^q - A^p$. $C(\varphi)$ equals to volume of parallelepiped constructed on vectors $\{u^i\} \cup \{x'^i\}$, which equals to volume of parallelepiped constructed on $\{x'^i\}$, which equals to $C(\psi)$ (regardless of $\varphi$ and $\psi$ types). That is $\varphi = \psi$. □

**Lemma 2.7.9.** *The number of basis vectors in lineal $A^p - B^q$ doesn't exceed the number of basis vectors of orthogonal complement of lineal $B^q - A^p$.*

*Proof.* Let dimension of lineal $A^p + B^q$ be equal to $u$ ($u+1$ of basis vectors), and dimension of lineal $A^p \cap B^q$ be equal to $m$ ($m + 1$ of basis vectors). Then dimension $A^p - B^q$ equals to $p - m = u - q$ (lineal is improper, same number of basis vectors), dimension of $B^q - A^p$ equals to $u - p$ (lineal is improper, same number of basis vectors), and dimension of its orthogonal complement equals to $n + 1 - (u - p) = n + p - u + 1$. Dimension of sum doesn't exceed the dimension of space: $u \leq n$. It means:

$$u + 1 \leq n + 1 < m + 1 + n + 1,$$
$$p - m \leq n + p - u + 1.$$

□

**Proposition 2.7.10.** *If the measure between lineals equals to zero or they are orthogonal, then the measure type is, generally speaking, ambiguous.*

**Example.** Consider $\mathbb{B}^2$ with signature $\{1, -1\}$. We can approximate vector $e = \{1 : 0 : 0\}$ by vectors $x = \{\cos\varphi : \sin\varphi : 0\}, \varphi \to 0$ or $y = \{\cosh\psi : 0 : \sinh\psi\}, \psi \to 0$ or $z = \{1 : \theta : \theta\}, \theta \to 0$. Measures $\varphi, \psi, \theta$ have zero limit and have the types $1, -1, 0$ respectively.

**Example.** Consider $\mathbb{B}^3$ with signature $\{1, -1, 1\}$. Then the plane constructed on vectors $a = \{1 : 0 : 0 : 0\}$, $b = \{0 : 1 : 0 : 0\}$ and $c = \{0 : 0 : 1 : 0\}$, is orthogonal to vector $d = \{0 : 0 : 0 : 1\}$. Vector $d$ can be approximated by vectors $x = \{0 : \sinh\varphi : 0 : \cosh\varphi\}, \varphi \to 0$ or $y = \{0 : 0 : \sin\psi : \cos\psi\}, \psi \to 0$ or $z = \{0 : \theta : \theta : 1\}, \theta \to 0$. All them tend to $d$ and angle types between them and plane are $-1, 1, 0$ respectively.

**Algorithm 2.8** (Measure between lineals).

1. Find orthonormal bases $\{a^i\}, \{b^j\}$ of the lineals $A^p - B^q$ and $B^q - A^p$. Consider that the number of vectors $\{a^i\}$ doesn't exceed the number of vectors $\{b^j\}$. Otherwise interchange the roles of lineals.

2. Decompose vectors $a^i = a'^i + a''^i$, where the vectors $\{a'^i\}$ are projections of vectors $\{a^i\}$ on $B^q - A^p$, and vectors $\{a''^i\}$ on orthogonal complement to $B^q - A^p$ with dimension not less then the number of vectors $\{a^i\}$.

3. Find determinants $w', w''$ of the status matrices $W', W''$ of vector families $\{a'^i\}, \{a''^i\}$.

   (a) If $w' = 1, w'' = 0$, then $\varphi = 0$, $\psi$ is undefined and $k$ is ambiguous.

   (b) If $w' = 0, w'' = 1$, then $\varphi$ is undefined, $\psi = 0$ and $k$ is ambiguous.

   (c) If $w' + w'' = 1$, then $\varphi = \tan^{-1}\sqrt{\frac{w''}{w'}}$, $\psi = \tan^{-1}\sqrt{\frac{w'}{w''}}$ and $k = 1$.

   (d) If $w' = 1, w'' \neq 0$, then $\varphi = \sqrt{w''}$, $\psi$ unmeasurable and $k = 0$.

   (e) If $w' \neq 0, w'' = 1$, then $\varphi$ unmeasurable, $\psi = \sqrt{w'}$ and $k = 0$.

   (f) If $w' - w'' = 1$, then $\varphi = \tanh^{-1}\sqrt{\frac{w''}{w'}}$, $\psi$ unmeasurable and $k = -1$.

   (g) If $w'' - w' = 1$, then $\varphi$ unmeasurable, $\psi = \tanh^{-1}\sqrt{\frac{w'}{w''}}$ and $k = -1$.

   (h) If at least one of the two lineals is limit and $w' = w''$, then $\varphi = \psi = \infty$ and $k = -1$.

### 2.8. Volume

**Proposition 2.8.1.** *For each $\mathbb{B}^n$, presented as unite sphere in $\mathbb{R}^{n+1}$, the surface is orthogonal to radius.*

*Proof.* Consider $X, Y \in \mathbb{B}^n$ and the distance between $X$ and $Y$ is small. Show that $X \odot (X - Y) \to 0$, when $Y \to X$ in sense of distance between them:

$$X \odot (X - Y) = X \odot X - X \odot Y = 1 - C(d(X, Y)),$$

where $d(X, Y)$ is the distance between $X$ and $Y$. Then:

$$Y \to X,$$
$$d(X, Y) \to 0,$$
$$C(d(X, Y)) \to 1,$$
$$1 - C(d(X, Y)) \to 0.$$

□

Consider $A, B \in \mathbb{B}^1$. $A = (C(a) : S(a))$, $B = (C(b) : S(b))$. Calculate area in $\mathbb{R}^2$ of sector $\mathbb{B}^1$ between $A$ and $B$. In Euclidean polar coordonate system the argument

$$\tan\varphi_e = \frac{y}{x} = \frac{S(p)}{C(p)} = T(p),$$

where $p$ is natural argument in $\mathbb{B}^1$. Euclidean radius:

$$\rho = \sqrt{x^2 + y^2} = \sqrt{C^2(p) + S^2(p)}.$$

Next, having:

$$\begin{aligned} \varphi_e &= \arctan T(p), \\ d\varphi_e &= \frac{dp}{(1 + T^2(p))C^2(p)} = \frac{dp}{C^2(p) + S^2(p)}, \end{aligned}$$

calculate area:

$$\begin{aligned} s &= \frac{1}{2}\int_A^B \rho(\varphi_e)^2 d\varphi_e = \frac{1}{2}\int_a^b (C^2(p) + S^2(p))\frac{dp}{C^2(p) + S^2(p)} \\ &= \frac{1}{2}\int_a^b dp = \frac{1}{2}p\Big|_a^b = \frac{b-a}{2}. \end{aligned}$$

That is, $2s$ equals to distance $AB$.

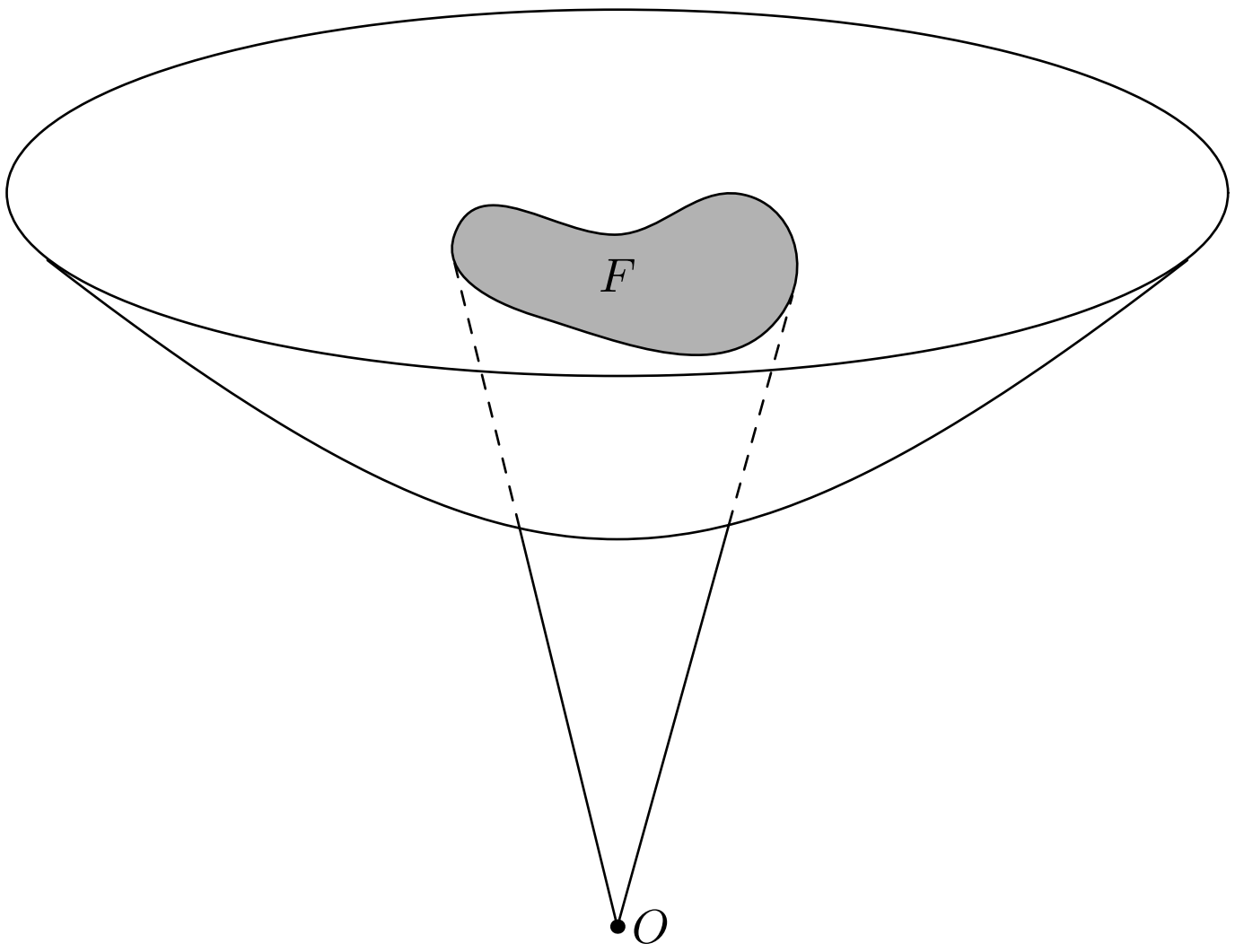

Figure 2.13: Volume measurement of figure $F \subset \mathbb{B}^n$ using cone in $\mathbb{R}^{n+1}$.

Define volume in space $\mathbb{B}^n$ using volume in metaspace $\mathbb{R}^{n+1}$. This approach has the advantage,

that volume in linear metaspace is usually easier to compute than in (not necessarily linear) space. In this case, domain of integration has one extra dimension. On the other hand, we can choose the integration order in several ways and obtain different models of space volume measuring.

**Proposition 2.8.2.** *If $F \subset \mathbb{B}^n$ is some figure with volume $v_{\mathbb{B}}$ (in sense of $\mathbb{B}^n$), and $v_{\mathbb{R}}$ is volume (in sense of $\mathbb{R}^{n+1}$) of cone with base in $F \subset \mathbb{B}^n$ and vertex in $O \notin \mathbb{B}^n$ origin of $\mathbb{R}^{n+1}$ (Figure 2.13), then*

$$v_{\mathbb{B}} = (n+1)v_{\mathbb{R}}$$

*Proof.* By proposition 2.8.1, having $\mathbb{B}^n$ orthogonal to radius, $F$ is also orthogonal to it. The radius equals to 1, because $\forall X \in \mathbb{B}^n, X \odot X = 1$.

For sufficiently small figure $F$ we can consider the sphere region, that contains it in $\mathbb{R}^{n+1}$, as linear. Then the cone volume $v_{\mathbb{R}}$ equals to base area (that is, to volume $v_{\mathbb{B}}$) multiplied by height (that is, by radius orthogonal to base and equal to 1) and multiplied by $\frac{1}{n+1}$. The equality has place.

Any large figure can be divided into sufficiently small parts, so that for each of them the equality has place. It remains to observe, that $v_{\mathbb{B}}$ equals to the sum of parts volumes and $v_{\mathbb{R}}$ equals to the sum of part cones volumes, and thus, the equality holds again. □

**Corollary.** *The volume $v_{\mathbb{B}}$ doesn't change on motions.*

*Proof.* The motion matrix determinant absolute value equals to $\pm 1$, it follows, that motions preserve volume in metaspace $\mathbb{R}^{n+1}$. Motions preserve also sphere $\mathbb{B}^n$. It means, that figure cones remain figure cones on motions. Because their volume $v_{\mathbb{R}}$ don't change, $v_{\mathbb{B}} = (n+1)v_{\mathbb{R}}$ also preserve their value on motions. □

### 2.8.1. Area of the right triangle on linear planes.

The easiest case of area calculus is the case of right triangle on linear plane ($k_1 = 0$). These planes spheres equations, $x_0^2 = 0$, represent two planes in metaspace. Consider a pyramid $OABC$ in metaspace $\mathbb{R}^3$ with base in right triangle $\triangle ABC$ and vectex in metaspace origin $O = (0, 0, 0)$ (Figure 2.14). Let triangle $\triangle ABC$ have point $A = E = (1 : 0 : 0)$ in origin of $\mathbb{B}^2$, catheti $a, b$ and hypotenuse $c$. Construct a secant parallel to lines $OA$ and $BC$ (the metaspace is linear), having distance $x$ from point $A$. This plane intersects the pyramid by rectangle with sides $y$ and $z$. By (2.48):

$$a = bT_2(\alpha),$$
$$y = xT_2(\alpha),$$
$$y = x\frac{a}{b}$$

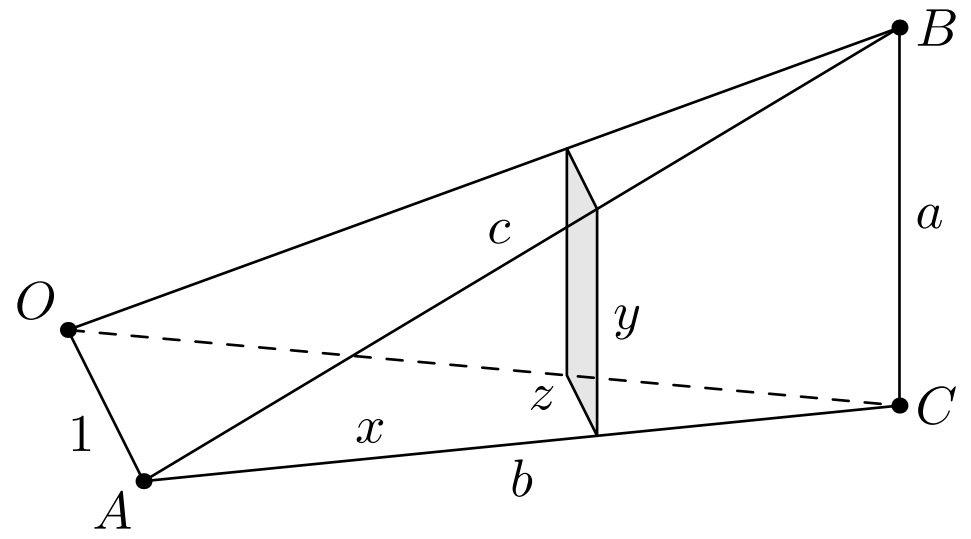

Figure 2.14: Finding right triangle area.

Similarly, having the radius equal to 1:

$$\frac{z}{b-x} = \frac{1}{b},$$
$$z = \frac{b-x}{b}$$

The rectangular section area equals to $yz$. Then pyramid $OABC$ volume equals to:

$$v = \int_0^b y(x)z(x)dx = \frac{a}{b}\frac{1}{b}\int_0^b x(b-x)dx = \frac{a}{b^2}\int_0^b (bx-x^2)dx$$
$$= \frac{a}{b^2}\left(\frac{bx^2}{2} - \frac{x^3}{3}\right)\Bigg|_0^b = \frac{a}{b^2}\left(\frac{b^3}{2} - \frac{b^3}{3}\right) = \frac{a}{b^2}\frac{b^3}{6} = \frac{ab}{6}$$

By proposition 2.8.2, triangle $\triangle ABC$ area equals to:

$$s = 3v = \frac{3ab}{6} = \frac{ab}{2} \tag{2.62}$$

As it follows from (2.62), the right triangle area on linear planes doesn't depend on angular type $k_2$.

### 2.8.2. Area of the right triangle on non-linear planes.

On non-linear planes ($k_1 \neq 0$), we will search triangle area equation as function of its angles. This way isn't applicable when $k_1 = 0$, because in this case plane $\mathbb{B}^2$ allows distance dilation and thus, the angles don't define lengths and areas. But in this case the area equals to semi product of catheti. We will generalize the result for linear case and justify this way.

Consider polar coordonate system $(\rho, \varphi)$ (Figure 2.15). As it was shown (proposition 2.8.1) in all homogeneous spaces the circular arc is orthogonal to radius. Let $E = (0,0)$ be origin, $P = (\rho, \varphi), Q = (\rho, \varphi + d\varphi)$ be two points on circle of radius $\rho$. Calculate distance $PQ = dl$. Drop from $Q$ to $EP$ perpendicular $QR \perp EP$. Then $\triangle EQR$ is right triangle. Let $ER = a$. Then:

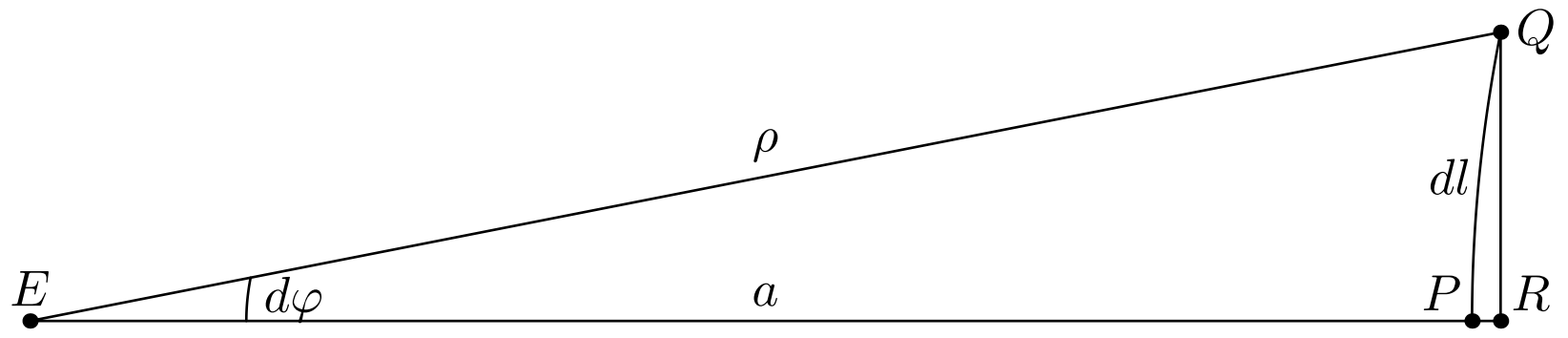

Figure 2.15: Area element deduction in polar coordonate system.

$$
\begin{aligned}
&QR \to dl,\\
&C(dl) \to 1,\\
&S(dl) = T(dl)C(dl) \to dl,\\
&C(d\varphi) \to 1,\\
&S(d\varphi) = T(d\varphi)C(d\varphi) \to d\varphi,\\
&T(a) = T(\rho)C(d\varphi), \quad a \to \rho.
\end{aligned}
$$

That is, when $d\varphi \to 0$, also $dl \to 0$ and $a \to \rho$. So, in polar coordonate system the circlar arc length by (2.47) is:

$$dl \approx S(dl) = S(\rho)S(d\varphi) \approx S(\rho)d\varphi.$$

Consider right triangle $\triangle ABC$ with catheti $a, b$, hypotenuse $c$, internal angle $\alpha$ and external angle $\beta'$. Place it so that $A = E$ is origin and $AB = e^1$ concides with coordonate vector. By proposition 2.8.2, this thiangle area is 3 times as much as volume of pyramid with this triangle base and vertex in $O = (0, 0, 0)$ origin of $\mathbb{R}^3$. The pyramid volume equals to sum of volumes of pyramids with quadrangular base $d\rho \times dl$ and height equal to sphere radius 1, each of them having the volume:

$$dv = \frac{1}{3}dld\rho \cdot 1 = \frac{1}{3}S(\rho)d\varphi d\rho.$$

It means, that area $s$ of triangle $\triangle ABC$ equals to:

$$s = \int_{\triangle ABC} S(\rho)d\varphi d\rho.$$

Moreover, by (2.45) triangle side $BC$ in polar coordinate system has equation:

$$T_1(b) = T_1(\rho)C_2(\varphi).$$

It follows:

$$s = \int_0^\alpha d\varphi \int_0^{T_1^{-1}\left(\frac{T_1(b)}{C_2(\varphi)}\right)} S_1(\rho)d\rho.$$

For further analysis we need a list of derivatives of generalized trigonometric and their inverse functions. It is easy to see that:

$$S'(x) = C(x), \tag{2.63}$$

$$C'(x) = -kS(x), \tag{2.64}$$

$$T'(x) = \left(\frac{S(x)}{C(x)}\right)' = \frac{S'(x)C(x) - S(x)C'(x)}{C^2(x)} = \frac{C^2(x) + kS^2(x)}{C^2(x)} = \frac{1}{C^2(x)}. \tag{2.65}$$

Inverse functions:

$$\begin{aligned} y &= S^{-1}(x), \\ x &= S(y), \\ dx &= C(y)dy, \\ \frac{dy}{dx} &= \frac{1}{C(y)} = \frac{1}{\sqrt{1 - kS^2(y)}} = \frac{1}{\sqrt{1 - kx^2}}; \end{aligned} \tag{2.66}$$

$$\begin{aligned} y &= C^{-1}(x), \\ x &= C(y), \\ dx &= -kS(y)dy, \\ \frac{dy}{dx} &= -\frac{1}{kS(y)} = -\frac{1}{k\sqrt{\frac{1-C^2(y)}{k}}} = -\frac{1}{\sqrt{k(1 - x^2)}}; \end{aligned} \tag{2.67}$$

$$\begin{aligned} y &= T^{-1}(x), \\ x &= T(y), \\ dx &= \frac{dy}{C^2(y)}, \\ \frac{dy}{dx} &= C^2(y) = \frac{1}{1 + kT^2(y)} = \frac{1}{1 + kx^2.} \end{aligned} \tag{2.68}$$

So:

$$\begin{aligned} s &= \int_0^\alpha d\varphi \int_0^{T_1^{-1}\left(\frac{T_1(b)}{C_2(\varphi)}\right)} S_1(\rho)d\rho = -\frac{1}{k_1}\int_0^\alpha \left( C_1(\rho)\Big|_0^{T_1^{-1}\left(\frac{T_1(b)}{C_2(\varphi)}\right)}\right) d\varphi \\ &= -\frac{1}{k_1}\int_0^\alpha \left( C_1\left(T_1^{-1}\left(\frac{T_1(b)}{C_2(\varphi)}\right)\right) - C_1(0)\right) d\varphi. \end{aligned}$$

Table 2.3: Derivatives of generalized trigonometric and their inverse functions

| | |
|---|---|
| $S'(x) = C(x)$ | $(S^{-1})'(x) = \dfrac{1}{\sqrt{1 - kx^2}}$ |
| $C'(x) = -kS(x)$ | $(C^{-1})'(x) = -\dfrac{1}{\sqrt{k(1 - x^2)}}$ |
| $T'(x) = \dfrac{1}{C^2(x)}$ | $(T^{-1})'(x) = \dfrac{1}{1 + kx^2}$ |

Next, observe that:

$$C(x) = \frac{1}{\sqrt{1 + kT^2(x)}},$$

obtain, using (2.54):

$$\begin{aligned} s &= \frac{1}{k_1}\int_0^\alpha d\varphi - \frac{1}{k_1}\int_0^\alpha \frac{d\varphi}{\sqrt{1 + k_1\frac{T_1^2(b)}{C_2^2(\varphi)}}} = \frac{\alpha}{k_1} - \frac{1}{k_1}\int_0^\alpha \frac{C_2(\varphi)d\varphi}{\sqrt{C_2^2(\varphi) + k_1T_1^2(b)}} \\ &= \frac{\alpha}{k_1} - \frac{1}{k_1}\int_0^\alpha \frac{d(S_2(\varphi))}{\sqrt{1 - k_2S_2^2(\varphi) + k_1T_1^2(b)}} = \frac{\alpha}{k_1} - \frac{1}{k_1}\int_0^\alpha \frac{d(S_2(\varphi))}{\sqrt{\frac{1}{C_1^2(b)} - k_2S_2^2(\varphi)}} \\ &= \frac{\alpha}{k_1} - \frac{1}{k_1}\int_0^\alpha \frac{C_1(b)d(S_2(\varphi))}{\sqrt{1 - k_2C_1^2(b)S_2^2(\varphi)}} = \frac{\alpha}{k_1} - \frac{1}{k_1}\int_0^\alpha \frac{d(C_1(b)S_2(\varphi))}{\sqrt{1 - k_2(C_1(b)S_2(\varphi))^2}} \\ &= \frac{\alpha}{k_1} - \frac{1}{k_1}S_2^{-1}(C_1(b)S_2(\varphi))\Big|_0^\alpha = \frac{\alpha}{k_1} - \frac{1}{k_1}S_2^{-1}(C_1(b)S_2(\alpha)) \\ &= \frac{\alpha}{k_1} - \frac{1}{k_1}S_2^{-1}(S_2(\beta')) = \frac{\alpha}{k_1} - \frac{\beta'}{k_1} = \frac{\alpha - \beta'}{k_1}. \end{aligned} \tag{2.69}$$

### 2.8.3. Area of the right triangle.

In order to deduce general form of right triangle area equation, use equations (2.69, 2.48, 2.51). For the moment, restrict ourselves to case $k_1 = \pm 1$. Because the function $T(x)$ is odd, $T(-x) = -T(x)$, and thus with our restriction:

$$T_2\left(\frac{\alpha - \beta'}{k_1}\right) = \frac{T_2(\alpha - \beta')}{k_1}.$$

Then:

$$T_2(s) = \frac{T_2(\alpha - \beta')}{k_1} = \frac{T_2(\alpha) - T_2(\beta')}{k_1(1 + k_2 T_2(\alpha) T_2(\beta'))} = \frac{\frac{T_{12}(a)}{S_1(b)} - \frac{S_{12}(a)}{T_1(b)}}{k_1\left(1 + k_2 \frac{T_{12}(a)}{S_1(b)} \frac{S_{12}(a)}{T_1(b)}\right)}$$
$$= \frac{T_{12}(a)T_1(b) - S_{12}(a)S_1(b)}{K_2 T_{12}(b) S_{12}(a) + K_1 T_1(b) S_1(b)}.$$

Next note, that when $k_1 = \pm 1, k_1^2 = 1$, thus $K_1 K_1 = k_1^2 k_2 = k_2$. Then:

$$S_2^2(s) = \frac{T_2^2(s)}{1 + k_2 T_2^2(s)} = \frac{T_2^2(s)}{1 + K_1 K_2 T_2^2(s)} = \frac{\left(\frac{T_{12}(a)T_1(b) - S_{12}(a)S_1(b)}{K_2 T_{12}(b) S_{12}(a) + K_1 T_1(b) S_1(b)}\right)^2}{1 + K_1 K_2 \left(\frac{T_{12}(a)T_1(b) - S_{12}(a)S_1(b)}{K_2 T_{12}(b) S_{12}(a) + K_1 T_1(b) S_1(b)}\right)^2}$$
$$= \frac{(T_{12}(a)T_1(b) - S_{12}(a)S_1(b))^2}{(K_2 T_{12}(a) S_{12}(a) + K_1 T_1(a) S_1(b))^2 + K_1 K_2 (T_{12}(a)T_1(b) - S_{12}(a)S_1(b))^2}$$
$$= \frac{(T_{12}(a)T_1(b) - S_{12}(a)S_1(b))^2}{K_2^2 T_{12}^2(a) S_{12}^2(a) + K_1^2 T_1^2(b) S_1^2(b) + K_1 K_2 T_{12}^2(a) T_1^2(b) + K_1 K_2 S_{12}^2(a) S_1^2(b)}$$
$$= \frac{(T_{12}(a)T_1(b) - S_{12}(a)S_1(b))^2}{K_2 S_{12}^2(a)(K_2 T_{12}^2(a) + K_1 S_1^2(b)) + K_1 T_1^2(b)(K_2 T_{12}^2(a) + K_1 S_1^2(b))}$$
$$= \frac{(T_{12}(a)T_1(b) - S_{12}(a)S_1(b))^2}{(K_2 T_{12}^2(a) + K_1 S_1^2(b))(K_1 T_1^2(b) + K_2 S_{12}^2(a))}$$
$$= \frac{(T_{12}(a)T_1(b) - S_{12}(a)S_1(b))^2}{(1 + K_2 T_{12}^2(a) - (1 - K_1 S_1^2(b)))(1 + K_1 T_1^2(b) - (1 - K_2 S_{12}^2(a)))}$$
$$= \frac{(T_{12}(a)T_1(b) - S_{12}(a)S_1(b))^2}{\left(\frac{1}{C_{12}^2(a)} - C_1^2(b)\right)\left(\frac{1}{C_1^2(b)} - C_{12}^2(a)\right)} = \frac{(T_{12}(a)T_1(b) - S_{12}(a)S_1(b))^2 C_{12}^2(a) C_1^2(b)}{(1 - C_{12}^2(a) C_1^2(b))^2}.$$

Now:

$$\begin{aligned} S_2(s) &= \frac{(T_{12}(a)T_1(b) - S_{12}(a)S_1(b)) C_{12}(a) C_1(b)}{(1 - C_{12}^2(a) C_1^2(b))} \\ &= \frac{S_{12}(a)S_1(b)(1 - C_{12}(a)C_1(b))}{(1 + C_{12}(a)C_1(b))(1 - C_{12}(a)C_1(b))} = \frac{S_{12}(a)S_1(b)}{1 + C_{12}(a)C_1(b)}. \end{aligned} \tag{2.70}$$

Equation (2.70) doesn't contain explicit parameters $k_1, k_2$ and is true when $k_1 = \pm 1$. Now let $k_1 = 0$. Equation becomes (omitting functions $C(x), S(x)$ indices):

$$s = S(s) = \frac{S(a)S(b)}{1 + C(a)C(b)} = \frac{ab}{1 + 1 \cdot 1} = \frac{ab}{2}.$$

This form coincides with (2.62) for Euclidean, Galilean and Minkowski planes, where $k_1 = 0$. So,

equations (2.70) is true for all 9 homogeneous planes.

### 2.8.4. Type of the area.

From the above follows, that the triangle area enters in equations as other geometric quantities, lengths and angles. It means, that the area is a measure, that has some type. From equations (2.69) follows that when $k_1 = \pm 1$ the area is commensurate with angles, so in these cases its type equals to angular type. Obviously it isn't true for planes where $k_1 = 0$. On Euclidean plane $\{0, 1\}$ and on Minkowski plane $\{0, -1\}$ exist transformations of distance dilation, it means that area type equals to $k = 0$, while angular type on these planes equal to $k_2 = 1$ and $k_2 = -1$ respectively. During deduction of (2.70) we assume that area type is $k = K_1 K_2$. This assumption leads to correct result including cases $k_1 = K_1 = K_1 K_2 = 0$ (in this case $s = S(s)$). Obviously, all areas of some space have the same type. So, the following theorem is proved:

**Theorem 2.8.3.** *Area type $k$ equals to product of coordinate vectors types:*

$$k = K_1 K_2.$$

*Remark.* Orthonormalized basis of 2-dimensional plane $\mathbb{B}^2$ consists of 3 vectors $\{e^0, e^1, e^2\}$. But vector $e^0$ type always equals to $K_0 = 1$. In this way, $k = K_1 K_2 = K_0 K_1 K_2$.

The theorem 2.8.3 explains why during composing the area integral, the polar coordinates (containing angles) were chosen. For considered cases area is commensurate with angles and it is easier to find exactly area, than some function on area.

**Proposition 2.8.4.** *Volume type $k = 0$ if and only if space signature $\{k_1, ..., k_n\}$ contains at least one zero element $k_i = 0$, $i = \overline{1, n}$.*

*Proof.* When some signature element $k_i = 0$, $i = \overline{1, n}$, then space allows transformation of dilation of $(i - 1)$–dimensional angles (or distances, if $i = 1$). This transformation beside the angles, changes the volume. Thus, the volume type $k = 0$.

On the other hand, when all elements $k_i \neq 0, i = \overline{1, n}$, then any measure defines all others, and with them, also the volume. The space is free of measure or volume dilation transformation. Thus, the volume type $k \neq 0$. □

**Proposition 2.8.5.** *If in space with signature $\{k_1, ..., k_n\}$ some element $k_i = 0, i = \overline{1, n}$, then the transformation of $(i - 1)$–dimensional angles dilation by $\alpha$ times changes volume by $\alpha^{n-i+1}$ times.*

*Proof.* In the origin $E = (1 : 0 : ... : 0)$ the volume element $dv = \prod_{j=1}^{n} dx^j$. The transformation of $(i - 1)$–dimensional angles dilation by $\alpha$ times changes all $dx^j$, $j = \overline{i, n}$ by $\alpha$ times. So, new

volume element:

$$dv' = \prod_{j=1}^{n} dx'^j = \left(\prod_{j=1}^{i-1} dx'^j\right)\left(\prod_{j=i}^{n} dx'^j\right) = \left(\prod_{j=1}^{i-1} dx^j\right)\left(\prod_{j=i}^{n} (\alpha dx^j)\right)$$

$$= \alpha^{n-i+1} \prod_{j=1}^{n} dx^j = \alpha^{n-i+1} dv$$

Because all points of homogeneous space are congruent, all volumes changes the same way as volume element in origin. □

Proposition 2.8.5 is interesting because the volume dilation factor exponent coincides with exponent of $k_j$ in expression:

$$\prod_{i=1}^{n} K_i = \prod_{i=1}^{n}\left(\prod_{j=1}^{i} k_j\right) = \prod_{j=1}^{n} k_j^{n-j+1}$$

At the end, one-dimensional volume of segment is its length. It means that for vectors $e^i$, the volume type is $k = K_i$. For two-dimensional case area type $k = K_i K_j$, where $K_i, K_j$ are the types of basis vectors $e^i, e^j$. In the same order of things, theorem 2.8.3 with propositions 2.8.4, 2.8.5 allow to formulate the following conjecture.

**Conjecture 2.8.6.** *The type of volume in space $\mathbb{B}^n$ with signature $\{k_1, ..., k_n\}$ equals to*

$$k = \prod_{i=1}^{n} K_i.$$

### 2.8.5. Other equations of area of the right triangle.

Starting from equation (2.70) and theorem 2.8.3, calculate:

$$
\begin{aligned}
C^2(s) &= 1 - kS^2(s) = 1 - \frac{K_1K_2S^2(a)S^2(b)}{(1+C(a)C(b))^2}\\
&= \frac{1 + 2C(a)C(b) + C^2(a)C^2(b) - K_1K_2S^2(a)S^2(b)}{(1+C(a)C(b))^2}\\
&= \frac{(C^2(a) + K_1S^2(a))(C^2(b) + K_2S^2(b))}{(1+C(a)C(b))^2}\\
&+ \frac{2C(a)C(b) + C^2(a)C^2(b) - K_1K_2S^2(a)S^2(b)}{(1+C(a)C(b))^2}\\
&= \frac{C^2(a)C^2(b) + K_1S^2(a)C^2(b) + K_2C^2(a)S^2(b) + K_1K_2S^2(a)S^2(b)}{(1+C(a)C(b))^2}\\
&+ \frac{2C(a)C(b) + C^2(a)C^2(b) - K_1K_2S^2(a)S^2(b)}{(1+C(a)C(b))^2}\\
&= \frac{C^2(a)(C^2(b) + K_2S^2(b)) + C^2(b)(C^2(a) + K_1S^2(a)) + 2C(a)C(b)}{(1+C(a)C(b))^2}\\
&= \frac{C^2(a) + 2C(a)C(b) + C^2(b)}{(1+C(a)C(b))^2} = \left(\frac{C(a)+C(b)}{1+C(a)C(b)}\right)^2,
\end{aligned}
$$

$$C(s) = \frac{C(a)+C(b)}{1+C(a)C(b)}. \tag{2.71}$$

And:

$$T(s) = \frac{S(s)}{C(s)} = \frac{S(a)S(b)}{1+C(a)C(b)}\frac{1+C(a)C(b)}{C(a)+C(b)} = \frac{S(a)S(b)}{C(a)+C(b)} \tag{2.72}$$

Next use:

$$T\left(\frac{x}{2}\right) = \frac{S(x)}{1+C(x)}.$$

Table 2.4: Right triangle area

| | |
|---|---|
| $S(s)$ | $= \dfrac{S(a)S(b)}{1+C(a)C(b)}$ |
| $C(s)$ | $= \dfrac{C(a)+C(b)}{1+C(a)C(b)}$ |
| $T(s)$ | $= \dfrac{S(a)S(b)}{C(a)+C(b)}$ |
| $T\left(\frac{s}{2}\right)=T\left(\frac{a}{2}\right)T\left(\frac{b}{2}\right)$ | |

It is easy to show it is true:

$$T(x) = T\left(2\frac{x}{2}\right) = \frac{2T\left(\frac{x}{2}\right)}{1-kT^2\left(\frac{x}{2}\right)} = \frac{2\frac{S(x)}{1+C(x)}}{1-k\left(\frac{S(x)}{1+C(x)}\right)^2} = \frac{2S(x)(1+C(x))}{(1+C(x))^2-kS^2(x)}$$
$$= \frac{2S(x)(1+C(x))}{1+2C(x)+C^2(x)-kS^2(x)} = \frac{2S(x)(1+C(x))}{2C(x)+C^2(x)+C^2(x)} = \frac{2S(x)(1+C(x))}{2C(x)(1+C(x))}$$
$$= \frac{S(x)}{C(x)} = T(x).$$

Calculate:

$$\begin{aligned} T\left(\frac{s}{2}\right) &= \frac{S(s)}{1+C(s)} = \frac{\frac{S(a)S(b)}{1+C(a)C(b)}}{1+\frac{C(a)+C(b)}{1+C(a)C(b)}} = \frac{S(a)S(b)}{1+C(a)C(b)+C(a)+C(b)} \\ &= \frac{S(a)S(b)}{(1+C(a))(1+C(b))} = \frac{S(a)}{1+C(a)}\frac{S(b)}{1+C(b)} = T\left(\frac{a}{2}\right)T\left(\frac{b}{2}\right). \end{aligned} \tag{2.73}$$

### 2.9. Conclusions of chapter 2

In this chapter the analytic geometry of homogeneous spaces was constructed. It includes:

1. The new notion of signature was introduced. With its aim, the model of homogeneous space with given signature was constructed. That permits the classification of homogeneous spaces based on the concept of signature. Also, for each known homogeneous space (spaces of constant curvature, Galilean, Minkowskii, De Sitter, Anti de Sitter among other spaces) its signature, and its place in the presented classification, were found [56, 57].

2. In dependence of the space signature, the parameterized form of some important axioms

were given. Based on them, the formulation and the proof of the theorems, parameterized by signature, is also possible in this unified manner. Generalized trigonometric functions were introduced by means of the signature. These functions make it possible to find universal form of trigonometric equations, common for each homogeneous space, also intruduced here [57, 59, 62].

3. By the new concept of the group of generalized orthogonal matrix, the isometry group of any homogeneous space was described. Conform to Felix Klein's concept of geometry, described by him in the well known Erlangen Program, the isometry group of the space determines the geometry of that space [63].

4. Via the new concept of signature, the type (elliptic, parabolic and hyperbolic) of geometric quantities (distances, angles, areas and volumes) was established [58, 60].

For distances, angles and areas the complete research is presented, for volumes only the parabolic type was rigorously established and the general case was conjectured.

# 3. APPLICATION OF THEORY

This chapter describes some applications of analytic geometry of homogeneous spaces to different areas of mathematics: algebraic geometry, topology and differential geometry.

## 3.1. Algebraic geometry

Theory described in previous chapter is universal and easy to use. However its usage is prevented by the fact geometric spaces are usually defined and classified in a different way than described here.

### 3.1.1. Signature of a space and a lineal.

In order to not loose the feeling of reality, describe algorithm that allows to find the specefication of all geometric spaces. The algorithm is applicable to any space in which the notions of points, lines, planes, subspaces, distances, angles and / or motions are defined. For example, it helps to find signatures of the 8 Thirston geometries [75] or signatures of 11 classes of three–dimensional algebra Lie in Bianchi classification [209, 100].

**Algorithm 3.1** (Space signature)**.**

1. Let $m$ be equal to the greatest number of general situated points, or equivalently, the least number of vertices of positive volumed polyhedron.
2. Calculate space dimension as $n = m - 1$.
3. For $i = \overline{1, n}$ (0-dimensional places are points, 1-dimensional planes are lines), do:
   (a) If among $(i - 1)$–dimensional planes there exist non-congruent, then the space definition or used terminology is inconsequent. The theory still can be used, but for its correct usage it is necessary either to change the terminology, or to define some space elements otherwise.
   (b) If measure between $(i - 1)$–dimensional planes is bounded, then $k_i = 1$.
   (c) If measure between $(i - 1)$–dimensional planes allows dilation, then $k_i = 0$.
   (d) Otherwise $k_i = -1$.
4. Use the theory knowing space dimension $n$ and signature $\{k_1, ..., k_n\}$.

The necessity of proper terminology, uniform among all spaces, is driven by desire to have a theory that isn't misleading and aims to analyse the space structure and to compare it with other spaces. At the same time, the inconsequent theory / terminology is not the contradictory one, but a theory / terminology that is improper with respect to theory notions. It assumes the following:

- All planes of any dimension are congruent, including points and lines.
- Theory admits the duality principle of $(m-1)$–dimensional planes and $(n-m)$–dimensional planes.

It should be acknowledged, that "proper notions" may change over the time. In order to understand what it is about, consider an example of inconsequent terminology. Minkowski space is successfully used in physics for Special Relativity Theory study. Unfortunately, from geometric point of view this space has no proper terminology. The notions "spacelike lines", "timelike lines" and "light lines" make sense in physics, but not in geometry. Corresponding geometric notions would be: "first sort lines", "second sort lines" and "third sort lines". No space motion maps a line of one sort to a line of another sort. There is no contradiction here, but there is inconsequence. What if somebody wishes to define a space with five sorts of lines? Depending on space, it may exist more different sorts of two–dimensional planes. The greater is plane dimension, the more different sorts of them can exist. Nobody defines different sorts of points. All points are congruent. Why lines should not be congruent? At the other hand, *relative position* of point pairs may be different.

By declaration of only first sort lines as lines, we need to exclude second and third sorts of lines from line set. At first sight it contradicts with familiar axiom, that it may be drawn a line through any two points. But this axiom may not be true in certain spaces. Moreover, just Euclidean plane, where all points are connectable, allows lines to be parallel (having no common point). By duality principle, it should exist the notion of unconnectable points (having no common line).

### 3.1.2. Signature of semi–Euclidean and semi–Riemannian spaces.

**Definition 3.1.1** (Semi–Euclidean and pseudo–Euclidean space)**.** A linear space $\mathbb{R}^n$ with vector product $\bar{x}, \bar{y} \in \mathbb{R}^n, \bar{x} = \{x_1, ..., x_n\}, \bar{y} = \{y_1, ..., y_n\}$ defined as:

$$\bar{x} \cdot \bar{y} = \sum_{i=1}^{m} x_i y_i - \sum_{i=m+1}^{n-d} x_i y_i \tag{3.1}$$

is called *semi–Euclidean* space ${}^d\mathbb{E}^m_n$ (or ${}^{(d)}\mathbb{E}^m_n$) of dimension $n \geq m+d$ with positive inertion index $m$ and defect $d$. Particularly, if $d = 0$, the space is called *pseudo–Euclidean* $\mathbb{E}^{m,n-m}$ with classical signature $(m, n-m)$.

**Definition 3.1.2** (semi–Riemannian and pseudo–Riemannian space)**.** (Not necessarily linear) space ${}^d\mathbb{V}^m_n$ (or ${}^{(d)}\mathbb{V}^m_n$) is called *semi–Riemannian*, if at each point its tangent space is semi–Euclidean ${}^d\mathbb{E}^m_n$. Particularly, if $d = 0$, the space is called *pseudo–Riemannian* $\mathbb{V}^{m,n-m}$ with classical signature $(m, n-m)$.

In order to determine the signature (in the sense of described theory) of semi–Euclidean space

${}^d\mathbb{E}^m_n$, rewrite its bilinear form of vector product as follows:

$$\bar{x}\cdot\bar{y}=\sum_{i=1}^{n}p_i x_i y_i;\quad p_i=\begin{cases}\pm1, & i=\overline{1,n-d},\\ 0, & i=\overline{n-d+1,n}.\end{cases}\tag{3.2}$$

*Remark.* In semi–Euclidean spaces, by convention, basis vectors with positive signs are written first, because reindexing of basis vectors doesn't change space properties. However, in homogeneous spaces we agreed to distinguish equivalent vectors (there is a motion that interchanges them up to sign) from interchangeable vectors (there is isomorphic space where corresponding vectors are interchanged). By theorem 2.4.17 on equivalence of coordinate vectors, their order is important. In any cases the first coefficient should be positive and all vectors with zero coefficient should be the last. Otherwise, it can be found contradictions in semi–Euclidean space definition in terms of analytic geometry of homogeneous spaces.

Construct homogeneous space $\mathbb{B}^n$ so, that $k_1=0, K_{1i}=p_i, i=\overline{2,n}$; $(p_1=1)$. Its signature is $\left\{0,\frac{p_2}{p_1},\frac{p_3}{p_2},\ldots,\frac{p_n}{p_{n-1}}\right\}$. It is always possible, because when $i=\overline{2,(n-d+1)}$, the denominator of fraction $\frac{p_i}{p_{i-1}}$ is different from zero and then the element $k_i$ has finite value, and when $i=\overline{(n-d+2),n}$, the element $k_i=\frac{p_i}{p_{i-1}}=\frac{0}{0}$ is undefined: 1, 0 or $-1$.

From $k_1=0$ in $\mathbb{B}^n$ follows that $K_0=1, K_i=0, i>0$. From sphere equation $x\odot x=x_0^2=1$ follows that $x_0=1$ or $x_0=-1$. Since $x\in\mathbb{B}^n$ and also $-x\in\mathbb{B}^n$, we can consider that $x_0=1$. We can consider $\mathbb{B}^n$ as hyperplane of $\mathbb{R}^{n+1}$ with equation $x_0=1$ and identify it with ${}^d\mathbb{E}^m_n$. In this construction, each vector $\bar{x}\in{}^d\mathbb{E}^m_n$; $\bar{x}=\{x_1,\ldots,x_n\}$ corresponds to some vector $x\in\mathbb{B}^n$; $x=\{1:x_1:\ldots:x_n\}$.

**Lemma 3.1.1.** *Nonzero distance $d(\bar{x},\bar{y})$ between vectors from semi–Euclidean space ${}^d\mathbb{E}^m_n$ is related to distance $d(x,y)$ between corresponding vectors from homogeneous space $\mathbb{B}^n$ with signature $\left\{0,\frac{p_2}{p_1},\frac{p_3}{p_2},\ldots,\frac{p_n}{p_{n-1}}\right\}$ by the following equality:*

$$d(x,y)=\begin{cases}d(\bar{x},\bar{y}), & d(\bar{x},\bar{y})\in\mathbb{R},\\ -i\,d(\bar{x},\bar{y}), & d(\bar{x},\bar{y})\in i\mathbb{R}.\end{cases}$$

*Proof.* If all coordinates of $\bar{x}$ from 1-th to $(n-d)$-th equal to corresponding coordinates of $\bar{y}$, then $d(\bar{x},\bar{y})=0$. Let among coordinates from 1-th to $(n-d)$-th of vectors $\bar{x},\bar{y}$ there are different ones.

Because vectors $x, y \in \mathbb{B}^n$ are points, have:

$$
\begin{aligned}
&x \odot x = 1, \\
&y \odot y = 1, \\
&x \odot y = x_0 y_0 + k_1 \sum_{i=1}^{n} K_{1i} x_i y_i = 1 \cdot 1 + 0 \sum_{i=1}^{n} K_{1i} x_i y_i = 1.
\end{aligned}
$$

Next note, that zero product:

$$(y - x) \odot_0 (y - x) = x \odot_0 x + y \odot_0 y - 2x \odot_0 y = 1 + 1 - 2 = 0.$$

It means, that for vector $y - x$ the natural product is:

$$(y - x) \odot_i (y - x) = \frac{1}{K_{1i}} (y - x) \odot_1 (y - x) = \pm (y - x) \odot_1 (y - x).$$

Having $p_1 = 1$, find:

$$
\begin{aligned}
(y - x) \odot (y - x) = \pm (y - x) \odot_1 (y - x) &= \pm \left( \frac{(y_0 - x_0)^2}{k_1} + \sum_{i=1}^{n} K_{1i} (y_i - x_i)^2 \right) \\
&= \pm \left( \frac{(1 - 1)^2}{0} + \sum_{i=1}^{n} \frac{p_i}{p_1} (y_i - x_i)^2 \right) = \pm \frac{1}{p_1} \sum_{i=1}^{n} p_i (y_i - x_i)^2 \\
&= \pm (\bar{y} - \bar{x}) \cdot (\bar{y} - \bar{x}) = \pm d^2(\bar{x}, \bar{y}).
\end{aligned}
$$

By algorithm 2.8 of measure computing between lineals, in order to find $d(x, y)$, we need to construct projections $y', y''$ of vector $y$ on vector $x$ and on its orthogonal complement in lineal constructed on $\{x, y\}$:

$$
\begin{aligned}
y' &= (x \odot y) x = 1 \cdot x = x, \\
y'' &= y - y' = y - x.
\end{aligned}
$$

So, having distance measure type $k_1 = 0$:

$$
\begin{aligned}
&C_1^2(d(x, y)) = y' \odot y' = x \odot x = 1, \\
&\quad d^2(x, y) = S_1^2(d(x, y)) = y'' \odot y'' = (y - x) \odot (y - x) = \pm d^2(\bar{x}, \bar{y}), \\
&\qquad d(x, y) = \begin{cases} d(\bar{x}, \bar{y}), & d(\bar{x}, \bar{y}) \in \mathbb{R}, \\ -i\, d(\bar{x}, \bar{y}), & d(\bar{x}, \bar{y}) \in i\mathbb{R}. \end{cases}
\end{aligned}
$$

□

*Remark.* In semi–Euclidean space the distance may be real or imaginary, and may exist different vectors with zero distance between them. In homogeneous space $\mathbb{B}^n$ corresponding to semi–Euclidean space ${}^d\mathbb{E}^m_n$ exists real measure between any vectors that has zero value only for equal vectors.

*Remark.* In case of semi–Euclidean space with defect $d > 1$, homogeneous space signature and geometry are ambiguous. The space allows dilation transformation of angles, and distances don't define all its metric.

**Corollary.** *Homogeneous space $\mathbb{B}^n$ with signature $\{k_1, ..., k_n\}$ is embedded in homogeneous meta–space $\mathbb{R}^{n+1}$ with signature $\{0, k_1, ..., k_n\}$.*

*Proof.* Let $x = \{x_0 : ... : x_n\}, y = \{y_0 : ... : y_n\} \in \mathbb{B}^n \subset \mathbb{R}^{n+1}$ with signature $\{k_1, ..., k_n\}$. Then:

$$x \odot y = \sum_{i=0}^{n} K_i x_i y_i.$$

Considering the space $\mathbb{R}^{n+1}$ as semi–Euclidean, for vectors $\bar{x}' = \{x'_1, ..., x'_{n+1}\}, \bar{y}' = \{y'_1, ..., y'_{n+1}\} \in \mathbb{R}^{n+1}$, $x'_{i+1} = x_i$, $y'_{i+1} = y_i$, $i = \overline{0, n}$ have distance bilinear form:

$$\bar{x} \cdot \bar{y} = \sum_{i=1}^{n+1} p_i x'_i y'_i = \sum_{i=0}^{n} p_{i+1} x_i y_i,$$

where $p_{i+1} = K_i, i = \overline{0, n}$.

Next, by lemma 3.1.1, considering semi–Euclidean space $\mathbb{R}^{n+1}$ as linear homogeneous with signature $\{k'_1 = 0, k'_2, ..., k'_{n+1}\}$, get:

$$K'_{1\,i+1} = p_{i+1} = K_i, \quad i = \overline{0, n},$$
$$k'_{i+1} = \frac{K'_{1\,i+1}}{K'_{1i}} = \frac{K_i}{K_{i-1}} = k_i, \quad i = \overline{1, n}.$$

In this way, the signature of homogeneous meta–space $\mathbb{R}^{n+1}$ is $\{0, k_1, ..., k_n\}$. □

**Lemma 3.1.2.** *Semi–Riemannian space ${}^d\mathbb{V}^m_n$ corresponds (in sense of metric) to homogeneous space with signature $\left\{k_1, \frac{p_2}{p_1}, ..., \frac{p_n}{p_{n-1}}\right\}$, where $k_1 = 0$ for linear, $k_1 = 1$ for positive curved and $k_1 = -1$ for negative curved semi–Riemannian space, and $p_i$, $i = \overline{1, n}$ are coefficients of its metric tensor bilinear form.*

*Proof.* Consider vectors $x, y \in \mathbb{B}^n \subset \mathbb{R}^{n+1}$, corresponding to vectors $\bar{x}, \bar{y} \in {}^d\mathbb{V}^m_n$. Using algorithm 2.8 of measure computing between lineals, find projections $x', x''$ of vector $x$ on $y$ and on its orthogonal complement in lineal constructed on $\{x, y\}$:

$$x' = (x \odot y)y,$$
$$x'' = x - x' = x - (x \odot y)y.$$

Let distance $d(x, y) \to 0$ be small. Then:

$$C^2(d(x, y)) = x' \odot x' = ((x \odot y)y) \odot ((x \odot y)y) = (x \odot y)^2(y \odot y) = (x \odot y)^2,$$
$$C(d(x, y)) = x \odot y.$$

By (2.3):

$$x \odot y = C(d(x, y)) = 1 + O(d^2(x, y)).$$

Next:

$$\begin{aligned}
S^2(d(x, y)) &= x'' \odot x'' = (x - (x \odot y)y) \odot (x - (x \odot y)y) \\
&= (x - (1 + O(d^2(x, y)))y) \odot (x - (1 + O(d^2(x, y)))y) \\
&= (x - y + O(d^2(x, y))y) \odot (x - y + O(d^2(x, y))y) \\
&= (x - y) \odot (x - y) + O(d^2(x, y))(x - y) \odot y + O(d^4(x, y))(y \odot y) \\
&= (x - y) \odot (x - y) + O(d^2(x, y))(x \odot y - y \odot y) + O(d^4(x, y)) \\
&= (x - y) \odot (x - y) + O(d^2(x, y))(1 + O(d^2(x, y)) - 1) + O(d^4(x, y)) \\
&= (x - y) \odot (x - y) + O(d^4(x, y)) + O(d^4(x, y)) \\
&= (x - y) \odot (x - y) + O(d^4(x, y)).
\end{aligned}$$

On the other hand, by (2.4):

$$\begin{aligned}
S^2(d(x, y)) &= (d(x, y) + O(d^3(x, y)))^2 = d^2(x, y) + d(x, y)O(d^3(x, y)) + O(d^6(x, y)) \\
&= d^2(x, y) + O(d^4(x, y)).
\end{aligned}$$

And then:

$$d^2(x, y) = (x - y) \odot (x - y) + O(d^4(x, y)).$$

Moreover, by proposition 2.8.1, when $d(x, y) \to 0$, vector $(x - y)$ is orthogonal to vectors $x$ and $y$. Because vectors $x$ and $y$ indices equal to 0, index of $(x - y)$ equals to at least 1. It means, that their natural product:

$$(x-y)\odot_i(x-y)=\frac{1}{K_{1i}}(x-y)\odot_1(x-y)=\pm(x-y)\odot_1(x-y).$$

If we choose the basis with center in $y$, then $y=\{1:0:\ldots:0\}$ and:

$$C(d(x,y))=x\odot y=\sum_{i=0}^{n}K_ix_iy_i=x_0\cdot 1+\sum_{i=1}^{n}K_ix_i\cdot 0=x_0.$$

Hence (again assume $p_1=1$):

$$\begin{aligned}
d^2(x,y)&=\pm((x-y)\odot_1(x-y))+O(d^4(x,y))\\
&=\pm\left(\frac{(x_0-y_0)^2}{k_1}+\sum_{i=1}^{n}K_{1i}(x_i-y_i)^2\right)+O(d^4(x,y))\\
&=\pm\left(\frac{(C(d(x,y))-1)^2}{k_1}+\sum_{i=1}^{n}\frac{p_i}{p_1}(x_i-y_i)^2\right)+O(d^4(x,y))\\
&=\pm\left(\frac{(1+O(d^2(x,y))-1)^2}{k_1}+\frac{1}{p_1}\sum_{i=1}^{n}p_i(x_i-y_i)^2\right)+O(d^4(x,y))\\
&=\pm\sum_{i=1}^{n}p_i(x_i-y_i)^2+O(d^4(x,y))=\pm d^2(\bar{x},\bar{y})+O(d^4(x,y)).
\end{aligned}$$

In this way, when $d(x,y)\to 0$:

$$\begin{aligned}
d^2(x,y)&=\pm d^2(\bar{x},\bar{y}),\\
d(x,y)&=\begin{cases}d(\bar{x},\bar{y}), & d(\bar{x},\bar{y})\in\mathbb{R},\\ -i\,d(\bar{x},\bar{y}), & d(\bar{x},\bar{y})\in i\mathbb{R}.\end{cases}
\end{aligned}$$

Finally, note that $k_1=1$ introduce elliptic distance measure (the space is positive curved), $k_1=0$ introduces parabolic distnce measure (linear space) and $k_1=-1$ introduces hyperbolic distance measure (negative curved space). □

*Remark.* Algorithm 2.8 of measure computing between lineals is based on linear algebra apparatus of metaspace and nowhere assumes this measure is small. For definition of metric in semi–Riemannian spaces in terms of metric tensor, the distance differentiability is essential, that enforce usage of differential geometry apparatus.

**Corollary.** *At small distances, non–linear homogeneous space with signature* $\{\pm 1,k_2,\ldots k_n\}$ *is best approximated by linear space with signature* $\{0,k_2,\ldots,k_n\}$.

We can use this bilinear form to find whole signature except $k_1$. Use this method to describe some special spaces by finding their signature.

**Elliptic, Euclidean and hyperbolic spaces.** Elliptic, linear (Euclidean) and hyperbolic (Bolyai–Lobachevsky) spaces are approximated by Euclidean one. We can find all the rest of elements of signature using Euclidean bilinear form. Let dimension be 3:

$$\begin{aligned} d^2(x,y) &= (y_1 - x_1)^2 + (y_2 - x_2)^2 + (y_3 - x_3)^2 \\ &= (y_1 - x_1)^2 + k_2(y_2 - x_2)^2 + k_2k_3(y_3 - x_3)^2 \end{aligned}$$

So, $k_2 = 1$ and $k_2k_3 = 1$, $k_3 = 1$. $k_1 = 1$ for elliptic, $k_1 = 0$ for Euclidean and $k_1 = -1$ for hyperbolic space.

**Minkowski, De Sitter and Anti de Sitter spaces.** The distance between $x$ and $y$ is defined (for timelike vectors) as:

$$\begin{aligned} d^2(x,y) &= (y_1 - x_1)^2 - (y_2 - x_2)^2 - (y_3 - x_3)^2 - (y_4 - x_4)^2 \\ &= (y_1 - x_1)^2 + k_2(y_2 - x_2)^2 + k_2k_3(y_3 - x_3)^2 + k_2k_3k_4(y_4 - x_4)^2, \end{aligned}$$

where coordinate 1 is time and coordinates 2, 3 and 4 are space ones. Then, $k_2 = -1$, $k_2k_3 = -1$, $k_3 = 1$ and $k_2k_3k_4 = -1$, $k_4 = 1$. Since Minkowski space is linear, $k_1 = 0$.

When introduce the curvature in space, its structure changes. For example, let $k_1 = 1$ (Anti de Sitter space). Then:

$$\begin{aligned} x \odot y &= x_0y_0 + k_1x_1y_1 + k_1k_2x_2y_2 + k_1k_2k_3x_3y_3 + k_1k_2k_3k_4x_4y_4 \\ &= x_0y_0 + x_1y_1 - x_2y_2 - x_3y_3 - x_4y_4 \end{aligned}$$

Thus time type $K_1$ becomes elliptic and space types $K_2, K_3, K_4$ become hyperbolic. When $k_1 = -1$ (De Sitter space), then:

$$\begin{aligned} x \odot y &= x_0y_0 + k_1x_1y_1 + k_1k_2x_2y_2 + k_1k_2k_3x_3y_3 + k_1k_2k_3k_4x_4y_4 \\ &= x_0y_0 - x_1y_1 + x_2y_2 + x_3y_3 + x_4y_4 \end{aligned}$$

Time type $K_1$ becomes hyperbolic and space types $K_2, K_3, K_4$ become elliptic.

**Galilean space.** The Galilean space–time metric is not defined as the metric of Minkowski space–time by means of distance bilinear form. Therefore its signature can't be obtained by previously described method. Still it is an interesting case. Find its signature from equations of its reference system change:

$$\begin{cases} t' &= t + \Delta t, \\ x_1' &= \alpha_1 x_1 + \beta_1 x_2 + \gamma_1 x_3 + v_1 t + \Delta x_1, \\ x_2' &= \alpha_2 x_1 + \beta_2 x_2 + \gamma_2 x_3 + v_2 t + \Delta x_2, \\ x_3' &= \alpha_3 x_1 + \beta_3 x_2 + \gamma_3 x_3 + v_3 t + \Delta x_3 \end{cases}$$

In this case the sub–matrix:

$$\begin{pmatrix} \alpha_1 & \beta_1 & \gamma_1 \\ \alpha_2 & \beta_2 & \gamma_2 \\ \alpha_3 & \beta_3 & \gamma_3 \end{pmatrix}$$

is elliptic GM-orthogonal. In the matrix form this equation system is rewritten as:

$$\begin{pmatrix} 1 \\ t' \\ x_1' \\ x_2' \\ x_3' \end{pmatrix} = \begin{pmatrix} 1 & 0 & 0 & 0 & 0 \\ \Delta t & 1 & 0 & 0 & 0 \\ \Delta x_1 & v_1 & \alpha_1 & \beta_1 & \gamma_1 \\ \Delta x_2 & v_2 & \alpha_2 & \beta_2 & \gamma_2 \\ \Delta x_3 & v_3 & \alpha_3 & \beta_3 & \gamma_3 \end{pmatrix} \begin{pmatrix} 1 \\ t \\ x_1 \\ x_2 \\ x_3 \end{pmatrix}$$

It follows that Galilean space–time has signature $\{0, 0, 1, 1\}$. Moreover, its time axis is non–interchangeable with any space axes. Hence, it is possible the only axes order: time axis, then space axes.

When introduce the curvature in Galilean space, its signature becomes $\{k \neq 0, 0, 1, 1\}$. The curvature affects only time type: $K_1 = k \neq 0$. Space types equal: $K_2 = k \cdot 0 = 0, K_3 = k \cdot 0 \cdot 1 = 0, K_4 = k \cdot 0 \cdot 1 \cdot 1 = 0$.

**Minkowski space with two–dimensional time.** Consider 4-dimensional space with distance bilinear form having 2 positive and 2 negative signs. This space sometimes is named Minkowski space with two–dimensional time:

$$\begin{aligned} d^2(x, y) &= (y_1 - x_1)^2 + (y_2 - x_2)^2 - (y_3 - x_3)^2 - (y_4 - x_4)^2 \\ &= (y_1 - x_1)^2 + k_2(y_2 - x_2)^2 + k_2 k_3 (y_3 - x_3)^2 + k_2 k_3 k_4 (y_4 - x_4)^2 \end{aligned}$$

It means that $k_2 = 1$, $k_2 k_3 = -1$, $k_3 = -1$ and $k_2 k_3 k_4 = -1$, $k_4 = 1$. As in any linear space, $k_1 = 0$.

**Spaces with degenerate distance bilinear form.** Consider linear 4-dimensional space ($k_1 = 0$) with degenerate distance bilinear form:

$$\begin{aligned} d^2(x, y) &= (y_1 - x_1)^2 + (y_2 - x_2)^2 + (y_3 - x_3)^2 \\ &= (y_1 - x_1)^2 + k_2(y_2 - x_2)^2 + k_2k_3(y_3 - x_3)^2 + k_2k_3k_4(y_4 - x_4)^2 \end{aligned}$$

It means that $k_2 = k_3 = 1$ and $k_4 = 0$.

From here results, that transformation:

$$\begin{cases} x'_1 = x_1 \\ x'_2 = x_2 \\ x'_3 = x_3 \\ x'_4 = \alpha x_1 + \beta x_2 + \gamma x_3 + x_4 \end{cases}$$

is the motion.

But the transformation:

$$\begin{cases} x'_1 = x_1 \\ x'_2 = x_2 \\ x'_3 = x_3 \\ x'_4 = \delta x_4 \end{cases}$$

is not the motion. Although it preserves the distances, it doesn't preserve the volume, except case $\delta = 1$ or $-1$. This is example of angle dilation transformation.

### 3.1.3. Signature of the subspaces product.

One more method of space definition is as subspaces product. One should be careful here. Geometric space is not only a structure of points. It is also structure of all subspaces. It is wrong to assume that since $\mathbb{R}^1$ is isometric to one–dimensional Euclidean space $\mathbb{E}^1$, then from $\mathbb{R}^1 \times \mathbb{R}^1 = \mathbb{R}^2$ follows that $\mathbb{E}^1 \times \mathbb{E}^1 = \mathbb{E}^2$ (in signature language, $\{0\} \times \{0\} = \{0, 1\}$). The problem is that algebraic product doesn't define the way to measure the angle between multiplied subspaces, which can be defined in different way, for instance, $\{0, 0\}$ or $\{0, -1\}$. In other words, from $K_1 = 0$, $K_2 = 0$ doesn't result that $k_2 = \frac{K_2}{K_1} = 1$ (it may happen that $k_2 = 0$ or $k_2 = -1$).

Even worse is the situation when multiplied subspaces $X^m$ and $Y^n$ have different signatures. There are two types of one-dimensional objects: $X^1 \times Y^0$ (isometric to $X^1$) and $X^0 \times Y^1$ (isometric to $Y^1$). But if $X^1$ and $Y^1$ have different signatures, these one-dimensional lines are not equivalent, the fact that leads to severe restrictions of possible angle types in lineal $X^1Y^1$, often not allowing Euclidean angle type.

**Example.** Let us construct homogeneous geometry on cylinder. The first idea in one's mind is $\mathbb{S}^1 \times \mathbb{E}^1$ $(\{1\} \times \{0\})$, where $\mathbb{S}^1$ is one–dimensional elliptic space. However, in this case some lines are circles, other lines are usual lines and yet other lines are right or left helices, which may not intersect, or may intersect in one point, or may intersect in infinity number of points. One example of homogeneous geometry on cylinder is positive curved Galilean space with signature $\{1, 0\}$.

Additionally, it is wrong to think that if from geometric point of view $\mathbb{E}^1$ is isometric to $\mathbb{H}^1$ (one–dimensional hyperbolic space), then structures like $\mathbb{H}^2 \times \mathbb{E}^1$ $(\{-1, 1\} \times \{0\})$ and $\mathbb{H}^2 \times \mathbb{H}^1$ $(\{-1, 1\} \times \{-1\})$ are also isometric. From topology point of view (see section 3.2.1), points on $\mathbb{E}^1$ are weak separable, while points on $\mathbb{H}^1$ are strong separable.

**Example.** We can construct spaces with signatures $\{-1, 1, 0\}$ and $\{-1, 1, -1\}$. The difference between them is the following. In the first space, on two–dimensional plane (with signature $\{-1, 1\}$), that doesn't contain some point, there is the only point unconnectable with it. In the second space, on two–dimensional plane (with signature also $\{-1, 1\}$), that doesn't contain some point, there are infinity number of points unconnectable with it.

**Example.** We need to find what 4-dimensional homogeneous spaces $\mathbb{B}^4$ allow two 2-dimensional absolutely orthogonal lineals (subspaces), so that each of them is isomorphic to hyperbolic plane $\mathbb{H}^2$. Choose the basis such that coordinate vectors $\{e^0, e^1, e^2\}$ to define the first lineal and $\{e^0, e^3, e^4\}$ to define the second lineal. Because the first lineal is isometric to $\mathbb{H}^2$, its signature is $\{-1, 1\}$, and main rotations are:

$$\mathfrak{R}_1(\varphi) = \begin{pmatrix} \cosh\varphi & \sinh\varphi & 0 & 0 & 0 \\ \sinh\varphi & \cosh\varphi & 0 & 0 & 0 \\ 0 & 0 & 1 & 0 & 0 \\ 0 & 0 & 0 & 1 & 0 \\ 0 & 0 & 0 & 0 & 1 \end{pmatrix}$$

$$\mathfrak{R}_2(\varphi) = \begin{pmatrix} 1 & 0 & 0 & 0 & 0 \\ 0 & \cos\varphi & -\sin\varphi & 0 & 0 \\ 0 & \sin\varphi & \cos\varphi & 0 & 0 \\ 0 & 0 & 0 & 1 & 0 \\ 0 & 0 & 0 & 0 & 1 \end{pmatrix}$$

The second lineal is also isometric to $\mathbb{H}^2$ with signature $\{-1, 1\}$. Its main rotations are:

$$\mathfrak{R}_{03}(\varphi) = \begin{pmatrix} \cosh\varphi & 0 & 0 & \sinh\varphi & 0 \\ 0 & 1 & 0 & 0 & 0 \\ 0 & 0 & 1 & 0 & 0 \\ \sinh\varphi & 0 & 0 & \cosh\varphi & 0 \\ 0 & 0 & 0 & 0 & 1 \end{pmatrix}$$

$$\mathfrak{R}_4(\varphi) = \begin{pmatrix} 1 & 0 & 0 & 0 & 0 \\ 0 & 1 & 0 & 0 & 0 \\ 0 & 0 & 1 & 0 & 0 \\ 0 & 0 & 0 & \cos\varphi & -\sin\varphi \\ 0 & 0 & 0 & \sin\varphi & \cos\varphi \end{pmatrix}$$

Motions $\mathfrak{R}_1(\varphi), \mathfrak{R}_2(\varphi), \mathfrak{R}_4(\varphi)$ are also main rotations of $\mathbb{B}^4$, but $\mathfrak{R}_{03}(\varphi)$, although is rotation, is not the main one. Based on known main rotation types, the possible signature of $\mathbb{B}^4$ is $\{-1, 1, k_3, 1\}$. Then, from the type of $\mathfrak{R}_{03}(\varphi)$ follows that $K_3 = -1 = k_1 k_2 k_3 = -1 \cdot 1 \cdot k_3 = -k_3$, then $k_3 = 1$ and whole signature of $\mathbb{B}^4$ is $\{-1, 1, 1, 1\}$. It is 4-dimensional hyperbolic space $\mathbb{B}^4 = \mathbb{H}^4$.

Any other combination of hyperbolic planes $\mathbb{H}^2$ doesn't give a homogeneous space. For example, algebraic product $\mathbb{H}^2 \times \mathbb{H}^2$ is not homogeneous space. By theorem 2.4.16, the number of freedom degrees in $\mathbb{H}^2$ equals to $\frac{2 \cdot 3}{2} = 3$. The number of freedom degrees in $\mathbb{H}^2 \times \mathbb{H}^2$ equals to sum of freedom degrees in subspaces, $3 + 3 = 6$. However, the 4-dimensional homogeneous space $\mathbb{H}^4$ have the number of freedom degrees equal to $\frac{4 \cdot 5}{2} = 10$. In this way, isometry group of $\mathbb{H}^2 \times \mathbb{H}^2$ is much smaller then isometry group of $\mathbb{H}^4$.

One more product, that often occures in geometry is the representation of group $Iso(\mathbb{B}^n)$ of space motions (discrete or continuous) as semidirect product of normal subgroup $\mathfrak{T}^n$ of translations and subgroup $\mathfrak{R}_E^{n-1}$ of rotations around the origin — its stabilizer:

$$Iso(\mathbb{B}^n) = \mathfrak{T}^n \rtimes \mathfrak{R}_E^{n-1}$$

Even if this representation may seem very different from above analyzed, in fact it is naturally fits into the theory using the language of GM-orthogonal matrices. In this case matrices represent not orthonormal basis, but space motions. Recall that translations $\mathfrak{T}_i(\varphi)$ are rotations $\mathfrak{R}_{0i}(\varphi)$, then all other rotations $\mathfrak{R}_{ij}(\varphi)$ don't contain 0 index among $i$ and $j$. In this case, the subgroups of matrices are "glued" not by 0-th column, but by all other columns.

### 3.1.4. Examples of crystallographic groups on homogeneous planes.

Consider a motion group $G$ in homogeneous space $\mathbb{B}^n$.

**Definition 3.1.3** (Point orbit and stabilizer, lattice and fundamental domain of a group)**.** The *orbit* of some point $P \in \mathbb{B}^n$ is the set of points $\{T_i = g_i P\}$, where $g_i \in G$ are all possible elements of motion group. The maximal subgroup $H < G$ that preserves some point $P \in \mathbb{B}^n$ is named *stabilizer* of the point $P$. The orbit of such invariant point is named the *lattice* of group $G$. Elements of a lattice are called *nodes.* The set of points $\{F_g\}$ containing one point from all point orbits is called *fundamental domain* of group $G$.

**Definition 3.1.4** (Discrete, uniform, crystallographic group)**.** If for some group $G$ and its lattice $\{T_p\}$ there exists $r > 0 \in \mathbb{R}$, such that minimal distance between two lattice nodes $\min_{T_i, T_j \in \{T_p\}} d(T_i, T_j) \geq r$, then the group is called *discrete.* If for $G$ and $\{T_p\}$ there exists $R > 0 \in \mathbb{R}$, such that maximal distance between any space point and some lattice node $\max_{P \in \mathbb{B}^n, T_i \in \{T_p\}} d(P, T_i) \leq R$, then the group is called *uniform.* If a motion group $G$ is discrete and uniform, it is called *crystallographic.*

Crystallographic groups [30, 118, 183] are usually studied in Euclidean space. For any dimension the number of groups is finite. In hyperbolic space there exist at least countable set of groups. In other spaces the interest for them is insufficiently large to be closely examined. In any case, their number is very large and quickly grows with dimension, while their construction and classification is not trivial.

Without elaborating the full theory of crystallographic groups, construct their examples on 9 homogeneous planes.

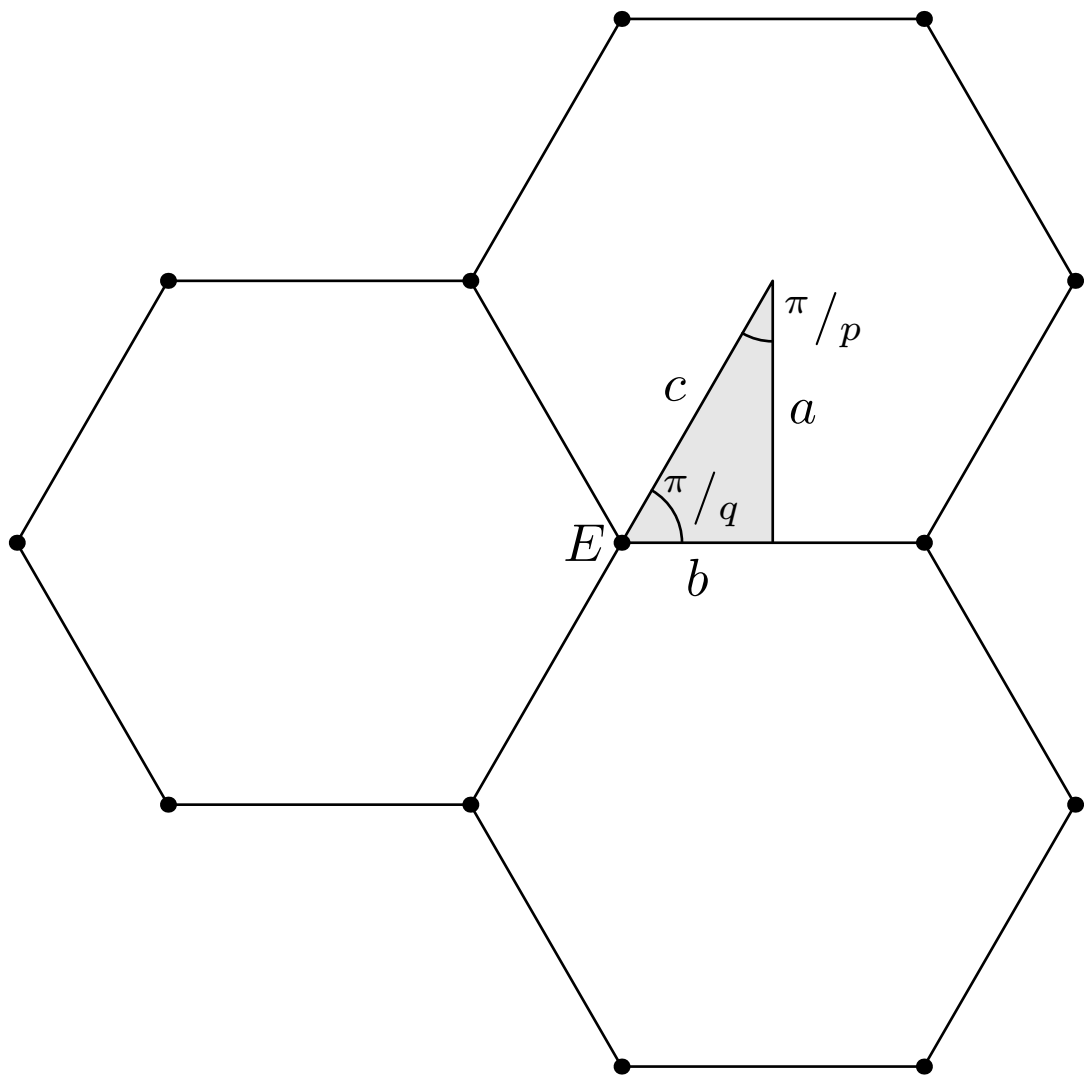

Figure 3.1: Fundamental domain of regular hexagon tiling on Euclidean plane.

**Elliptic $\{1, 1\}$, Euclidean $\{0, 1\}$ and hyperbolic $\{-1, 1\}$ planes.** As classical example, take symmetry group of regular $p$-gon tiling on plane, having $q$ edges in each vertex (Figure 3.1). The fundamental domain of this group is right triangle with vertices in tiling vertex, edge center and $p$-gon center. Because on these planes angular type $k_2 = 1$, the value of angle is bounded. The acute angles of the triangle equal to $\frac{\pi}{p}, \frac{\pi}{q}$ respectively. By (2.52, 2.53, 2.54), keeping in mind that $\beta = \frac{\pi}{2} - \beta'$, obtain catheti and hypotenuse lengths $a, b, c$:

$$C_1(a) = \frac{\cos\frac{\pi}{q}}{\sin\frac{\pi}{p}},\ C_1(b) = \frac{\cos\frac{\pi}{p}}{\sin\frac{\pi}{q}},\ C_1(c) = \cot\frac{\pi}{p}\cot\frac{\pi}{q}.$$

From here we can find the length type $k_1$. In case $k_1 = 1$:

$$\begin{aligned} C_1(c) = \cot\frac{\pi}{p}\cot\frac{\pi}{q} &< 1, \\ \tan\left(\frac{\pi}{2} - \frac{\pi}{p}\right) &< \tan\frac{\pi}{q}, \\ \frac{\pi}{2} - \frac{\pi}{q} &< \frac{\pi}{p}, \\ pq - 2p - 2q &< 0, \\ (p-2)(q-2) &< 4. \end{aligned}$$

Knowing that minimal value for $p, q$ equals to 3, obtain possible values for elliptic plane $(p, q) = \{(3,3), (4,3), (3,4), (5,3), (3,5)\}$, corresponding to symmetry groups of tetrahedron, cube, octahedron, dodecahedraon and icosahedron respectively. The groups of cube and octahedron are isomorphic, as well as of dodecahedron and icosahedron.

Similarly, for Euclidean plane:

$$(p-2)(q-2) = 4$$

It gives square, hexagon and triangle tilings $(p, q) = \{(4,4), (6,3), (3,6)\}$.

Finally, for hyperbolic plane:

$$(p-2)(q-2) > 4$$

It gives infinitely many tilings.

For all these cases, the generators of proper motion groups are:

$$\mathfrak{T}=\begin{pmatrix}C_1(2d) & -k_1S_1(2d) & 0\\ S_1(2d) & C_1(2d) & 0\\ 0 & 0 & 1\end{pmatrix},\ d=\begin{cases}b, & q\in 2\mathbb{Z},\\ b+c, & q\in 2\mathbb{Z}+1, p\in 2\mathbb{Z},\\ a+b+c, & q\in 2\mathbb{Z}+1, p\in 2\mathbb{Z}+1,\end{cases}$$

$$\mathfrak{R}=\begin{pmatrix}1 & 0 & 0\\ 0 & \cos 2\varphi & -\sin 2\varphi\\ 0 & \sin 2\varphi & \cos 2\varphi\end{pmatrix},\ \varphi=\frac{\pi}{q},$$

$$a=C_1^{-1}\left(\frac{\cos\frac{\pi}{q}}{\sin\frac{\pi}{p}}\right),\ b=C_1^{-1}\left(\frac{\cos\frac{\pi}{p}}{\sin\frac{\pi}{q}}\right),\ c=C_1^{-1}\left(\cot\frac{\pi}{p}\cot\frac{\pi}{q}\right),$$

$$k_1=\begin{cases}1, & (p-2)(q-2)<4,\\ 0, & (p-2)(q-2)=4,\\ -1, & (p-2)(q-2)>4.\end{cases}$$

The lattice node can be chosen as $E=(1:0:0)$.

*Remark.* Much more elegant is to define these groups solely by rotations on angles $\frac{2\pi}{p}$ and $\frac{2\pi}{q}$, with distance between centres equal to $c$. However, for futher analysis we need to take main rotations as group generators.

*Remark.* For Euclidean plane it is impossible to calculate $a, b, c$ from equations above. It is possible to take some positive value $d$ instead. As Euclidean plane allows lenght dilation transformation, all groups corresponding to the same value of $\varphi$ and all possible values of $d$, are isomorphic.

**Positively curved Galilean $\{1, 0\}$ and Minkowski $\{1, -1\}$ planes.** In order to describe the groups on these planes, use the duality principle.

**Lemma 3.1.3.** *Crystallographic group of space with signature $\{k_1, ..., k_n\}$, generated by main rotations $\mathfrak{R}_1(\varphi_1), ..., \mathfrak{R}_n(\varphi_n)$, is isomorphic to crystallographic group of dual space with signature $\{k_n, ..., k_1\}$, generated by main rotations $\mathfrak{R}'_1(\varphi_n), ..., \mathfrak{R}'_n(\varphi_1)$.*

*Proof.* Construct matrix transformation $\eta$ as "transposition" around secondary diagonal:

$$\eta\left(\begin{pmatrix}m_{00} & \cdots & m_{0n}\\ \vdots & \ddots & \vdots\\ m_{n0} & \cdots & m_{nn}\end{pmatrix}\right)=\begin{pmatrix}m_{nn} & \cdots & m_{0n}\\ \vdots & \ddots & \vdots\\ m_{n0} & \cdots & m_{00}\end{pmatrix}$$

It is easy to see, that:

$$\eta(\mathfrak{R}_i(\varphi))=\mathfrak{R}'_{n-i}(\varphi)$$

Additionally:

$$\eta(\mathfrak{R}_i(\varphi)\mathfrak{R}_j(\psi)) = \eta(\mathfrak{R}_j(\psi))\eta(\mathfrak{R}_i(\varphi))$$

Check the last equality for $n = 2, \quad \mathfrak{R}_1(\varphi), \mathfrak{R}_2(\psi)$. This case is trivially generalized to any dimension. The case when matrices have no common elements different from $\pm 1, 0$ also trivially verifies.

$$\eta\left(\begin{pmatrix} C_1(\varphi) & -k_1S_1(\varphi) & 0 \\ S_1(\varphi) & C_1(\varphi) & 0 \\ 0 & 0 & 1 \end{pmatrix}\begin{pmatrix} 1 & 0 & 0 \\ 0 & C_2(\psi) & -k_2S_2(\psi) \\ 0 & S_2(\psi) & C_2(\psi) \end{pmatrix}\right)$$
$$=\eta\left(\begin{pmatrix} C_1(\varphi) & -k_1S_1(\varphi)C_2(\psi) & k_1k_2S_1(\varphi)S_2(\psi) \\ S_1(\varphi) & C_1(\varphi)C_2(\psi) & -k_2C_1(\varphi)S_2(\psi) \\ 0 & S_2(\psi) & C_2(\psi) \end{pmatrix}\right)$$
$$=\begin{pmatrix} C_2(\psi) & -k_2S_2(\psi)C_1(\varphi) & k_2k_1S_2(\psi)S_1(\varphi) \\ S_2(\psi) & C_2(\psi)C_1(\varphi) & -k_1C_2(\psi)S_1(\varphi) \\ 0 & S_1(\varphi) & C_1(\varphi) \end{pmatrix}$$

On the other hand:

$$\eta\left(\begin{pmatrix} 1 & 0 & 0 \\ 0 & C_2(\psi) & -k_2S_2(\psi) \\ 0 & S_2(\psi) & C_2(\psi) \end{pmatrix}\right)\eta\left(\begin{pmatrix} C_1(\varphi) & -k_1S_1(\varphi) & 0 \\ S_1(\varphi) & C_1(\varphi) & 0 \\ 0 & 0 & 1 \end{pmatrix}\right)$$
$$=\begin{pmatrix} C_2(\psi) & -k_2S_2(\psi) & 0 \\ S_2(\psi) & C_2(\psi) & 0 \\ 0 & 0 & 1 \end{pmatrix}\begin{pmatrix} 1 & 0 & 0 \\ 0 & C_1(\varphi) & -k_1S_1(\varphi) \\ 0 & S_1(\varphi) & C_1(\varphi) \end{pmatrix}$$
$$=\begin{pmatrix} C_2(\psi) & -k_2S_2(\psi)C_1(\varphi) & k_2k_1S_2(\psi)S_1(\varphi) \\ S_2(\psi) & C_2(\psi)C_1(\varphi) & -k_1C_2(\psi)S_1(\varphi) \\ 0 & S_1(\varphi) & C_1(\varphi) \end{pmatrix}$$

It means that $\eta(\mathfrak{R}_1(\varphi)\mathfrak{R}_2(\psi)) = \eta(\mathfrak{R}_2(\psi))\eta(\mathfrak{R}_1(\varphi))$.

Next:

$$\eta\left(\begin{pmatrix}1 & 0 & 0\\ 0 & C_2(\psi) & -k_2S_2(\psi)\\ 0 & S_2(\psi) & C_2(\psi)\end{pmatrix}\begin{pmatrix}C_1(\varphi) & -k_1S_1(\varphi) & 0\\ S_1(\varphi) & C_1(\varphi) & 0\\ 0 & 0 & 1\end{pmatrix}\right)$$
$$=\eta\left(\begin{pmatrix}C_1(\varphi) & -k_1S_1(\varphi) & 0\\ C_2(\psi)S_1(\varphi) & C_2(\psi)C_1(\varphi) & -k_2S_2(\psi)\\ S_2(\psi)S_1(\varphi) & S_2(\psi)C_1(\varphi) & C_2(\psi)\end{pmatrix}\right)$$
$$=\begin{pmatrix}C_2(\psi) & -k_2S_2(\psi) & 0\\ S_2(\psi)C_1(\varphi) & C_2(\psi)C_1(\varphi) & -k_1S_1(\varphi)\\ S_2(\psi)S_1(\varphi) & C_2(\psi)S_1(\varphi) & C_1(\varphi)\end{pmatrix}$$

On the other hand:

$$\eta\left(\begin{pmatrix}C_1(\varphi) & -k_1S_1(\varphi) & 0\\ S_1(\varphi) & C_1(\varphi) & 0\\ 0 & 0 & 1\end{pmatrix}\right)\eta\left(\begin{pmatrix}1 & 0 & 0\\ 0 & C_2(\psi) & -k_2S_2(\psi)\\ 0 & S_2(\psi) & C_2(\psi)\end{pmatrix}\right)$$
$$=\begin{pmatrix}1 & 0 & 0\\ 0 & C_1(\varphi) & -k_1S_1(\varphi)\\ 0 & S_1(\varphi) & C_1(\varphi)\end{pmatrix}\begin{pmatrix}C_2(\psi) & -k_2S_2(\psi) & 0\\ S_2(\psi) & C_2(\psi) & 0\\ 0 & 0 & 1\end{pmatrix}$$
$$=\begin{pmatrix}C_2(\psi) & -k_2S_2(\psi) & 0\\ S_2(\psi)C_1(\varphi) & C_2(\psi)C_1(\varphi) & -k_1S_1(\varphi)\\ S_2(\psi)S_1(\varphi) & C_2(\psi)S_1(\varphi) & C_1(\varphi)\end{pmatrix}$$

In this way, $\eta(\mathfrak{R}_2(\psi)\mathfrak{R}_1(\varphi)) = \eta(\mathfrak{R}_1(\varphi))\eta(\mathfrak{R}_2(\psi))$.

Because any GM-orthogonal matrix can be obtained as product of diagonal matrix $Diag(\pm 1)$ and main rotations, we can construct the isomorphism $\xi$ as:

$$\xi(\mathfrak{M}) = \eta(\mathfrak{M}^{-1}),$$
$$\xi(\mathfrak{M}\mathfrak{W}) = \eta((\mathfrak{M}\mathfrak{W})^{-1}) = \eta(\mathfrak{W}^{-1}\mathfrak{M}^{-1}) = \eta(\mathfrak{M}^{-1})\eta(\mathfrak{W}^{-1}) = \xi(\mathfrak{M})\xi(\mathfrak{W}).$$

On this mapping, the matrices signature of the second space, and with them also space signature, equals to $\{k_n, ..., k_1\}$. It means that the second space is really dual to the first one with signature $\{k_1, ..., k_n\}$.

If complete isometry group of one space is isomorphic to complete isometry group of the second space, then its crystallographic subgroup is isomorphic to crystallographic subgroup of the second space. Also the first crystallographic group generators are isomorphic to the second group generators and the first group relations are equivalent to the second group relations. □

*Remark.* Although crystallographic groups of dual spaces are isomorphic, the dual spaces are, generally speaking, not isometric, as well as group lattices.

Applying duality principle to elliptic groups obtain elliptic groups. In this procedure symmetry groups of regular polyhedrons corresspond to symmetry groups of dual polyhedrons ($(3,3) \leftrightarrow (3,3)$; $(4,3) \leftrightarrow (3,4)$; $(5,3) \leftrightarrow (3,5)$). Next, applying duality principle to Euclidean ($\{0,1\}$) groups obtain groups on positively curved Galilean plane ($\{1,0\}$), applying duality to hyperbolic ($\{-1,1\}$) groups obtain groups on positively curved Minkowski plane ($\{1,-1\}$). In all cases generators are:

$$\mathfrak{T} = \begin{pmatrix} \cos 2d & -\sin 2d & 0 \\ \sin 2d & \cos 2d & 0 \\ 0 & 0 & 1 \end{pmatrix}, \; d = \frac{\pi}{q},$$

$$\mathfrak{R} = \begin{pmatrix} 1 & 0 & 0 \\ 0 & C_2(2\varphi) & -k_2 S_2(2\varphi) \\ 0 & S_2(2\varphi) & C_2(2\varphi) \end{pmatrix}, \; \varphi = \begin{cases} b, & q \in 2\mathbb{Z}, \\ b + c, & q \in 2\mathbb{Z}+1, p \in 2\mathbb{Z}, \\ a + b + c, & q \in 2\mathbb{Z}+1, p \in 2\mathbb{Z}+1, \end{cases}$$

$$a = C_2^{-1}\left(\frac{\cos\frac{\pi}{q}}{\sin\frac{\pi}{p}}\right), \; b = C_2^{-1}\left(\frac{\cos\frac{\pi}{p}}{\sin\frac{\pi}{q}}\right), \; c = C_2^{-1}\left(\cot\frac{\pi}{p}\cot\frac{\pi}{q}\right),$$

$$k_2 = \begin{cases} 1, & (p-2)(q-2) < 4, \\ 0, & (p-2)(q-2) = 4, \\ -1, & (p-2)(q-2) > 4. \end{cases}$$

As lattice node we can choose, as earlier, point $E = (1:0:0)$.

*Remark.* As earlier, in case $k_2 = 0$, instead of $\varphi$ we can take any positive value.

**Galilean plane** $\{0,0\}$**.** Search the crystallographic group on linear Galilean plane as symmetry group of the lattice $T_{ij} = E + ix + jy$; $E = (1:0:0), x = \{0:a:0\}, y = \{0:0:b\}$; $i, j \in \mathbb{Z}$. Evidently, the following motions are its generators:

$$\mathfrak{T} = \begin{pmatrix} 1 & 0 & 0 \\ d & 1 & 0 \\ 0 & 0 & 1 \end{pmatrix}, \; d = a \in \mathbb{R},$$

$$\mathfrak{R} = \begin{pmatrix} 1 & 0 & 0 \\ 0 & 1 & 0 \\ 0 & \varphi & 0 \end{pmatrix}, \; \varphi = \frac{b}{a} \in \mathbb{R}.$$

As usually, take the lattice node $E = (1:0:0)$.

**Minkowski plane** $\{0, -1\}$. For linear Minkowski plane also search crystallographic group as symmetry group of lattice $T_{ij} = E + ix + jy$; $E = (1 : 0 : 0), x = \{0 : a : 0\}, y = \{0 : 0 : b\}$. Evidently, translation along the vector $x$ by $a$ maps the lattice to itself. Consider rotation $\mathfrak{R}$ by angle $\varphi$. Let two lattice points, $P = (1 : a : 0), Q = (1 : 0 : b)$ be mapped to some lattice points $\mathfrak{R}P = (1 : ua : vb), \mathfrak{R}Q = (1 : ra : tb)$; $u, v, r, t \in \mathbb{Z}$:

$$\mathfrak{R}P = \begin{pmatrix} 1 & 0 & 0 \\ 0 & \cosh\varphi & \sinh\varphi \\ 0 & \sinh\varphi & \cosh\varphi \end{pmatrix} \begin{pmatrix} 1 \\ a \\ 0 \end{pmatrix} = \begin{pmatrix} 1 \\ a\cosh\varphi \\ a\sinh\varphi \end{pmatrix} = \begin{pmatrix} 1 \\ ua \\ vb \end{pmatrix}$$

From here results:

$$\cosh\varphi = u,$$
$$\sinh\varphi = v\frac{b}{a} = \sqrt{u^2 - 1}.$$

Next:

$$\mathfrak{R}Q = \begin{pmatrix} 1 & 0 & 0 \\ 0 & \cosh\varphi & \sinh\varphi \\ 0 & \sinh\varphi & \cosh\varphi \end{pmatrix} \begin{pmatrix} 1 \\ 0 \\ b \end{pmatrix} = \begin{pmatrix} 1 \\ b\sinh\varphi \\ b\cosh\varphi \end{pmatrix} = \begin{pmatrix} 1 \\ ra \\ tb \end{pmatrix}$$

From here results:

$$\cosh\varphi = t,$$
$$\sinh\varphi = r\frac{a}{b} = \sqrt{t^2 - 1}.$$

Summing up all the above, obtain:

$$u = \cosh\varphi,$$
$$v = \frac{a}{b}\sqrt{u^2 - 1},$$
$$r = \frac{b}{a}\sqrt{u^2 - 1} = \frac{a^2}{b^2}v,$$
$$t = u.$$

This way, in order $u, v, r, t$ to be all integers, we need to choose right values $\varphi, \frac{b}{a}$:

$$\varphi = \cosh^{-1} u,$$
$$\frac{b}{a} = \frac{\sqrt{u^2 - 1}}{v}.$$

For integer value of $r$ require that:

$$r = \frac{b}{a}\sqrt{u^2 - 1} = \frac{u^2 - 1}{v} \in \mathbb{Z}$$

It is possible for many different values of $u, v \in \mathbb{Z}$, for instance $v = u \pm 1$. But the group is crystallographic even when $r \in \mathbb{Q}$ (that is, when $u, v \in \mathbb{Z}$). In this case the angle $\varphi$ will be multiplied by $v$. When $v = u \pm 1$:

$$\frac{b}{a} = \frac{\sqrt{u^2 - 1}}{u \pm 1} = \sqrt{\frac{u \pm 1}{u \mp 1}}.$$

At the end, the group generators are:

$$\mathfrak{T} = \begin{pmatrix} 1 & 0 & 0 \\ d & 1 & 0 \\ 0 & 0 & 1 \end{pmatrix}, \; d = a \in \mathbb{R}, \; \frac{b}{a} = \sqrt{\frac{u \pm 1}{u \mp 1}}, u \in \mathbb{Z},$$
$$\mathfrak{R} = \begin{pmatrix} 1 & 0 & 0 \\ 0 & \cosh\varphi & \sinh\varphi \\ 0 & \sinh\varphi & \cosh\varphi \end{pmatrix}, \; \varphi = \cosh^{-1} u, u \in \mathbb{Z}.$$

And usual lattice node is $E = (1 : 0 : 0)$.

**Negatively curved Galilean plane $\{-1, 0\}$.** Because negatively curved Galilean plane $\{-1, 0\}$ is dual to linear Minkowski plane $\{0, -1\}$, by lemma 3.1.3 on group duality obtain crystallographic group generators:

$$\mathfrak{T} = \begin{pmatrix} \cosh d & \sinh d & 0 \\ \sinh d & \cosh d & 0 \\ 0 & 0 & 1 \end{pmatrix}, \; d = \cosh^{-1} u, u \in \mathbb{Z},$$
$$\mathfrak{R} = \begin{pmatrix} 1 & 0 & 0 \\ 0 & 1 & 0 \\ 0 & \varphi & 1 \end{pmatrix}, \; \varphi = \sqrt{\frac{u \pm 1}{u \mp 1}}, u \in \mathbb{Z}.$$

Again, the lattice node is $E = (1 : 0 : 0)$.

**Negatively curved Minkowski plane** $\{-1,-1\}$. There is isomorphism between positive ($\{1,-1\}$) and negative ($\{-1,-1\}$) curved Minkowski planes: $\eta((x_0 : x_1 : x_2)) = (x_0 : x_2 : x_1)$. This isomorphism preserves crystallographic groups, their generators *and lattices.* The only difficulty is the fact that one generator is not more a main rotation.

In order to construct generators of crystallographic group on plane $\{-1,-1\}$ as main rotations, follow all concatenation of isomorphisms up to described groups. Let $\mathfrak{T} = \mathfrak{R}_1(d), \mathfrak{R} = \mathfrak{R}_2(\varphi)$. Firstly, corresponding motions on the plane $\{1,-1\}$ are $\mathfrak{R}'_{02}(d), \mathfrak{R}'_2(\varphi)$. Secondly, corresponding motions on the plane $\{-1,1\}$ are $\mathfrak{R}''_{02}(d), \mathfrak{R}''_1(\varphi)$. The second one is exactly generator $\mathfrak{T}''$ — translation along vector $x^1$. The first is translation along vector $x^2$. This translation is easily expressed in terms of groups on hyperbolic plane if the number of edges in a vertex is even. Thus, for simplicity consider that $q \in 2\mathbb{Z}$. Then the generators are:

$$\mathfrak{T} = \begin{pmatrix} \cosh 2d & \sinh 2d & 0 \\ \sinh 2d & \cosh 2d & 0 \\ 0 & 0 & 1 \end{pmatrix}, \; d = b,$$

$$\mathfrak{R} = \begin{pmatrix} 1 & 0 & 0 \\ 0 & \cosh 2\varphi & \sinh 2\varphi \\ 0 & \sinh 2\varphi & \cosh 2\varphi \end{pmatrix}, \; \varphi = \begin{cases} b, & q \in 4\mathbb{Z}, \\ a + c, & q \in 4\mathbb{Z} + 2, \end{cases}$$

$$a = \cosh^{-1}\left(\frac{\cos\frac{\pi}{q}}{\sin\frac{\pi}{p}}\right), \; b = \cosh^{-1}\left(\frac{\cos\frac{\pi}{p}}{\sin\frac{\pi}{q}}\right), \; c = \cosh^{-1}\left(\cot\frac{\pi}{p}\cot\frac{\pi}{q}\right),$$

$$(p-2)(q-2) > 4.$$

As always, the lattice node is the point $E = (1 : 0 : 0)$.

## 3.2. Topology

### 3.2.1. Separability of the points on a line.

As it was shown, elliptic, parabolic and hperbolic lines are all different. In what exactly this difference consists from topology point of view? In points separability in these lines.

Traditionally, points separability is defined as follows. We say, points in some line are *separable*, if among any three different points $A, B, C$ one (let it be $B$) divides the line into two half–lines, and remaining two points $A, C$ lie in different half–lines. In this case we call $B$ the *middle* point. Otherwise, we call points *non–separable.*

*Remark.* Obviously, points in elliptic line are non-separable. Points in parabolic and hyperbolic line are separable.

In order to draw the difference between the last two cases of separable points, give more precise definition. This definition gives finer distinction and is based only on points connectability notion. However, it is given in terms of line metaplane.

**Definition 3.2.1** (Separable and non–separable points)**.** We call points on a line *non–separable*, if all points on this line are connectable with any point on metaplane. We call points on a line *separable*, if for any three points $A, B, C$ in this line and some point $D$ on metaplane, that is connectable with $A, C$ and unconnectable with $B$, the angle $\angle ADC$ is unmeasurable (Figure 3.2).

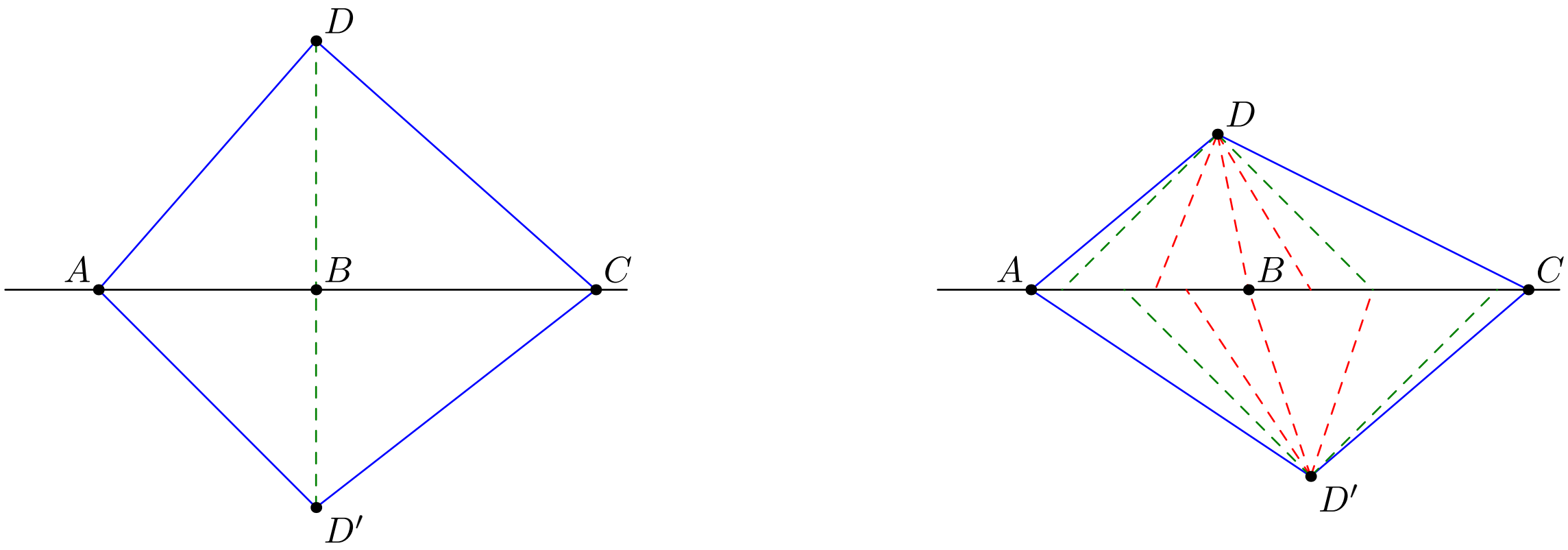

Figure 3.2: Points separability on a plane: weak (left) and strong (right).

*Remark.* For separable points $A, B, C$ a single point ($B$) has described property. For other points ($A, C$) and some unconnectable with them points $D_A$, $D_C$, the angles $\angle BD_AC$ and $\angle AD_CB$ are measurable.

**Definition 3.2.2** (Middle point)**.** In case of separable points we say, the point $B$ *lies between* points $A$ and $C$. The point $B$ we call *middle* one.

*Remark.* In case of non–separable points among any three points there is no middle, and it is impossible to talk about position of some point between others two.

**Definition 3.2.3** (Weak and strong separable points)**.** We call the points in some line *weak separable* (Figure 3.2, left), if on this line metaplane any point $D$, being unconnectable with middle point $B$ and connectable with other two points $A, C$, is connectable with all points from neighborhood of $B$. We call the points in some line *strong separable* (Figure 3.2, right), if in the same conditions any point $D$ is unconnectable not only with $B$, but also with all points from some its neighborhood.

In these definitions, the points in elliptic line are non-separable, the points in parabolic line are weak separable, and the points in hyperbolic line are strong separable.

### 3.2.2. Neighborhood notion generalization.

In previous section we mentioned the notion of neighborhood on a line. This notion in general case, however, needs to be refined. And the cause is, that the classical definition of neighborhood does not satisfy the duality principle.

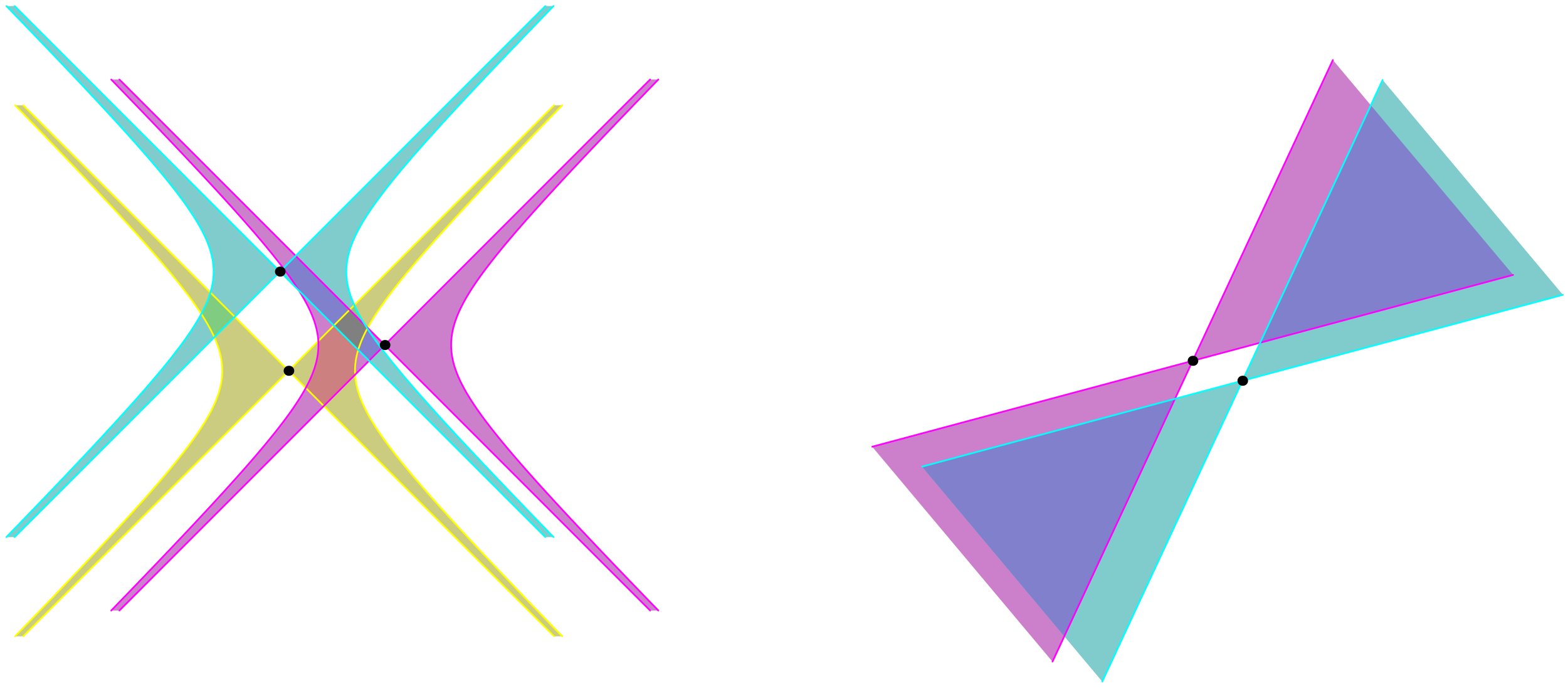

Figure 3.3: On Minkowski plane all points neighborhoods are intersected (left). On Euclidean plane all points anti–neighborhoods are intersected (right).

Classical neighborhood of some point is defined as a set of all points situated at the distance less than some fixed value $r$ away from this point. Let construct the definition of "anti–neighborhood", dual to neighborhood definition. *Anti–neighborhood* of some line $a$ (Figure 3.3, right) consists of all lines $b$ on a plane, so that the angle between lines $a$ and $b$ is less than some fixed angle $\rho$. Reasoning further in set operation language, consider a line as set of points. In this way, the anti–neighborhood of point $A$ (having some fixed direction $AO$ in this point) consists of all points $B$, so that the angle $\angle BAO < \rho$. Unlikely this definition have some practical application, comparable with application of neighborhood notion.

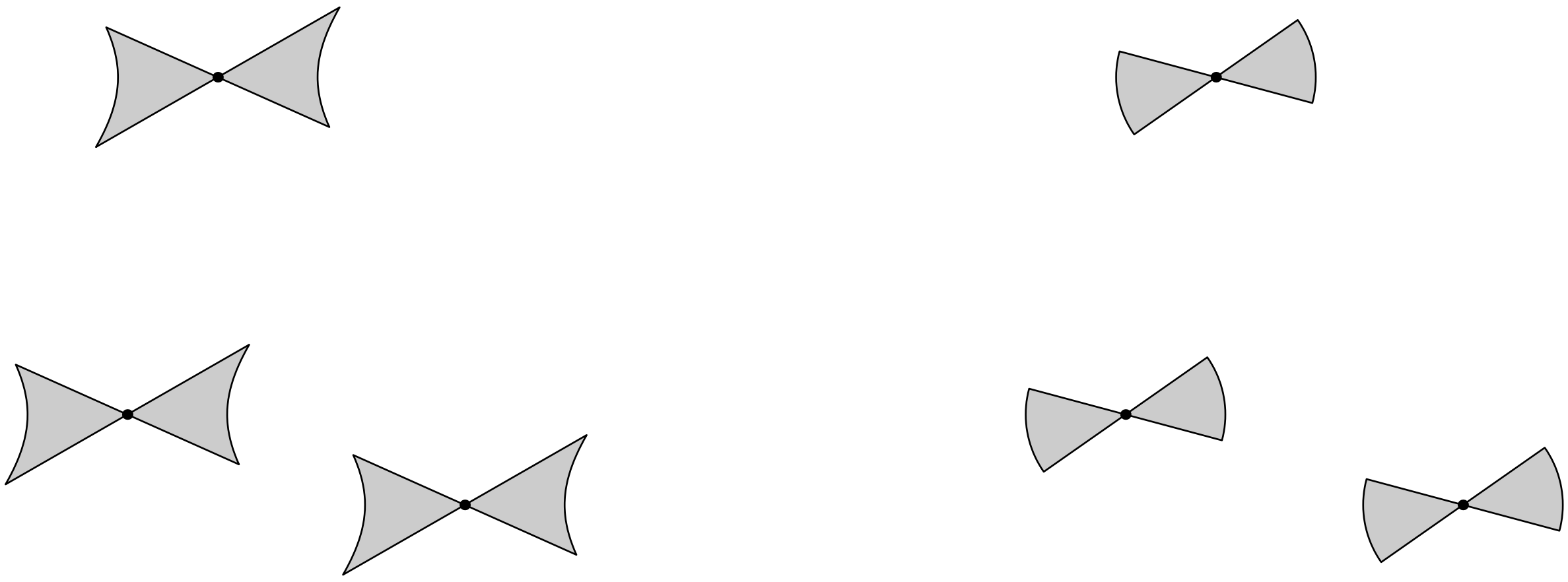

Figure 3.4: Each two points have non–intersecting generalized neighborhoods. Shown Minkowski (left) and Euclidean (right) planes.

Let generalize the neighborhood notion.

**Definition 3.2.4** (Generalized neighborhood)**.** We call *neighborhood* of point $A$ with defined in it basis $\{a^0 = A, a^1, ..., a^n\}$ in some homogeneous space $\mathbb{B}^n$ (Figure 3.4) as set of points, lines and planes of different dimensions so that $m$-dimensional measure between any (point, line or plane) of them and coordinate one (point, line or plane) of the same dimension $m$ is less then $r_m$, $m = \overline{0, n-1}$.

*Remark.* This definition is self-dual, that is, the dual definition of this neighborhood coincides with this definition of neighborhood.

*Remark.* The classical definition of neighborhood represent particular case of this definition, when all angles measures are finite. In this case $r_0 = r, r_i = \frac{\pi}{2}, i = \overline{1, n-1}$. This definition is essential in spaces with unbounded measure of some angles.

### 3.2.3. Hausdorff spaces.

One of the main reasons, the homogeneous spaces are viewed with suspection, except elliptic, Euclidean and hyperbolic, may be the fact they are not Hausdorff, unlike mentioned ones.

Recall what Hausdorff space means:

**Definition 3.2.5** (Separable or Hausdorff space)**.** Topologic space is called *separable*, or *Hausdorff*, if any two its points have non–intersected neighborhoods.

Property to be Hausdorff is important for "good" topologic spaces. For instance, in figure 3.3 on left side, there are presented neighborhoods of three points on Minkowski plane. Obviously, all these neighborhoods are intersected regardless of how large is the distance between the points and how small is the radius.

Does it mean, that homogeneous spaces are worse then others? Of course, no. The Hausdorff space definition depends on neighborhood definition, which suffers of being "prejudiced". If one changes the neighborhood to generalized neighborhood in Hausdorff space definition, everything becomes allrigt. In the figure 3.4, there are presented generalized neighborhoods on planes Minkowski (left) and Euclidean (right). Evidently, any two points on both these planes have non–intersected generalized neighborhoods, and therefore are generalized Hausdorff.

On the other hand, Euclidean plane is not "anti Hausdorff" in sense, that anti–neighborhoods of any two points (actually, of two lines) are intersected (Figure 3.3, right). But it doesn't prevent Euclidean plane to be respectable object of reserch in topology and differential geometry.

### 3.2.4. Examples of homogeneous manifolds.

As it already was noted when homogeneous space definition was given, globality of their properties makes difference between spaces and manifolds. Let's start with several definitions.

**Definition 3.2.6** (Topologic, differential and Reimannian manifold)**.** Topologic space $\boldsymbol{M}$, in which each point $\boldsymbol{P} \in \boldsymbol{M}$ have the neighborhood $\boldsymbol{U}(\boldsymbol{P})$, homeomorphic to subset of $\mathbb{E}^n$ ($n$-dimensional

Euclidean space), $\varphi : U(P) \to V \subset \mathbb{E}^n, \forall P \in M$, is called $n$-dimensional *topologic manifold.* If for each point exists smooth mapping between different homeomorphims $\varphi, \psi$, that is, if $\varphi \circ \psi^{-1} \in C^k, \forall P \in M$, the manifold is called *differentiable.* If at each point of manifold, for the tangent vectors $p, q$ is defined scalar product $p \cdot q \in \mathbb{R}$, whose value smoothly changes from one point to another, the manifold is called *Riemannian.*

This definition can be generalized to *homogeneous* topologic, differentiable and Riemannian manifolds, if homeomorphism between $M$ and Euclidean space $\mathbb{E}^n$ condition changes to homeomorphism between $M$ and homogeneous space $\mathbb{B}^n$ condition.

*Remark.* In topologic manifold definition, in requirement of its homeomorphism to $\mathbb{E}^n$, the condition that $M$ is Hausdorff is implicitely expressed. The Hausdorff condition can't be weaken, however it may be generalized, as it was shown in section 3.2.3.

Since there exist isomorphic homogeneous spaces with different signature, and manifold $M$ is defined up to homeomorphism, the signature of $\mathbb{B}^n$ may be ambiguous. Additionally, as there exist automorphisms in homogeneous spaces, that interchange subspaces with non–subspaces lineals, generally it is impossible to speak about the difference between "linear" and "non-linear" geodesics. All geodesics play the role of lineals. Finally, because usually isometry groups of manifolds are discrete, it is often impossible to detect whether two geodesics are of the "same sort" or not (they are still distinguishable in homogeneous tangent space).

Construct examples of homogeneous manifolds.

**Example.** The spheres of homogeneous spaces with signatures $\{1, -1\}$ and $\{-1, -1\}$ have equations, respectively:

$$\begin{aligned} B_1 &: x_0^2 + x_1^2 - x_2^2 = 1, \\ B_2 &: x_0^2 - x_1^2 + x_2^2 = 1. \end{aligned}$$

These two spheres intersect by limit lines $x_2 = x_1$ and $x_2 = -x_1$. Compose manifold:

$$M = \begin{cases} B_1, & |x_2| \geq |x_1|, \\ B_2, & |x_2| < |x_1|. \end{cases}$$

Two parts of the manifold are glued smoothly. In order to see it, we can represent the manifold $M$ as:

$$M = \begin{cases} x_0 = \sqrt{1 - x_1^2 + x_2^2} - 1, & |x_2| \geq |x_1|, \\ x_0 = 1 - \sqrt{1 + x_1^2 - x_2^2}, & |x_2| < |x_1|. \end{cases}$$

More than that, $M$ is completely covered by two charts of homogeneous spaces. So, it is Riemannian one (actually, semi–Riemannian).

The point $E = (1 : 0 : 0)$ has the property, that beside two limit geodesics with parabolic type, all the rest geodesics have hyperbolic type. This property is impossible on homogeneous spaces with hyperbolic angular type. However, the only point has this property. Additionally, manifold $M$ allows automorphism $\mathfrak{F} : M \to M$, defined as $\mathfrak{F}((x_0 : x_1 : x_2)) = (x_0 : x_2 : x_1)$. This automorphism is present on linear Minkowski plane with signature $\{0, -1\}$, but not on curved Minkowski planes with signatures $\{1, -1\}$ and $\{-1, -1\}$. Finally, manifold $M$ is simply connected and orientable, unlike both spaces $\{1, -1\}, \{-1, -1\}$ which are biconnected and and non–orientable.

**Example.** We can choose another parts of planes to construct the manifold $M$:

$$M = \begin{cases} B_1, & |x_2| < |x_1|, \\ B_2, & |x_2| \geq |x_1|. \end{cases}$$

In this case, all geodesics (except two limit), that pass through the point $E$, have elliptic type. As in previous example, the manifold allows automorphism that interchanges geodesics $x^1$ and $x^2$. The manifold is non–orientable and multiply–connected.

In constructed manifolds, as "raw material" it were taken positively and negatively curved Minkowski planes, which *are not* Hausdorff. They are, however, generalized Hausdorff. We essentially used generalization of point $E$ neighborhood notion here. More precisely, one part of manifold is neighborhood of system $\{E, x^1\}$ and another is neighborhood of system $\{E, x^2\}$.

### 3.3. Differential Geometry

#### 3.3.1. Geodesic as the shortest or the longest path.

**Definition 3.3.1** (Geodesic line). *Geodesic line* in space (generally, in manyfold) between two points is the line, that connects these two points and has minimal length.

In its essence, geodesic line is described as "straight line in curved space", that is, geodesic line is viewed as the "most straight" line. Our experience shows, that stretched thread (of the shortest length), free or on some surface, is rectified (becomes the most straight). However, geodesic line length isn't view as its *property*, but as its *definition.* Let show, that this definition is improper.

**Example.** Consider Minkowski plane with signature $\{0, -1\}$. Let $A, B$ be two connectable points on this plane, and point $C$ lies on line $AB$ between points $A$ and $B$. In this case the length $|AB| = |AC| + |CB|$. Let points $A, B$ be connected also with some curve having all points mutually connectable. Move the point $C$ from line to curve (Figure 3.5). From triangle inequation (proposition 2.3.2) follows, that $|AB| > |AC| + |CB|$. Each segment can be further splitted in subsegments by points in the curve. The sum of subsegment lengths is always less then segment length:

$$|AC| > |AE| + |EC|,$$
$$|CB| > |CF| + |FB|,$$
$$|AB| > |AC| + |CB| > |AE| + |EC| + |CF| + |FB|.$$

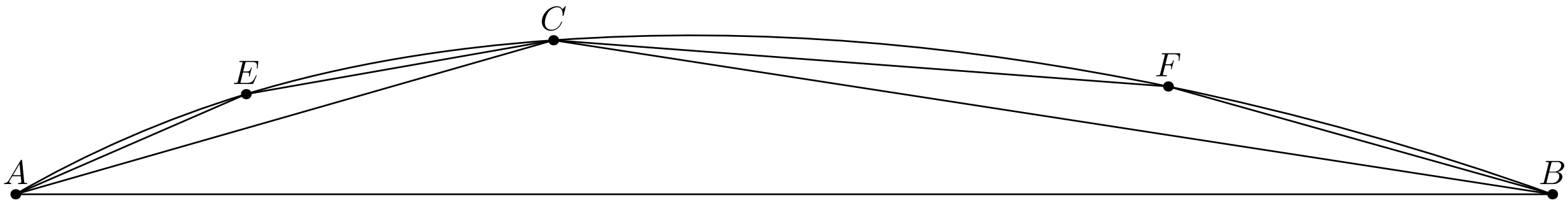

Figure 3.5: On Minkowski plane straight line between any two points is longer then curve.

It is possible to approximate arbitrarily closely the curve by polygonal path. The polygonal path length, and thus the curve length, is always less then the segment $|AB|$ length. That is why, the segment, that connects points $A$, $B$, is not the shortest, but *the longest* path (among all paths whose points are all connectable).

It is also incorrect to define geodesic line on Minkowski plane as the longest curve.

**Proposition 3.3.1.** *Generally speaking, in homogeneous space the geodesic line between points $A$, $B$ can't be defined as shortest or longest path.*

*Proof.* The triangle inequation (proposition 2.3.2) states, that the inequality sign between *lengths* ($k_1$) depends on *angular* type ($k_2$). It means, that the property of being the shortest or the longest path is not intrinsic property of a curve in which the length is measured, but the property of the space in which this curve is embedded. That is why, the same curve can be the shortest in one space and the longest in another space. □

**Example.** Consider space $\mathbb{B}^4$ with signature $\{0, 1, -1, 0\}$. Its subspace $\mathbb{B}^1$, constructed on basis vectors $\{e^0, e^1\}$, is straight line, and thus is geodesic. This line belongs to two–dimensional planes constructed on $\{e^0, e^1, e^2\}$, $\{e^0, e^1, e^3\}$, $\{e^0, e^1, e^4\}$. These planes have angular types, respectively $K_{12} = 1, K_{13} = -1, K_{14} = 0$. It means, that on the first plane any curve is *longer* then the line segment $\{e^0, e^1\}$, on the second plane any curve is *shorter* then the line segment, while on the third plane any curve has *the same length* as the line segment.

### 3.4. Conclusions of chapter 3

This chapter was focused on the application of the constructed theory to different areas of mathematics: algebraic geometry, topology and differential geometry, for whose objects it is possible to apply the notion of signature. This permits to formulate the following conclusions:

1. Based on the new signature concept, the notion of duality of homogeneous spaces was formalized. This leads to the theorem on the isomorphism of crystallographic groups of dual spaces, which was presented here [61].

2. Using the first element of the signature, the topological and metric distrinction among elliptic, parabolic and hyperbolic lines was given. This distinction leads to generalization of the neighborhood notion and Hausdorff space notions and makes homogeneous spaces the first class citizens among the metric spaces. This opens the door for homogeneous manifolds study [61].

## 4. GENERAL CONCLUSIONS AND RECOMMENDATIONS

The theory of homogeneous spaces raised from:

- the researches related to the problem of the independence of the Euclid's parallel axiom;
- the tendentions in creation of a common theory for different geometries introduced in XIX century by Arthur Cayley and Felix Klein.

The Erlangen Program was aimed to elaborate the overview and the classification of geometric spaces. However this objective was not achieved for rigid geometries of dimension greater than 2, with exception of spaces of constant curvature and spaces needed in physics. That is why the following problem is actual: **to investigate the homogeneous spaces via linear methods applying the concept of signature**. In the Chapter 2, the analytic geometry of homogeneous spaces was constructed. In the Chapter 3 the application of the constructed theory to different areas of mathematics was developed: algebraic geometry, topology and differential geometry, for whose objects it is possible to apply the notion of signature. This permits to formulate the following conclusions:

1. The new notion of signature was introduced. With its aim, the model of homogeneous space with given signature was constructed. That permits the classification of homogeneous spaces based on the concept of signature. Also, for each known homogeneous space (spaces of constant curvature, Galilean, Minkowskii, De Sitter, Anti de Sitter among other spaces) its signature, and its place in the presented classification, were found [56, 57].

2. In dependence on the space signature, the parameterized form of some important axioms were given. Based on them, the formulation and the proof of theorems, parameterized by signature, is also possible in this unified manner. Generalized trigonometric functions were introduced by means of the signature. These functions make it possible to find universal form of trigonometric equalities, common for each homogeneous space, also intruduced here [57, 59, 62].

3. By the new concept of the group of generalized orthogonal matrix, the isometry group of any homogeneous space was described. In accordance with Felix Klein's concept of geometry, described by him in the well known Erlangen Program, the isometry group of a space determines the geometry of that space [63, 154].

4. Via the new concept of signature, the type (elliptic, parabolic and hyperbolic) of geometric quantities (distances, angles, areas, volumes) was established [58, 60].

5. Based on the new signature concept, the notion of duality of homogeneous spaces was formalized. This leads to the theorem on the isomorphism of crystallographic groups of dual spaces, which was presented here [61].

6. Using the first element of the signature, the topological and metric distinction among elliptic, parabolic and hyperbolic lines was given. This distinction leads to generalization of the neighborhood notion and Hausdorff space notions and makes homogeneous spaces the first class citizens among the metric spaces. This opens the door for homogeneous manifolds study [61].

7. That permits to affirm that the the following problem is completely solved: **to investigate the homogeneous spaces via linear methods applying the concept of signature**.

Different concrete homogeneous spaces were successfully used in mathematics and physics: Euclidean, hyperbolic, Minkowskii, De Sitter. In different courses of mathematics the important role have applications of the geometry of homogeneous spaces. It follows that the analytic geometry of homogeneous spaces has a potential to enlarge and deepen these applications. Our recommandation is to use the analytic geometry of homogeneous spaces

**Recommandations:** The obtained results and developed methods can be used for:

- further study of different concrete homogeneous spaces;
- research and theoretic modelling of different physical phenomena;
- facultative and optional courses. For this objective, the software application GeomSpace was developed.

## APPENDIX

### Appendix 1. Theory of Relativity and Quantum Physics

This section presents some ideas regarding possible application of analytic theory of homogeneous spaces to different areas of physics.

#### A1.1. The shape of space–time.

As it was noted at page 128, if the angle between space and time component of space–time is hyperbolic one (that is, corresponds to Minkowski space model), then the space and the time are more interrelated than it may seem at first glance.

**Proposition A1.1.** *In Minkowski space (not necessarily linear) with signature* $\{k_1, -1, 1, 1\}$*,* ($k_1 = \{-1, 0, 1\}$)*, the time type* $k_t = K_1$ *uniquely determines the space type* $k_s = K_2$*, and the space type uniquely determines the time type.*

*Proof.* Regardless of $k_1$ value, the values of time and space types are opposite:

$$k_s = K_2 = K_1 k_2 = k_t(-1) = -k_t$$

Therefore, there are three possible cases:

1. $k_t = 1, k_s = -1$ (Anti de Sitter space),
2. $k_t = 0, k_s = 0$ (Minkowski space),
3. $k_t = -1, k_s = 1$ (De Sitter space).

□

Let examine each possible case in physical space–time. If $k_t = 1$, then the time is elliptic. It follows that points in the time axis are not separable (see section 3.2.1), or, among any three points in time axis, no one lies between other two, separating the past from the future with respect to it. It means, that in elliptic time axis doesn't exist the notion of past and future. Any moment of time repeats with interval $\pi$ (in natural elliptic time scale). It seems to not correspond to physical time. Exactly, it conflicts with cause and effect principle, which states, that no effect occurs before its cause. In this way, the cause and effect principle draws well ordered time with past and future, which is not the case of elliptic time.

If $k_s = 1$, then the space is elliptic, and thus it has finite volume equal to $\pi^2$ (in natural scale of elliptic space). It seems to disagree with universe expansion, and especially with accelerated expansion (only if expansion of elliptic space doesn't mean increase of its scale, that is, decrease of space curvature, or slowdown of time course).

From the above results the only case: $k_s = k_t = 0$. This is in a good agreement with existing measurements of space curvature, which is not very different from zero. If the reasoning above is correct, then the space curvature exactly equals to zero (on universe scale, excluding the curvature of space–time by gravity). As for the time, it is simple stated without checking, because at the moment there is no possibility to measure time curvature.

Particularly, if these reasonongs are right, then the hypothesis that physical space represents Pioncaré Sphere (Pioncaré Sphere is three–dimensional manyfold, that can be obtained from dodecahedron by glueing opposite faces after turning them by $\frac{\pi}{5}$ with respect to the opposite ones), proposed in 2003 [36], is not true, because in this case the space should be elliptic.

### A1.2. On the way toward "Theory of Everything".

As known, there is no universal proved and working theory yet, that can combine all fundamental interactions (gravitational, electromagnetic, strong and weak nuclear) as well as to explain existence of all known elementary particles [77]. The main disagreement exists between theory of relativity, that describes macrocosmos, and quantum mechanic, that describes microcosmos. One of the reasons of disagreement consists of the fact, that from quantum mechanic point of view there is minimal discreteness size of space and of time, but from theory of relativity point of view, even if distance between two events (points in space–time) is greater than this minimal size, the curve length between these events can be arbitrarily small, hence less then minimal size (see section 3.3.1).

One way to reconcile these two directions is to discretize the space–time, that is, not only to tessellate the space–time into elementary "bricks", but also to allow only discrete motions in this space–time (by the way, this method is successfully used for practical computing). The desadvantage consists of violation of invariance with respect to Lorentz transformations. Formulate the space–time discretization problem with Lorentz invariance conservation from geometry point of view:

*In Minkowski space with signature* $\{0, -1, 1, 1\}$ *to find crystallographic group with lattice size sufficiently small.*

In order to solve this problem, we need firstly to deduce all not isomorphic crystallographic groups of Minkowski space, then to choose from them one, that has appropriate lattice size. In section 3.1.4 there were constructed examples of such a groups for two–dimensional Minkowski space (with signature $\{0, -1\}$, p. 139). It was shown that the number of these groups is infinite. The same schema of crystallographic group deduction can be easily adopted also for 4-dimensional space. That is, there are infinitely many 4-dimensional Minkowski groups (the number of crystallographic groups in Euclidean space is finite for any dimension). As in case of crystallographic groups in hyperbolic space, it may happen, that, although their number is infinite, there exists minimal possible discreteness.

*Remark.* Since in suggested method of group deduction, the functions of hyperbolic angle $\cosh\varphi$ and $\sinh\varphi$ are expressed as radicals from integer numbers, in any group with space component dimension 2 or more (whole dimension is 3 or more), there are always lattice nodes lying in some limit vector. At first glance it may seem that the group loose its crystallographic property, because traditionally the distance on limit vectors is considered 0. In reality, it isn't the case, because even in limit vectors, the distance equals to 0 if and only if the points coincide. More than that, on motions by group generators the lattice nodes always map only to lattice nodes, thus any group motions map the lattice nodes only to lattice nodes, preventing their arbitrary approaching.

Second, we need to show, what means "minimal possible discreteness" for Minkowski space. Because the space is linear (see section A1.1), there exist length and time interval dilation transformations. From here results, that it can be reached any necessary level of length and time discreteness.

It is impossible to reach any level of rotations discreteness (including the Lorentz transformations). But the discreteness limitation plays different role on elliptic rotations (actually rotations) and on hyperbolic rotations (Lorentz transformations). Elliptic measure is bounded and its maximal value equals to $\frac{\pi}{2}$. That is, if minimal elliptic rotation has the angle $\varphi$, then only finite number of minimal rotations fit in maximal rotation. Opposite to that, hyperbolic measure has infinite maximal value. That is, if minimal hyperbolic rotation has the angle $\psi$, then infinitely many minimal rotations fit in maximal rotation. And exactly hyperbolic rotations create difficulties.

So, by finding sutable crystallographic group of 4-dimensional Minkowski space, the problem of Lorentz invariance is automatically solved. The problem of sufficient discreteness may be solved in one of the following ways:

- *It may happen*, that truly necessary discreteness is related only to lengths and time intervals. In this case any crystallographic group is sutable.
- If hyperbolic angle discreteness is really necessary (the elliptic angle discreteness isn't necessary), *it may happen* to find crystallographic group with its sufficient discreteness.
- If there is no such a group, then *possibly* the speed is what needs to be discrete, not the hyperbolic angle (the speed is hyperbolic tangent of this angle). If $\psi$ is the minimal angle, then speed increment $\Delta v$:

$$\begin{aligned}\tanh\psi &= \tanh((n+1)\psi - n\psi) = \frac{\tanh(n+1)\psi - \tanh n\psi}{1-\tanh(n+1)\psi\tanh n\psi},\\ \Delta v &= v_{n+1} - v_n = c\tanh(n+1)\psi - c\tanh(n\psi)\\ &= c\tanh\psi(1-\tanh(n+1)\psi\tanh n\psi).\end{aligned}$$

Then, on large speeds ($n \to \infty$):

$$\begin{aligned} \tanh n\psi &\to 1, \\ \tanh(n+1)\psi &\to 1, \\ \tanh(n+1)\psi \tanh n\psi &\to 1, \\ \Delta v = c \tanh \psi (1 - \tanh(n+1)\psi \tanh n\psi) &\to 0. \end{aligned}$$

In this way, the larger is the speed, the more its discreteness is possible. If the speed is small, then its arbitrary discreteness can be achieved by changind the crystallographic group of Minkowski space by crystallographic group of Galilean space, which allows angle dilation transformation, and thus arbitrary discreteness of speed.

### A1.3. Looking at space–time with photon's eye.

As known, in physical space–time, where the relativity principle of Einstein is valid, all inertial reference system are equivalent. In other words, in Minkowski space all orthonormal bases (complete orthonormal vector families) are congruent. Regardless of the relative speed of one body to another, the space–time looks equally from both positions.

Actually, it is related only to bodies, whose speed is less then light speed. We would like to say, that no body has the speed equal to light speed. In reality this is not true. Photons speed exactly equals to light speed. And the space–time looks differently from photon's position, than from position of bodies with smaller than light speed. As it was shown in section 2.6.6, the photon's space–time basis *is not a limiting case* of Minkowski space basis.

We can theoretically describe what does the space–time looks like from photon's position, even if it is experementally impossible. For this purpose it should be noted, that along the light ray, the space coordinate is uniquely defined by time coordinate and vice–versa. These two coordinates are not independent. We can change two coordinate vectors, the time one and the space one along the photon direction, by one limit vector (two previous vectors are its decomposition vectors). So, the photon space represents the limit lineal in Minkowski space. Make use of algorithm 2.5 to find this limit lineal signature.

The Minkowski space signature is $\{0, -1, 1, 1\}$. Its basis consists of vectors $\{x_0, x_t, x_{s_1}, x_{s_2}, x_{s_3}\}$, where $x_0$ is meta–space vector orthogonal to space, $x_t$ is time vector and $x_{s_1}, x_{s_2}, x_{s_3}$ are space vectors. The photon space basis consists of vectors $\{x_0, x_{s_1}, x_{s_2}, x_l\}$, where $x_l = x_t + x_{s_3}$ is limit vector and vectors $x_{s_1}, x_{s_2}$ are orthogonal to photon direction. Among the vectors $x_0, x_{s_1}, x_{s_2}, x_l$ there are two vector equivalence groups: $\{x_0\}$ and $\{x_{s_1}, x_{s_2}, x_l\}$. In the first group there is no motions between the vectors (there is one vector), hence there is no signature elements. In the second group there are two indexed and one limit vectors. The motion type between indexed vectors $x_{s_1}, x_{s_2}$ equals to 1 and between indexed and limit vectors equals to 0. So, two parts of signatures are $\{\}, \{1, 0\}$. The photon space signature is composed from these parts by putting the number 0 between them. Finally, complete photon space signature is $\{\{\}, 0, \{1, 0\}\} = \{0, 1, 0\}$.

Firstly note, that this signature has dimension 3 rather than 4 (by the first corollary from proposition 2.6.8). It is so, because one space coordinate is "merged" with time coordinate. In this way, photons (and generally, any body having the speed of light) fundamentally differ from objects having less then light speed. Really, if some object has less then light speed, there are reference systems in which it is at rest. For the photon, there is no such reference system. The only possiblity to fix its space coordinate is to fix also its time coordinate, that is, to stop the time run.

Secondly note, that in photon space there is no time, instead there is limit vector (of space and time "merge"), that essentially differs from space one and from time one. While all space vectors are equivalent, this limit vector is non–interchangeable with space vectors. The rotation from space vector to limit vector (its physical meaning is difficult to explain) is parabolic. However, unlike the Galilean space, where the rotation from time to space is also parabolic, in photon space the parabolic rotation is performed in the opposed direction (for non–interchangeable vector pairs it matters).

Thirdly note, that parabolic "length" along the limit vector (distance or time interval between photons of the same ray) *is not invariant* of Minkowski space. It means, that in different reference systems, moving relative to each another, this distance or this interval may be different, but always the ratio of distance and time interval equals to light speed. Still, in any reference system, the ratio of two limit lengths–intervals between two photon pairs is the same.

When introduce the curvature in Minkowski space, $\{k_1 = \pm 1, -1, 1, 1\}$, then, applying the algorithm of lineal signature finding, obtain its signature $\{k_1, 1, 0\}$. So, the curvature will affect only space component of photon space, leaving its limit component unchangeable.

## Appendix 2. Geometrical Optics

In this section the "non Euclidean physics" is analyzed — how may look physical phenomena, if the space would not be Euclidean. Only space component of space–time is analyzed and only geometrical optics (see e.g. [44]). The light nature isn't considered in this chapter.

### A2.1. Light ray as line model.

As known, the light propagates only by straight line. In some cases (for example, in case of gravitational lenses) the apparent light trajectory may be a curve. In fact the space becomes curved, and the light tries to find the most straight path, moving along a geodesic (see section 3.3.1). In some cases (for instance, near the black holes) there exist areas of space, that are never illuminated. Anyway, as our senses are not able to detect the space curvature, we take the light ray as straight line standard.

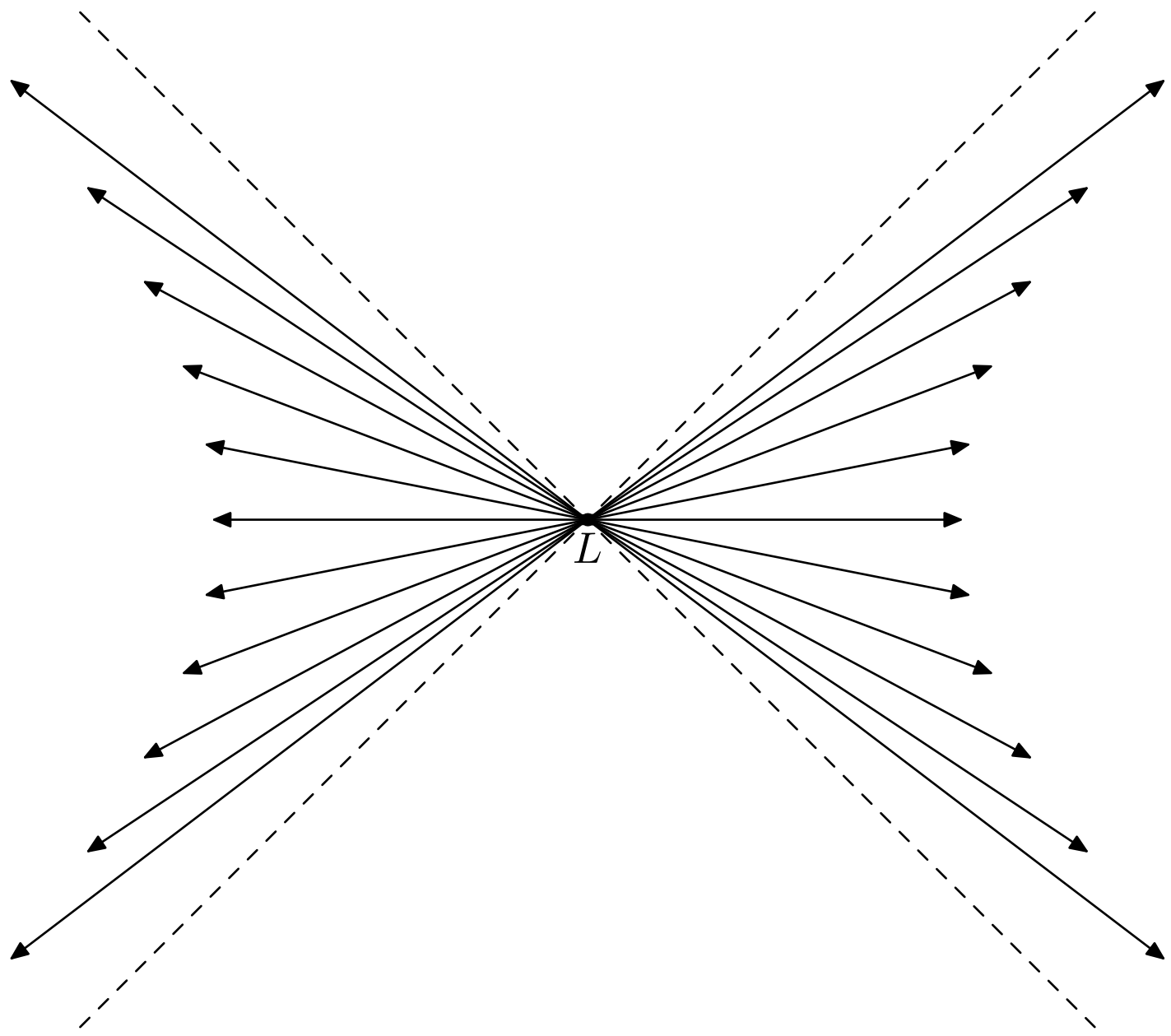

Figure A2.1: The light propagates only by straight line. When $k_2 = -1$, then upper and lower areas between limit directions are not illuminated.

When constructing the illumination model we take this proposition as axiom: the light propagates only by straight line.

**Proposition A2.1.** *In three–dimensional homogeneous space with orthonormal basis* $\{e^0, e^1, e^2, e^3\}$, *the light ray vector has index* 1.

*Proof.* All lines in space $\mathbb{B}^3$ are images of subspace $\mathbb{B}^1$ on all possible motins. For each line we can choose space basis in such a way, that this line to coincide with $\mathbb{B}^1$, constructed on coordinate vectors $\{e^0, e^1\}$. So, the light ray vector index always equals to 1. By lemma 2.4.3, this index doesn't change on motions. □

**Corollary.** *The light doesn't propagate in directions of vectors with indices* 2, 3, *if they are not equivalent with vectors with index* 1.

In other words, all space points, that are unconnectable with light source position, remain unilluminated (Figure A2.1). Of course, it doesn't mean, that there are points unable to illuminate. Any point can be illuminated by right positioning of light source.

*Remark.* As it was shown, the properties of one–dimensional lineals, different from lines, may be very different from line properties. Additionally, there are no space motions that map a line to a lineal, that is not a line.

### A2.2. Transparency.

Speaking about transparency, we should distinguish between at least two different phenomena:

- Transparency related to property of substance to *filter* source light of certain wave length. This transparency is related to light nature, which isn't discussed here, as we agreed, depends only on substance properties and gives the same result regardless of transparent material layer thickness.

- Transparency related to property of substance to *absorb* the intensity of transmitted through it light. On such transparency, the light intensity after material layer depends on its thickness.

Of course, both phenomena can occur simultaneously. We consider only the last case.

Let intensity of source light ray be equal to $I_0$. Falling on the surface, one part of the light is absorbed, other is dispersed:

$$I_0 = I_a + I_d \tag{A2.1}$$

Here, $I_a$ is intensity of absorbed light and $I_d$ — of dispersed light.

Let, after passing of distance $l$ through transparent layer, the ray intensity becomes $I_a(l) = I_a\alpha$, $0 \le \alpha \le 1$. Then, after passing once more the distance $l$, the intensity becomes $I_a(2l) = I_a(l)\alpha = I_a\alpha^2$. So, the value $\alpha$ is a function on distance with property $\alpha(2l) = (\alpha(l))^2$. This property has only one function $\alpha(l) = e^{-pl}$ (the "−" sign means, that after passing positive distance, its value becomes less then unite). Parameter $p$ depends on substance ability to absorb the ligt intensity by length unite. The larger is $p$, the less transparent is the substance. Let introduce substance transparency $\lambda = \frac{1}{p}$. Then:

$$I_a(l) = I_a e^{-\frac{l}{\lambda}} \tag{A2.2}$$

Consider all substances more or less "transparent", that is, let equation (A2.2) describe light intensity transmitted through any substance. It seems reasonable, having the inter–molecular distances, where can penetrate light photons. Even absolutely opaque for visible light substances are rather transparent for certain particles, for instance, neutrino. When parameter $\lambda$ is close to zero, the intensity quickly decreases to nearly zero even on small thickness. When $\lambda$ tends to infinity, the intensity decreases slowly.

What happens, if the light falls to transparent layer obliquely, instead of orthogonally (Figure A2.2)? In this case changes the distance $l$. If the layer thickness equals to $l_n$, and the angle

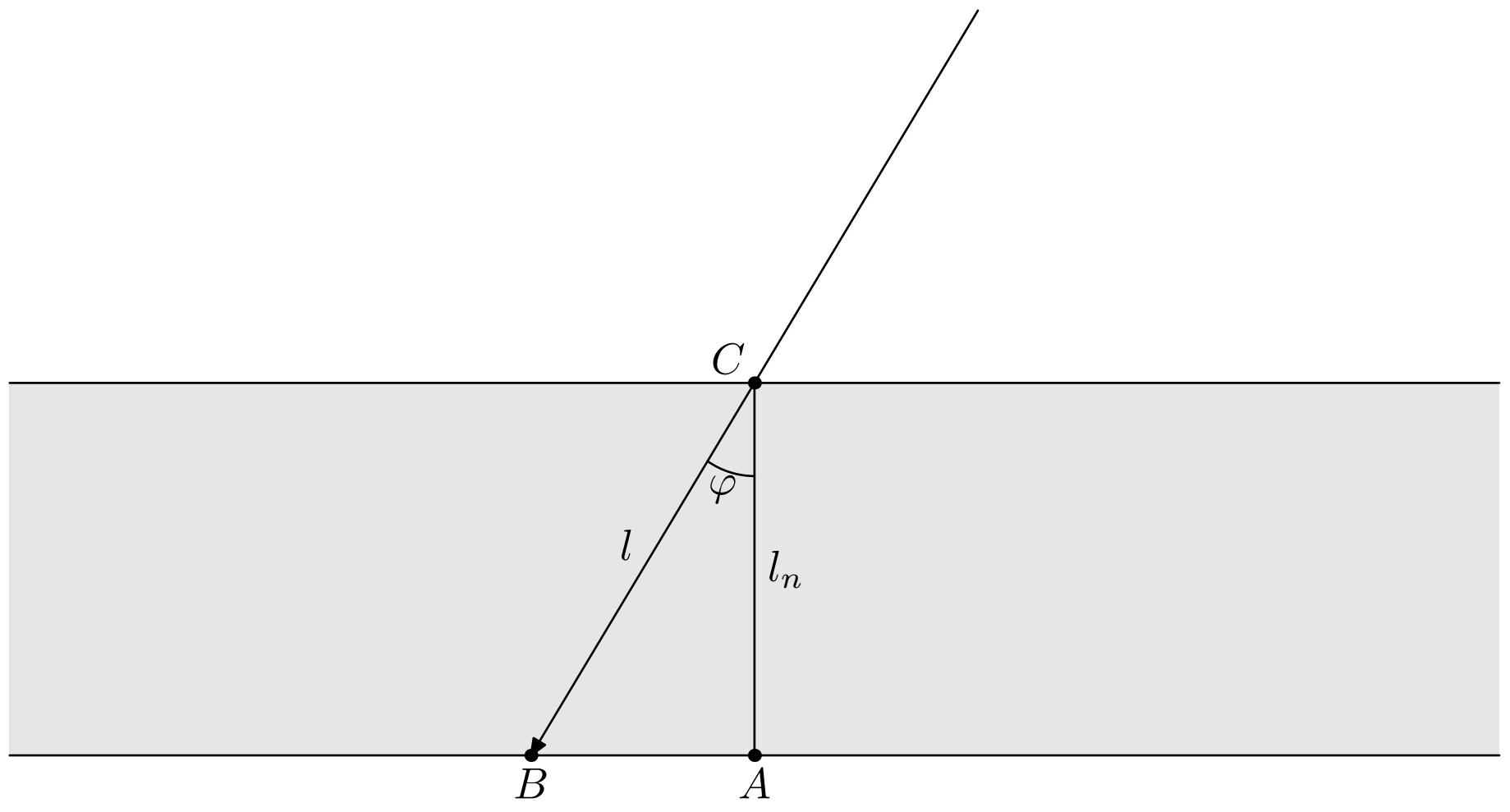

Figure A2.2: Transmission of oblique light through transparent material.

between light ray and surface normal equals to $\varphi$, then by (2.45):

$$T_1(l) = \frac{T_1(l_n)}{C_i(\varphi)} \tag{A2.3}$$

In case of signature $\{-1, 1, k_3\}$, the equation takes the form:

$$\tanh l = \frac{\tanh l_n}{\cos \varphi}$$

Since the functions from equation have the property $\tanh x < 1$, $\cos x < 1$, it may happen, that $\tanh l_n = \cos \varphi$. Then:

$$\begin{aligned} \tanh l_n &= \cos \varphi, \\ \tanh l &= \frac{\tanh l_n}{\cos \varphi} = 1, \\ l &= \infty \end{aligned}$$

In this case the light ray transmitted into transparent layer with the angle of parallelism $\varphi_0 \geq \arccos \tanh l_n$ is parallel or divergent with its lower border. The intensity equals:

$$I_a(\varphi_0) = I_a e^{-\frac{\infty}{\lambda}} = 0$$

That is, starting from the angle of parallelism, the substance is opaque.

In the rest, when thickness $l_n$ is small, $T_1(l_n) \approx l_n, T_1(l) \approx l$:

$$l = \frac{l_n}{C_i(\varphi)}$$

Because depending on the space and direction, hypotenuse may be either greater then, or less then cathetus (proposition 2.3.2), in case of angular type $k_i = 1$ thickness infinitely grows with rotation, and in case $k_i = -1$ decreases to about zero. Intensity of absorbed oblique light equals:

$$I_a(l_n, \varphi) = I_a e^{-\frac{l}{\lambda}} = I_a e^{-\frac{l_n}{\lambda C_i(\varphi)}} \tag{A2.4}$$

In other words, transparency becomes:

$$\lambda(\varphi) = \lambda C_i(\varphi) \tag{A2.5}$$

Rotating in spaces with angular type $k_2 = -1$, the layer becomes infinitely transparent on angle of incidence growth. The light is practically completely absorbed, only small portion is dispersed. As source light intensity $I_0$ doesn't change on surface rotation, the sum of intensities of dispersed $I_d(\varphi)$ and absorbed $I_a(\varphi)$ light is constant given fixed angle of incidence on the surface. By (A2.4), obtain:

$$I_0 = I_d(\varphi) + I_a(\varphi),$$
$$I_a(\varphi) = I_0 e^{-\frac{l_n}{\lambda C_i(\varphi)}},$$
$$I_d(\varphi) = I_0 - I_a(\varphi) = I_0 \left(1 - e^{-\frac{l_n}{\lambda C_i(\varphi)}}\right)$$

**A2.3. Reflection.**

Although it may seem meaningless, the reflected light ray is also a light ray. In other words, its vector has index 1. Let show, that it is always so.

**Lemma A2.2.** *Reflected light ray has index* 1.

*Proof.* Regardless of the fact, whether the mirror is a plane or a lineal, different from plane, the basis can be chosen in the way, that the mirror is a lineal constructed on coordinate vectors. In this case, if light ray vector $x$ has coordinates $\{x_0 : x_1 : x_2 : x_3\}$, then reflected vector $x'$ has coordinates $\{x_0 : \pm x_1 : \pm x_2 : \pm x_3\}$, and the only coordinate of $x'$ has opposite sign than corresponding coordinate of $x$. Then the first norm of vector $x'$ equals:

$$x' \odot_1 x' = \frac{x_0^2}{K_1} + \sum_{i=1}^{3} K_{1i}(\pm x_i)^2 = \frac{x_0^2}{K_1} + \sum_{i=1}^{3} K_{1i} x_i^2 = x \odot_1 x > 0$$

Hence, the reflected vector $x'$ index equals to 1, like initial vector $x$ index. $\square$

On nonspecular (diffuse) reflection, besides the main vector $x'$, the reflected ray spreads from reflection point also in other directions. As in case of light source, the reflected ray vector has index 1, this fact limits the possible directions (in all cases angle $\psi$ between the diffuse light ray and $x'$ is measurable). The intesity of reflected ray in direction $x'$ equals to $I_d$, and the intensity of diffuse light in direction, that makes the angle $\psi$ with vector $x'$ is, by model of diffuse reflection:

$$I_d(\psi) = I_d C_i(\psi) \tag{A2.6}$$

In all cases, since on radiation, the energy transmission occurs from the light source toward the surface, and not backward, the intensity can't be negative, even when the angular type is $k_i = 1$ and $\frac{\pi}{2} < \psi < \pi$, $C_i(\psi) < 0$. In real world, the surfaces, whose perpendicular is turned more than by $\frac{\pi}{2}$ from ray vector, are shadowed ($I_d(\psi) = 0$). Reasoning in similar fashon, from the low of energy conservation results that the diffuse light intensity $I_d$ can't exceed the source light intensity $I_0$, even on angular type $k_i = -1$, when $C_i(\psi)$ infinitely grows. Let introduce function ($min < max$):

$$clamp_{[min,max]}(x) = \begin{cases} min, & x < min; \\ x, & min \le x \le max; \\ max, & x > max \end{cases} \tag{A2.7}$$

Equation (A2.6) becomes:

$$\begin{aligned} I_d(\psi) = clamp_{[0,I_0]}(I_d C_i(\psi)) = clamp_{[0,I_0]}\left(I_0\left(1 - e^{-\frac{l_n}{\lambda C_i(\varphi)}}\right)C_i(\psi)\right) \\ = I_0\, clamp_{[0,1]}\left(\left(1 - e^{-\frac{l_n}{\lambda C_i(\varphi)}}\right)C_i(\psi)\right) \end{aligned} \tag{A2.8}$$

**A2.4. Luminosity.**

Consider, that the light propagates *evenly* in all directions. Because there are spaces, where the sphere has infinite surface, it is impossible to compute the luminosity based on light source intensity (in this case on any bounded portion the luminosity is finite, while the light source intensity should be infinite). Instead, it is possible to compute how luminosity of finite surface decreases on its departure from the light source.

As it was shown on area element deduction in polar coordonate system (section 2.8.2), the length of arc $l$ of a circle with radius $r$, bounded by the angle of measure $\varphi$, equals:

$$l = S_1(r)\varphi.$$

Arguing similarly, obtain surface $s$ of the sphere portion with radius $r$, bounded by solid angle with measure $\vartheta$, equal:

$$s = S_1^2(r)\vartheta \tag{A2.9}$$

The solid angle $\vartheta$ is proportional to the "ray quantity" $\sigma_0$, within this angle. Let fix the area $s$ and vary the distance $r$. Then luminosity $\sigma(r)$ on distance $r$ away from light source is:

$$\sigma(r) = \sigma_0\vartheta = \frac{\sigma_0 s}{S_1^2(r)} \tag{A2.10}$$

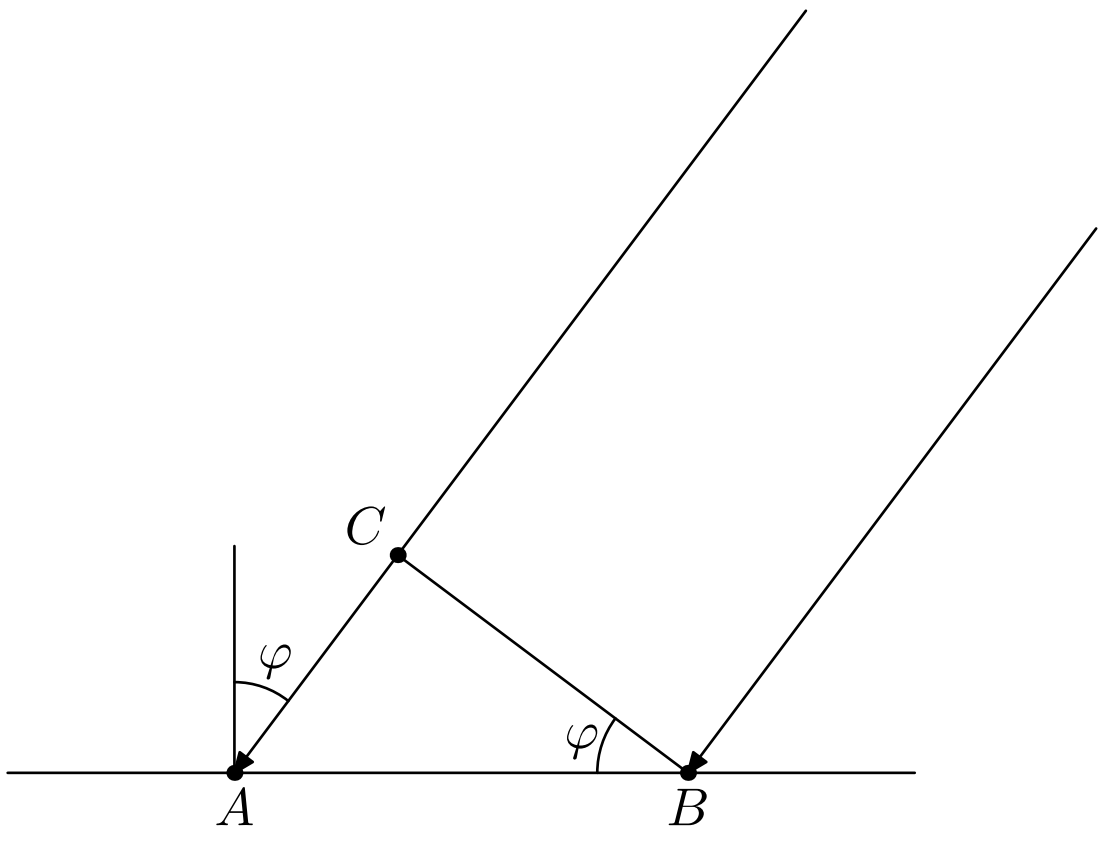

Figure A2.3: Luminosity on oblique light.

When compute the surface luminosity, assume that this surface is orthogonal to the ray. When this is not the case, in equation (A2.10) instead of surface $s$ (region $AB$ in figure A2.3), enters its normal component $s_n$ (region $BC$ in figure A2.3).

If the area is small, and the angle between light ray and surface perpendicular equals to $\varphi$, then the luminosity, given (2.45) and $T(x) \approx x$ equals:

$$\sigma(r,\varphi) = \frac{\sigma_0 s_n}{S_1^2(r)} = \frac{\sigma_0 s\, C_i(\varphi)}{S_1^2(r)} \tag{A2.11}$$

If the surface region is rather large, then, generally speaking, each ray has different angle of incindce to the surface.

As in case of diffuse illumination, the value $C_i(\varphi), i = 2,3$ may be less then one or greater then one. Even in one homogeneous space, on surface inclination in different directions, this value may be both less then and greater then one. This property follows from triangle inequation (proposition 2.3.2), thich states, that the hypotenuse of right triangle may be shorter, or longer then its cathetus, or equal to it.

As in case of diffuse illumination, the limunosity can not be negative and can't exceed $\sigma_0$. That is why, equation (A2.11) should be rewritten as:

$$\sigma(r,\varphi) = clamp_{[0,\sigma_0]}\left(\frac{\sigma_0 s\, C_i(\varphi)}{S_1^2(r)}\right) = \sigma_0\, clamp_{[0,1]}\left(\frac{s\, C_i(\varphi)}{S_1^2(r)}\right) \tag{A2.12}$$

**A2.5. Refraction.**

The refraction phenomena (Figure A2.4) is related to the difference between the light speed of propagation in different medium. The ratio of light speed in some medium and light speed in vacuum, is called index of refraction of a substance $\eta = \frac{c_0}{c}$, where $c$ is the light speed in medium propagation and $c_0$ is the light speed in vacuum.

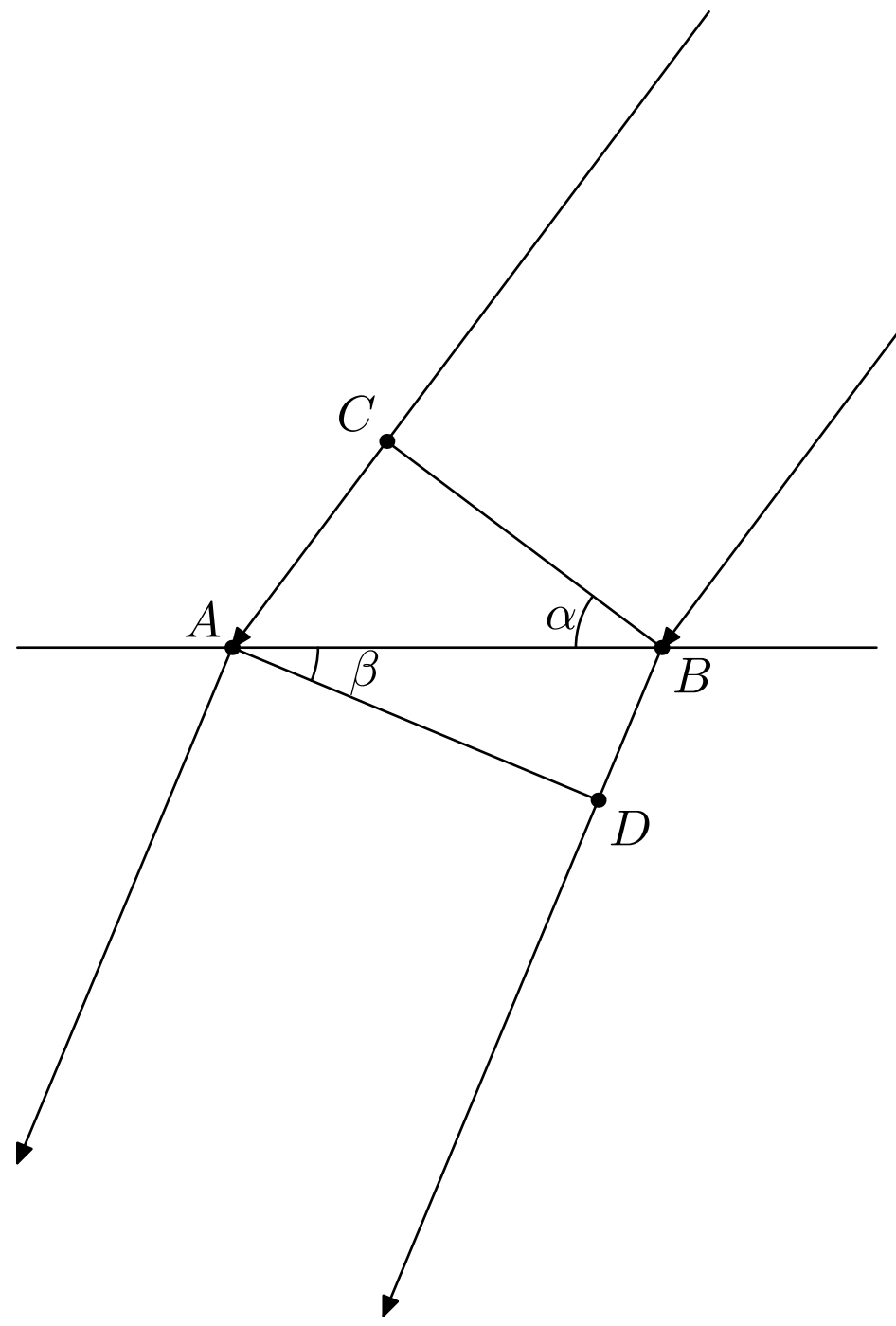

Figure A2.4: Light refraction.

Let some light ray angle of incidence equal to $\alpha$, and transmission angle equal to $\beta$. Consider two right triangles: $\triangle ABC$, $\triangle ABD$ with common hypotenuse $AB$. At the time the left margin of the ray passes the distance $AC$, its right margin passes the distance $BD$. It means, that catheti $AC$ and $BD$ are proportional to $c_0$ and $c$ respectively:

$$\frac{AC}{BD} = \frac{c_0}{c} = \eta$$

At other hand, by (2.47), considering small $AB$ ($S(x) \approx x$):

$$\frac{S_1(AC)}{S_i(\alpha)} = S_{1i}(AB) = \frac{S_1(BD)}{S_i(\beta)},$$
$$\frac{S_i(\alpha)}{S_i(\beta)} = \frac{S_1(AC)}{S_1(BD)} \approx \frac{AC}{BD} = \eta \tag{A2.13}$$

Because the function $S(x)$ is increasing, and index of refraction of substance is always greater then one, the angle of incidence is always greater then the angle of transmission:

$$\frac{S_i(\alpha)}{S_i(\beta)} = \eta > 1,$$
$$S_i(\alpha) > S_i(\beta),$$
$$\alpha > \beta \tag{A2.14}$$

## Appendix 3. GeomSpace Project

The GeomSpace project (Figure A3.5) was created as testing platform for theoretic results in analytic geometry of homogeneous spaces. Today it is independent project, with aim of which it is possible to demonstrate the theory achievements. The project address is http://sourceforge.net/projects/geomspace/. As it is expected from this kind of collaboration of theory with practice, these projects are closely related. At one hand, new theoretic results become available in GeomSpace. At other hand, if some algorithm doesn't work in GeomSpace, it is revised also in theory.

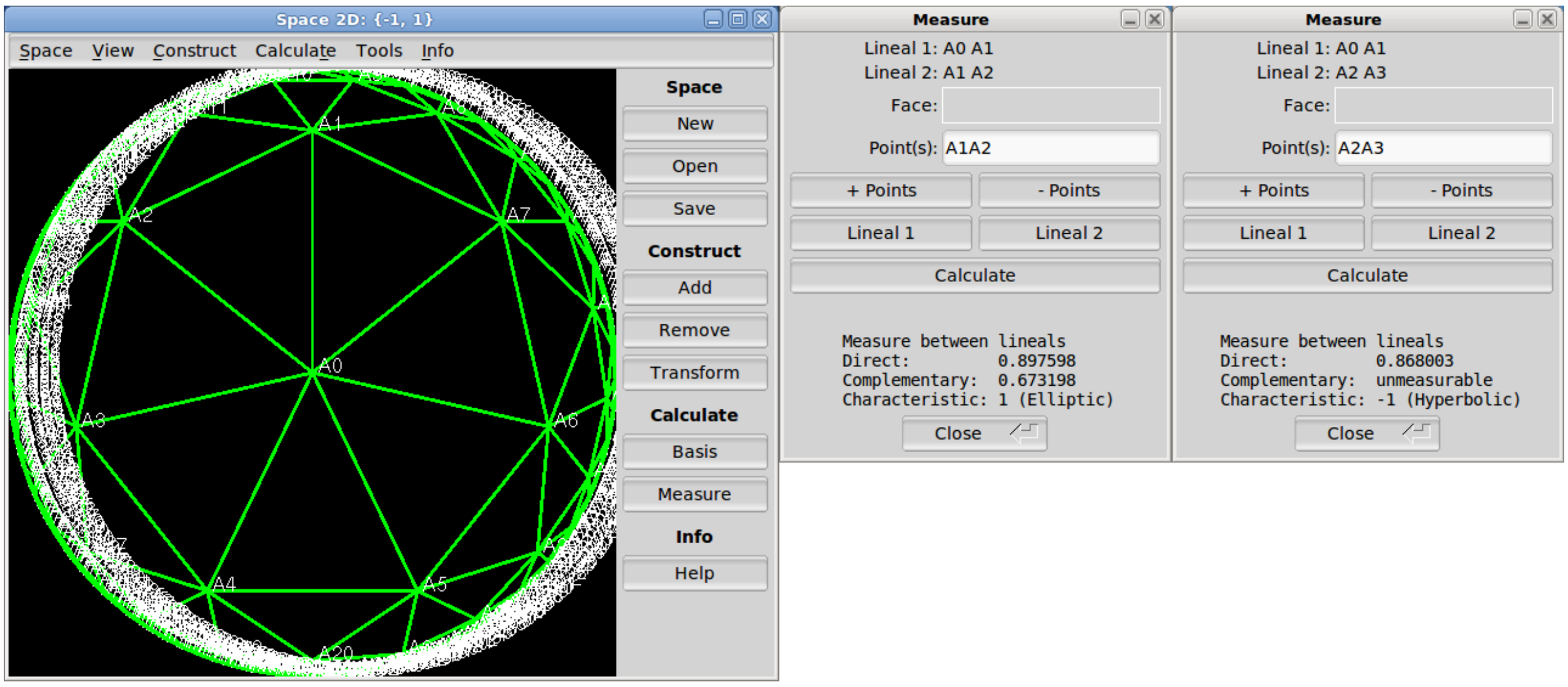

Figure A3.5: Using GeomSpace to find the angle between intersected lines and the distance between divergent lines on hyperbolic plane.

GeomSpace is interactive geometry system for homogeneous spaces. It is unique in several aspects:

1. It focuses on the space, rather than on space model;
2. It is unlimited in space dimension;
3. It operates optimal for spaces: Euclidean, Riemannian, Lobachevsky (hyperbolic), Galilean (including curved ones), Minkowski (including curved ones) and others. Actually, GeomSpace uses the same algorithms for all them;
4. It makes use of OpenGL hardware acceleration advantages not only for Euclidean motions, but also for all possible motions.

Additionally, the advantages of GeomSpace over similar software of interactive geometry include:

1. Small size;
2. Light and portable or common standards and dependent libraries (*OpenGL, FLTK*);
3. Really crossplatform (compiles without source change under Linux, Windows, OS X, FreeBSD, Solaris, including 32 and 64 bit versions);
4. Open source licence GPL.

GeomSpace users can (now or in future):

1. Choose any homogeneous space;
2. Construct, modify or remove figures in it;
3. Move selected figures in the space;
4. Find point coordinates and lineal orthonormal basis;
5. Find distance between points or a point and a lineal; distances, angles or inclinations between lineals (all sort of measures between figures);
6. Construct the lineals sum, intersection or difference;
7. Find volumes of bounded figures;
8. Exchange constructed models in GeomSpace format.

GeomSpace comes with ready to use set of models, constructed in different spaces, so that one can observe the properties of space motions as well as experiment with constructions or calculus. At the project website the latest models are always available.

**Use case.** Consider the problem:

**Given:** Two lines $a$ and $b$ on $\mathbb{E}^2$ (Euclidean plane).

**To find:**

1. Determine whether the lines $a$ and $b$ are intersected or parallel
2. (a) If $a$ and $b$ are intersected, find the angle between them.
   (b) If $a$ and $b$ are parallel, find the distance between them.

Usually, in order to resolve this problem, one needs 3 different algorithms:

1. Algorithm of line parallelism detection.
2. Algorithm of angle computing between intersected lines (not applicable if the lines are parallel).
3. Algorithm of distance computing between parallel lines (not applicable if the lines intersect each other).

Now, let change the problem issue:

**Given:** Two lines $a$ and $b$ on $\mathbb{S}^2$ (elliptic plane).

**To find:** The angle between $a$ and $b$.

Since on elliptic plane all lines intersect, the problem greatly simplifies. Nevertheless, the algorithm of the angle computing between intersected lines, that is used in previous problem, isn't applicable here. We need completely new algorithm for elliptic plane.

Change the problem issue again:

**Given:** Two lines $a$ and $b$ on $\mathbb{H}^2$ (hyperbolic plane).

**To find:**

1. Determine whether the lines $a$ and $b$ are intersected, parallel or divergent.
2. (a) If $a$ and $b$ are intersected, find the angle between them.
   (b) If $a$ and $b$ are parallel, find the inclination between them.
   (c) If $a$ and $b$ are divergent, find the distance between them.

Not only appears new relative position of lines (divergent), but all algorithms, including the algorithm of relative position determination, need to be revised for hyperbolic plane case. As earlier, no one algorithm is applicable in improper cases.

GeomSpace users have two advantages. Firstly, the only algorithm is necessary: computation of measure between lineals algorithm, which works correctly in all cases and at the same time determines relative position of lineals. Secondly, this algorithm works in all homogeneous spaces. So, the problem can be issued as follows:

**Given:** Two one–dimensional lineals $a$ and $b$ on a homogeneous plane $\mathbb{B}^2$.

**To find:**

1. Determine relative position between $a$ and $b$.
2. Find the most appropriate measure between $a$ and $b$.

In order to resolve this problem the algorithm 2.8 of measure computation between lineals can be used. As result, it gives the following information:

- Measure type (elliptic, parabolic or hyperbolic), and
- Measure value

Having Euclidean plane signature $\{0, 1\}$, the possible measure types are 0 (parabolic) when the lines are parallel, or 1 (elliptic) when the lines are intersected. Elliptic plane has signature $\{1, 1\}$, which leaves place for only one measure signature: 1 (elliptic). Hyperbolic plane has signature $\{-1, 1\}$, therefore the possible types are 1 (elliptic) when the lines are intersected, 0 (parabolic) when the lines are parallel, and $-1$ (hyperbolic) when the lines are divergent. In all cases one universal algorithm is quite enough. In Figure A3.5 there are shown how GeomSpace computes the measure between intersected and divergent lines of $\{3, 7\}$ tiling on hyperbolic plane.

Note, that algorithm should not know a priori, what relative positions between two lineals are possible on a plane. It works correctly for all them. The lineal signature (algorithm 2.5) together with type and value of measure between them gives all necessary information about the lineals relative position.

It is important, because there exist much more complex situations. For example, on pseudo–hyperbolic plane (or, negative curved Minkowski plane) with signature $\{-1, -1\}$ there are 10 possible relative positions of lineals (use term “line” when the lineal is a line, and “lineal” when it isn’t a line):

1. Two intersected elliptic lineals with hyperbolic angle between them;
2. Elliptic lineal intersects parabolic limit lineal with infinite hyperbolic angle between them;
3. Elliptic lineal intersects hyperbolic line with hyperbolic complementary angle between them (direct angle isn’t measurable in this case);

4. Two intersected parabolic limit lineals with infinite hyperbolic angle between them;
5. Two parallel parabolic limit lineals with parabolic distance between them;
6. Parabolic limit lineal intersects hyperbolic line with infinite hyperbolic angle between them;
7. Parabolic limit lineal is parallel to hyperbolic line with parabolic inclination between them;
8. Two intersected hyperbolic lines with hyperbolic angle between them;
9. Two parallel hyperbolic lines with parabolic inclination between them;
10. Two divergent hyperbolic lines with elliptic distance between them.

## List of Figures

**List of Definitions**

## List of Algorithms

## List of Tables

## Declaration on liability

The undersigned, declare under sole responsibility that the material presented in the doctoral thesis, is the result of personal research and development. I realize that otherwise would be liable in accordance with applicable law.

Popa Alexandru

---

April 30, 2017

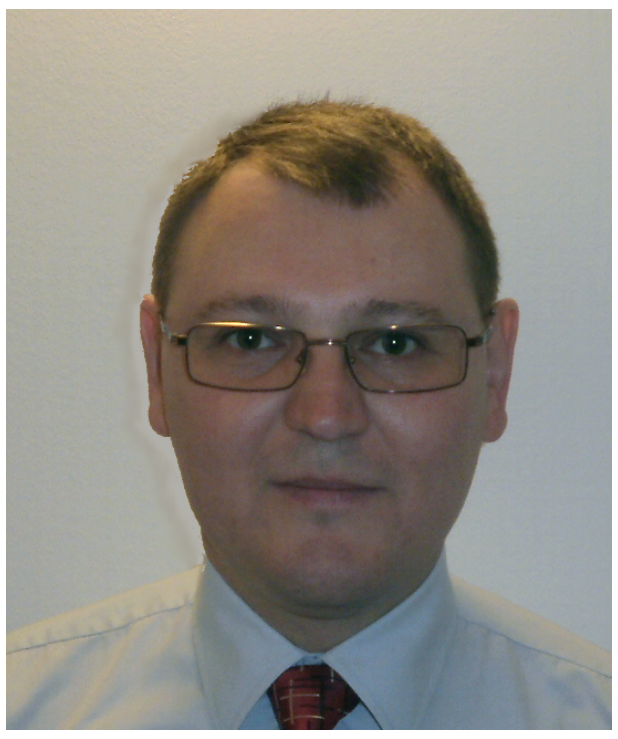

# Curriculum Vitae

Name: Popa Alexandru

Birth: July 22, 1977, Chișinău, Moldova

Nationality: Moldova, Romania

Studies:

**1995 — 2000** State University of Moldova (Bachelor degree in Theoretic Mathematics, 9.83 mark)

**2000 — 2001** West University of Timisoara (Master degree in Computational Mathematics, 9.5 mark)

Field of scientific interests: I have the geometrical intuition for many–dimensional and non–Euclidean spaces. Research the properties of figures and of spaces that can be generalized to any homogeneous geometry and extrapolated to higher dimensional spaces.

Proffessional activity:

**2001 — 2003** West University Vasile Goldiș, Proffessor Assistant

**2004 — 2005** Linux Magazin SRL, Technical Redactor

**2005 — 2011** Siemens VDO / Continental Automotive SRL, Software Ingineer

**2011 — present** SSI Schaefer SRL, Software Ingineer

Projects: **GeomSpace** http://sourceforge.net/projects/geomspace/

Scientific activity:

**1999** Chişinău, Moldova — Students Scientific Conference, ed. IV,

**2000** Chişinău, Moldova — Students Scientific Conference, ed. V,

**2003** Timișoara, România — The 5th International Workshop on Symbolic and Numeric Algorithms for Scientific Computation SYNASC2003,

**2003** Arad, România — Scientific Communication Session "Aradean Academic Days", ed. XIII,

**2009** Alba Iulia, România — International Conference on Theory and Applications in Mathematics and Informatics,

**2010** Moscow, Russia — International conference "Metric geometry of surfaces and polyhedra", dedicated to 100th anniversary of N. V. Efimov,

**2010** Moscow, Russia —The International Conference "Geometry, Topology, Algebra and Number Theory, Applications" dedicated to the 120th anniversary of B. N. Delone,

**2014** Chișinău, Moldova — The Third Conference of Mathematical Society of Moldova IMCS-50,

**2015** Tula, Russia — International conference, dedicated to 85th anniversary of professor S. S. Ryshkov. Algebra, Number Theory and Discrete Geometry: Modern Problems and Applications,

**2015** Iași, România — The 8th Congress of Romanian Mathematicians,

**2016** Chişinău, Moldova — International Conference Mathematics & Information Technologies: Research and Education, MITRE — 2016.

Publications: **2** monographies, **13** scientific publications, **10** technical publications, **1** software product

Awards:

**2000** Chişinău — Student Eminent Medal, State University of Moldowa

**2000** Chişinău — Diploma for Academic Performance and Extracurricular Activities, Soros Foundation Moldova

Affiliation: Academy of Sciences of Moldova, Institute of Mathematics and Computer Science

Languages:
1. Russian — Native speaker
2. Romanian — Advanced
3. English — Fluent
4. Classical Greek — Beginner
5. Ancient Hebrew — Beginner

Contact:

**Address:** Str. Cronicar Ion Neculce nr. 11 bl. L40 ap. 14, 300109, Timișoara, România

**Phone:** +40 742 967 251

**Email:** alpopa@gmail.com